\documentclass[12pt,twoside, english]{amsbook}

\usepackage{amsthm,thmtools,xcolor}

\usepackage{mathrsfs}
\usepackage{mathtools}
\usepackage{amssymb}
\usepackage{verbatim}
\usepackage{amsmath}
\usepackage{comment}
\usepackage{amsthm,thmtools,xcolor}
\usepackage[colorlinks,linkcolor=blue,citecolor=blue, pdfstartview=FitH]
{hyperref}
\usepackage{verbatim}
\usepackage{amsmath}
\usepackage{bm}
\usepackage{a4wide}
\usepackage{hyperref}

\usepackage{babel,latexsym,amsmath,amsthm,eucal}
\usepackage[utf8]{inputenc}

\usepackage[T1]{fontenc}
\usepackage{times}
\usepackage{amssymb,latexsym}

\usepackage{amssymb,latexsym}
\usepackage{enumerate}

\newcommand{\dbar}{\ensuremath{\overline\partial }}

\makeatletter
\newcommand{\sumprime}{\if@display\sideset{}{'}\sum%
            \else\sum'\fi}
\makeatother

\newcommand{\del}[2]{\ensuremath{\frac{\partial  #1}{\partial  #2}}}

\newcommand{\ip}[2]{\ensuremath{\left\langle #1, #2\right\rangle}}

\newcommand{\s}[2]{\ensuremath{   \sum_{#1}^{#2}    }}

\newcommand{\vect}[1]{\ensuremath{  \mathbf{#1}   }}

\newcommand{\rbrac}[1]{\ensuremath{  \left( #1 \right)  }}

\newcommand{\fr}[2]{\ensuremath{  \frac{#1}{#2}  }}

\newcommand{\bt}[1]{\textbf{#1}}

\newcommand{\id}[1]{\ensuremath{\bt{1}_{#1}}}

\newcommand{\ita}[1]{\emph{#1}}

\newcommand{\set}[1]{\ensuremath{\mathbb{#1}}}

\newcommand{\ld}{\ensuremath{\lambda}}

\newcommand{\Ld}{\ensuremath{\Lambda}}

\newcommand{\cbrac}[1]{\ensuremath{\left\{#1\right\}}}

\newcommand{\sub}{\ensuremath{\subseteq}}

\newcommand{\norm}[1]{\ensuremath{\left| \left |#1 \right|\right|}}

\newcommand{\supp}[1]{\ensuremath{\text{supp}(#1)}}

\newcommand{\abs}[1]{\ensuremath{\left| #1 \right|}}

\DeclareSymbolFont{rsfs}{U}{rsfs}{m}{n}
\DeclareSymbolFontAlphabet{\mathscrsfs}{rsfs}

\newcommand{\fancy}[1]{\ensuremath{\mathscrsfs{#1}}}
\newcommand{\ml}[1]{\ensuremath{#1^{\uparrow_{\perp}}}}

\usepackage{import}

\begin{document}

\numberwithin{equation}{section}

\newtheorem{theorem}{Theorem}[section]
\newtheorem{proposition}[theorem]{Proposition}
\newtheorem{conjecture}[theorem]{Conjecture}
\def\theconjecture{\unskip}
\newtheorem{corollary}[theorem]{Corollary}
\newtheorem{lemma}[theorem]{Lemma}
\newtheorem{observation}[theorem]{Observation}
\newtheorem{definition}{Definition}
\numberwithin{definition}{section} 
\newtheorem{remark}{Remark}
\def\theremark{\unskip}
\newtheorem{kl}{Key Lemma}
\def\thekl{\unskip}
\newtheorem{question}{Question}
\def\thequestion{\unskip}
\def\theexample{\unskip}
\newtheorem{problem}{Problem}
\newtheorem{example}{Example}[section]

%



\title[Thesis]{A Hilbert Bundles Description of Complex Brunn-Minkowski Theory\\[0.4cm]
\begin{small}
\end{small}
}

\author{Tai Terje Huu Nguyen}
\date{\today}

\maketitle
\tableofcontents

\chapter*{Preface}
This thesis is submitted in partial fulfillment of the requirements for the degree of Philosophiae Doctor (Ph.D.) at the Norwegian University of Science and Technology. 

A special thanks goes to my advisor, Xu Wang, who has contributed immensely in making these past few years especially enjoyable. I would also like to thank my co-advisor, Bo Berndtsson, who has been, and who continues to be, a great source of inspiration. It has been a great pleasure, and an honor, to work with both of you, and I hope our collaboration continues in the future.

Finally, my greatest gratitude goes to my parents and my family, to whom I owe everything I am today.\\[0.3cm]

Tai Terje Huu Nguyen, \\[0.2cm]
Trondheim, December 2022.

\chapter*{Introduction}
The theme of the following thesis is two-fold. One of the folds concerns Hilbert bundles, or more specifically, trivial Hilbert bundles with smooth weighted hermitian metrics (\cite{Tai1}), while the other fold concerns complex Brunn-Minkowski theory (see for instance \cite{Bo0}, \cite{Bo2}, \cite{Bo3}). The complex Brunn-Minkowski theory is a theory introduced by Berndtsson, and deals with the positivity of curvatures of certain vector bundles of infinite rank associated with Bergman (type of) spaces (\cite{Bo2}). Parts of the theory may be viewed as complex analogues of the real-variables (functional version of the) Brunn-Minkowski inequality from convex analysis (for convex functions) (\cite{Bo2}, \cite{Bo3}). An aspect (of the theory) that we find to be particularly beautiful and interesting, is that rather deep results may be obtained by using relatively simple convexity arguments. One of our favourite examples here is a ''complex Brunn-Minkowski proof'' of the Suita conjecture, and building on the same ideas, an extension theorem of Ohsawa-Takegoshi type with optimal estimates, first demonstrated by Berndtsson and Lempert in \cite{BoLempOT}. With the language and formalism that we shall develop, the optimal estimate in the extension theorem can be seen to follow from general monotonicity properties of (weighted) norms of holomorphic sections of certain trivial Hilbert quotient bundles. To explain this, we first take a dual point of view and consider sections of the dual bundle of the quotient bundle. By a quotient bundle curvature formula, which we shall detail (see \cite{Tai1}), if the initial Hilbert bundle is (semi)positive (in the sense of Griffiths), then so is the quotient bundle. By complex Brunn-Minkowski theory, this positivity now implies (or rather means) that all holomorphic sections of the dual bundle (of the quotient bundle) satisfy that their squared norm is logarithmically plurisubharmonic as a function of the base parameter. Choosing as base parameter the real part of the left half-plane in the complex plane this yields that the norms (squared) are logarithmically convex. In some special situations, one may verify that as a function of the base parameter, the sums of logarithms of the norms and some linear terms are bounded from above as the parameter tends to negative infinity. By the convexity from before, these sums must then be monotonically increasing functions of the base parameter. This latter monotonicity finally gives the optimal estimates, after turning back to the norms on the quotient bundle (times some exponential terms). Very interestingly (in our opinion), one finds that analogous arguments may be used to address also (the) strong openness (conjecture of Demailly) (\cite{TaiXuHilb}).

This thesis consists of a collection of 4 papers, each of which corresponds to its own chapter. In chapter 1, the first paper (\cite{Tai1}), we introduce what we call a ''Hilbert bundles approach to complex Brunn-Minkowski theory''. We develop and introduce here a formalism, or theory, of Hilbert bundles, more specifically, what we call trivial Hilbert bundles with smooth weighted hermitian metrics, and show how to use this formalism to give an abstract approach to complex Brunn-Minkowski theory in a general setting. Our approach is rooted in and inspired by a new approach to the variation of Bergman kernels in \cite{Bo0} (which might be said to be the starting point of complex Brunn-Minkowski theory) that does not rely on the classical use of the deep regularity properties of the $\dbar$-Neumann operator associated with smoothly bounded strictly pseudoconvex domains. Instead, we use a smoothness theorem of Berndtsson for the Bergman projection map. We generalize this smoothness theorem to the setting of trivial Hilbert bundles with certain special weights, and discuss a few applications to complex Brunn-Minkowski theory, as well as some novel results akin to the variations of Bergman kernels. At this point we are only able to obtain these results presuming certain positivity properties of the curvatures of certain subbundles. These presumed positivity assumptions are addressed in a more concrete setting in the sequel paper and chapter 2 of the thesis.

Chapter 2, the second paper (\cite{Tai2}), is a continuation of chapter 1 and the first paper. We address here the presumed positivity assumptions in (some of) the main results from the first chapter in more concrete settings. In particular, we generalize the positivity of direct image bundles from \cite{Bo1}, in the special case of trivial fibrations, to the case of $(n,q)$-forms with $q\geq 0$ (and not just $q=0$), and a more general class of Kähler manifolds. This class is the class of so-called quasi-complete Kähler manifolds, a notion which generalizes that of complete Kähler manifolds, and which is due to Xu Wang and Bo Yong Chen (\cite{XuDiag}). To give a fairly general statement, we introduce a type of Hörmander estimate condition. The main idea here, similar to Berndtsson's approach in \cite{Bo1}, is to use $L^2$-methods to control the second fundamental form term in a formula for the curvature of the subbundle (\cite{Tai1}), viewing it as the norm of the $L^2$ (type of) minimal solution to a certain linear equation. The actual computations and end-results, however, differ from those that appear in \cite{Bo1}. We also give a novel plurisubharmonicity result for holomorphic sections of certain trivial Hilbert quotient bundles in the same setting, which we prove using similar ideas. Very briefly stated, and using the formalism and notation introduced in the first paper (\cite{Tai1}) and chapter 1, this result says that under certain conditions on the right-hand side of an equation of the form $\dbar(\cdot)=v$ (that is, on $v$) that relates to an appropriate smooth weighted hermitian metric $h$, $-f(\norm{u_{m}}^2_{h})$ admits plurisubharmonicity properties for all strictly increasing concave and twice differentiable functions $f$ on $[0,\infty)$, where $u_m$ denotes the $h$-minimal solution to the equation. It might be that this latter result is related to variations of so-called generalized Green functions, but the details are far from clear to us. We hope maybe to return to say more about this in a future paper. Finally, by considering special kinds of hermitian metrics on our Hilbert bundles, we prove versions of the above results where the metric on the line bundle may be possibly singular. The possibly non-smoothness of the line bundle metric is a main novelty, and we deal with it using a Hörmander type of theorem for singular weights, on quasi-complete Kähler manifolds. To obtain  the latter, we use regularization of quasi-plurisubharmonic currents on compact manifolds, due to Demailly. The special form of the hermitian metrics that we consider plays a rather essential role here. We conclude the paper with a few comments on the theme of the sequel paper and chapter 3 of the thesis. 

The first three chapters of the thesis are parts of a series of papers with the common title \ita{A Hilbert bundles approach to complex Brunn-Minkowski theory}. Chapter 3, and the third paper, is the third installment in the series (and, of course, chapters 1 and 2 are respectively the first and second installments), and focuses on specific applications to two important topics in complex differential geometry and analysis of several variables: $L^2$-holomorphic extension and (strong) openness. Building on the ideas developed throughout the former two parts of the series, we give a complex Brunn-Minkowski proof of an Ohsawa-Takegoshi extension theorem with sharp effective $L^2$-estimates, by generalizing the Berndtsson-Lempert method (\cite{BoLempOT}) to the global setting of complex manifolds, and the strong openness conjecture of Demailly (or the strong openess theorem of Guan-Zhou (\cite{GZ0})). The crux of the argument in either case is a general monotonicity property, more specifically, a certain decresingness property, of (squared norms of) holomorphic sections of certain trivial Hilbert quotient bundles. This paper, and the third part of the series, is joint work with Xu Wang. 

The final chapter of the thesis, and fourth paper, is also joint work with Xu Wang. We address a question by Ohsawa in a remark in his 2017-paper \cite{Ohsawa17} on whether a complex Brunn-Minkowski proof, or rather a proof using the Berndtsson-Lempert method, of (either of) two Ohsawa-Takegoshi type of extension theorems (Theorems 0.1 and 4.1  in \cite{Ohsawa17}) can be had. Using so-called weak geodesics from pluripotential theory to construct an appropriate weight function to be used in the Berndtsson-Lempert method, we prove a general Ohsawa-Takegoshi type of extension theorem that answers Ohsawa's question in the affirmative (for both extension theorems). The paper has been published in Arkiv För Matematik. After its publication, the first (listed) author of the paper discovered that the main theorem of the paper may be used, in an interesting application of it, to give a new simple and short proof of a well-known first order Bergman kernel asymptotics. This proof is included in the version of the paper that appears in the thesis that follows.

\chapter{Paper 1}

\begin{center}
\normalsize{\bt{A HILBERT BUNDLES APPROACH TO COMPLEX BRUNN-MINKOWSKI THEORY, I}}\\[0.5cm]
\small{TAI TERJE HUU NGUYEN}

\begin{abstract}
We present an abstract and novel approach to complex Brunn-Minkowski theory (see \cite{Bo0}, \cite{Bo1}, \cite{Bo2}, \cite{Bo3} for some first references) in a general setting using a theory (or formalism) of Hilbert bundles. Our approach is based on and inspired by a new approach of Berndtsson to the variation of Bergman kernels in \cite{Bo0} that does not rely on the (classical use of the) deep Hamilton-Kohn regularity theory for the $\dbar$-Neumann operator associated with smoothly bounded strictly pseudoconvex domains (see the proof of Lemma 2.1 in \cite{Bo0}). 

This is the first installment in a series of papers, and we establish in here the language and formalism, as well as some rudimentary results, that we use in subsequent papers.
\end{abstract}

\end{center}

%
%
%
%
%
%
%
%
%


\section{Introduction}
In a milestone paper in 2006, entitled \ita{Subharmonicity properties of the Bergman kernel and some other functions associated with pseudoconvex domains} (\cite{Bo0}), Berndtsson proved a rather remarkable result on plurisubharmonicity properties of Bergman kernels depending on a parameter, generalizing an earlier result of Maitani and Yamaguchi (\cite{MaiYam}). This is in the literature often known as the \ita{variation of Bergman kernels}, and might be said to be the starting point of the \ita{complex Brunn-Minkowski theory} (see \cite{Bo0}, \cite{Bo1}, \cite{Bo2}, \cite{Bo3} for some first references). It is also the starting point of our story. The reader is directed to \cite{Bo0}, specifically to Theorem 1.1, for the most general statement(s). Here, we will focus on the following (special) situation (see, however, also sections 2 and 3 in \cite{Bo0}): 

Let $U\times \Omega $ be a domain in $\set{C}^{m}\times \set{C}^{n}$, with $\Omega$ bounded, and let $\phi$ be a strictly plurisubharmonic function on $U\times \Omega$ which is smooth up to the boundary. We use $t$ to denote a generic point in $U$. Let $L^2_t$, respectively $A^2_t$, denote the weighted $L^2$-space, respectively the weighted Bergman space, on $\Omega$ with weight function $e^{-\phi(t,\cdot)}$. Following \cite{Bo0}, we write $\phi^{t}$ for the restriction of $\phi$ to $\{t\}\times \Omega$. Thus, $L^2_t$ and $A^2_t$ are the (complex) Hilbert spaces

\begin{align}
L^2_{t}&=\left\{f\in L^2_{\text{loc}}(\Omega):\norm{f}^2_{t}:=\int_{\Omega}\abs{f}^2e^{-\phi^{t}}<\infty\right\}\quad\quad \text{ and }\\
A^2_t&=L^2_t\cap \fancy{O}(\Omega)
\end{align}
with respect to the weighted $L^2$-inner product $(\cdot,\cdot)_{t}$ given by 
\begin{align}
(f,g)_{t}&:=\int_{\Omega}f\bar{g}e^{-\phi^{t}}.
\end{align}Integration is with respect to the Lebesgue measure, and $\fancy{O}(\Omega)$ denotes the family of holomorphic functions on $\Omega$. Let $k^{t}$ denote the Bergman kernel of $A^2_t$, and $K^{t}$ the Bergman kernel of $A^2_{t}$ \ita{on the diagonal}. Thus, $k^{t}$ is the map
\begin{align}
k^{t}:\Omega\times \Omega\to \set{C}, (z,w)\mapsto k^{t}_{w}(z),
\end{align}where $k^{t}_{z}$ is the unique element in $A^2_t$ satisfying that
\begin{align}
f(z)&=\int_{\Omega}f\overline{k^{t}_{z}}e^{-\phi^{t}}
\end{align}for all $f\in A^2_t$, and $K^{t}$ is the map
\begin{align}
K^{t}:\Omega\to \set{C}, z\mapsto k^{t}_{z}(z).
\end{align}
The (version of the) variation of Bergman kernels (that will be important to us) asserts the following:

\begin{theorem}[Variation of Bergman kernels (special case), \cite{Bo0}]\label{intro: variation of Bergman kernels}
With the above notation and set-up, suppose that $\Omega$ is pseudoconvex. Then for each $z\in\Omega$, the map $U\ni t\mapsto \log K^{t}(z)$ is plurisubharmonic.
\end{theorem}
 Let us comment briefly on the proof of the theorem; we denote the map $t\mapsto \log K^{t}(z)$ in Theorem \ref{intro: variation of Bergman kernels} by $\log K(z)$, and the map $t\mapsto K^{t}(z)$ by $K(z)$. The proof may be split into three key steps:
\begin{enumerate}[(1)]
\item Show that $K(z)$ is smooth on $U$.
\item Compute $i\partial_t\dbar_t \log K(z)$; we use $t$-subscripts to emphasize that we are considering the $\partial$ and $\dbar$ operators on $U$ (as opposed to $\Omega$).
\item Verify that $i\partial_t\dbar_t \log K(z)\geq 0$.
\end{enumerate}
We shall in this paper mainly be concerned with steps (1) and (2), at least to begin with, and here merely mention that in step (3),  Hörmander’s $L^2$-estimates for (solutions of) the $\dbar$-equation (see \cite{Bo0}, \cite{Bo1}) are utilized. The pseudoconvexity of $\Omega$, and the (strict) plurisubharmonicity of $\phi$, are then assumptions of natural relevance and importance. Interestingly, however, and as we shall see, as far as steps (1) and (2) are concerned, the assumptions on $\Omega$ and $\phi$ may be (much) weaker. 

In \cite{Bo0}, the smoothness of the Bergman kernel $K(z)$ in step (1) is established by use of the deep Hamilton-Kohn regularity theory for the $\dbar$-Neumann operator associated with smoothly bounded strictly pseudoconvex domains, under the additional assumption that $\Omega$ be smoothly bounded and strictly pseudoconvex; see (the proof of) Lemma 2.1 in \cite{Bo0}. We will refer to this (way of showing that $K(z)$ is smooth) as \ita{the classical approach} to the variation of Bergman kernels. In addition to the classical approach, there is recently an alternative and novel approach owing to Berndtsson, more pertinent in the following paper, and that we now turn to. In this new approach, the role of the deep regularity properties of the $\dbar$-Neumann operator in the classical approach is instead played by a smoothness theorem that is taught to us by Berndtsson. Said theorem admits an integral part in the present paper, and we have entitled it \bt{Berndtsson’s regularity theorem}. For the statement of it, recall that the Bergman projection associated with $A^2_t$ is the canonical orthogonal projection map $L^2_t\to A^2_t$ with respect to the (complete) inner product $(\cdot,\cdot)_t$. Note that in our current situation, $L^2:=L^2(\Omega)$ and $L^2_{t}$ are equal as vector spaces for each $t\in U$, so this map is also the canonical orthogonal projection map $L^2\to A^2$ with respect to $(\cdot,\cdot)_t$; $A^2$ being the usual Bergman space $L^2\cap \fancy{O}(\Omega)$ on $\Omega$. The point here is of course that $L^2$ and $A^2$ do not depend on (the) (parameter) $t$. Let $P^{t}$ denote the Bergman projection associated with $A^2_t$, and define the map $P$ by
\begin{align}
P:U\times L^2\to A^2, (t,f)\mapsto P^{t}f.
\end{align}Using the above notation, Berndtsson's regularity theorem may now be stated as follows:

\begin{theorem}[Berndtsson's regularity theorem]\label{intro: Berndtsson's regularity theorem, special case}
Let $U\times \Omega$ be a domain in $\set{C}^{m}\times \set{C}^{n}$, and let $\phi$ be a real-valued function on $U\times \Omega$. Assume that $\Omega$ is bounded, that $\phi$ is smooth up to the boundary, and let for each $t\in U$, $A^2_t=A^2(\Omega, \phi^{t})$, $P^{t}$, and $P$ be defined as above. Then $P$ is smooth (in the (real) Fréchet sense; see \cite{Cartan}).
\end{theorem}
Note that we in Theorem \ref{intro: Berndtsson's regularity theorem, special case} are not assuming that $\Omega$ be pseudoconvex, nor even that $\phi$ be plurisubharmonic. The situation in there is therefore already notably different from that in the classical approach. Indeed,  for instance, since $\Omega$ is only assumed to be bounded, the deep regularity properties of the $\dbar$-Neumann operator also never come into play. That we are allowed weaker assumptions on $\Omega$ and $\phi$, suggests that some natural generalizations are plausible, and becomes the starting point for the theme of the present paper. We shall say a little bit more about this theme in a moment, but first we use Berndtsson's regularity theorem (Theorem \ref{intro: Berndtsson's regularity theorem, special case}) to carry out step (1) above. We do this by using Theorem \ref{intro: Berndtsson's regularity theorem, special case} to prove the following smoothness lemma for the Bergman kernel of $A^2_t$:

\begin{lemma}\label{intro: smoothness of the Bergman kernel}
Let $z,w\in \Omega$. Then under the assumptions in Theorem \ref{intro: Berndtsson's regularity theorem, special case}, we have that $k^{t}_z(w)$ is smooth in $t$.
\end{lemma}

\begin{proof}
We shall give two proofs of the lemma. The first mimics in parts the proof of Lemma 2.1 in \cite{Bo0}. Said proof uses the specific setting that we are currently in. The second proof works also in a more general setting, and we shall meet it again in section \ref{sec: special metric and Berndtsson's regularity theorem}. The second proof also gives a slightly more general statement (in that we need not fix $w$ to begin with), and we include it here to serve as a nice introduction to what to come.\\

\subsubsection*{First proof}
It follows by definition of $P^{t}$ that for all $g\in L^2_t$,
\begin{align}
P^{t}g(z)&=\int_{\Omega}g\overline{k^{t}_{z}} e^{-\phi^{t}}.\label{eq: definition of Bergman projection, in proof}
\end{align}
Let $w:=(w_1,\ldots, w_{n})$ and $\chi_{j}$ be a rotation-invariant smooth real-valued non-negative function on $\set{C}$ supported in a small disc $\set{D}_{\epsilon_{j}}(w_{j})$ centred around $w_{j}$; the notation here means that the disc has center $w_{j}$ and radius $\epsilon_j$. Assume also for simplicity, by way of normalization, that $\displaystyle \int_{\set{D}_{\epsilon_j}(w_{j})}\chi_{j}\equiv 1$. Let $\epsilon:=(\epsilon_1,\ldots, \epsilon_{n})$ and $\Delta^{n}_{\epsilon}(w):=\prod_{j=1}^{n}\set{D}_{\epsilon_j}(w_{j})$. That is, $\Delta^{n}_{\epsilon}(w)$ is the polydisc centred at $w$ with polyradius $\epsilon$. Define $\chi:=\prod_{j=1}^{n}\chi_{j}$ as a function on $\Omega$ supported in $\Delta^{n}_{\epsilon}(w)$, and consider $g:=e^{\phi^{t}}\chi$. By using polar coordinates, the rotation-invariance of each $\chi_{j}$, and the mean-value property for $k^{t}_{z}$, we then get from \eqref{eq: definition of Bergman projection, in proof} that
\begin{align}
P^{t}g(z)&=c\overline{k^{t}_{z}(w)},
\end{align}where $c=1$ by the normalization of each of the $\chi_{j}$'s. Hence, the assertion follows by Berndtsson' regularity theorem, Theorem \ref{intro: Berndtsson's regularity theorem, special case} (and by composition of smooth functions; see the upcoming remark below).

\subsubsection*{Second proof}
We show that $k^{t}_z$ is smooth in $t$ as a map $U\to A^2$. Since evaluation at $w$ is smooth (see \cite{Cartan}), this gives the assertion. Let $f\in A^2_t$, and denote by $\xi_{z}$ the evaluation at $z$ as a continuous linear functional on $A^2_t$. By the Riesz representation theorem there exists a unique (Riesz representative) $R^{t}(\xi_{z})\in A^2_t$ (of $\xi_z$ with respect to the $(\cdot,\cdot)_t$-metric) such that for all $f\in A^2_t$, 
\begin{align}
\xi_{z}(f)=f(z)=\displaystyle \int_{\Omega}f\overline{R^{t}(\xi_z)}e^{-\phi^{t}}.\label{eq: Riesz representation, in proof, introduction}
\end{align} It follows that $R^{t}(\xi_z)=k^{t}_{z}$, so it suffices to show that $R^{t}(\xi_z)$ is smooth in $t$. Let us denote the usual metric on $L^2$ by $(\cdot,\cdot)_{L^2}$, and take $R^{L^2}(\xi_z)$ to be the Riesz representative of $\xi_z$ with respect to this metric. Then for all $f\in A^2\simeq A^2_t$ (here $\simeq$ means vector space isomorphism), we have that
\begin{align}
\xi_{z}(f)=f(z)=\int_{\Omega}f\overline{R^{L^2}(\xi_z)}=\int_{\Omega}f\overline{e^{\phi^{t}}R^{L^2}(\xi_z)}e^{-\phi^{t}}=\int_{\Omega}f\overline{P^{t}(e^{\phi^{t}}R^{L^2}(\xi_z))}e^{-\phi^{t}},
\end{align}where the last equality follows by orthogonality.
Thus by \eqref{eq: Riesz representation, in proof, introduction}, $R^{t}(\xi_z)=P^{t}(e^{\phi^{t}}R^{L^2}(\xi_z))$, so the assertion follows (as in the first proof) by Berndtsson's regularity theorem, Theorem \ref{intro: Berndtsson's regularity theorem, special case} (and by composition of smooth functions; see the upcoming remark).
\end{proof}

\begin{remark}
In the two previous proofs we have used the fact that if $u\in L^2$, then $t\mapsto e^{\pm \phi^{t}}u$ is smooth on $U$. This is also used in the proof of Theorem \ref{intro: Berndtsson's regularity theorem, special case} below, for the smoothness of $F$ (see the proof below). This fact is not obvious (at least not to us), but certainly believable. It can be shown by use of the mean-value theorem and the uniform continuity of continuous functions on compact sets; recall that in our setting, $\overline{\Omega}$ is compact (and $\phi$ is smooth up to the boundary).
\end{remark}

It clearly follows from Lemma \ref{intro: smoothness of the Bergman kernel} that $K(z)$ is smooth. We have therefore carried out step (1) above using Berndtsson's regularity theorem (Theorem \ref{intro: Berndtsson's regularity theorem, special case}), and the latter gives a new approach to the variation of Bergman kernels different from the classical approach. Let us also briefly discuss the proof of Berndtsson's regularity (Theorem \ref{intro: Berndtsson's regularity theorem, special case}); we shall give a more detailed proof of a generalized version of it later (see Theorem \ref{intro: Berndtsson's regularity theorem, general case} and also Theorem \ref{thm: Berndtsson's regularity theorem, general version}). The main tool used in the proof is the implicit mapping theorem (see \cite{Cartan}).
 \begin{proof}[Proof of Theorem \ref{intro: Berndtsson's regularity theorem, special case}] Fix $(t_0,f_0)\in U\times L^2$, and assume for simplicity that  $\phi^{t_0}=0$. It suffices to show that $P$ is smooth near $(t_0, f_0)$. Consider the smooth map (see the previous remark)
\begin{align}
F:U\times L^2\times A^2\to A^2, (t,f,g)\mapsto P_{A^2}(e^{-\phi^{t}}(f-g)),
\end{align}where $P_{A^2}$ denotes the Bergman projection associated with $A^2$. It can be shown that we have $F(t,f,g)=0$ if and only if $g=P^{t}f$, and that the partial derivative of $F$ in the $g$-direction (the $A^2$-direction) at the point $(t_0,f_0, g_0)$, for any $g_0$, is invertible (it is equal to $-\id{A^2}$, where we let $\id{X}$ denote the identity map on a set $X$). Hence, by the implicit mapping theorem, there exist an open neighbourhood $O$ of $(t_0,f_0)$ and a unique smooth $A^2$-valued function $G$ on $O$, such that $F(t,f,g)=0$ if and only if $g=G(t,f)$, for all $(t,f)\in O$. By comparison we have $P=G$ on $O$, so $P$ is smooth(!).
\end{proof}
We now observe that only elementary tools are used in the proof (of Theorem \ref{intro: Berndtsson's regularity theorem, special case}), the implicit mapping theorem being the main one. Looking at the classical approach, in which the deep regularity properties of the $\dbar$-Neumann operator are rather important, we find this to be quite remarkable. This aspect of (the proof of) the theorem has moreover the noteworthy effects of making the new approach to the variation of Bergman kernels (much) less rigid, and (in our opinion) somewhat ''more accessible'', than the classical approach. For example, we see that as far as steps (1) and (2) above (in the 3-steps proof of Theorem \ref{intro: variation of Bergman kernels}) are concerned (leaving aside step (3) for the time being), we may completely dispose of the pseudoconvexity and plurisubharmonicity assumptions on $\Omega$ and $\phi$. Thus, if we are only interested in the smoothness of the Bergman kernel $K(z)$ (and, not yet, its (logarithmic) plurisubharmonicity), we may consider a more general setting in which the assumptions in Theorem 1.1 are (much) weaker. We also see that the deep regularity properties of the $\dbar$-Neumann operator have actually no play in the smoothness of $K(z)$, something which is not at all visible from the classical approach. These observations suggest some further natural generalizations, and we will in this paper explore such a generalization where we replace $L^2$ and $A^2$ above with more general Hilbert spaces. This naturally leads to the abstract and (more) general setting of (what we call)  \ita{(trivial) Hilbert bundles}, and we introduce and develop a theory (or formalism) for these. We then use this theory to generalize the new approach to the variation of Bergman kernels discussed above to the abstract setting. In particular, we shall give a generalization of Berndtsson's regularity theorem (Theorem \ref{intro: Berndtsson's regularity theorem, special case}) to the setting of (trivial) Hilbert bundles. Using the latter as a key, we discuss an abstract and novel approach to complex Brunn-Minkowski theory in this (more) general setting; what we call \ita{a Hilbert bundles approach to complex Brunn-Minkowski theory}, and which is precisely the theme of the present paper. To obtain the aforementioned generalization of Berndtsson's regularity theorem (Theorem \ref{intro: Berndtsson's regularity theorem, special case}), whom we also entitle \bt{Berndtsson's regularity theorem}, and which is one of the main results in the paper (aside from the theory of Hilbert bundles and the abstract approach itself), we introduce two special kinds of (what we call) \ita{hermitian metrics on (trivial) Hilbert bundles}. This allows us to adapt the original proof of Berndtsson's regularity theorem (Theorem \ref{intro: Berndtsson's regularity theorem, special case}) to the (more) general setting and obtain a more general version using an analogous argument as before. We also show how to use this (more general) Berndtsson's regularity theorem to prove a generalization, or a natural analogue, of Lemma \ref{intro: smoothness of the Bergman kernel} in the (more) general setting. Finally, we discuss steps (2) and (3) (in the proof of Theorem \ref{intro: variation of Bergman kernels}) from above in the (more) general context of (trivial) Hilbert bundles, as well as some novel results, akin to Theorem \ref{intro: variation of Bergman kernels}, or a variant of it.

With the terminology and formalism that we shall be developing throughout the rest of the paper, the aforementioned generalization of Berndtsson's regularity theorem (Theorem \ref{intro: Berndtsson's regularity theorem, special case}) may be stated as follows:

\begin{theorem}[Berndtsson's regularity theorem, Hilbert bundles version]\label{intro: Berndtsson's regularity theorem, general case}
Let $\fancy{H}$ be a trivial Hilbert bundle over an $m$-dimensional complex manifold, and let $h$ be a smooth weighted hermitian metric on $\fancy{H}$ induced by a hermitian metric $h^0$ of zero variation. Let $\fancy{H}_0$ be a subbundle of $\fancy{H}$, and assume that $\fancy{H}_0$ is closed in $\fancy{H}$ with respect to $h$ and $h_0$. Then the total orthogonal projection map $P=P(h)$ associated with  $h$ (of $\fancy{H}_0$ in $\fancy{H}$) is smooth.
\end{theorem}
This concludes our introduction. In the remainder of the paper, and starting with the next section, we detail our aforementioned Hilbert bundles approach to complex Brunn-Minkowski theory.

\subsection{Acknowledgements}
It is a pleasure to thank Bo Berndtsson for being a great source of inspiration and for teaching us the implicit mapping theorem argument, and Xu Wang for invaluable discussions.

\section{Hilbert bundles and hermitian metrics on Hilbert bundles}
\subsection{Hilbert bundles}
In this section we begin our Hilbert bundles approach to complex Brunn-Minkowski theory, and we commence by discussing a few generalities on (what we call) \ita{Hilbert bundles} and \ita{hermitian metrics} (on Hilbert bundles). An immediate natural generalization of what we discussed in the introduction is to replace the families $\{L^2_t\}_{t}$ and $\{A^2_t\}_{t}$ with more general families of Hilbert spaces, say $\{\fancy{H}_t\}_{t}$ and $\{(\fancy{H}_0)_{t}\}_{t}$. We may also let $t$ vary over a more general index set $\fancy{B}$. We shall call a family of Hilbert spaces such as $\{\fancy{H}_t\}_{t}$, which depends on a parameter ($t$), a \ita{(general) Hilbert bundle}. While this precisely defines what we shall mean by a (general) Hilbert bundle, to make our definition resemble somewhat more that of a vector bundle, we will employ the following formulation:

\begin{definition}\label{def: Hilbert bundles}
A \bt{(general) Hilbert bundle} is an ordered triple $(\fancy{H},\pi,\fancy{B})$ consisting of two (non-empty) sets $\fancy{H}$ and $\fancy{B}$, and a surjective map $\pi:\fancy{H}\to \fancy{B}$, with the property that $\pi^{-1}(b)$ is a (complex) Hilbert space for each $b\in \fancy{B}$.
\end{definition}
In particular, all finitely-ranked vector bundles, and all infinitely-ranked vector bundles (see  \cite{Lang}, \cite{Abraham}) whose fibers are Hilbert spaces, classify as Hilbert bundles according to our definition. Given a Hilbert bundle $(\fancy{H},\pi, \fancy{B})$,  we will refer to $\fancy{H}, \fancy{B}$, and $\pi$ as respectively \bt{the total space, the base space}, and \bt{the projection map}, of the (Hilbert) bundle. Also, given a Hilbert bundle $(\fancy{H},\pi, \fancy{B})$, and $b\in \fancy{B}$, we will refer to the Hilbert space $\pi^{-1}(b)$ as \bt{the fiber of $\fancy{H}$ over $b$}, and generically denote it by $\fancy{H}_b$. Finally, we will typically refer to a Hilbert bundle by referring to its total space, and to a Hilbert bundle whose base space is $\fancy{B}$ as a Hilbert bundle \bt{over} $\fancy{B}$.

Of particular interest to us will be what we call \ita{trivial} Hilbert bundles. One way to define these is to first define what one might call \ita{Hilbert bundle isomorphisms}, meant to be ''the Hilbert bundles analogues'' of vector bundle isomorphisms, and then mimic the definition of trivial vector bundles (see \cite{Loring}). For simplicity, however, we have settled for the following definition that shall suffice for all our purposes:

\begin{definition}\label{def: trivial Hilbert bundles, simplified}
Let $H$ be a (complex) Hilbert space (or more generally, a Hilbert manifold; see \cite{Lang}, \cite{Abraham}) and $\fancy{B}$ a (non-empty) set. The \bt{trivial Hilbert bundle over $\fancy{B}$ with typical fiber $H$} is (the Hilbert bundle over $\fancy{B}$ given by) the triple $(\fancy{H},\pi,\fancy{B})$, where $\fancy{H}=\fancy{B}\times H$ and $\pi:\fancy{H}\to \fancy{B}$ is the natural projection map $(b,u)\mapsto b$. 
\end{definition}

Let us reserve the symbol $B$ throughout the paper to denote a fixed $m$-dimensional complex manifold; this to avoid having to repeatedly write the phrase ''over an $m$-dimensional complex manifold''. We shall mostly focus on trivial Hilbert bundles \ita{over $B$} with Hilbert spaces (contra more general Hilbert manifolds) as their typical fibers. If $\fancy{H}$ is a trivial Hilbert bundle over $B$ with typical fiber $H$, we identify for simplicity $H$ with $\fancy{H}$, and also say that $H$ is a trivial Hilbert bundle over $B$. 

\begin{example}\label{ex: first ex of trivial Hilbert bundles}
Let $L^2$, $A^2$, and $U$ be defined as in the introduction (that is, $U\times \Omega$ is a domain in $\set{C}^{m}\times \set{C}^{n}$ with $U$ open and $\Omega$ bounded, and $L^2$ and $A^2$ are respectively the $L^2$-space and Bergman space on $\Omega$). Then $L^2$ and $A^2$ are trivial Hilbert bundles over $U$. Since we are dealing with domains in $\set{C}^{n}$ and $\set{C}^{m}$, we like to think of this example as being in \ita{the local setting}. There is a natural \ita{global analogue} of it in the setting of complex manifold as follows: Let $X$ be an $n$-dimensional complex manifold, and let $L$ be a hermitian holomorphic line bundle over $X$ with metric $e^{-\phi}$. For two smooth $L$-valued sections $u$ and $v$ of the canonical bundle of $X$, here denoted $K_X$, define $h^0$ by
\begin{align}
h^{0}(u,v)&:=\int_{X}i^{n^2}u\wedge \overline{v}e^{-\phi},
\end{align} and let $H$ be the completion of the space of smooth $L$-valued sections $u$ of $K_{X}$ such that $h^{0}(u,u):=\norm{u}^2_{h^{0}}<\infty$, with respect to $h^0$. Let $H_0$ denote the subspace of $H$ consisting of holomorphic sections, and equip both $H$ and $H_0$ with $h^0$ as a metric (or Hilbert inner product). Then $H$ and $H_0$ are trivial Hilbert bundles over $B$. 
\end{example}
Let $\fancy{H}$ be a (general) Hilbert bundle over $\fancy{B}$, and let $\fancy{U}\sub \fancy{B}$ be non-empty. Adopting terminology from the theory of vector bundles, we define \bt{a (local) section of $\fancy{H}$ over $\fancy{U}$} to be a map $u:\fancy{U}\to \fancy{H}$ such that 
\begin{align}
u(b):=u^{b}\in \fancy{H}_b
\end{align}for all $b\in \fancy{U}$. In the case that $\fancy{H}$ is trivial over $B$ with typical fiber $H$, we may, and will, identify (local) sections of $\fancy{H}$ with $H$-valued maps defined on subsets of $B$, in the natural way. When the occasion calls for it, we use $\Gamma(\fancy{U};\fancy{H})$ to denote the family of all local sections of $\fancy{H}$ over $\fancy{U}$; in the global case that $\fancy{U}=\fancy{B}$, we simply write $\Gamma(\fancy{H})$. 

We shall also need the notions of \ita{(Hilbert) subbundles}, \ita{(Hilbert) dual bundles}, and \ita{(Hilbert) quotient bundles} in our theory, which are analogous to the corresponding notions from the theory of vector bundles. For these, we use the following definition:

\begin{definition}\label{def: subbundles, dual bundles, and quotient bundles}
Let $\fancy{H}$ be a (general) Hilbert bundle over $\fancy{B}$ (some non-empty set, but we shall stop pointing this out). Then:
\begin{enumerate}[(i)]
\item A \bt{(Hilbert) subbundle of $\fancy{H}$} is a Hilbert bundle $\fancy{H}_0\sub \fancy{H}$ over $\fancy{B}$ such that $(\fancy{H}_0)_{b}$ is a closed subspace of $\fancy{H}_b$ for each $b\in \fancy{B}$.
\item The \bt{(Hilbert) dual bundle of $\fancy{H}$}, here generically denoted $\fancy{H}^*$, is the Hilbert bundle over $\fancy{B}$ such that $(\fancy{H}^*)_{b}=(\fancy{H}_b)^*$, the dual space of $\fancy{H}_b$, for each $b\in \fancy{B}$.
\item If $\fancy{H}_0$ is a (Hilbert) subbundle of $\fancy{H}$, we refer to $\fancy{H}/\fancy{H}_0$ as a \bt{(Hilbert) quotient bundle} and define it to be the Hilbert bundle over $\fancy{B}$ such that $(\fancy{H}/\fancy{H}_0)_{b}=\fancy{H}_b/(\fancy{H}_0)_b$ for each $b\in \fancy{B}$.
\end{enumerate}
\end{definition}
Let $\fancy{H}$ be a Hilbert bundle over $\fancy{B}$. We will use the notation $\fancy{H}_0\leq \fancy{H}$ to denote that $\fancy{H}_0$ is a subbundle of $\fancy{H}$, and in the case that $\fancy{H}_0\leq \fancy{H}$ and $\fancy{H}_0$ is fixed, or understood from the context, we will frequently denote the quotient bundle $\fancy{H}/\fancy{H}_0$ by $\fancy{Q}$ (for ''quotient''). In the case of trivial Hilbert bundles, we shall typically use the symbols $H,H_0$, and $Q$ respectively in place of $\fancy{H},\fancy{H}_0$ and $\fancy{Q}$.

\begin{example}\label{ex: second ex of trivial Hulbert bundles}
Let $L^2, A^2, H,$ and $H_0$ be as in Example \ref{ex: first ex of trivial Hilbert bundles}. Then $A^2\leq L^2$ and $H_0\leq H$, as trivial Hilbert bundles. 
\end{example}
\subsection{Hermitian metrics on Hilbert bundles}
Next, we equip our Hilbert bundles with natural hermitian structures by introducing on these what we call \ita{hermitian metrics}. In principle, the definition of a hermitian metric on a Hilbert bundle is straight-forward. There is, however, a slight subtlety that appears when we also take a ''dual point of view''. To explain this, suppose that $h=\{h_{b}\}_{b\in \fancy{B}}$ is collection of complete inner products on the fibers of $\fancy{H}$, indexed of course such that $h_b$ is a complete inner product on $\fancy{H}_b$. Consider now the dual bundle of $\fancy{H}$, $\fancy{H}^*$, and its fiber over $b$, $(\fancy{H}^*)_b=(\fancy{H}_b)^*:=\fancy{H}_b^*$. The fiber $\fancy{H}_b^*$ is by definition the dual space of $\fancy{H}_b$, and thus consists of complex-linear continuous functionals on $\fancy{H}_b$. The point to stress here is that the continuity in this context is with respect to the fixed inner product on $\fancy{H}_b$ with which it is equipped from the outset (being a Hilbert space). Let us denote this inner product by $g_b$, as it in general may differ from $h_b$. In the situation at hand, we may also consider the space of complex-linear functionals on $\fancy{H}_b$ which are continuous with respect to $h_b$, and an inconvenience arises if this latter space and $\fancy{H}_b^*$ do not coincide as vector spaces. We shall say that $h$ satisfies \bt{the duality condition (of continuity)} if the two former vector spaces \ita{do coincide}. Thus, $h$ satisfies the duality condition of continuity if it is the case that any complex-linear functional on $\fancy{H}_b$ is continuous with respect to $g_b$ if and only if it is continuous with respect to $h_b$. In particular, if $g_b$ and $h_b$ are equivalent, this is of course the case. To avoid the aforementioned inconvenience, we incorporate the duality condition directly as part of our definition of hermitian metrics on Hilbert bundles. The latter is as follows:

\begin{definition}\label{def: hermitian metric on Hilbert bundles}
Let $\fancy{H}$ be a (general) Hilbert bundle over $\fancy{B}$. A \bt{hermitian metric on $\fancy{H}$} is an assignment $h:=\{h_b\}_{b\in \fancy{B}}$ to each $b\in \fancy{B}$ a complete inner product $h_b$ on $\fancy{H}_b$ which satisfies the duality condition (of continuity).
\end{definition}
If $\fancy{H}$ is a Hilbert bundle over $\fancy{B}$, and $h$ is a hermitian metric on $\fancy{H}$, then we define for each pair of (local) sections $u$ and $v$ of $\fancy{H}$, $h(u,v)$ to be the (local) function on $\fancy{B}$ with values in $\set{C}$ given by
\begin{align}
h(u,v)(b)&:=h_{b}(u^{b},v^{b})
\end{align}(for all $b$ in the intersection of the domain of $u$ and $v$). We also put 
\begin{align}
\norm{u}_{h}^2:=(\norm{u}_{h})^2:=h(u,u),
\end{align} and refer to $\norm{u}_{h}$ as \bt{the $h$-norm of $u$}, or the \bt{norm of $u$ with respect to $h$}. The hermitian metric $h$ is completely determined by specifying $h(u,v)$ for all pairs of local sections $u$ and $v$ of $\fancy{H}$.

\subsection{Adjoints, minimal lifts, and total orthogonal projections}

Suppose that $h$ is a hermitian metric on $\fancy{H}$, and let $\fancy{H}_0\leq \fancy{H}$. Then $h$ obviously induces a natural hermitian metric on $\fancy{H}_0$, simply by restriction. Similarly, $h$  induces a natural hermitian metric on the dual bundle $\fancy{H}^*$, and also on the quotient bundle $\fancy{Q}$ under the additional natural assumption that $\fancy{H}_0$ be what we call \ita{closed in $\fancy{H}$ with respect to $h$}. By the latter is meant the following: We shall say that $\fancy{H}_0$ is \bt{closed in $\fancy{H}$ with respect to $h$}, or \bt{$h$-closed in $\fancy{H}$}, and denote this by writing
\begin{align}
\overline{\fancy{H}_0}^{(\fancy{H},h)}=\fancy{H}_0,
\end{align}if $(\fancy{H}_0)_b$ is closed in $\fancy{H}_b$ with respect to $h_b$ for  each $b\in \fancy{B}$. 

To discuss the naturally induced hermitian metrics on $\fancy{H}^*$ and $\fancy{Q}$, we introduce next what we call \ita{adjoint operators} and \ita{minimal lifting operators}. Let us begin with \ita{adjoint operators}. Let $u^*$ be a (local) section of $\fancy{H}^*$, $u$ be a (local) section of $\fancy{H}$, and the intersection of the domains of $u^*$ and $u$ be denoted by $\fancy{U}$. We define 
\begin{align}
\ip{u^*}{u}=u^*(u)
\end{align} to be the complex-valued function on $\fancy{U}$ given by
\begin{align}
\ip{u^*}{u}(b):=(u^{*}(b))u^{b}:=\ip{u^*(b)}{u^{b}}.
\end{align}
Let $D_{u^*}$ denote the domain of definition of $u^*$. By the Riesz representation theorem (and the duality condition satisfied by $h$(!)), there exists a unique local section $\sharp u^*$ of $\fancy{H}$ over $D_{u^*}$ such that 
\begin{align}
\ip{u^*}{u}&=h(u,\sharp u^*)
\end{align}for all local sections $u$ of $\fancy{H}$. We give special names to the (local) section $\sharp u^*$ of $\fancy{H}$, and the operator that sends $u^*$ to $\sharp u^*$:

\begin{definition}\label{def: adjoint operators}
With the above notation and set-up, we define $\sharp u^*$ (or $\sharp^{h}u^*$ if we want to emphasize $h$) to be \bt{the adjoint of $u^*$ (with respect to $h$)}. We also define the \bt{adjoint operator (for $\fancy{H}^*$) with respect to $h$} to be the complex-antilinear operator that sends any local section $u^*$ of $\fancy{H}^*$ to its adjoint (with respect to $h$). We will generically denote the latter by $\sharp$ (or $\sharp^{h}$ if we want to emphasize $h$).
\end{definition}
Note that the adjoint operator $\sharp$ in Definition \ref{def: adjoint operators} has an obvious inverse operator. We will generically denote this inverse operator by $\sharp^{-1}$ (or $(\sharp^{h})^{-1}$ if $h$ needs to emphasized), and its value at a ''point'' $u$ (a local section of $\fancy{H}$) by $\sharp^{-1}u$ (or $(\sharp^{h})^{-1}u$ if $h$ needs to emphasized). The adjoint operator $\sharp$ may be used to define a natural hermitian metric on $\fancy{H}^*$ induced from $h$. We will denote this natural hermitian metric on $\fancy{H}^*$ induced by $h$, by $h^{-1}$, and refer to it as \bt{the dual (hermitian) (metric) of $h$}. We define the dual metric $h^{-1}$ by letting
\begin{align}
h^{-1}(u^*,v^*):=h(\sharp v^*, \sharp u^*)
\end{align}for all local sections $u^*$ and $v^*$ of $\fancy{H}^*$.

We next consider the induced metric on the quotient bundle. For this, we shall always assume that \begin{align}\overline{\fancy{H}_0}^{(\fancy{H},h)}=\fancy{H}_0.
\end{align} It then follows that for each $b\in \fancy{B}$, we have an orthogonal direct sum decomposition
\begin{align}
\fancy{H}_b=(\fancy{H}_0)_b\oplus (\fancy{H}_0)_{b}^{\perp},
\end{align}with respect to $h_{b}$, where
\begin{align}
(\fancy{H}_0)_b^{\perp}&=\cbrac{a_b\in \fancy{H}_b:h_{b}(a_{0,b},a_{b})\equiv 0\;\forall\;a_{0,b}\in (\fancy{H}_0)_b}
\end{align}is (what we may call) the $h_b$-orthogonal complement of $(\fancy{H}_0)_b$ in $\fancy{H}_b$. We shall usually write all these decompositions, as $b$ varies, as a single decomposition
\begin{align}
\fancy{H}&=\fancy{H}_0\oplus \fancy{H}_0^{\perp}.
\end{align}
Let $u$ be a (local) section of $\fancy{H}$ over $\fancy{U}$. The above decomposition induces a decomposition of $u$ of the form 
\begin{align}
u=u_0\oplus u_{0}^{\perp},
\end{align} where $u_0$ is a (local) section of $\fancy{H}_0$ over $\fancy{U}$, and where 
\begin{align}
h(v_0,u_{0}^{\perp})\equiv 0
\end{align} for all local sections $v_0$ of $\fancy{H}_0$. In general, if $v$ is any local section of $\fancy{H}$ satisfying 
\begin{align}
h(a_0,v)\equiv 0
\end{align} for all local sections $a_0$ of $\fancy{H}_0$, we shall say that $v$ is a local section of $\fancy{H}^{\perp}_0$, and also write $v\perp \fancy{H}_0$; if we need to emphasize $h$ we may write $\perp_{h}$ for $\perp$. Let $[u]$ be a (local) section of $\fancy{Q}$ over $\fancy{U}$. By a \bt{representative of $[u]$}, we shall mean a (local) section $u$ of $\fancy{H}$ over $\fancy{U}$ such that for all $b\in \fancy{U}$, \begin{align}
[u^{b}]=[u]^{b},
\end{align} where $[a_{b}]$ denotes the equivalence class of $a_b\in \fancy{H}_b$ in $\fancy{Q}_{b}=(\fancy{H}_b)/(\fancy{H}_0)_b$. It follows that $[u]$ admits a canonical representative, namely $u_0^{\perp}$ where $u=u_0\oplus u_0^{\perp}$ is \ita{any} representative of $[u]$. This suggests the following ''quotient analogue'' to Definition \ref{def: adjoint operators}:

\begin{definition}\label{def: minimal lifts}
With the above notation and set-up, we define \bt{the minimal lift(ing) of $[u]$ (with respect to $h$)} to be the canonical representative $u_{0}^{\perp}$ of $[u]$, where $u=u_{0}\oplus u^{\perp}_0$ is \ita{any} representative of $[u]$, and denote it by $\ml{[u]}$ or $m_{l}[u]$. We also define \bt{the minimal lift(ing) operator for $\fancy{Q}$ (with respect to $h$)} to be the complex-linear operator that sends any local section of $\fancy{Q}$ to its minimal lift(ing) (with respect to $h$), and denote it by $m_{l}$. (If we want to emphasize $h$, we may write $\perp^{h}$ instead of just $\perp$, and $m_{l}^{h}$ instead of just $m_{l}$.)
\end{definition}Similar to the case for the adjoint operator, there is an inverse operator of $m_{l}$ on its image. We will denote this inverse operator by $m_{l}^{-1}$, and its value at a ''point'' $u$ by $u_{\downarrow}$ or $m_{l}^{-1}u$. (As before, we may replace $\downarrow$ with $\downarrow^{h}$, and $m_{l}^{-1}$ with $(m_{l}^{h})^{-1}$ if we want to emphasize $h$.) Using the minimal lifting operator in Definition \ref{def: minimal lifts}, we may now define a natural hermitian metric on $\fancy{Q}$ induced from $h$. We generically denote this metric by $h^{\fancy{Q}}$, and we define it by letting
\begin{align}
h^{\fancy{Q}}([u],[v])&:=h(\ml{[u]},\ml{[v]})\label{eq: definition of quotient metric}
\end{align}for all (local) sections $[u]$ and $[v]$ of $\fancy{Q}$. We shall refer to the metric $h^{\fancy{Q}}$ as \bt{the quotient (hermitian) metric (on $\fancy{Q}$) induced by $h$}.

We conclude this section by defining the ''Hilbert bundles analogue'' of the map $P$ in Berndtsson's regularity theorem (Theorem \ref{intro: Berndtsson's regularity theorem, special case}) from the introduction, and with a preliminary result that we shall later show relates to Lemma \ref{intro: smoothness of the Bergman kernel}. We begin with the definition, which is the following:

\begin{definition}\label{def: total orthogonal projection map}
Let $\fancy{H}$ be a (general) Hilbert bundle over $\fancy{B}$, and suppose that $h$ is a hermitian metric on $\fancy{H}$. Assume that $\fancy{H}_0\leq \fancy{H}$ with $\overline{(\fancy{H}_0)}^{(\fancy{H},h)}=\fancy{H}_0$. Then we define \bt{the total orthogonal projection map $P=P(h)$ associated with $h$ (of $\fancy{H}_0$ in $\fancy{H}$) } to be the complex-linear operator that sends any local section $u$ of $\fancy{H}$ to the local section $Pu$ of $\fancy{H}_0$ defined by
\begin{align}
(Pu)(b)&:=P^{b}u^{b},
\end{align}where $P^{b}$ denotes the canonical orthogonal projection $\fancy{H}_b\to (\fancy{H}_0)_b$ with respect to $h_b$. 
\end{definition}

\begin{example}
Let $L^2_t, A^2_t$, and $(\cdot,\cdot)_t$ be defined as in the introduction. Then $\fancy{H}:=\{L^2_t\}_{t\in U}$ and $\fancy{H}_0:=\{A^2_t\}_{t\in U}$ are (trivial) Hilbert bundles over $U$ with $\fancy{H}_0\leq \fancy{H}$, and we may view the collection $h:=\{(\cdot,\cdot)_t\}_{t\in U}$ as a (smooth) hermitian metric on these such that $\overline{\fancy{H}_0}^{(\fancy{H},h)}=\fancy{H}_0$. The total orthogonal projection map associated with $h$ of $\fancy{H}_0$ in $\fancy{H}$ is precisely the map $P$ that appears in Berndtsson's regularity theorem (Theorem \ref{intro: Berndtsson's regularity theorem, special case}).
\end{example}

We next give the preliminary result relating to Lemma \ref{intro: smoothness of the Bergman kernel}; we shall say more about the relation later (see Corollary \ref{cor: generaization of smoothness of Bergman kernels}). To prepare for the result, let $\fancy{H}, \fancy{H}_0$, $h$, and $P$ be as in Definition \ref{def: total orthogonal projection map}. In this situation, there are (at least) two dual bundles that we may consider, namely $\fancy{H}^*$ and $(\fancy{H}_0)^*:=\fancy{H}_0^*$. Let $\sharp$ denote the adjoint operator for $\fancy{H}^*$, and $\sharp_0$ the adjoint operator for $\fancy{H}_0^*$, both with respect to $h$ (to be more precise, the first adjoint operator is with respect to $h$, and the second is with respect to the restriction of $h$ as a hermitian metric on $\fancy{H}_0$).  We may ask for the relationship between the two adjoint operators, if any, and have then the following quite useful result:

\begin{proposition}\label{prop: relation between adjoints of sub and exterior bundles}
With the above notation and set-up, we have 
\begin{align}
\sharp_{0}=P\sharp.
\end{align}
\end{proposition}
It might be good to comment on the precise meaning of the identity in Proposition \ref{prop: relation between adjoints of sub and exterior bundles}. By definition, $\sharp_0$ acts on (local) sections of $\fancy{H}_0^*$, while $\sharp$ acts on (local) sections of $\fancy{H}^*$. The identity means the following: If $u^*$ is a (local) section of $\fancy{H}_0^*$, then
by the Hahn-Banach theorem, $u^*$ extends to a (local) section of $\fancy{H}^*$. The $\sharp$-operator can act on such an extension, and the identity means that if we first let $\sharp$ act on such an extension and then apply $P$ to the result, we get $\sharp_0 u^*$. For this to make sense, $P$ acting on $\sharp$ acting on such an extension needs to be independent of the choice of extension of $u^*$. This is indeed the case. To see this, we let $\xi^*$ and $\eta^*$ be two extensions of $u^*$ as (local) sections of $\fancy{H}^*$, and verify that $P\sharp \xi^*=P\sharp \eta^*$. It suffices of course (due to the projection) to show that for all (local) sections $u_0$ of $\fancy{H}_0$, 
\begin{align}
h(u_0, P\sharp \xi^*-P\sharp \eta^*)=0.
\end{align}But this follows by definition of $\sharp$, orthogonality, and linearity. Indeed, we get
\begin{align}
h(u_0, P\sharp \xi^*-P\sharp \eta^*)&=h(u_0, P(\sharp \xi^*-\sharp \eta^*))=\ip{\xi^*-\eta^*}{u_0}=\ip{u^*-u^*}{u_0}=0,
\end{align}with the two last equalities following from the fact that $\xi^*$ and $\eta^*$ are both extensions of $u^*$. Conversely, if $u^*$ is a (local) section of $\fancy{H}^*$, and for each $b$ in a subset of the domain of definition of $u^*$, $u^*(b)$ happens to also be defined on $(H_0)_b$, then we may naturally view $u^*$ (restricted to this subset) as a (local) section of $\fancy{H}_0^*$. The identity then means that $\sharp_0$ acting on $u^*$ viewed as such a (local) section of $\fancy{H}_0^*$, is equal to $P\sharp u^*$. Our discussion actually also gives the proof of the proposition:
\begin{proof}[Proof of Proposition \ref{prop: relation between adjoints of sub and exterior bundles}]
Let $u^*$ be a local section of $\fancy{H}_0^*$, and let $\xi^*$ be any extension of $u^*$ to a local section of $\fancy{H}^*$. Let $u_0$ be a local section of $\fancy{H}_0$. Then, by definition of $\sharp_0$,
\begin{align}
\ip{u^*}{u_0}&=h(u_0,\sharp_0 u^*).
\end{align}On the other hand, since $\xi^*$ is an extension of $u^*$, we also have, by definition of $\sharp$,
\begin{align}
\ip{u^*}{u_0}&=\ip{\xi^*}{u_0}=h(u_0, \sharp \xi^*),
\end{align} By orthogonality the latter is equal to $h(u_0,P\sharp \xi^*)$, so by comparison it follows that we have $\sharp_0 u^*=P\sharp \xi^*$ (we must not believe that $h(u_0, \sharp \xi^*)=h(u_0, \sharp_0 u^*)$, which does hold true, gives that $\sharp_0 u^*=\sharp \xi^*$(!)).
\end{proof}

\section{Special hermitian metrics and Berndtsson's regularity theorem}\label{sec: special metric and Berndtsson's regularity theorem}

\subsection{Hermitian metrics of zero variation}
In the previous section we looked at some generalities of Hilbert bundles and hermitian metrics. Most of what we discussed in there applies to general Hilbert bundles. In this section we shall mostly restrict our attention to trivial Hilbert bundles. Let $H$ be a trivial Hilbert bundle over $B$. Then the Hilbert bundle $H$ admits a natural smooth and holomorphic structure. Let us start by briefly recapitulating this. Let $U\sub B$ be open, $u\in \Gamma(U;H)$, and $b\in U$. We consider what it means for $u$ to be smooth at $b$. Since smoothness is a local property, we may, by passing to local holomorphic coordinates near $b$, restrict to the case that $b$ is point in some open subset of $\set{C}^{m}$. We may then view $u$ as an $H$-valued map on this open subset, and we say that $u$ is \bt{smooth at $b$} if it is smooth in the real-Fréchet sense (see \cite{Cartan}) at $b$ when viewed as such a map. We say that $u$ is \bt{smooth (on $U$)} if it is smooth at each point in $U$. Let us denote by $u'(b)$ the real-Fréchet derivative of $u$ at $b$, viewed as such a map. The following definition asserts what it means for $u$ to be \ita{holomorphic} at $b$:

\begin{definition}\label{def: holomorphic sections, trivial bundles}
Using the above notation, we say that $u$ is \bt{holomorphic at $b$} if $u'(b)$ (which apriori is only real-linear) is complex-linear. We also say that $u$ is \bt{holomorphic (on $U$)} if it  is holomorphic at each point in $U$.
\end{definition}
Let us denote the subfamily of $\Gamma(U;H)$ consisting of smooth, respectively holomorphic, sections, by $\fancy{C}^{\infty}(U;H)$, respectively $\fancy{O}(U;H)$; in the global case that $U=B$, we simply write $\fancy{C}^{\infty}(H)$, respectively $\fancy{O}(H)$. Let us also write $\fancy{C}^{\infty}(U)$, respectively $\fancy{O}(U)$, for $\fancy{C}^{\infty}(U;\set{C})$, respectively $\fancy{O}(U;\set{C})$. Then $\fancy{C}^{\infty}(U;H)$ is a $\fancy{C}^{\infty}(U)$-module, and $\fancy{O}(U;H)$ is an $\fancy{O}(U)$-module. The smooth structure of $H$ gives rise to a notion of \ita{smooth hermitian metrics} on $H$:

\begin{definition}\label{def: smoothness of hermitian metrics}
Suppose that $h$ is a hermitian metric on $H$. Then we shall say that $h$ is \bt{smooth} if $h(u,v)$ is a smooth complex-valued function for all smooth local sections $u$ and $v$ of $H$.
\end{definition}
In general, if $\fancy{H}$ is a general Hilbert bundle over some complex, possibly infinite-dimensional manifold $\fancy{B}$, and there is an appropriate notion of smooth/holomorphic (local) sections of $\fancy{H}$ (for example, the families of such smooth/holomorphic sections should have natural module structures), we shall say that $\fancy{H}$ \bt{admits a smooth/holomorphic structure}. In the case that $\fancy{H}$ admits a smooth structure, and $h$ is a hermitian metric on $\fancy{H}$, we may define $h$ to be smooth in exactly the same way as in Definition \ref{def: smoothness of hermitian metrics}. Namely, by saying that $h$ is smooth if $h(u,v)$ is smooth for all smooth (local) sections $u$ and $v$ of $\fancy{H}$. Now, suppose that $\fancy{H}_0\leq \fancy{H}$. If $\fancy{H}$ admits a smooth/holomorphic structure, then so does $\fancy{H}_0$. Indeed, given a (local) section $u^*$ of $\fancy{H}_0$, we may view it as a (local) section of $\fancy{H}$, and then simply define $u^*$ to be smooth/holomorphic as a (local) section of $\fancy{H}_0$ if it is smooth/holomorphic when viewed as a (local) section of $\fancy{H}$. It is natural to ask whether the smooth/holomorphic structure of $\fancy{H}$ also induces a smooth/holomorphic structure on the dual bundle $\fancy{H}^*$ and or the quotient bundle $\fancy{Q}:=\fancy{H}/\fancy{H}_0$. This seems to be a rather hard question in general, and we will say just a little bit about the trivial case. Of course, since in this setting both the dual bundle and quotient bundle are trivial Hilbert bundles, they automatically admit smooth and holomorphic structures. The point is whether and or how these are induced from those of $\fancy{H}$.

Our discussion shall be by means of what we call \ita{hermitian metrics of zero variations}, whom we introduce next. Let $h^0$ denote the fixed inner product on $H$ (as a Hilbert space). Then $h^0$ induces a natural hermitian metric $\mathfrak{h}^0$ on $H$ (as a trivial Hilbert bundle) defined simply by letting, for all local sections $u$ and $v$ of $H$, and all $b$ in the intersection of the domains of definition of $u$ and $v$,
\begin{align}
\mathfrak{h}^0(u,v)(b)&:=h^{0}(u^{b}, v^{b}).
\end{align}In particular, if $u$ and $v$ are elements of $H$ (viewed as constant sections), then of course we get 
\begin{align}
\mathfrak{h}^{0}(u,v)(b)=\mathfrak{h}^{0}(u,v)(b')
\end{align} for all $b,b'\in B$. As a short-hand form, we shall express this by writing 
\begin{align}
\mathfrak{h}^{0}_{b}=\mathfrak{h}^{0}_{b'}
\end{align} for all $b,b'\in B$. We think of this as $\mathfrak{h}^{0}$ not varying with respect to the base parameter, and hence being of \ita{zero variation}. This motivates the following definition:

\begin{definition}\label{def: hermitian metrics of zero variation}
Let $H$ be a trivial Hilbert bundle over $B$. We shall say that a hermitian metric $h$ on $H$ is of, or has, \bt{zero variation}, if $h_{b}=h_{b'}$ for all $b,b'\in B$.
\end{definition} 
Despite their simplicity, hermitian metrics of zero variation turn out to be very important objects in our theory. Among other things, they play an essential role in the statement and the proof of Berndtsson's regularity theorem below (Theorem \ref{thm: Berndtsson's regularity theorem, general version}). Another aspect that is worth mentioning has to do with the existence of \ita{Chern connections}; we discuss this more closely in section \ref{sec: Positivity and Variation}, when these objects are introduced properly. Using calculus (see \cite{Cartan}), we get the following simple, but useful, characterization of hermitian metrics of zero variation:

\begin{lemma}\label{lemma: characterization of hermitian metrics of zero variation}
Let $H$ be a trivial Hilbert bundle over $B$, and let $h$ be a hermitian metric on $H$. Then $h$ is of zero variation if and only if $h$ is smooth and satisfies
\begin{align}
\partial h(u,v)&=h(\partial u,v)+h(u,\dbar v)\label{eq: identity for characterization of hermitian metrics of zero variation}
\end{align}for all smooth local sections $u$ and $v$ of $H$.
\end{lemma}
\begin{proof}
Suppose first that $h$ is of zero variation. Then, by calculus (see \cite{Cartan}), $h$ is smooth, and we have
\begin{align}
h(u,v)'&=h(u',v)+h(u,v').
\end{align}The identity \eqref{eq: identity for characterization of hermitian metrics of zero variation} thus follows by definition of partial derivatives (see \cite{Cartan}) and the sesquilinearity of $h$. Conversely, suppose that $h$ is smooth and that \eqref{eq: identity for characterization of hermitian metrics of zero variation} holds. Note that $h$ is of zero variation if and only if $\partial h(u,v)=0$ for all $u,v\in H$ (viewed as constant sections). This is of course a special case of \eqref{eq: identity for characterization of hermitian metrics of zero variation}, so $h$ is of zero variation.
\end{proof}
Suppose that $E$ is a finitely-ranked holomorphic vector bundle over $U$, and that $u$ is a smooth (local) section of $E$ over $U$. Let $\dbar=\dbar_t$, where we choose generic local holomorphic coordinates $t:=(t_1,\ldots, t_m)$ on $U$; recall also Example \ref{ex: first ex of trivial Hilbert bundles}. Then $u$ is holomorphic if and only if $u$ satisfies the \bt{Cauchy-Riemann equation} $\dbar u=0$. Not too surprisingly, the same characterization holds also in the possibly infinitely-ranked case that $u\in \fancy{C}^{\infty}(U;H)$:

\begin{corollary}\label{cor: characterization of holomorphic sections in the trivial case using Cauchy-Riemann equation}
Let $u\in \fancy{C}^{\infty}(U;H)$. Then, using the above notation, $u\in \fancy{O}(U;H)$ if and only if $\dbar u=0$.
\end{corollary}

\begin{proof}
Let $h$ be any hermitian metric of zero variation on $H$, and let $v\in H$ (viewed as a constant section, but we shall stop pointing this out). From the proof of Lemma \ref{lemma: characterization of hermitian metrics of zero variation}, 
\begin{align}
h(u,v)'&=h(u',v),
\end{align}from which it follows that $u'$ is complex-linear if and only if $h(u,v)'$ is complex-linear for all $v\in H$. The latter is equivalent to $\dbar h(u,v)=0$ for all $v\in H$, which by Lemma \ref{lemma: characterization of hermitian metrics of zero variation} is the same as $h(\dbar u,v)=0$ for all $v\in H$. The latter again is equivalent to $\dbar u=0$, which proves the assertion.
\end{proof}
We shall use the following notation: Let $\set{K},\set{F}\in \{\set{R}, \set{C}\}$. If $E$ and $F$ are two normed spaces over $\set{F}$, with $\set{K}\sub \set{F}$, then we denote by $\mathcal{L}_{\set{K}}(E;F)$ the space of continuous $\set{K}$-linear maps $E\to F$, and by $\text{Iso}_{\set{K}}(E;F)$, the space of $\set{K}$-linear isomorphisms $E\to F$. Consider now the dual bundle of $H$, $H^*$, and let $\sharp$ denote the adjoint operator with respect to the hermitian metric $\mathfrak{h}^0$ of zero variation. Then we have the following characterization of smooth and holomorphic sections of $H^*$:
\begin{proposition}\label{prop: characherization of smooth/holom structure of dual bundle, trivial case, using trivial hermitian metric of zero variation}
With the above notation and set-up, the following holds:
\begin{enumerate}[(i)]
\item Let $u^*\in \Gamma(U;H^*)$. Then $u^*\in \fancy{C}^{\infty}(U;H^*)\iff \sharp u^*\in \fancy{C}^{\infty}(U;H)$.
\item Let $u^*\in \fancy{C}^{\infty}(U;H^*)$. Then $u^*$ is holomorphic if and only if $\sharp u^*$ is antiholomorphic; that is, if $(\sharp u^*)'$ is complex-antilinear.
\end{enumerate}
\end{proposition}
\begin{proof}
$\phantom{-}$
\begin{enumerate}[(i)]
\item Let $\tilde{\sharp}$ be defined on $B$ by $\tilde{\sharp}(b):=\sharp^{b}$, where $\sharp^{b}$ is defined by $\sharp^{b}u^*:=(\sharp u^*)^{b}$ for all $u^*\in H^*$. We claim that $\tilde{\sharp}$ is a map $B\to \text{Iso}_{\set{R}}(H^*;H)$. It is clear that $\tilde{\sharp}(b)$ is real-linear for each $b\in B$. It is also continuous for we have for all $u^*\in H^*$, that:
\begin{align}
\norm{\tilde{\sharp}(b)u^*}&=\norm{(\sharp u^*)^{b}}=\norm{\sharp u^*}_{\mathfrak{h}^0}(b)=\norm{u^*}_{(\mathfrak{h}^0)^{-1}}.
\end{align}Thus, $\tilde{\sharp}(b)$ is continuous with norm 1 for all $b\in B$. Finally, it also has an inverse, namely $(\tilde{\sharp}(b))^{-1}:=(\tilde{\sharp}^{-1})^{b}$, where $(\tilde{\sharp}^{-1})^{b}$ is defined similar to $\sharp^{b}$: $(\tilde{\sharp}^{-1})^{b}u:=(\sharp^{-1}u)^{b}$ for all $u\in H$. Since $\mathfrak{h}^0$ is of zero variation, it is easy to check that $\sharp^{b}=\sharp^{b'}$ for all $b,b'\in B$ (see also the proof of Proposition \ref{prop: smoothness and holomorphicity of dual bundles in trivial setting under additional assumption on hermitian metrics of zero variation} below). In particular, $\tilde{\sharp}$ is constant in the sense that $\tilde{\sharp}(b)=\tilde{\sharp}(b')$ for all $b,b'\in B$, and hence smooth. It follows that for all local sections $u^*$ of $H^*$, $\sharp u^*=\tilde{\sharp}\circ u^*$, which shows that $\sharp u^*$ is smooth if $u^*$ is smooth. For the converse, let $\mathfrak{i}$ denote the inversion operator $\text{Iso}_{\set{R}}(H^*;H)\to \text{Iso}_{\set{R}}(H;H^*), T\mapsto T^{-1}$. By Calculus (see \cite{Cartan}), $\mathfrak{i}$ is smooth, so it follows that $\tilde{\sharp}^{-1}:=\mathfrak{i}\circ \tilde{\sharp}$ is smooth as a map $B\to \text{Iso}_{\set{R}}(H,H^*)$. Hence, if $\sharp u^*$ is smooth, then $u^*=\tilde{\sharp}^{-1}\circ (\sharp u^*)$ is smooth. This proves (i).
\item By the proof of (i), it follows that the dual of $\mathfrak{h}^0$, $(\mathfrak{h}^{0})^{-1}$, is of zero variation. Thus, by the proof of Corollary \ref{cor: characterization of holomorphic sections in the trivial case using Cauchy-Riemann equation}, $u^*$ is holomorphic if and only if
\begin{align}
\dbar (\mathfrak{h}^0)^{-1}(u^*,v^*)=0
\end{align}for all $v^*\in H^*$. By definition of $(\mathfrak{h}^{0})^{-1}$, we therefore have that $u^*$ is holomorphic if and only if $\partial \sharp^*u=0$, which precisely means that $\sharp u^*$ is antiholomorphic; one may also observe that $\partial\sharp u^*=\sharp(\dbar u^*)$, (where the right-hand side has a natural and obvious meaning).
\end{enumerate}
\end{proof}
If we try to prove a similar statement as Proposition \ref{prop: characherization of smooth/holom structure of dual bundle, trivial case, using trivial hermitian metric of zero variation} where $\mathfrak{h}^0$ is replaced by a more general hermitian metric of zero variation, we seem to run into some problems. We will prove a similar statement under an additional assumption. Let $h$ be a hermitian metric on $H$, and let $\sharp$ be the adjoint operator with respect to $h$. We shall say that $\sharp$ is \bt{smooth} if $\sharp u^*$ is (a) smooth (local) (section of $H$) for all smooth (local) (sections) $u^*$ (of $H^*$). Similarly, we shall say that $\sharp^{-1}$ is \bt{smooth} if $\sharp^{-1}u$ is (a) smooth (local section of $H^*$) whenever $u$ is (a) smooth (local section of $H$). We may also consider similar definitions in the non-trivial setting, but let us not get into that here. The proof of Proposition \ref{prop: characherization of smooth/holom structure of dual bundle, trivial case, using trivial hermitian metric of zero variation} gives first the following proposition:

\begin{proposition}\label{prop: smoothness and holomorphicity of dual bundles in trivial setting under additional assumption on hermitian metrics of zero variation}
Suppose that $h$ is a hermitian metric of zero variation on $H$, and let $\sharp$ denote the adjoint operator with respect to $h$. Assume that if we define $\tilde{\sharp}$ as a map on $B$ by $\tilde{\sharp}(b):=\sharp^{b}$ with $\sharp^{b}u^*:=(\sharp u^*)^{b}$ for all $u^*\in H^*$, $\tilde{\sharp}$ is a map $B\to \text{Iso}_{\set{R}}(H^*;H)$. Then $\sharp$ and $\sharp^{-1}$ are smooth.
\end{proposition}

\begin{proof}
By the proof of Proposition \ref{prop: characherization of smooth/holom structure of dual bundle, trivial case, using trivial hermitian metric of zero variation}, it suffices to show that $\tilde{\sharp}$ is smooth, but this follows since $\tilde{\sharp}(b)=\tilde{\sharp}(b')$ for all $b,b'\in B$ as $h$ is of zero variation. Indeed, given $u^*\in H^*$ and $u\in H$, 
\begin{align}
\ip{u^*}{u}(b)&=\ip{u^*}{u}(b')\implies h(u, \sharp^{b}u^*)=h(u, \sharp^{\tilde{b}}u^*)\implies (\tilde{\sharp}(b)-\tilde{\sharp}(b'))u^*=0.
\end{align}
\end{proof}
Using Proposition \ref{prop: smoothness and holomorphicity of dual bundles in trivial setting under additional assumption on hermitian metrics of zero variation} (and its proof), we then get the following result that is similar to Proposition \ref{prop: characherization of smooth/holom structure of dual bundle, trivial case, using trivial hermitian metric of zero variation}:

\begin{corollary}\label{cor: smoothness and holomorphicity of dual bundles in terms of hermitian metrics of zero variation with additional property}
With the same assumptions, and using the same notation, as in Proposition \ref{prop: smoothness and holomorphicity of dual bundles in trivial setting under additional assumption on hermitian metrics of zero variation}, the following holds:
\begin{enumerate}[(i)]
\item Given $u^*\in \Gamma(U;H^*)$, $u^*\in \fancy{C}^{\infty}(U;H^*)\iff \sharp u^*\in \fancy{C}^{\infty}(U;H)$.
\item Given $u^*\in \fancy{C}^{\infty}(U;H)$, $u^*$ is holomorphic if and only if $\sharp u^*$ is antiholomorphic.
\end{enumerate}
\end{corollary}

\begin{proof}
$\phantom{-}$
\begin{enumerate}[(i)]
\item Suppose that $u^*$ is smooth. Then $\sharp u^*$ is smooth since $\sharp$ is smooth by Proposition \ref{prop: smoothness and holomorphicity of dual bundles in trivial setting under additional assumption on hermitian metrics of zero variation}. Conversely if $\sharp u^*$ is smooth so is $u^*=\sharp^{-1}(\sharp u^*)$, again by the proposition.
\item By the proof of Proposition \ref{prop: smoothness and holomorphicity of dual bundles in trivial setting under additional assumption on hermitian metrics of zero variation}, since $h$ is of zero variation, $h^{-1}$, its dual metric, is also of zero variation. Thus the assertion follows from the same argument as in the proof of Proposition \ref{prop: characherization of smooth/holom structure of dual bundle, trivial case, using trivial hermitian metric of zero variation} (ii).
\end{enumerate}
\end{proof}
Similar characterizations hold for the smooth and holomorphic structure on the quotient bundle $Q:=H/H_0$. Since the arguments are very similar to those in the case of the dual bundle, we will here content ourselves with stating the simplest result. Let $m_{l}$ denote the minimal lifting operator with respect to $\mathfrak{h}^0$. Then we have the following proposition that can be thought of as the ''quotient analogue'' to Proposition \ref{prop: characherization of smooth/holom structure of dual bundle, trivial case, using trivial hermitian metric of zero variation}:

\begin{proposition}\label{prop: characterization of smooth/holomorphic structure of trivial quotient bundle in terms of trivial minimal lifting operator}
With the above notation, the following holds:
\begin{enumerate}[(i)]
\item Given $[u]\in \Gamma(U;Q)$, we have $[u]\in \fancy{C}^{\infty}(U;Q)\iff m_{l}[u]\in \fancy{C}^{\infty}(U;H)$.
\item If $[u]\in \fancy{C}^{\infty}(U;Q)$, then $[u]\in \fancy{O}(U;Q)\iff m_{l}[u]\in \fancy{O}(U;H)$.
\end{enumerate}
\end{proposition}
We shall have use for, and come back to, the above characterizations later. For example in the proof of Corollary \ref{cor: generaization of smoothness of Bergman kernels} below, and in section \ref{sec: Positivity and Variation}, when we discuss \ita{Chern connections}. In there, we shall also discuss more general versions of the above results (see Propositions \ref{prop: Chern connection for dual bundle, trivial case} and \ref{prop: Chern connection for quotient bundle, trivial case}).

\subsection{Weighted hermitian metrics}

Suppose now that $h$ is a smooth hermitian metric on $H$, not necessarily of zero variation, and that $\overline{H_0}^{(H,h)}=H_0$. Let $P=P(h)$ denote the total orthogonal projection map associated with $h$ (of $H_0$ in $H$). Similar to the case of the adjoint operator (and its inverse), we shall say that $P$ is \bt{smooth} if $Pu$ is (a) smooth (local) (section of $H_0$) whenever $u$ is (a) smooth (local) (section of $H$). We will see several instances throughout the story where the smoothness of $P$ is of natural interest and importance; Proposition \ref{prop: relation between adjoints of sub and exterior bundles} already provides an example: Using the notation in there, if $\sharp$ is smooth, that $P$ is smooth implies that $\sharp_0$ is smooth. As we shall see, this may be viewed as a direct counterpart to Lemma \ref{intro: smoothness of the Bergman kernel} in the introduction. Another example is furnished by the smoothness of minimal lifting operators. First, we need to properly state what we mean by ''the smoothness of minimal lifting operators''; this will of course be similar to the smoothness of adjoint operators and total orthogonal projection maps: Let $m_{l}$ denote the minimal lifting operator with respect to $h$. We shall say that $m_{l}$ is \bt{smooth} if $m_{l}[u]$ is smooth whenever $[u]$ is smooth. Similarly, we  shall say that $m_{l}^{-1}$ is \bt{smooth} if $m_{l}^{-1}u$ is smooth whenever $u$ is smooth. With this, we then have that $m_{l}$ is smooth if $P$ is smooth: Suppose that $[u]$ is smooth, and let $[u]$ be a smooth representative of $[u]$ (an example is $m_{l}^{0}[u]$ where $m_{l}^{0}$ denotes the minimal lifting operator with respect to $\mathfrak{h}^0$). Then $m_{l}[u]=u-Pu$, from which it follows that $m_{l}$ is smooth given that $P$ is. We also see from this that the converse is true. That is, if $m_{l}$ is smooth, then so is $P$: Let $u$ be smooth, and let $\pi_{Q}:H\to Q$ be the natural projection. Then $\pi_{Q}u:=[u]$ is smooth, so $Pu=u-m_{l}[u]$ is smooth. 

We shall show that if $h$ is what we call a \ita{smooth weighted hermitian metric} \ita{induced from a hermitian metric of zero variation}, then $P$ is smooth. This is Berndtsson's regularity theorem (Theorem \ref{thm: Berndtsson's regularity theorem, general version} below), one of our main results, and the generalization of Theorem \ref{intro: Berndtsson's regularity theorem, special case} from the introduction to our current setting of (trivial) Hilbert bundles. Our notion of \ita{weighted hermitian metrics} is motivated by the following elementary, but important, example:

\begin{example}\label{example: motivation for weighted hermitian metrics}
Let $H$ be a trivial Hilbert bundle over $B$, and suppose that $h$ is a smooth hermitian metric on $H$. Let $f>0$ be a smooth function on $B$. We write $f:=e^{-\phi}$ for some $\phi$ (namely $\phi=-\log(f)$); this is of course inspired by the inner product $(\cdot,\cdot)_t$ in the introduction. If we let $u$ be a smooth (local) section of $H$ over $U$, some open subset of $B$, then by definition, $fu$ is the smooth (local) section of $H$ over $U$ given by $(fu)(b):=f(b)u^{b}$. It follows that $\mathfrak{h}$ defined by
\begin{align}
\mathfrak{h}(v_1,v_2):=h(e^{-\phi}v_1,v_2)
\end{align}for all local sections $v_1$ and $v_2$ of $H$, is a smooth hermitian metric on $H$. Moreover, $h$ can be recovered from $\mathfrak{h}$ by
\begin{align}
h(v_1,v_2)&=\mathfrak{h}(e^{\phi}v_1,v_2).
\end{align}
Observe that for each $b\in B$, (multiplication by) $e^{-\phi(b)}$ gives a complex-linear automorphism $H\to H$. We may therefore view $f$ also as a map $f:B\to \text{Aut}(H)$, where $\text{Aut}(H)=\text{Iso}_{\set{C}}(H;H)$ is the space of complex-linear isomorphisms $H\to H$. To distinguish $f$ as such a map from $f$ as a map $B\to \set{R}$, let us for the moment denote the former map $B\to \text{Aut}(H)$ by $\tilde{f}$. Thus, $\tilde{f}=\iota \circ f$, where $\iota: \set{R}\to \text{Aut}(H)$ is the canonical inclusion map. Since $f$ and $\iota$ are both smooth, so is $\tilde{f}$. Note that $f$ also induces a map $\hat{f}:B\times H\to H, (b,v)\mapsto \tilde{f}(b)v$, which in turn for a fixed local section $u$ of $H$ over $U$ induces a map $U\to H, b\mapsto \tilde{f}(b)u^{b}$. This is precisely the map $fu$ from earlier. By direct inspection, we have that $fu$ is smooth. This follows simply because $\fancy{C}^{\infty}(U;H)$ is a $\fancy{C}^{\infty}(U)$-module. Now, given that $\tilde{f}$ is smooth, the smoothness of $fu$ can also be established as follows: Let $\hat{u}$ denote the map $U\to B\times H, b \mapsto (b, u^{b})$. Then $\hat{u}$ is smooth if and only if $u$ is, and we have $fu=\hat{f}\circ \hat{u}$. We claim that $\hat{f}$ is smooth when $\tilde{f}$ is. To see this, we note that $\hat{f}$ can be written $\hat{f}=m_2\circ m_1$, where $m_2:=\text{ev}$ denotes the evaluation map $\text{Aut}(H)\times H\to H, (T,v)\mapsto Tv$, and $m_1=m_{1,\tilde{f}}$ is the map $B\times H\to \text{Aut}(H)\times H, (b,v)\mapsto (\tilde{f}(b),v)$. By calculus (see \cite{Cartan}) $m_{2}$ is continuous, thus (being bilinear) smooth. It is also clear that $m_1$ is smooth when $\tilde{f}$ is. Therefore $\hat{f}$ is smooth as required. Finally, let $\mathfrak{i}:\text{Aut}(H)\to \text{Aut}(H)$ denote the inversion map. We know (again, by calculus) that $\mathfrak{i}$ is smooth. From this it follows that $f^{-1}u:=b\mapsto e^{\phi(b)}u^{b}$ is smooth. The argument is analogous to before; using similar notation, the first step is that $\tilde{f}^{-1}$ is smooth, which follows from $\tilde{f}$ being smooth, since $\tilde{f}^{-1}=\mathfrak{i}\circ \tilde{f}$. The conclusion of all this is that given that $\tilde{f}$ is smooth, we have that $fu$ and $f^{-1}u$ are smooth whenever $u$ is.
\end{example}
The rather elaborate discussion in Example \ref{example: motivation for weighted hermitian metrics} has its purposes(s). Indeed, the property that the map $B\times H\to H$ induced from $f$ ($\tilde{f}$ in the example), is smooth, is quite essential in the proof of Berndtsson's regularity theorem below (Theorem \ref{thm: Berndtsson's regularity theorem, general version}); see also the upcoming remark. The discussion in Example \ref{example: motivation for weighted hermitian metrics} leads to the following definition:

\begin{definition}\label{def: weighted hermitian metrics}
Let $H$ be a trivial Hilbert bundle over $B$, and suppose that $h$ is a (smooth) hermitian metric on $\fancy{H}$. By a \bt{(smooth) weight (map) for $h$} we shall mean a (smooth) map 
\begin{align}
w:B\to \text{Aut}(H), b\mapsto w^{b}
\end{align} such that defining $\mathfrak{h}$ by
\begin{align}
\mathfrak{h}(u,v)&:=h(wu,v),
\end{align}for all (smooth) local sections $u$ and $v$ of $H$, where $wu$ is the (smooth) local section of $H$ given by 
\begin{align}
(wu)^{b}:=w^{b}u^{b},
\end{align} $\mathfrak{h}$ is a (smooth) hermitian metric on $\fancy{H}$. In this case, we refer to the (smooth) map \begin{align}
w^{-1}:B\to \text{Aut}(H), b\mapsto (w^{b})^{-1}:=(w^{-1})^{b}
\end{align} as the \bt{inverse (smooth) weight (map) of $w$}, and to $\mathfrak{h}$ as a \bt{(smooth) weighted hermitian metric on $H$ (induced by $h$)}. We also say in this case that $w$ is a \bt{(smooth) weight \underline{of} $\mathfrak{h}$}.
\end{definition}

\begin{remark}
If $w$ is a smooth weight for $h$ of $\mathfrak{h}$, then $w^{-1}$ is a smooth weight for $\mathfrak{h}$ of $h$; this follows by simply replacing $u$ in the definition of $\mathfrak{h}$ by $w^{-1}u$ so that we get
\begin{align}
h(u,v)&=\mathfrak{h}(w^{-1}u,v).
\end{align}Note also that by our discussion in Example \ref{example: motivation for weighted hermitian metrics}, if $w$ is smooth, then $w^{-1}$ is automatically smooth.
\end{remark}

\begin{example}\label{ex: example of trivial hilbert bundles with weighted hermitian metrics}
Let $L^2,A^2,U$, and $\phi$ be as in the introduction, with $\phi$ smooth up to the boundary and $\Omega$ bounded. We have seen that $L^2$ and $A^2$ are trivial Hilbert bundles over $U$ and that $A^2\leq L^2$ as such bundles. The fixed inner product on $L^2$ (and $A^2$) is the standard $L^2$-inner product on $\Omega$. Let us denote the hermitian metric of zero variation induced by it by $h^{0}$. Then (multiplication by) $e^{-\phi^{t}}$, viewed as a map $U\to \text{Aut}(L^2)$ given by $t\mapsto e^{-\phi^{t}}$, defines a smooth weight for $h^{0}$. Let us denote the smooth weighted hermitian metric induced by $h^0$ by $h$. Then for all $f,g\in L^2$, $h_{t}(f,g)=(f,g)_{t}$, where the right-hand side inner product is the one from the introduction.

Analogous to Example \ref{ex: first ex of trivial Hilbert bundles}, there is also a global analogue of the above example. Let $H,H_0$, and $h^0$ be as in Example \ref{ex: first ex of trivial Hilbert bundles}. We have seen that $H$ and $H_0$ are trivial Hilbert bundles over $B$ with $H_0\leq H$, and that $h^0$ is a hermitian metric of zero variation on $H$ (and $H_0$). For each $b\in B$, suppose that $\psi^{b}:=\{\psi^{b}\}$ is a collection of locally defined real-valued functions such that $e^{-(\phi+\psi^{t})}$ can be integrated over $X$ when twisted with some $i^{n^2}u\wedge \bar{v}$, for $u,v\in H$, and let $\psi:=\{\psi\}$, where $\psi(b,z):=\psi^{b}(z)$ for $b\in B$ and $z\in X$. Suppose also for simplicity that $b\mapsto \psi$ is smooth and that $e^{-\psi^{b}}$ is bounded for each $b$. Then (multiplication by) $e^{-\psi}$ defines a smooth weight for $h^{0}$. If we let $h$ denote the induced smooth weighted hermitian metric with smooth weight $e^{-\psi}$, then for all $u,v\in H$ and all $b\in B$, we have
\begin{align}
h_b(u,v)&=\int_{X}i^{n^2}u\wedge \bar{v}e^{-\theta^{b}},
\end{align} where we put $\theta:=\phi+\psi$ and $\theta^{b}:=\phi+\psi^{b}$. Our most important examples in practice are akin to this example.
\end{example}
Our next result is a cute proposition on adjoint operators of weighted hermitian metrics:

\begin{proposition}\label{prop: relation between adjoint of metric and adjoint of weighted metrics}
Let $H$ be a trivial Hilbert bundle over $B$, and suppose that $h$ is a hermitian metric on $H$. Suppose also that $\mathfrak{h}$ is a weighted hermitian metric on $H$ induced by $h$ with weight $w$. Let $\sharp_{h}$ denote the adjoint operator for $H^*$ with respect to $h$, and $\sharp_{\mathfrak{h}}$ the adjoint operator for $H^*$ with respect to $\mathfrak{h}$. Then 
\begin{align}
\sharp_{\mathfrak{h}}&=w^{-1}\sharp_{h}.
\end{align}
\end{proposition}

\begin{proof}
Let $u^*$ be local section of $H^*$. Then for all $u\in H$,
\begin{align}
\ip{u^*}{u}=h(u, \sharp_{h}u^*)=\mathfrak{h}(u, \sharp_{\mathfrak{h}}u^*).
\end{align}Conjugating, we therefore have
\begin{align}
\mathfrak{h}(\sharp_{\mathfrak{h}}u^*,u)=h(\sharp_h u^*,u).
\end{align}By definition of $\mathfrak{h}$, the left-hand side of this identity is equal to $h(w \sharp_{\mathfrak{h}}u^*,u)$, which implies that $w\sharp_{\mathfrak{h}}=\sharp_{h}$. Applying $w^{-1}$ then gives the desired relation.
\end{proof}
In particular, if $\sharp_{h}$ and $w$ are smooth, then $\sharp_{\mathfrak{h}}$ is also smooth. Combining this with Proposition \ref{prop: relation between adjoints of sub and exterior bundles} and Corollary \ref{cor: smoothness and holomorphicity of dual bundles in terms of hermitian metrics of zero variation with additional property}, we get a natural generalization, or analogue, of Lemma \ref{intro: smoothness of the Bergman kernel} in the introduction. 

\subsection{Berndtsson's regularity theorem}
We end this section with Berndtsson's regularity theorem and a discussion on this. First up is Berndtsson's regularity theorem, which generalizes Theorem \ref{intro: Berndtsson's regularity theorem, special case} from the introduction:

\begin{theorem}[Berndtsson's regularity theorem]\label{thm: Berndtsson's regularity theorem, general version}
Let $H$ be a trivial Hilbert bundle over $B$, $h^0$ a hermitian metric of zero variation on $H$, and  $h$ a smooth weighted hermitian metric on $H$ induced by $h^0$. Suppose that $H_0\leq H$ is closed in $H$ with respect to $h$ and $h^0$, and let $P=P(h)$ denote the total orthogonal projection map associated with $h$. Then $P$ is smooth.
\end{theorem}

\begin{proof}
We may consider $P$ as a map $B\times H\to H_0$ given by sending $(b,u)\in B\times H$ to $P^{b}u$ in $H_0$. It suffices to show that $P$ as such a map is smooth, since it then follows that $P$ is smooth (in our usual sense) by composition. It might be useful here to refer back to (our discussion in) Example \ref{example: motivation for weighted hermitian metrics}. Fix $ a \in H$ and $b_0\in B$. Since smoothness is a local property, it suffices to show that $P$ as a map $B\times H\to H_0$ is smooth near $(b_0, a)$. Suppose that the smooth weight for $h^0$ of $h$ is $w$, and consider the auxiliary smooth map
\begin{align}
F^{b_0}:B\times H\times H_0\to H_0, (b,u,v)\mapsto P^{b_0}((w^{-1})^{b_0}(w^{b}(v-u))).
\end{align}We claim that $F^{b_0}(b,u,v)=0$ if and only if $v=P^{b}(u)$. Assume that the claim is true, and note that $F^{b_0}(b_0, a, P^{b_0}a)=0$. By direct computation, for all $v,\Delta v\in H_0$, we have
\begin{align}
F^{b_0}(b_0, a, v+\Delta v)-F^{b_0}(b_0, a, v)&=P^{b_0}(v+\Delta v -a)-P^{b_0}(v-a)=\Delta v,
\end{align}since $P^{b_0}$ is linear and restricts to the identity on $H_0$. It follows from this that the partial derivative (see \cite{Cartan}) of $F^{b_0}$ at $(b_0, a,v)$ in the $v$-direction (in the $H_0$-direction), is the identity map. In particular, it is invertible. Thus, by the implicit mapping theorem (see \cite{Cartan}), there is an open neighbourhood $U$ of $(b_0, a)$ and a unique smooth map $f$ on $U$ with values in $H_0$, such that 
\begin{align}
F^{b_0}(b,u,v)=0
\end{align} if and only if 
\begin{align}
v=f(b,u)
\end{align} for all $(b,u)\in U$. By the claim it follows then that $f=P$ on $U$. Hence, $P$ is smooth. To complete the proof, it therefore remains to show that the claim is true. We now do this. By definition, \begin{align}
F^{b_0}(b,u,v)=0
\end{align} if and only if 
\begin{align}P^{b_0}((w^{-1})^{b_0}(w^{b}(v-u)))=0.
\end{align} That is, if and only if 
\begin{align}
h_{b_0}((w^{-1})^{b_0}(w^{b}(v-u)), v_0)=0
\end{align} for all $v_0\in H_0$, which by definition of $w^{-1}$ is equivalent to
\begin{align}
(h^{0})_{b_0}(w^{b}(v-u), v_0)=0
\end{align} for all $v_0\in H_0$. Now, since $h^{0}$ is of zero variation, $(h^0)_{b_0}=(h^{0})_{b}$, so this is in turn equivalent to
\begin{align}
(h^{0})_{b}(w^{b}(v-u), v_0)=0
\end{align} for all $v_0\in H_0$. By definition of $w$ and $h$, the latter is equivalent to
\begin{align}h_b(v-u,v_0)=0
\end{align} for all $v_0\in H_0$, which in turn is equivalent to
\begin{align}P^{b}(v-u)=0.\label{eq: last line in proof of Berndtsson's reg theorem, general}
\end{align} Finally, since $v\in H_0$ already, $P^{b}v=v$, so linearity gives that \eqref{eq: last line in proof of Berndtsson's reg theorem, general} holds if and only if 
\begin{align}v=P^{b}u.
\end{align} This proves the claim, which completes the proof.
\end{proof}
A very important corollary to Berndtsson's regularity theorem (and previous results) is now the following:

\begin{corollary}\label{cor: generaization of smoothness of Bergman kernels}
With the same assumptions, and the same notation, as in Theorem \ref{thm: Berndtsson's regularity theorem, general version}, let $\sharp^{0}$ denote the adjoint operator for $H^*$ with respect to $h^0$, and suppose that $B\ni b\mapsto (\sharp^{0})^{b}$ defined by $(\sharp^{0})^{b}u:=(\sharp^{0}u)^{b}$, for all $u\in H$, gives a map $B\to \text{Iso}_{\set{R}}(H^*;H)$. Let $\sharp$ denote the adjoint operator with respect to $h$ for $H^*$, and $\sharp_0$ the adjoint operator with respect to (the restriction of) $h$ (as a hermitian metric on $H_0$) for $(H_0)^*$. Then $\sharp_0$ is smooth.
\end{corollary}
As will be seen from the proof, we may replace the assumption that $b\mapsto (\sharp^0)^{b}$ is a map $B\to \text{Iso}_{\set{R}}(H^*;H)$ with the assumption that $\sharp^0$ is smooth. We then get a slightly shorter and more general result. 
\begin{proof}[Proof of Corollary \ref{cor: generaization of smoothness of Bergman kernels}]
Let $w$ denote the smooth weight map for $h^0$ of $h$. By Corollary \ref{cor: smoothness and holomorphicity of dual bundles in terms of hermitian metrics of zero variation with additional property}, $\sharp^{0}$ is smooth and by Proposition \ref{prop: relation between adjoint of metric and adjoint of weighted metrics}, we have $\sharp=w^{-1}\sharp^{0}$. Thus, by Proposition \ref{prop: relation between adjoints of sub and exterior bundles}, we have $\sharp_{0}=P\sharp=P(w^{-1}\sharp^{0})$. Hence it follows by Theorem \ref{thm: Berndtsson's regularity theorem, general version} that $\sharp_0$ is smooth.
\end{proof}

Corollary \ref{cor: generaization of smoothness of Bergman kernels} may be seen as a generalization, or a natural analogue, of Lemma \ref{intro: smoothness of the Bergman kernel} from the introduction. To explain this, let $L^2$, $A^2$, $U$, $\phi$, $h^0$, and $h=\{h_t\}_{t\in U}$ be as in Example \ref{ex: example of trivial hilbert bundles with weighted hermitian metrics}. We know that $L^2$ and $A^2$ are trivial Hilbert bundles over $U$. The dual bundle of $A^2$ is $(A^2)^*$. Let $\sharp_{A^2}$ denote the adjoint operator for it with respect to the smooth weighted hermitian metric $h$. By Corollary \ref{cor: generaization of smoothness of Bergman kernels}, given any smooth local section $\xi^{*}$ of $(A^2)^*$, $\sharp_{A^2}\xi^*$ is smooth. Fix $z\in \Omega$ and consider $\xi:=\xi_{z}$ to be the evaluation map at $z$ for elements in $A^2$. That is, given any $f\in A^2$, we let $\ip{\xi_{z}}{f}:=f(z)$. It is clear that $\xi$ is smooth as a section of $(A^2)^*$. Using the same notation as in the introduction, let $k^{t}$ denote the Bergman kernel of $A^2_t$. Then $\sharp_{A^2}\xi(t)$ is precisely $k^{t}_z$. Thus we see that Lemma \ref{intro: smoothness of the Bergman kernel} is a special case of the corollary. An interesting generalization of this, which concerns derivatives of the Bergman kernel of $A^2_t$, is the following:

\begin{example}\label{ex: derivatives of Bergman kernels}
With the same notation as above, let $\xi^{\alpha}_{z}$ be defined as an element of $(A^2)^*$ by $\ip{\xi^{\alpha}_z}{f}:=\partial^{\alpha}f(z)$, where $\alpha$ is a multiindex. We may of course replace $f$ with a local section of $A^2$. Considering specifically $k^{t}_{w}$, we find (by the reproducing property of the Bergman kernel) that
\begin{align}
\ip{\xi^{\alpha}_{z}}{k^{t}_{w}}&=h_{t}(k^{t}_{w},\sharp_{A^2}\xi^{\alpha}_{z})=\overline{h_{t}(\sharp_{A^2}\xi^{\alpha}_{z},k^{t}_{w})}=\overline{\sharp_{A^2}\xi^{\alpha}_{z}(w)}.
\end{align}On the other hand, by definition of $\xi^{\alpha}_{z}$, this is also equal to $\partial^{\alpha}k^{t}_{w}(z)$. Thus, 
\begin{align}
\sharp_{A^2}\xi^{\alpha}_{z}(w)&=\overline{\partial^{\alpha}k^{t}_{w}(z)}=\dbar^{\alpha}k^{t}_{z}(w),
\end{align}where the last equality follows by the conjugate symmetric property of the Bergman kernel (that is, $\overline{k^{t}(z,w)}=k^{t}(w,z)$). It follows that $\sharp_{A^2}\xi^{\alpha}_{z}=\dbar^{\alpha}k^{t}_{z}$, and hence by Corollary \ref{cor: generaization of smoothness of Bergman kernels} that $\dbar^{\alpha}k^{t}_{z}$ is smooth in $t$. The case $\alpha=0$ gives back Lemma \ref{intro: smoothness of the Bergman kernel}. 
\end{example}

\section{Positivity of Hilbert bundles and plurisubharmonic variation}\label{sec: Positivity and Variation}

\subsection{Connections and curvatures on holomorphic Hilbert bundles}
In the previous section we showed how to carry out step (1) (in proving Theorem \ref{intro: variation of Bergman kernels}) in the introduction in the (more) general setting of (trivial) Hilbert bundles; recall Corollary \ref{cor: generaization of smoothness of Bergman kernels}, which is the natural counterpart of Lemma \ref{intro: smoothness of the Bergman kernel} from the introduction in this setting. In this section we will explore step (2) in the more general setting. The idea is to use what we call \ita{Chern connections} (Definition \ref{def: Chern connections on holomorphic Hilbert bundles}) to study, or rather compute, the variations of (squared norms of) smooth sections. This may not make too much sense at the moment; after all, we have yet to elaborate on the terminology. Hopefully, things will make more sense as we go along. Essentially, we are going to use calculus, in much the same way as we have been doing up till now, to compute $\partial\dbar \norm{u}_{h}^2$, where $u$ is a smooth local section of some (trivial) Hilbert bundle equipped with a smooth hermitian metric $h$. 

To introduce the notion of \ita{Chern connections} (and, later, also \ita{Chern curvatures}) on (what we call holomorphic) Hilbert bundles, we first extend the notion of local sections of Hilbert bundles to what we may call \ita{Hilbert-bundle-valued differential forms}:

\begin{definition}
Let $\fancy{H}$ be a (not necessarily trivial) Hilbert bundle over $B$. Then we define a(n) \bt{(local) $\fancy{H}$-valued differential $k$-form on (an open subset $U$ of) $B$} to be an assignment to each $b\in (U\sub)B$ an $\fancy{H}_b$-valued alternating complex-$k$-multilinear map on $T_{b}B\otimes \set{C}$, the (complexification of the) tangent space at $b$ of $B$.
\end{definition}
Let $\fancy{H}$ be as in the definition, and let $U\sub B$ be open. Locally, with respect to local holomorphic coordinates, say $t:=(t_1,\ldots, t_{m})$, on $U$, defined say on $V\sub U$, an $\fancy{H}$-valued differential form $u$ on $U$ takes the form
\begin{align}
u&=\s{|I|=p,|J|=q}{'}u_{I\bar{J}}dt_{I}\wedge d\bar{t}_{J},\label{eq: identity for Hilbert bundle valued diff-forms, locally}
\end{align}where $I,J$ are multiindices of respective lengths $p$ and $q$ such that $p+q=k$, $dt_{I}=dt_{i_1}\wedge\cdots\wedge dt_{i_p}, d\bar{t}_{J}=d\bar{t}_{j_1}\wedge\cdots\wedge d\bar{t}_{j_q}$, the symbol "$'$" in the summation sign designates that we sum over strictly ascending (increasing) multiindices, and $u_{I\bar{J}}$ are local sections of $\fancy{H}$ over $V$. If we write $u$ above as 
\begin{align}
\s{p+q=k}{}u^{p,q},
\end{align} we shall say that each $u^{p,q}$ of \bt{bidegree} or \bt{type} $(p,q)$, and we shall call $k$ \bt{the (total) degree of $u$}. In the case that $\fancy{H}$ admits a smooth structure, $u$ above is smooth if each coefficient  $u_{I\bar{J}}$ is. We will denote the family of smooth $\fancy{H}$-valued differential $k$-forms, respectively $(p,q)$-forms, on $U$, by $\fancy{C}^{\infty}(U;\Ld^{k}(T^*B)\otimes \fancy{H})$, respectively $\fancy{C}^{\infty}(U;\Ld^{(p,q)}(T^*B)\otimes \fancy{H})$; in the global case that $U=B$, we simply write $\fancy{C}^{\infty}(\Ld^{k}(T^*B)\otimes \fancy{H})$, respectively $\fancy{C}^{\infty}(\Ld^{(p,q)}(T^*B)\otimes \fancy{H})$. Let now $h$ be a hermitian metric on $\fancy{H}$. We extend $h$ to a real-bilinear (and not quite sesquilinear) form acting on pairs of (local) $\fancy{H}$-valued differential forms, possibly of different degrees. We denote the extension by the same symbol $h$ as the metric, and define it as follows: Consider first pairs of local $\fancy{H}$-valued differential forms $u$ and $v$ of the special forms $u=\hat{u}\otimes \xi, v=\hat{v}\otimes \eta$, where $\hat{u}$ and $\hat{v}$ are local sections of $\fancy{H}$, and $\xi$ and $\eta$ are (local) differential forms, possibly of different degrees. For $\fancy{H}$-valued differential forms of this special form, we let $h$ be defined by
\begin{align}
h(u,v)&:=h(\hat{u},\hat{v})\xi\wedge \overline{\eta}.
\end{align}We then define $h$ in general from this by invoking real-bilinearity and appealing to \eqref{eq: identity for Hilbert bundle valued diff-forms, locally}. We are now ready to define \ita{Chern connections} on \ita{holomorphic} Hilbert bundles (over $B$); let us agree to call $\fancy{H}$ a \bt{holomorphic Hilbert bundle} if it admits a holomorphic structure such that there is a well-defined $\dbar$ operator on $B$ acting on smooth local sections of $\fancy{H}$ which vanishes precisely on the holomorphic ones:

\begin{definition}\label{def: Chern connections on holomorphic Hilbert bundles}
Let $\fancy{H}$ be a holomorphic Hilbert bundle over $B$, and let $h$ be a smooth hermitian metric on $\fancy{H}$. A \bt{Chern connection on $\fancy{H}$ with respect to $h$}, or simply \bt{a Chern connection for/of $h$}, is a complex-linear map $D$ with the following properties:
\begin{enumerate}[(i)]
\item $D$ acts on an element $u\in \fancy{C}^{\infty}(U;\fancy{H})$ with value $Du\in \fancy{C}^{\infty}(U;\Ld^{1}(T^*B)\otimes \fancy{H})$, where $U$ is an open subset of $B$.
\item $D$ splits according to type or bidegree of image points as 
\begin{align}
D=\delta+\dbar,
\end{align} where $\delta u\in \fancy{C}^{\infty}(U;\Ld^{(1,0)}(T^*B)\otimes \fancy{H})$ and $\dbar u\in \fancy{C}^{\infty}(U;\Ld^{(0,1)}(T^*B)\otimes \fancy{H})$, for all $u\in \fancy{C}^{\infty}(U;\fancy{H})$, and $U\sub B$ open.
\item For all smooth local sections $u$ and $v$ of $\fancy{H}$, 
\begin{align}
d h(u,v)&=h(Du,v)+h(u,Dv).
\end{align}
We shall refer to $\delta$ as \bt{the (1,0)-part of $D$} and to $\dbar$ as \bt{the (0,1)-part of $D$}. Often we will also refer to (iii) as \bt{(the) metric compatibility of $D$}.
\end{enumerate}
\end{definition}
If there exists a Chern connection $D$ on $\fancy{H}$ with respect to $h$, where both are as in Definition \ref{def: Chern connections on holomorphic Hilbert bundles}, then we shall also say that $\fancy{H}$ \bt{admits} $D$ as Chern connection (with respect to $h$). In the case of holomorphic hermitian vector bundles of finite rank, it is well-known that the Chern connection always exists uniquely. In general, it is not clear to us whether we have existence, but we may ensure uniqueness by adding an additional assumption on $\fancy{H}$ (which immediately holds in the case that $\fancy{H}$ is trivial):

\begin{proposition}\label{prop: uniqueness of Chern connection}
Let $\fancy{H}$ be a holomorphic Hilbert bundle over $B$, and $h$ a smooth hermitian metric on $\fancy{H}$. Suppose that the following holds: Given $b\in \fancy{B}$ and $a_{b}\in \fancy{H}_b$, there exists a smooth (local) section $v$ of $\fancy{H}$ such that $v^{b}=a_{b}$. 
Then $\fancy{H}$ admits at most one Chern connection with respect to $h$.
\end{proposition}

\begin{proof}
Suppose that $D=\delta+\dbar$ and $D'=\delta'+\dbar$ are two Chern connections for $h$, and let $u$ be a smooth (local) section of $\fancy{H}$, say over $U$. Fix $b\in U$ and $a_b\in \fancy{H}_b$. It suffices by the non-degeneracy of $h_{b}$ to show that $h_{b}(((\delta-\delta)'u)^{b},a_{b})=0$. By assumption there exists a smooth (local) section $v$ (near $b$) of $\fancy{H}$ such that $v^{b}=a_{b}$. Using the metric compatibility of each of the connections and subtracting, we get 
\begin{align}
0&=\partial h(u,v)-\partial h(u,v)=h((\delta-\delta')u,v).
\end{align}Evaluating this at $b$ gives then 
\begin{align}
0&=h_{b}(((\delta-\delta')u)^{b},a_{b}),
\end{align}which proves the assertion.
\end{proof}
While it is not clear to us whether Chern connections exist in general, there is an important case in which the existence of Chern conenction can be established, namely when our Hilbert bundle is trivial. Let $H$ be a trivial Hilbert bundle over $B$. Then $H$ is automatically holomorphic. In the previous section, we discussed two particular types of hermitian metrics on trivial Hilbert bundles. Let us first consider the simplest of these two: Those of zero variation. Observe that an equivalent way of stating Lemma \ref{lemma: characterization of hermitian metrics of zero variation} is that a hermitian metric $h$ on $H$ is of zero variation if and only if it is smooth and admits the trivial Chern connection $d$. Hence it follows that the Chern connection on $H$ with respect to any hermitian metric of zero variation exists and is equal to the trivial connection $d$. Since the connection is trivial, this does not seem to give much, but when this is combined with the next result, we actually get a way of manufacturing a lot of smooth hermitian metrics admitting non-trivial Chern connections. Let us first give the result and then explain this a bit more closely. The result concerns the second type of hermitian metrics with which we have been concerned in the previous section, namely, the weighted ones, and is the following lemma:

\begin{lemma}\label{lemma: Chern connection from weighted hermitian metrics}
Let $H$ be a trivial Hilbert bundle over $B$, and suppose that $h$ is a smooth hermitian metric on $H$ that admits the Chern connection $D:=\delta+\dbar$. Assume that $\mathfrak{h}$ is a smooth weighted hermitian metric on $H$ induced by $h$ with smooth weight $w$. Then $H$ admits the Chern connection $\mathfrak{D}:=\mathfrak{d}+\dbar$ with respect to $\mathfrak{h}$, where \begin{align}
\mathfrak{d}=w^{-1}\circ \delta \circ w:=w^{-1}\delta w.
\end{align}
\end{lemma}
\begin{proof}
Let $u$ and $v$ be smooth local sections of $\fancy{H}$. By definition of $\mathfrak{h}$ and metric compatibility of $D$, we have
\begin{align}
\partial \mathfrak{h}(u,v)&=\partial h(wu,v)=h(\delta wu,v)+h(wu,\dbar v)=\mathfrak{h}(w^{-1}\delta w u,v)+\mathfrak{h}(u,\dbar v).
\end{align}This shows that $\mathfrak{D}$ satisfies metric compatibility with $\mathfrak{h}$, so the assertion follows by uniqueness (Proposition \ref{prop: uniqueness of Chern connection}).
\end{proof}
Using Lemma \ref{lemma: Chern connection from weighted hermitian metrics}, we may (at least in theory) manufacture numerous smooth hermitian metrics on $H$ admitting non-trivial Chern connections. The process is as follows: We begin with a smooth hermitian metric on $H$ that admits Chern connection.  We then consider a new smooth hermitian metric on $H$ given as a smooth weighted hermitian metric induced from the former metric. By the lemma, this new metric also admits Chern connection, and there is even an explicit  relationship between the ((1,0)-part of the) new and old Chern connections given by twisting with the smooth weight and its inverse.  We now repeat the above process with a new smooth weight to get yet another smooth hermitian metric that admits Chern connection, and so on. Of course, this process would be void if we cannot from the outset find a smooth hermitian metric on $H$ that admits Chern connection. This is where hermitian metrics of zero variation enter. We know that these always admit the trivial connection, so we may commence the above process using these. 

\begin{example}\label{ex: Chern connection example trivial Hilbert bundles}
Let $H,h^0$, and $h$ be as in Example \ref{ex: example of trivial hilbert bundles with weighted hermitian metrics}. Then $H$ is a trivial Hilbert bundle over $B$, and $h^0$ is a (smooth) hermitian metric of zero variation on $H$. The metric $h$ is a smooth weighted hermitian metric on $H$ induced by $h^0$ with smooth weight $w=e^{-\psi}$. By Lemma \ref{lemma: Chern connection from weighted hermitian metrics}, the Chern connection on $H$ with respect to $h$ exists and is given by $D:=\delta_{t}+\dbar_{t}$ (where we as usual use $t$-subscripts to emphasize that we are considering differential operators on $B$ and not $X$, taking $t$ as generic local holomorphic coordinates on $B$), (and) where
\begin{align}
\delta_{t}&=w^{-1}\partial_t w=e^{\psi}\partial_{t} (e^{-\psi})=\partial_{t}-\partial_{t}\psi\wedge;
\end{align}a formula that might look familiar to many.
\end{example}
Suppose that $h$ is a smooth hermitian metric on $H$ that admits the Chern connection $D:=\delta+\dbar$. Suppose also that $H_0\leq H$ and that $\overline{H_0}^{(H,h)}=H_0$, and let $P=P(h)$ be the total orthogonal projection map associated with $h$. It is natural to ask whether $H_0$ admits a (necessarily unique) Chern connection with respect to $h$. Let us suppose that the Chern connection does exist and write it as $D_0:=\delta_0+\dbar$. Let $u$ and $v$ be smooth (local) sections of $H_0$. By metric compatibility of $D_0$, we have
\begin{align}
\partial h(u,v)&=h(\delta_0 u, v)+h(u,\dbar v).
\end{align}On the other, by metric compatibility of $D$ we also have
\begin{align}
\partial h(u,v)&=h(\delta u,v)+h(u,\dbar v)=h(P\delta u,v)+h(u,\dbar v);
\end{align}the last equality following by orthogonality. Comparing, we then get that 
\begin{align}
\delta_0 u=P\delta u,
\end{align} which gives that 
\begin{align}
\delta_0=P\delta.
\end{align} Note that we have not assumed any smoothness on $P$ from the outset. Nevertheless, we see that the existence of $\delta_0$ forces $P\delta$ to be smooth (when restricted to sections of $H_0$). Conversely, we see that if $P\delta$ is smooth and we define $\delta_0:=P\delta$, then $D_0:=\delta_0+\dbar$ must be the Chern connection of $H_0$. Our discussion therefore gives the following result:

\begin{proposition}\label{prop: Chern connection for subbundles, trivial case}
With the above notation and set-up, the Chern connection of $H_0$ with respect to $h$ exists if and only if $P\delta$ (restricted to sections of $H_0$) is smooth, in which case it is given by $D_0=\delta_0+\dbar$, with $\delta_0=P\delta$.
\end{proposition}

In particular, by Berndtsson's regularity theorem (Theorem \ref{thm: Berndtsson's regularity theorem, general version}) and our discussion above (in which Lemma \ref{lemma: Chern connection from weighted hermitian metrics} plays an important role), it follows that if $h$ is a smooth weighted hermitian metric on $H$ induced by a hermitian metric of zero variation, then the Chern connection on $H_0$ with respect to $h$ exists (and we also have a formula for it). Let $h^{-1}$ be the dual hermitian metric of $h$ on the dual bundle $^*$. We may similarly ask whether $D$ induces a (necessarily unique) Chern connection on $H^*$ with respect to $h^{-1}$. From the outset, this is already more subtle than the case of the subbundle. Indeed, it is not even clear whether $h^{-1}$ is smooth, so the question may not even make sense without additional assumptions. Adding the assumptions that the associated adjoint operator and its inverse be smooth, we get the following result:

\begin{proposition}\label{prop: Chern connection for dual bundle, trivial case}
With the set-up and notation above, let $\sharp$ denote the adjoint operator with respect to $h$,  and $\sharp^{-1}$ its inverse. Assume that both $\sharp$ and $\sharp^{-1}$ are smooth. Then the Chern connection for $H^*$ with respect to $h^{-1}$ exists and is given by
\begin{align}
D^{\vee}&=\delta^{\vee}+\dbar^{\vee},
\end{align}where 
\begin{align}
\delta^{\vee}=\sharp^{-1}\dbar\sharp
\end{align} and 
\begin{align}
\dbar^{\vee}=\sharp^{-1}\delta\sharp.
\end{align}
\end{proposition}

\begin{proof}
Let $u^*$ and $v^*$ be smooth (local) sections of $H^*$. By metric compatibility of $D$ and the definition of $h^{-1}$,
\begin{align}
\partial h^{-1}(u^*,v^*)&=h^{-1}(\sharp^{-1}\dbar \sharp u^*,v^*)+h^{-1}(u^*,\sharp^{-1}\delta \sharp v^*).
\end{align}This shows by our assumptions that $D^{\vee}$ as defined satisfies metric compatibility with $h^{-1}$. To complete the proof it remains only to check that $\dbar u^*=0$ precisely when $ \dbar^{\vee}u^*=0$. In fact, we actually have 
\begin{align}
\dbar =\sharp^{-1}\delta \sharp.\label{eq: dbar identity for adjoint operators}
\end{align} To see this, let $u^*$ be smooth and let $v\in H$. By the chain rule, 
\begin{align}
\dbar \ip{u^*}{v}=\ip{\dbar u^*}{v}.
\end{align} Since, by definition of $\sharp$, the latter is also equal to $h(v,\sharp \dbar^*u)$, we get that 
\begin{align}
\dbar \ip{u^*}{v}=h(v,\sharp \dbar u^*).
\end{align} On the other hand, we also have by metric compatibility of $D$, and the definition of $\sharp$ again (applied to $u^*$ now instead of $\dbar u^*$),  that 
\begin{align}
\dbar \ip{u^*}{v}=\dbar h(v,\sharp u^*)=h(v,\delta \sharp u^*).
\end{align} Hence we have 
\begin{align}
h(v,\sharp \dbar u^*)=h(v,\delta \sharp u^*),
\end{align} which gives $\sharp \dbar u^*=\sharp \delta u^*$, and hence \eqref{eq: dbar identity for adjoint operators}.
\end{proof}
Note that in the case that $h$ is of zero variation, we recover Corollary \ref{cor: smoothness and holomorphicity of dual bundles in terms of hermitian metrics of zero variation with additional property} part (ii) from Proposition \ref{prop: Chern connection for dual bundle, trivial case}. The proposition therefore gives a more general description of holomorphic sections of dual bundles: In the case of hermitian metrics of zero variation, a smooth (local) section $u^*$ of the dual bundle is holomorphic precisely when the the adjoint of $u^*$ (assuming additionally that the adjoint operator and its inverse are smooth) with respect to the metric vanishes under the action of the $\partial$-operator. If $h$ is on the other hand more generally only assumed to be smooth with a (necessarily unique) Chern connection, then $u^*$ is holomorphic precisely when its adjoint with respect to $h$  (still assuming the adjoint operator and its inverse are smooth) vanishes under the action of the $(1,0)$-part of the Chern connection. This is indeed more general since the $(1,0)$-part of the trivial connection, which is the Chern connection of any hermitian metric of zero variation, is precisely equal to $\partial$.

We have  looked at the Chern connections on $H_0$ and $H^*$, with respect to the natural hermitian metrics (induced by $h$) on these. Let us lastly discuss, in a similar fashion, the case of the quotient bundle. Thus, let $Q:=H/H_0$ be the quotient bundle of $H_0$ in $H$. We ask whether the Chern connection of $Q$ with respect to $h^{Q}$, the quotient hermitian metric on $Q$ induced by $h$, exists, and have then the following ''quotient analogue'' of Proposition \ref{prop: Chern connection for dual bundle, trivial case}:

\begin{proposition}\label{prop: Chern connection for quotient bundle, trivial case}
With the above set-up and notation, let $P=P(h)$ denote the total orthogonal projection map associated with $h$, and let $m_l$ denote the minimal lifting operator with respect to $h$. Assume that $P$ is smooth. Then the Chern connection of $Q$ with respect to $h^{Q}$ exists and is given by 
\begin{align}
D^{Q}&=\delta^{Q}+\dbar^{Q},
\end{align}where 
\begin{align}
\delta^{Q}=m_{l}^{-1}\delta m_{l}
\end{align}
 and 
 \begin{align}
 \dbar^{Q}=m_{l}^{-1}\dbar^{\perp}m_{l},
 \end{align} and we use the notation $\dbar^{\perp}=(1-P)\dbar$.
\end{proposition}

\begin{proof}
Recall that $P$ is smooth if and only if $m_l$ is smooth. Thus, $m_{l}$ is smooth under our assumptions. Let $u^{\perp}_0$ be a smooth (local) section of $H_0^{\perp_{h}}$, and let $\pi_{Q}:H\to Q$ denote the canonical projection map. Then $[u]:=\pi_{Q}\circ u^{\perp}_0$ is smooth and satisfies $m_{l}[u]=u_0^{\perp}$. Thus, $m_{l}^{-1}u_0^{\perp}=[u]$ is smooth. Therefore $m_{l}^{-1}$ is also smooth. We claim that $\dbar^{Q}=\dbar$ as an operator acting on smooth local sections of $Q$. Let $[u]$ be a smooth (local) section of $Q$. It suffices to show that $[u]$ is holomorphic if and only if $\dbar^{Q}[u]=0$. Let $h^0$ denote the hermitian metric of zero variation induced by the fixed inner product on $H$, and let $m_l^0$ denote the minimal lifting operator with respect to $h^0$. By Proposition \ref{prop: characterization of smooth/holomorphic structure of trivial quotient bundle in terms of trivial minimal lifting operator}, $[u]$ is holomorphic if and only if $m_{l}^{0}[u]$ is holomorphic, which by Corollary \ref{cor: characterization of holomorphic sections in the trivial case using Cauchy-Riemann equation} is if and only if $\dbar m_{l}^{0}[u]=0$. By orthogonal decomposition we may write 
\begin{align}
m_{l}^{0}[u]&=Pm_{l}^{0}[u]\oplus m_{l}(m_{l}^{0}[u])=Pm_{l}^{0}\oplus m_{l}[u],
\end{align}so we find that $[u]$ is holomorphic if and only if
\begin{align}
-\dbar P(m_{l}^{0}[u])&=\dbar m_{l}[u].
\end{align}The key observation now is that the left-hand side of this identity is $H_0$-valued. It therefore follows that the identity holds if and only if the $H^{\perp_{h}}$-valued part of the right-hand side vanishes, which is precisely means $\dbar^{Q}[u]=0$. To complete the proof it remains to explain that $\delta m_{l}$ lies in the image of $m_{l}$ so that we can apply $m_{l}^{-1}$ to it; that is, $\delta^{Q}$ is well-defined, and verify that $D^{Q}$ as defined satisfies metric compatibility with $h^{Q}$. Indeed, if we can do this, then the assertion follows by uniqueness (Proposition \ref{prop: uniqueness of Chern connection}). To see that $\delta m_{l}$ lies in the image of $m_{l}$, it suffices to check that $h(\delta m_{l}[u], v_0)=0$ for all $v_0\in H_0$. By orthogonality, $h(m_{l}[u],v_0)=0$ for all $v_0\in H_0$. Applying $\partial$ to this, and using the the definition of $D$, we get that $0=h(\delta m_{l}[u],v_0)=0$ for all $v_0\in H_0$ as required. It remains to verify metric compatibility, but this follows immediately from the definition of $h^{Q}, D$,$D^{Q}$, and orthogonality ($h(m_{l}[u], \dbar m_{l}[v])=h(m_{l}[u],\dbar^{\perp}m_{l}[v])$).
\end{proof}
Analogous to how Proposition \ref{prop: Chern connection for dual bundle, trivial case} may be seen to generalize Corollary \ref{cor: smoothness and holomorphicity of dual bundles in terms of hermitian metrics of zero variation with additional property}, so can Proposition \ref{prop: Chern connection for quotient bundle, trivial case} be seen to generalize Proposition \ref{prop: characterization of smooth/holomorphic structure of trivial quotient bundle in terms of trivial minimal lifting operator}. Indeed, note that 
\begin{align}
\dbar m_{l}[u]=\dbar^{\perp}m_{l}[u]
\end{align} for all smooth (local) sections $[u]$ of $Q$ in the case that $h$ is of zero variation.

Let $\fancy{H}$ be a holomorphic Hilbert bundle over $B$, and suppose that $h$ is a smooth hermitian metric on $\fancy{H}$ that admits the Chern connection $D$. We shall next define what we call the \ita{Chern curvature of $D$}. To to do this, we first extend $D$ from acting on smooth local sections of $\fancy{H}$ to acting on smooth local $\fancy{H}$-valued differential forms of total degree greater than 0; we denote as usual the extension using the same symbol $D$ as the connection. Let $f$ be a smooth (local) complex-valued function on (some open subset of) $B$. By metric compatibility, if $u$ is a smooth (local) section of $\fancy{H}$, we have that
\begin{align}
D(fu)&=df\wedge u+fDu.\label{eq: Leibniz rule}
\end{align}We shall often refer to \eqref{eq: Leibniz rule} as \ita{the Leibniz rule for $D$}; for obvious reasons. To extend $D$, we first define it on smooth (local) $\fancy{H}$-valued differential forms of the form $u=\omega\wedge \hat{u}$, where $\omega$ is a smooth (local) differential form of degree $\text{deg}(\omega)$, and $\hat{u}$ is a  smooth (local) section of $\fancy{H}$, by
\begin{align}
D(\omega\wedge \hat{u})&:=d\omega\wedge \hat{u}+(-1)^{\text{deg}(\omega)}\omega\wedge D\hat{u}.\label{eq: extended Leibniz rule}
\end{align}
We then define $D$ in general from this by requiring it to be complex-linear. The extension of $D$ satisfies ''an extended'' metric compatibility of the form:
\begin{align}
D h(u,v)&=h(Du,v)+(-1)^{\text{deg}(u)}h(u,Dv),\label{eq: extended metric compatibility}
\end{align}where $u$ and $v$ are smooth (local) $\fancy{H}$-valued differential forms, possibly of different degrees, and the total degree of $u$ is $\text{deg}(u)$. We may now define (the) \ita{Chern curvature} (of $D$) as follows:

\begin{definition}\label{def: Chern curvature}
Suppose that $D$ is a Chern connection. Then we define its \bt{Chern curvature} to be $D^2:\equiv D\circ D$.
\end{definition}
We shall generically denote the Chern curvature of $D$ by $\Theta$, or $\Theta^{D}$ if we want to emphasize the connection. In the case that $\fancy{H}$ is a finitely-ranked vector bundle it is well-known and easy to see that $\Theta$ is a zeroth-order differential operator, a tensor. The same is true, and just as easy to see, also in the our current setting. That is, if $f$ is a smooth (local) complex-valued function and $u$ is a smooth (local) section of $\fancy{H}$, then 
\begin{align}
\Theta (fu)&=D(df\wedge u+fDu)=d^2f\wedge u-df\wedge Du+df\wedge Du+f\Theta u=f\Theta u.
\end{align}
Since $D$ is a differential operator of order 1 by the \eqref{eq: Leibniz rule}, one might expect that the curvature is a second-order differential operator, but somewhat surprisingly it thus turns out that the curvature is a tensor. Let us write $D:=\delta+\dbar$ and denote its (Chern) curvature by $\Theta$. One may easily check that $\delta^2=0$ under the assumptions in Proposition \ref{prop: uniqueness of Chern connection}. As a result, $\Theta$ is in that setting always of type (1,1) and equal to $\delta\dbar+\dbar\delta$.  More generally, it is true without any further assumptions, that for all smooth (local) sections $u$ and $v$ of $\fancy{H}$, we have

\begin{align}
h(\delta^2 u,v)&=0
\end{align}and
\begin{align}
h(\Theta u,v)&=h((\delta\dbar+\dbar \delta)u,v).
\end{align}

To see this we may use \eqref{eq: extended metric compatibility}. A similar argument, using instead $d^2$, shows that $\Theta$ in general also satisfies the following antisymmetric property:
\begin{align}
h(\Theta u,v)&=-h(u,\Theta v)
\end{align}for all smooth (local) sections $u$ and $v$ of $\fancy{H}$. We will express this property by writing 
\begin{align}
\Theta^*=-\Theta.
\end{align}

\subsection{The subbundle and quotient bundle curvature formulae}

Now, suppose that $\fancy{H}_0\leq \fancy{H}$, that $\overline{\fancy{H}_0}^{(\fancy{H},h)}=\fancy{H}_0$, and that $D_0:=\delta_0+\dbar$ is a Chern connection for $\fancy{H}_0$ with respect to $h$. Let $u$ and $v$ be smooth (local) sections of $\fancy{H}_0$. By metric compatibility of $D$ and $D_0$, we have
\begin{align}
0&=h((\delta-\delta_0)u,v).
\end{align}This shows that the operator $\delta-\delta_0=D-D_0$ satisfies a certain orthogonality property. We shall give this operator a special name:

\begin{definition}\label{def: second fundamental form}
Let $\fancy{H}$ be a holomorphic Hilbert bundle over $B$, and suppose that $h$ is a smooth hermitian metric on $\fancy{H}$ that admits a Chern connection $D$. Suppose also that $\fancy{H}_0\leq \fancy{H}$, that $\overline{\fancy{H}_0}^{(\fancy{H},h)}=\fancy{H}_0$, and that $\fancy{H}_0$ admits a Chern connection $D_0$ with respect to $h$. Then we define the \bt{second fundamental form associated with $D$ and $D_0$} to be $D-D_0$.
\end{definition}
For fixed connections $D$ and $D_0$ (such as above) we shall generically denote the second fundamental form associated with $D$ and $D_0$ by $S_{f}$. Using the orthogonality property of $S_f$, we may now prove the following generalization of a subbundle curvature formula of Griffiths's from the setting of finitely ranked holomorphic hermitian vector bundles (\cite{Griffiths});  henceforth, referred to as \ita{the subbundle curvature formula}:

\begin{theorem}[The subbundle curvature formula]\label{thm: the subbundle curvature formula}
Suppose that $D,D_0$ and $S_f$ are as above, and denote the Chern curvature of $D$ and $D_0$ by respectively $\Theta$ and $\Theta_0$. Then we have the following subbundle curvature formula:
\begin{align}
\Theta_0&=\Theta-S_f^*S_f,
\end{align}by which we shall mean more precisely the following: For all smooth (local) sections $u$ and $v$ of $\fancy{H}_0$, we have:
\begin{align}
h(\Theta_0 u,v)&=h(\Theta u,v)-h(S_f u, S_f v).
\end{align}
\end{theorem}

\begin{proof}
Let $u$ and $v$ be smooth (local) sections of $\fancy{H}_0$, and let us write $D:=\delta+\dbar$ and $D_0:=\delta_0+\dbar$. Note that while we actually do not know in the current setting whether it is true that $\delta^2_{(0)}=0$, and hence that $\Theta_{(0)}=\delta_{(0)}\dbar+\dbar \delta_{(0)}$, it \ita{is} true by our discussion above that 
\begin{align}
h(\delta_{(0)}^2u,v)&=0,
\end{align}and hence that
\begin{align}
h(\Theta_{(0)}u,v)&=h((\delta_{(0)}\dbar+\dbar \delta_{(0)})u,v).
\end{align}Thus we find; for simplicity we write $S_f$ as $S$,
\begin{align}
h((\Theta_0-\Theta)u,v)&=-h((S \dbar+\dbar S)u,v).
\end{align}By the orthogonality property of $S$, $h(S\dbar u,v)=0$, so we are left with
\begin{align}
h((\Theta_0-\Theta)u,v)&=-h(\dbar Su,v).
\end{align}By the orthogonality property of $S$ again, $h(S u,v)=0$. Thus applying $\dbar$ to this we get (recall \eqref{eq: extended Leibniz rule}) 
\begin{align}
h(\dbar S u,v)=h(Su, \delta v)=h(S u, S v),
\end{align} where the last equality follows by orthogonality. Substituting this we get the desired curvature formula.
\end{proof}
Note that in the case that $H$ is trivial, it follows by Proposition \ref{prop: Chern connection for subbundles, trivial case} that $D_0=PD$, where $P=P(h)$ is the total orthogonal projection map associated with $h$. Then it is immediate that the image of $S=D-D_0$ is $H_0^{\perp}$-valued. In the more general case we do not know that we have $D_0=PD$, but we still get an orthogonality property of $S$ and the subbundle curvature formula above nonetheless. Using a similar argument as in the proof of the subbundle curvature formula (Theorem \ref{thm: the subbundle curvature formula} above), we may also prove a corresponding formula for the curvature of the quotient bundle $Q:=H/H_0$ in the case that $\fancy{H}:=H$ and $\fancy{H}_0:=H_0$ are trivial, and $P$ is smooth; note that in this case, Proposition \ref{prop: Chern connection for quotient bundle, trivial case} gives an explicit formula for the Chern connection of $Q$ with respect to the quotient hermitian metric $h^{Q}$ on $Q$ induced by $h$. We shall refer to this curvature formula as the \ita{quotient bundle curvature formula}. It is given by the following theorem:

\begin{theorem}[The quotient bundle curvature formula]\label{thm: the quotient bundle curvature formula}
Let $H$ be a trivial Hilbert bundle over $B$, and suppose that $h$ is a smooth hermitian metric on $H$ that admits the Chern connection $D$. Let $H_0\leq H$ with $\overline{H_0}^{(H,h)}=H_0$, and suppose that $P=P(h)$, the total orthogonal projection map associated with $h$, is smooth. Let $Q:=H/H_0$, and let $D^{Q}$ denote the Chern connection of $Q$ with respect to $h^{Q}$, the quotient hermitian metric on $Q$ induced by $h$. Let $m_{l}$ denote the minimal lifting operator with respect to $h$, and let $\dbar_0:=P\dbar$. Denote the Chern curvature of $D$ and $D^Q$ by respectively $\Theta$ and $\Theta^{Q}$. Then the following quotient bundle curvature formula holds:
\begin{align}
m_{l}\Theta^{Q}&=\Theta m_{l}-(\dbar_0 m_{l})^*(\dbar_0 m_{l}),
\end{align}by which we shall mean more precisely the following: For all smooth (local) sections $[u]$ and $[v]$ of $Q$, we have:

\begin{align}
h^{Q}(\Theta^{Q}[u],[v])&=h(\Theta m_{l}[u], m_{l}[v])-h(\dbar_0 m_{l}[u], \dbar_0 m_{l}[v]).
\end{align}

\end{theorem}

\begin{proof}
Let $[u]$ and $[v]$ be smooth (local) sections of $Q$, and let us also write $D:=\delta+\dbar$. By Proposition \ref{prop: Chern connection for quotient bundle, trivial case}, we have $D^{Q}=\delta^{Q}+\dbar^{Q}$ with $\delta^{Q}=m_{l}^{-1}\delta m_{l}$ and $\dbar^{Q}=m_{l}^{-1}\dbar^{\perp}m_{l}$, where $\dbar^{\perp}=(1-P)\dbar$. From our discussion above we also have:

\begin{align}
h^{Q}(\Theta^{Q}[u], [v])&=h^{Q}((\delta^{Q}\dbar^{Q}+\dbar^{Q}\delta^{Q})[u], [v])
\end{align}and 
\begin{align}
h(\Theta m_{l}[u], m_{l}[v])&=h((\delta\dbar+\dbar\delta)m_{l}[u], m_{l}[v]).
\end{align}By definition of $h^{Q}$ and $D^{Q}$, we thus get 
\begin{align}
h^{Q}(\Theta^{Q}[u], [v])&=h((\delta \dbar^{\perp}+\dbar^{\perp}\delta) m_{l}[u], m_{l}[v]).
\end{align}Hence, 
\begin{align}
h^{Q}(\Theta ^{Q}[u], [v])-h(\Theta m_{l}[u], m_{l}[v])&=-h((\delta \dbar_0+\dbar_0 \delta)m_{l}[u], m_{l}[v]).
\end{align}By orthogonality, $h(\dbar_0 \delta m_{l}[u], m_{l}[v])=0$, so it remains to deal with $h(\delta \dbar_0 m_{l}[u], m_{l}[v])$. By orthogonality again, $h(\dbar_0 m_{l}[u], m_{l}[v])=0$. Applying $\partial$ to this, and using the definition of $D$, we get 

\begin{align}
h(\delta \dbar_0 m_{l}[u], m_{l}[v])&=h(\dbar_0 m_{l}[u], \dbar m_{l}[v])=(\dbar_0 m_{l}[u], \dbar_{0}m_{l}[v]),
\end{align}where the last equality follows by orthogonality. Substituting this we get the desired curvature formula.
\end{proof}
We have discussed a subbundle curvature formula and a quotient bundle formula. It therefore stands to reason to ask for a formula for the curvature of the dual bundle of $\fancy{H}$. We comment on the case, similar to that of the quotient bundle, that $\fancy{H}:=H$ is trivial. Let $\sharp$ denote the adjoint operator with respect to $h$, and suppose that $\sharp$ and $\sharp^{-1}$ are smooth. Then we have an explicit formula for the Chern connection $D^{\vee}$ of the dual bundle $H^*$ with respect to the dual hermitian metric $h^{-1}$ of $h$, by Proposition \ref{prop: Chern connection for dual bundle, trivial case}. Namely, if we write $D^{\vee}=\delta^{\vee}+\dbar^{\vee}$, then we have $\delta^{\vee}=\sharp^{-1}\dbar \sharp$ and $\dbar^{\vee}=\sharp^{-1}\delta \sharp$. Let $\Theta^{\vee}$ denote the Chern curvature of $D^{\vee}$, and $\Theta$ the Chern curvature of $D:=\delta+\dbar$, the Chern connection of $H$ with respect to $h$. By our discussion above, we have in this situation that
\begin{align}
\Theta^{\vee}&=\delta^{\vee}\dbar^{\vee}+\dbar^{\vee}\delta^{\vee}=\sharp^{-1}\Theta \sharp.
\end{align}Let $u^*$ and $v^*$ be two smooth (local) sections of $H^*$. Thus, by definition of $h^{-1}$ and the antisymmetric property $\Theta=-\Theta^*$ of the curvature $\Theta$,
\begin{align}
h^{-1}(\Theta^{\vee} u^*, v^*)=h(\sharp v^*,\Theta^{\vee}\sharp u^*)=-h(\Theta \sharp u^*, \sharp v^*).
\end{align}As a short-hand form we shall write this as
\begin{align}
\Theta^{\vee}&=-\Theta \sharp \label{eq: relation between dual and non-dual curvatures}
\end{align}This identity says that $\Theta^{\vee}$, the curvature for $H^*$, and $\Theta$, the curvature for $H$, in some sense are \ita{negative} of each other. The precise sense in which the two are negative of each other is part of our upcoming discussion(s).

\subsection{Positivity notions and variational formulae}
We are now ready to explore step (2) in proving Theorem \ref{intro: variation of Bergman kernels} from the introduction in our more general setting. Recall, using the same notation as in the introduction, that this is to compute $i\partial_{t}\dbar_{t}\log K(z)$ given that $K(z)$ is smooth (from step (1)). Looking back at Corollary \ref{cor: generaization of smoothness of Bergman kernels} and Example \ref{ex: derivatives of Bergman kernels}, the counterpart of this in our current setting is then the following: Let $\fancy{H}$ be a Hilbert bundle over $B$ with a smooth structure, and suppose that $h$ is a smooth hermitian metric on $\fancy{H}$. Assume also that $\fancy{H}_0\leq \fancy{H}$ with $\overline{\fancy{H}_0}^{(\fancy{H},h)}=\fancy{H}_0$, and let $\sharp_0$ denote the adjoint operator for $\fancy{H}_0^*$ with respect to $h$. Assuming that $\sharp_0$ is smooth, which, as we have seen, corresponds to $K(z)$ being smooth (recall Example \ref{ex: derivatives of Bergman kernels} and its prequel remarks), the counterpart to computing $\partial_{t}\dbar_{t}\log K(z)$ is now to compute $i\partial \dbar \log \norm{\sharp_0 u^*}_{h}^2$, say for $\sharp_0 u^*$ not vanishing. Discarding the logarithm (and the constant $i$), this motivates us to compute and study ''second-order derivatives'' of the squared norms of smooth (local) sections of $\fancy{H}$ of the form $\partial \dbar \norm{u}^2_{h}$. We shall think of and refer to these second-order derivatives as \ita{variations of second order}. The key observation here is to note that $\norm{u}^2_h=h(u,u)$, so these variations can be computed by means of Chern connections (presuming existence). Let $f$ be a two times differentiable function of real variables. We introduce the following notation: The first order variations $\fancy{V}^{1}_{f(h)}$ and $\fancy{V}^{\bar{1}}_{f(h)}$ are the operators acting on any smooth (local) section $u$ of $\fancy{H}$ given by:

\begin{align}
\fancy{V}^{1}_{f(h)}(u)&:=\partial f(\norm{u}^2_{h})\quad\quad \text{ and }\\
\fancy{V}^{\bar{1}}_{f(h)}(u)&:=\dbar f(\norm{u}^2_{h}),
\end{align}and the second order variations $\fancy{V}^{1,\bar{1}}_{f(h)}$ and $\fancy{V}^{\bar{1},1}_{f(h)}$ are the operators acting on any smooth (local) section $u$ of $\fancy{H}$ given by:
\begin{align}
\fancy{V}^{1,\bar{1}}_{f(h)}(u)&:=\partial \fancy{V}^{\bar{1}}_{f(h)}(u)\quad\quad \text{ and }\\
\fancy{V}^{\bar{1},1}_{f(h)}(u)&:=\dbar \fancy{V}^{1}_{f(h)}(u).
\end{align}
We simply write $f(h)$ as $h$ in the case that $f$ is the function 1. With the above notation it follows immediately that we have the second-order variational formula:
\begin{align}
\fancy{V}^{1,\bar{1}}_{f(h)}&=f''(\norm{\cdot}_{h}^2)\fancy{V}^{1}_{h}\wedge \fancy{V}^{\bar{1}}_{h}+f'(\norm{\cdot}^2_{h})\fancy{V}^{1,\bar{1}}_{h}.\label{eq: general variational formula for f}
\end{align}
By a direct computation we furthermore have the following basic, but important, second-order variational formulae; the rest of what we do in the paper will be based on these formulae (and on formula \ref{eq: general variational formula for f}):

\begin{lemma}\label{lemma: second order variational formula}
Let $\fancy{H}$ be a holomorphic Hilbert bundle over $B$, and suppose that $h$ is a smooth hermitian metric on $\fancy{H}$ that admits a Chern connection $D:=\delta+\dbar$. Let $\Theta$ denote the Chern curvature of $D$. Then we have the following second-order variational formula:
\begin{align}
\fancy{V}^{1,\bar{1}}_{h}&=\delta \dbar+\delta^*\delta-\dbar^*\dbar+(\dbar \delta)^*,
\end{align}by which we more precisely mean the following: For all smooth (local) sections $u$ of $\fancy{H}$, we have:
\begin{align}
\fancy{V}^{1,\bar{1}}_{h}(u)&=h(\delta \dbar u,u)+\norm{\delta u}^2_{h}-\norm{\dbar u}^2_{h}+h(u,\dbar \delta u).
\end{align}
In particular, we have:
\begin{enumerate}[(i)]
\item \begin{align}
\fancy{V}^{1,\bar{1}}_{h}|_{\text{ker}(\dbar)}=-\Theta+\delta^*\delta,
\end{align}by which we more precisely the following: For all holomorphic (local) sections $u$ of $\fancy{H}$, we have
\begin{align}
\fancy{V}^{1,\bar{1}}_{h}(u)&=-h(\Theta u,u)+\norm{\delta u}^2_{h}.
\end{align}
\item \begin{align}
\fancy{V}^{1,\bar{1}}_{\text{ker}(\delta)}&=\Theta-\dbar^*\dbar,
\end{align}by which we more precisely mean the following: For all smooth (local) sections $u$ of $\fancy{H}$ such that $\delta u=0$, we have
\begin{align}
\fancy{V}^{1,\bar{1}}_{h}(u)&=h(\Theta u, u)-\norm{\dbar u}^2_h.
\end{align}
\end{enumerate}
\end{lemma}
\begin{proof}
The first formula follows from \eqref{eq: extended Leibniz rule} and the definition of $D$. Then (ii) follows since $h(\Theta u,u)=h((\delta\dbar+\dbar \delta )u,u)$ for all smooth (local) sections $u$ of $\fancy{H}$. Finally, (i) follows in the same as (ii) by using that $\Theta=-\Theta^*$.
\end{proof}
The connection between complex Brunn-Minkowski theory and the above variational formulae in Lemma \ref{lemma: second order variational formula}, is that if $i\fancy{V}^{1,\bar{1}}_{h}(u)\geq 0$ as a (1,1)-form, then $\norm{u}^2_h$ is plurisubharmonic on the domain of definition of $u$. We also see from the above formulae that there is an intimate relationship between variations of smooth (local) sections $u$ of $\fancy{H}$, say such that $\delta u=0$, and the positivity of $ih(\Theta u,u)$. Indeed, from (ii) in Lemma \ref{lemma: second order variational formula}, it follows immediately that if $ih(\Theta u,u)\geq 0$ for such a $u$, then $\norm{u}^2_h$ is plurisubharmonic on its domain of definition. This motivates the following definition of \ita{Griffiths (semi)positivity of $\Theta$}:

\begin{definition}\label{def: Griffiths semipositivity}
Let $\fancy{H}$ be a holomorphic Hilbert bundle over $B$, suppose that $h$ is a smooth hermitian metric on $\fancy{H}$ that admits the Chern connection $D$, and let $\Theta$ denote the Chern curvature of $D$. Then we shall say that $\Theta$ is \bt{semipositive in the sense of Griffiths}, or \bt{Griffiths semipositive}, and write $i\Theta\geq_{G}0$, if
\begin{align}
h(i\Theta u,u)&\geq 0
\end{align}for all smooth (local) sections $u$ of $\fancy{H}$. We  shall also say that $\Theta$ is \bt{positive in the sense of Griffiths}, or \bt{Griffiths positive}, and write $i\Theta>_{G}0$, if the inequality above is strict.
\end{definition}
The Griffiths (semi)positivity of $\Theta$ (in Definition \ref{def: Griffiths semipositivity}) may also be described as follows: Let us choose $t:=(t_1,\ldots, t_m)$ as generic holomorphic local coordinates on $B$, and let us write $\Theta$ with respect to these locally as $\Theta=\s{j,k=1}{m}\Theta_{j\bar{k}}dt_{j}\wedge d\bar{t}_{j}$; $[\Theta_{j\bar{k}}]$ is a matrix whose entries at each point (in some open subset of $B$) are complex-linear maps that act on (local) smooth sections of $\fancy{H}$. Then $\Theta$ is Griffiths semipositive if at each point $b$ 
\begin{align}
\s{j,k=1}{m}h(\Theta_{j\bar{k}}u, u)(b)\xi_{j}\bar{\xi}_{k}\geq 0
\end{align}for all smooth (local) sections $u$ of $\fancy{H}$, and all $\xi:=(\xi_1,\ldots, \xi_m)\in \set{C}^{n}$; $\Theta$ is Griffiths positive if the above holds with equality only when $\xi=\vect{0}$ or $u(b)=0$. This description of Griffiths (semi)positivity leads to a stronger notion of (semi)positivity for $\Theta$, namely what we call \ita{Nakano (semi)positivity}:

\begin{definition}\label{def: Nakano semipositivity}
We shall say that $\Theta$ is \bt{semipositive in the sense of Nakano}, or \bt{Nakano semipositive}, and write $i\Theta \geq_{N}0$, if at each point $b\in B$, 
\begin{align}
\s{j,k=1}{m}h(\Theta_{j\bar{k}}u_{j}, u_{k})(b)\geq 0
\end{align}for all smooth (local) sections $u_1,\ldots, u_{m}$ of $\fancy{H}$. We shall also say that $\Theta$ is \bt{positive in the sense of Nakano}, or \bt{Nakano positive}, and write $i\Theta>_{N}0$, if the above inequality is strict unless $u_{j}(b)=0$ for all $j$.
\end{definition}
By choosing $u_{j}:=\xi_{j}u$ we see that Nakano (semi)positivity implies Griffiths (semi)positive. We define Griffiths and Nakano (semi)\ita{negativity} of $\Theta$ of course in a similar fashion, using the opposite inequalities of those above, and the notation $i\Theta\leq_{G}0$ (or $i\Theta<_{G}0$) and $i\Theta\leq_{N}0$ (or $i\Theta<_{N}0$). We may now interpret the earlier identity \eqref{eq: relation between dual and non-dual curvatures} as $\Theta$ and $\Theta^{\vee}$ being ''negative of each other'' in the precise sense that 
\begin{align}
i\Theta_{G}\geq 0 \iff i\Theta^{\vee}_{G}\leq 0.
\end{align}
What we have said earlier subsequent to Lemma \ref{lemma: second order variational formula} may be formulated using the notion of Griffiths (semi)positivity as the following result:

\begin{corollary}\label{cor: simple first cor of variational formulae}
Under the assumptions of Lemma \ref{lemma: second order variational formula}, and using the same notation as in there, suppose that $i\Theta_{G}\geq 0$. Then $\norm{u}^2_{h}$ is plurisubharmonic on its domain of definition for any smooth (local) section $u$ of $\fancy{H}$ such that $\delta u=0$.
\end{corollary}
As mentioned, the corollary follows directly from Lemma \ref{lemma: second order variational formula}. Note that the lemma, however, says more. Indeed, due to the term $-\norm{\dbar u}^2_{h}$, we may allow for some negativity of the curvature. Note also that if we \ita{first fix} $u$, then the variational formula shows that it is sufficient that $ih(\Theta u,u)-i\norm{\dbar u}^2_{h}\geq 0$ as a (1,1)-form for just the particular chosen $u$. A stronger statement, which is also somehow more in line with the variation of Bergman kernels, would be that even $\log\norm{u}^2_{h}$ is plurisubharmonic. This does indeed hold under the same assumptions, and can be seen from the following general result:

\begin{theorem}\label{thm: a general psh property of smooth local functions with delta u=0}
Under the assumptions of Lemma \ref{lemma: second order variational formula}, and using the same notation as in there and above, suppose that $i\Theta\geq_{G}0$ and that $f=f(x)$ is strictly increasing and concave. Assume moreover that for all $x>0$, $f$ satisfies that
\begin{align}
1+\fr{f''(x)}{f'(x)}x\geq 0.
\end{align}Let $\epsilon>0$. Then $f(\norm{u}^2_{h}+\epsilon)$ is plurisubharmonic for all smooth (local) sections $u$ of $\fancy{H}$ such that $\delta u=0$.
\end{theorem}

\begin{proof}
Let $u$ be a section such as asserted. For simplicity we put $\epsilon:=0$ assuming that $\norm{u}^2_{h}\neq 0$. Since $\delta u=0$, it follows by Lemma \ref{lemma: second order variational formula}, (ii), that $\fancy{V}^{1,\bar{1}}_{h}(u)=h(\Theta u, u)-\norm{\dbar u}^2_{h}$. By a direct computation (using the definition of $D$) we also get $\fancy{V}^{1}_{h}(u)=h(u,\dbar u)$ and $ \fancy{V}^{\bar{1}}_{h}(u)=h(\dbar u, u)$. Hence by \eqref{eq: general variational formula for f},
\begin{align}
\fancy{V}^{1,\bar{1}}_{f(h)}(u)&=f''(\norm{u}^2_{h})\abs{h(u,\dbar u)}^2+f'(\norm{u}^2_{h})\left(h(\Theta u, u)-\norm{\dbar u}^2_{h}\right),
\end{align} where we use the notation $\abs{h(u,\dbar u)}^2:=h(u,\dbar u)\wedge h(\dbar u,u)$. By the Cauchy-Schwarz inequality, we have that
\begin{align}
i\abs{h(u,\dbar u)}^2&\leq -i \norm{u}^2_{h}\norm{\dbar u}^2_{h},
\end{align}and since $f$ is concave, $f''\leq 0$, so it follows that 
\begin{align}
if''(\norm{u}^2_{h})\abs{h(u,\dbar u)}^2&\geq -i f''(\norm{u}^2_{h})\norm{u}^2_{h}\norm{\dbar u}^2_{h}.
\end{align}Substituting this into the above, we therefore have
\begin{align}
i\fancy{V}^{1,\bar{1}}_{f(h)}(u)&\geq i f'(\norm{u}^2_{h})h(\Theta u, u)-i\norm{\dbar u}^2_{h}\left(f'(\norm{u}^2_{h})+f''(\norm{u}^2_{h})\norm{u}^2_{h}\right).
\end{align}Since $f$ is assumed to be strictly increasing, $f'>0$, so we may divide by $f'(\norm{u}^2_{h})$. The assertion now follows from noting that $i\fancy{V}^{1,\bar{1}}_{f(h)}(u)$ gives that $f(\norm{u}^2_h)$ is plurisubharmonic, and by considering $x:=\norm{u}^2_h$.
\end{proof}
Note that if we let $x>0$ and $f(x):=\log(x)$, then $f$ is strictly increasing and concave (on the positive real line) with $1+\fr{f''(x)}{f'(x)}x=0$, so the result on logarithmic plurisubharmonicity does indeed follow from Theorem \ref{thm: a general psh property of smooth local functions with delta u=0}. Note also that we can get rid of $\epsilon$ altogether by considering $f$ above to be defined on $[0,\infty)$ and satisfy that $1+\fr{f''(x)}{f'(x)}x\geq 0$ there; of course, $f$ should still be strictly increasing concave and twice differentiable. The previous version is recovered by considering $f_{\epsilon}(x):=f(x+\epsilon)$ for each $\epsilon>0$.

Above we have used (ii) in Lemma \ref{lemma: second order variational formula}. Using (i) in Lemma \ref{lemma: second order variational formula} instead, gives a similar result for holomorphic (local) sections of $\fancy{H}$ to which we shall come back below. Note also that the actual proof of the theorem gives a more general result that what is stated: If we define $x:=\norm{u}^2_{h}+\epsilon$ and let $x$ be a function on some open subset $U$ of $B$ on which $u$ is defined, then $f(x)$ is plurisubharmonic on the set of all $b\in U$ such that $1+\fr{f''(x)}{f'(x)}x\geq 0$. We also see that regardless of the positivity assumption on $\Theta$, we in general have the inequality
\begin{align}
i\fr{\fancy{V}^{1,\bar{1}}_{f(h)}(u)}{f'(\norm{u}^2_{h}+\epsilon)}&\geq i h(\Theta u,u)-i\norm{\dbar u}^2_{h}\left(1+\fr{f''(\norm{u}^2_{h}+\epsilon)}{f'(\norm{u}^2_{h}+\epsilon)}\norm{u}^2_{h}\right).
\end{align}

To use Theorem \ref{thm: a general psh property of smooth local functions with delta u=0}, we need to exploit the condition $\delta u=0$. One scenario in which we can do this is as follows: With the same set-up as in Lemma \ref{lemma: second order variational formula}, consider the dual bundle $\fancy{H}^*$ of $\fancy{H}$, and assume that it is holomorphic. Let $u$ be a smooth (local) section of $\fancy{H}$, and define the (local) section $u^*$ of $\fancy{H}^*$ by
\begin{align}
\ip{u^*}{v}:=h(v,u)
\end{align}for all (local) sections $v$ of $\fancy{H}$. We then get the following lemma:

\begin{lemma}\label{lemma: lemma for delta u=0, general case}
With the above notation and set-up, assume furthermore that the following holds: If $h(v,\delta u)=0$ for all $v$ holomorphic, then $\delta u=0$. Then $\delta u=0$ if $u^*$ is holomorphic.
\end{lemma}

\begin{proof}
If $u^*$ and $v$ are holomorphic, $\ip{u^*}{v}$ is holomorphic by composition. Thus, $h(v,u)$ is holomorphic so $\dbar h(u,v)=h(u,\delta u)=0$ since $v$ is holomorphic, and by definition of $D$. By our assumptions this implies that $\delta u=0$.
\end{proof}
In particular the lemma holds when $H$ is trivial. Indeed, then $h(v,\delta u)=0$ for all \ita{constant} $v$ already gives $\delta u=0$. In that case we also automatically have $h(v,\delta u)=0$ for all $v$ holomorphic if $u^*$ is holomorphic.  Considering Example \ref{ex: derivatives of Bergman kernels} and replacing the Hilbert bundle above with a suitable \ita{subbundle} of it, we get the following generalization of, or counterpart to, the variation of Bergman kernels from the introduction (Theorem \ref{intro: variation of Bergman kernels}):

\begin{theorem}\label{thm:counterpart to variation of Bergman kernels in general setting of trivial Hilbert bundles}
Let $H$ be a trivial Hilbert bundle over $B$, and suppose that $h$ is a smooth hermitian metric on $H$ that admits the Chern connection $D=\delta+\dbar$. Let $h^{-1}$ denote the dual hermitian metric of $h$ for $H^*$, and $\sharp$ the adjoint operator for $H^*$ with respect to $h$. Suppose that $\sharp$ is smooth, that $H_0\leq H$ with $\overline{H_0}^{(H,h)}=H_0$, and that the total orthogonal projection map $P=P(h)$ associated with $h$ is smooth. Let $D_0$ denote the Chern connection on $H_0$ with respect to $h$, and $\Theta_0$ its Chern curvature. Suppose that $i\Theta_0\geq_{G} 0$. Then for all holomorphic local sections $u^*$ of $H_0$, and all strictly increasing and concave second-order differentiable functions $f$ on $[0,\infty)$, $f(\norm{u^*}_{h^{-1}}^2)$ is plurisubharmonic on the domain of definition of $u^*$ intersected with the set of all $b\in B$ such that writing $x=x(b)$, we have 
\begin{align}
1+\fr{f''(x)}{f'(x)}x\geq 0.
\end{align}In particular, $\log(\norm{u^*}^2_{h^{-1}})$ is purisubharmonic on the whole domain of definition of $u^*$. 
\end{theorem}

\begin{proof}
Let $u^*$ be a holomorphic (local) section of $H_0^*$, and let $\sharp_0$ denote the adjoint operator for $H_0^*$ with respect to $h$. Then if we let $u:=\sharp_0 u^*$, it follows that $\ip{u^*}{v}=h(v,u)$ for all (local) sections $v$ of $H_0^*$, and moreover that $\norm{u^*}^2_{h^{-1}}=\norm{u}^2_{h}$. Since $\sharp$ is assumed smooth, $\sharp_0$ is smooth by Proposition \ref{prop: relation between adjoints of sub and exterior bundles} and the assumption that $P$ is smooth. Hence $u$ is smooth. Let us write $D_0:=\delta_0+\dbar$ ($D_0$ exists by Proposition \ref{prop: Chern connection for subbundles, trivial case}, and we moreover have $\delta_0=P\delta$). The assertion now follows by applying Theorem \ref{thm: a general psh property of smooth local functions with delta u=0} to $H_0$ instead of $H$, since $\delta_0 u=0$ as $u^*$ is holomorphic (if $u^*$ is holomorphic, then $\ip{u^*}{v}=h(v,u)$ is holomorphic whenever $v$ is holomorphic, by composition, so we get $\dbar h(v,u)=h(v,\delta_0 u)=0$ by definition of $D_0$).
\end{proof}
In particular, the conclusion of Theorem \ref{thm:counterpart to variation of Bergman kernels in general setting of trivial Hilbert bundles} holds by Berndtsson's regularity theorem (Theorem \ref{thm: Berndtsson's regularity theorem, general version}), under the assumptions that $h$ is a smooth weighted hermitian metric induced by a hermitian metric of zero variation $h^0$, that $H_0$ is closed in $H$ with respect to both these metrics, and that the adjoint operator with respect to $h^0$ is smooth (recall Proposition \ref{prop: relation between adjoint of metric and adjoint of weighted metrics}). 

We next consider the case of holomorphic (local) sections, that is, (i) in Lemma \ref{lemma: second order variational formula} above. An immediate analogue in this case to Theorem \ref{thm: a general psh property of smooth local functions with delta u=0} is the following:

\begin{theorem}\label{thm: a general variational formula and psh properties, general, for holomorphic sections}
Let $\fancy{H}$ be a holomorphic Hilbert bundle over $B$, and suppose that $h$ is a smooth hermitian metric on $\fancy{H}$ that admits a Chern connection $D:=\delta+\dbar$ with Chern curvature $\Theta$. Let $f=f(x)$ be a twice differentiable function on $[0,\infty)$ which is strictly increasing and concave, and let $u$ be a holomorphic (local) section of $\fancy{H}$. Then 
\begin{align}
i\fancy{V}^{1,\bar{1}}_{f(h)}(u)&\geq h(-i\Theta u,u)+\norm{\delta u}^2_{h}\left(f'(\norm{u}^2_h)+f''(\norm{u}^2_h)\norm{u}^2_h\right).
\end{align}In particular, if $i\Theta \leq_{G} 0$, then $f(\norm{u}^2_h)$ is plurisuharmonic where
\begin{align}
1+\fr{f''(\norm{u}^2_h)}{f'(\norm{u}^2_{h})}\norm{u}^2_{h}\geq 0.
\end{align}
\end{theorem}
\begin{proof}
The first assertion follows by \eqref{eq: general variational formula for f} and Lemma \ref{lemma: second order variational formula} (i) together with the Cauchy-Schwarz inequality (similar as in the proof of Theorem \ref{thm: a general psh property of smooth local functions with delta u=0} since $f$ is assumed concave). The second assertion follows then by dividing the first inequality by $f'(\norm{u}^2_h)>0$ and the definition of Griffith (semi)negativity of $\Theta$.
\end{proof}
Using the same notation and assumptions as in the theorem, Theorem \ref{thm: a general variational formula and psh properties, general, for holomorphic sections} says that if $i\Theta_{G}\leq 0$, then $f(\norm{u}^2_{h})$ is plurisubharmonic if $f$ is such that we always have $1+\fr{f''(x)}{f'(x)}x$ for $x>0$. We may ask whether the converse is true. That is, if $f(\norm{u}^2_{h})$ is plurisubharmonic for all holomorphic $u$ and such $f$'s, whether it follows that $i\Theta\leq_{G}0$. In the case that $\fancy{H}$ is trivial, this is indeed the case. To prove this, we shall have use for the following preliminary result, which may also be of independent interest in itself:

\begin{lemma}\label{lemma: preliminary result on prescribed holomorphic sections with vanishing (1,0)-parts of Chern connections, on trivial bundles}
Let $H$ be a trivial Hilbert bundle over $B$, and suppose that $h$ is a smooth hermitian metric on $H$ which admits the Chern connection $D:=\delta+\dbar$. Let $a\in H$ and $b\in B$. Then there exist an open neighbourhood $U$ of $b$ and a holomorphic local section $u$ of $H$ over $U$ such that $u^{b}=a$ and $\delta u(b)=0$.
\end{lemma}
The lemma says that we can find local holomorphic sections of $H$ with prescribed values and vanishing (1,0)-part of the Chern connection at a point.

\begin{proof}
Let $t:=(t_1,\ldots, t_m)$ be holomorphic local coordinates near $b$ such that $t(b)=0$. We let $U$ be the domain of definition of $t$. Let us locally, with respect to the $t$-coordinates, write $\delta:=\s{j=1}{m}dt_{j}\wedge \delta_{j}$. The idea is to choose $u$ to be of the form 
\begin{align}
u:=a+\s{j=1}{m}f_{j}\delta_{j}(a)
\end{align} where $f_{j}$ are holomorphic functions of $t$ to be chosen such that $u^{b}=a$ and $\delta u (b)=0$. Using \eqref{eq: Leibniz rule} and computing we find,
\begin{align}
\delta u&=\delta a+\s{j=1}{m}\rbrac{\partial f_{j}\wedge \delta_{j}(a)+f_{j}\delta(\delta_{j}(a))}=\s{k=1}{m}\delta_{k}(a)dt_{k}+\s{j=1}{m}\left(\s{k=1}{m}\del{f_{j}}{t_{k}}\delta_j(a)dt_{k}+f_{j}\delta (\delta_{j}(a))\right).
\end{align}Hence, we may take $f_{j}:=-t_{j}$.
\end{proof}
Using Lemma \ref{lemma: preliminary result on prescribed holomorphic sections with vanishing (1,0)-parts of Chern connections, on trivial bundles}, we now get the following characterization for Griffiths (semi)negativity of trivial Hilbert bundles:

\begin{theorem}\label{thm: characterization of Griffiths seminegativity of trivial Hilbert bundles}
Let $H$ be a trivial Hilbert bundle over $B$, and suppose that $h$ is a smooth hermitian metric on $H$ which admits the Chern connection $D$ with curvature $\Theta$. Let $f$ be a twice differentiable function on $[0,\infty)$ which is strictly increasing and concave, and which satisfies that
\begin{align}
1+\fr{f''(x)}{f'(x)}x\geq 0
\end{align}for all $x\geq 0$. Then $i\Theta\leq_{G}0$ if and only if $f(\norm{u}^2_{h})$ is plurisubharmonic for all (local) holomorphic sections $u$ of $H$.
\end{theorem}
\begin{proof}
Suppose that $i\Theta_{G}\leq 0$, and let $u$ a holomorphic (local) section of $H$. Then $f(\norm{u}^2_h)$ is plurisubharmonic by Theorem \ref{thm: a general variational formula and psh properties, general, for holomorphic sections}. Conversely, suppose that $f(\norm{u}^2_h)$ is plurisubharmonic for all holomorphic (local) sections $u$ of $H$. By (the first part of) Theorem \ref{thm: a general variational formula and psh properties, general, for holomorphic sections}, regardless of any assumption on $\Theta$, we have 
\begin{align}
i\fancy{V}^{1,\bar{1}}_{f(h)}(u)&\geq h(-i\Theta u,u)+i\norm{\delta u}^2_{h}\left(F(\norm{u}^2_{h})\right),\label{eq: an equation in the proof of characterization of Griffiths seminegativity}
\end{align}where $F$ is some function depending on (derivatives) of $f$. Fix $b\in B$ and $a\in H$. It suffices to show that $ih(\Theta a,a)(b)\leq 0$. We now apply Lemma \ref{lemma: preliminary result on prescribed holomorphic sections with vanishing (1,0)-parts of Chern connections, on trivial bundles}. By the lemma we can choose $u$ holomorphic near $b$ such that $u^{b}=a$ and $\delta u(b)=0$. Evaluating at $b$, \eqref{eq: an equation in the proof of characterization of Griffiths seminegativity} therefore gives $i\fancy{V}^{1,\bar{1}}_{f(h)}(u)(b)=h(-i\Theta u,u)(b)=h(-i\Theta a,a)(b)$, which by our assumption is no less than 0. It follows that $h(i\Theta a,a)(b)\leq 0$, which completes the proof.
\end{proof}
In particular, Theorem \ref{thm: characterization of Griffiths seminegativity of trivial Hilbert bundles} applies to $f(x):=\log(x+\epsilon)$ for any $\epsilon>0$. This recovers in our possibly infinitely-ranked setting the well-known characterization of Griffiths (semi)negativity of finitely-ranked holomorphic vector bundles asserting that the vector bundle is (semi)negative in the sense of Griffiths if and only if every non-vanishing holomorphic section is logarithmically plurisubharmonic. Also, by applying the theorem to the dual bundle of $H$, under appropriate conditions (for example, when the adjoint operator and its inverse, with respect to the metric $h$, are smooth), we get that $H$ is Griffiths (semi)positive (with respect to $h$) if and only if $f(\norm{u^*}_{h^{-1}}^2)$ is plurisubharmonic for all non-vanishing holomorphic (local) sections of the dual bundle, and where $h^{-1}$ of course is the dual metric of $h$. 

We conclude the paper with two novel results of similar sorts to Theorem \ref{thm:counterpart to variation of Bergman kernels in general setting of trivial Hilbert bundles}. We obtain these from our discussions above. The first is the following result on plurisubharmonic variational properties of (certain image points of) total orthogonal projection maps:

\begin{theorem}\label{thm: first novel result on psh properties of projection map}
Let $\fancy{H}$ be a holomorphic Hilbert bundle over $B$, and suppose that $h$ is a smooth hermitian metric on $H$ which admits the Chern connection $D:=\delta+\dbar$. Assume that $\fancy{H}_0\leq \fancy{H}$, that $\overline{H_0}^{(\fancy{H},h)}=\fancy{H}_0$, and that $P=P(h)$, the total orthogonal projection map associated with $h$ (of $\fancy{H}_0$ in $\fancy{H}$) is smooth. Suppose also that for $b\in B$, given any $a_{b,0}\in (\fancy{H}_{0})_b$, we can find a holomorphic (local) section $v_0$ of $\fancy{H}_0$ near $b\in B$, such that $(v_{0})^{b}=a_{b,0}$.
Let the Chern connection of $\fancy{H}_0$ be $D_0:=\delta_0+\dbar$, and let $\Theta$ and $\Theta_{0}$ denote respectively the Chern curvatures of $D$ and $D_0$. Let $f$ be a twice differentiable function on $[0,\infty)$ which is strictly increasing and concave. Let $u$ be a smooth (local) section of $\fancy{H}$ and suppose that $\delta_0 u=0$. Put $u_0:=Pu$. Then if

\begin{align}
i h(\Theta_0 u_0,u_0)-i\norm{\dbar u_0}^2_{h}\left(1+\fr{f''(\norm{u_0}^2_{h}}{f'(\norm{u_0}^2_{h}}\norm{u_0}^2_{h}\right)&\geq 0,
\end{align}
 $f(\norm{u_0}^2_{h})$ is plurisubharmonic.
In particular, if $i\Theta_0\geq_{G}0$, then $f(\norm{u_0}^2_{h})$ is plurisubharmonic wherever
\begin{align}
1+\fr{f''(\norm{u_0}^2_{h})}{f'(\norm{u_0}^2_{h}}\norm{u_0}^2_{h}\geq 0.
\end{align}
\end{theorem}

\begin{proof}
From our discussion above, and (the proof of) Theorem \ref{thm: a general psh property of smooth local functions with delta u=0} applied to $\fancy{H}_0$, it suffices to show that we have $\delta_0 u_0=0$. Our assumption is that $\delta_0 u=0$. Let $u_{m}:=u-u_0$. Then $u_{m}\perp_{h} v_0$ for all (local) sections $v_0$ of $\fancy{H}_0$, and $u_{m}$ is smooth. We claim that $\delta_0 u_{m}=0$. Let $b$ be in the domain of definition of $\delta_0 u_m$. It suffices to show that $\delta_0 u_{m}(b)=0$. Let $a_{b,0}\in (\fancy{H}_0)_b$ be arbitrary. It suffices to show that $h_{b}(\delta_0 u_{m}(b), a_{b,0})=0$. By our assumption we can find a holomorphic $v_0$ such that $v_0(b)=a_{b,0}$. By orthogonality we have $h(u_{m},v_0)=0$. Applying $\partial$ to this and using the definition of $D$, and orthogonality again, we get 
\begin{align}
0&=h(\delta_0 u_m, v_0).
\end{align}Evaluating at $b$ then gives our claim. Since $\delta_0 u_m=0$, it is clear that $\delta_0 u_0=0$ from the assumption that $\delta_0 u=0$. This completes the proof.
\end{proof}
The assumption that we can find holomorphic $v_0$'s such as in the statement of the theorem is rather strong. In the case that $\fancy{H}$ is trivial, however, it is immediately satisfied. Consider the case that $\fancy{H}:=H$ is trivial, and let $Q:=H/H_0$, where we write $H_0$ for $\fancy{H}_0$. Let $m_{l}$ denote the minimal lifting operator with respect to $h$, and $[u]$ be a smooth (local) section of $Q$. Assume that $u$ is a representative of $[u]$. Then $[u]=u-u_0$, where $u_0=Pu$, and we use the same notation as in Theorem \ref{thm: first novel result on psh properties of projection map}. The following interesting consequence of Theorem \ref{thm: first novel result on psh properties of projection map}, which is our second result, therefore gives in this setting variational properties of $[u]$:

\begin{corollary}\label{cor: second novel result on psh variations of minimal solutions}
With the same assumptions, and the same notation, as in Theorem \ref{thm: first novel result on psh properties of projection map}, suppose furthermore that $\norm{u}^2_{h}$ is pluri\ita{super}harmonic, and let $u_{m}:=u-u_0$. Then $-\norm{u_{m}}^2_{h}$ is plurisubharmonic.
\end{corollary}

\begin{proof}
By orthogonality we have $\norm{u_m}^2_{h}=\norm{u}^2_h-\norm{u_0}^2_h$, so $-\norm{u_m}^2_{h}=-\norm{u}^2_h+\norm{u_0}^2_{h}$. Hence the assertion follows from Theorem \ref{thm: first novel result on psh properties of projection map} taking $f(x):=x$.
\end{proof}

Let $T$ be a linear operator on $H$ and suppose that $v$ lies in the image of $T$, and that $u$ is a solution to $T(\cdot)=v$. Assume also that $\text{ker}(T):=H_0$ is a closed subspace of $H$. By \bt{the minimal solution} to the equation $T(\cdot)=v$ (with respect to $h$) we shall mean the solution $\hat{u}$ such that if $\tilde{u}$ is any other solution, we have $\norm{\hat{u}}^2_h\leq \norm{\tilde{u}}^2_h$. With the notation above, it follows that $u_{m}$ is the minimal solution to the equation $T(\cdot)=T(u)$. Thus, Corollary \ref{cor: second novel result on psh variations of minimal solutions} can also be seen to give plurisubharmonicity properties of (the negative of of the squared norms of) minimal solutions. 

\begin{example}
Let $H, H_0$, and $h$ be as in Example \ref{ex: Chern connection example trivial Hilbert bundles}. Then we know that $\delta =\partial_t-\partial_{t}\psi $. Let $u=u(z)$ be an element in $H$ which is smooth and suppose also that for all $z\in \supp{u}$, $\partial_{t}\psi(t,z)=0$. Then $\delta_{0}u=0$ and Theorem \ref{thm: first novel result on psh properties of projection map} applies. That is, (for example, if $h(i\Theta_0 u_0, u_0)\geq 0$, $\log\norm{u-u_{m}}^2_h=\log\norm{u_0}^2_{h}$ is plurisubharmonic, where $u_{m}$ is the minimal solution to the equation $\dbar_{z}(\cdot)=\dbar_{z}(u)$ (with respect to $h$). In this situation, we also have that $\partial_{t}\dbar_{t}\norm{u}^2_h=\partial_{t}(h(u,\delta u))=0$ so $\norm{u}^2_h$ is also plurisuperharmonic. Therefore, Corollary \ref{cor: second novel result on psh variations of minimal solutions} also applies and $-\norm{u}^2_{m}$ is plurisubharmonic.
\end{example}

%
%
%
%
%
%
%
%
%
%
%
\chapter{Paper 2}

\begin{center}
\normalsize{\bt{A HILBERT BUNDLES APPROACH TO COMPLEX BRUNN-MINKOWSKI THEORY, II}}\\[0.5cm]
\small{TAI TERJE HUU NGUYEN}

\begin{abstract}
The following is a sequel to \cite{Tai1}. We address in a particular setting the presumed positivity assumptions on the curvature(s) of the subbundle(s) in the main results in section 4 of \cite{Tai1}, and discuss some related results. In particular, we prove a version of the positivity of direct images from \cite{Bo1} in the case of a trivial fibration, for special types of weights where the metric on the line bundle may possibly be singular, and where the complex manifold $X$, over which we have the line bundle, belongs to a large class of Kähler manifolds. This class includes all the complete Kähler manifolds as well as the Zariski open sets in these. Our result is further a slight generalization in that we treat $(n,q)$-forms, $n$ being the complex dimension of $X$, where we allow that $q\geq 0$ (in \cite{Bo1}, $q=0$). Using similar ideas we also prove a  plurisubharmonicity result for certain minimal solutions to the $\dbar$-equation in the same setting. 
\end{abstract}

\end{center}

%
%
%
%
%
%
%

%

\section{Introduction}
The following is a sequel to \cite{Tai1}. We begin with a brief recollection of some of the main results proved in that paper, sticking to the following setting; we shall be in this setting also throughout the rest of this paper: Let $H$ be a trivial Hilbert bundle over an $m$-dimensional complex manifold, and suppose that $h$ is a smooth hermitian metric on $H$ which admits the Chern connection $D$. Following the notation in \cite{Tai1}, we shall generically write $D$ as $D=\delta+\dbar$. Suppose also that $H_0$ is a subbundle of $H$, closed in $H$ with respect to $h$, and let $P=P(h)$ denote the total orthogonal projection map associated with $h$ (of $H_0$ in $H$). One of the main results in \cite{Tai1}, Berndtsson's regularity theorem, asserts that $P$ is smooth if $h$ is a smooth weighted hermitian metric induced by a hermitian metric of zero variation with respect to which $H_0$ is also closed in $H$. Suppose that $P$ is smooth. Then, as was shown in \cite{Tai1}, the Chern connection on $H_0$ with respect to (the restriction of) $h$ (as a hermitian metric on $H_0$) exists and is given by $D_0:=\delta_0+\dbar$ with $\delta_0=P\delta$. Let $\Theta$ and $\Theta_0$ denote the Chern curvatures of $D$ and $D_0$ respectively. From \cite{Tai1}, we have the following three results on plurisubharmonic variations of smooth sections of $H$ and $H_0$ (with certain vanishing properties):

\begin{theorem}[Theorem 4.15 in \cite{Tai1}]\label{thm: characterization of Griffiths seminegativity of trivial Hilbert bundles}
With the above notation,
 and set-up, let $f$ be a twice differentiable function on $[0,\infty)$ which is strictly increasing and concave, and which satisfies that
\begin{align}
1+\fr{f''(x)}{f'(x)}x\geq 0
\end{align}for all $x\geq 0$. Then $i\Theta\leq_{G}0$ (or $i\Theta_{0}\leq_{G}0$) if and only if $f(\norm{u}^2_{h})$ is plurisubharmonic for all (local) holomorphic sections $u$ of $H$ (or $H_0$).
\end{theorem}

\begin{theorem}[Theorem 4.16 in \cite{Tai1}]\label{thm: first novel result on psh properties of projection map, II}
With the above notation and set-up, let $f$ be a twice differentiable function on $[0,\infty)$ which is strictly increasing and concave. Suppose that $u$ is a smooth (local) section of $H$ satisfying $\delta_0 u=0$, and put $u_0:=Pu$. Then if

\begin{align}
i h(\Theta_0 u_0,u_0)-i\norm{\dbar u_0}^2_{h}\left(1+\fr{f''(\norm{u_0}^2_{h}}{f'(\norm{u_0}^2_{h}}\norm{u_0}^2_{h}\right)&\geq 0,
\end{align}
 $f(\norm{u_0}^2_{h})$ is plurisubharmonic.
In particular, if $i\Theta_0\geq_{G}0$, then $f(\norm{u_0}^2_{h})$ is plurisubharmonic wherever \begin{align}
1+\fr{f''(\norm{u_0}^2_{h})}{f'(\norm{u_0}^2_{h}}\norm{u_0}^2_{h}\geq 0.\label{eq: the former}
\end{align}
\end{theorem}

\begin{theorem}[Corollary 4.17 in \cite{Tai1}]\label{cor: second novel result on psh variations of minimal solutions, II}
With the same assumptions and notation as in Theorem \ref{thm: first novel result on psh properties of projection map, II}, suppose furthermore that $\norm{u}^2_h$ is pluri\ita{super}harmonic. Then $-\norm{u-u_0}^2_{h}$ is plurisubharmonic.
\end{theorem}

As explained in \cite{Tai1}, Theorem \ref{thm: characterization of Griffiths seminegativity of trivial Hilbert bundles} in the case of $H_0$ may be seen as a counterpart to the variation of Bergman kernels from \cite{Bo0} in the (more) general setting of (trivial) Hilbert bundles,
while Theorem \ref{cor: second novel result on psh variations of minimal solutions, II} may be seen as a result on plurisubharmonic variation of minimal solutions to certain linear equations. Theorem \ref{thm: first novel result on psh properties of projection map, II} evidently gives variational properties of the total orthogonal projection map $P$, and is also naturally related to the the variation of Bergman kernels. 

A common factor for the above results is that we presume some kind of positivity on the subbundle curvature $\Theta_0$. In the original statement of the variation of the Bergman kernels from \cite{Bo0}, or rather, the vector bundles analogue of it (see Theorem 1.2, or Theorem 1.1 in the case of trivial fibrations, in \cite{Bo1}), this presumption is in fact the conclusion of the result. More specifically, the result asserts that a certain vector bundle associated with Bergman type of spaces is (semi)positive in the sense of Nakano. In this paper we shall address the positivity of the subbundle curvature $\Theta_0$ in a particular setting. The result that we give subsumes in the case of trivial fibrations the previous variation of Bergman kernels as well as its vector bundles analogue. In the vector bundles analogue, the sections of the relevant vector bundle are $L$-valued sections of the canonical bundle of some $n$-dimensional complex compact Kähler manifold, where $L$ is a(n) (ample) holomorphic hermitian line bundle over the manifold. Let us denote the manifold by $X$, and its canonical bundle by $K_{X}$. As the complex dimension of $X$ is $n$, these sections may alternatively be viewed as $L$-valued $(n,q)$-forms where $q=0$, and our result also subsumes a generalization of this where we allow for $q$ to be greater than 0, and $X$ to belong to a larger class of Kähler manifolds (namely, the class of so-called \ita{quasi-complete Kähler} manifolds); see Theorem \ref{thm: Nakano semipositivity of quasi-complete Kahler, smooth case}. In particular, all complete Kähler manifolds belong to this class, but the class is even closed under finite intersections with complements of analytic subsets. Using similar arguments, we also prove a result on plurisubharmonicity properties of holomorphic sections of the quotient bundle $Q:=H/H_0$, which is motivated by, and in the same spirit as, Theorem \ref{cor: second novel result on psh variations of minimal solutions, II} (see Theorems \ref{thm: second main result on psh properties of holomorphic sections of Q under certain assumptions} and \ref{thm: special case of psh for quotient sections, psh of minimal solutions in the setting of quasi-complete kahler manifolds, smooth case}). Finally, by considering certain special metrics on $H$, we prove versions of the former two results where the metric on $L$ may be possibly singular (see Theorems \ref{thm: Nakano semipositivity for singular metrics of special type on qck Kahler manifolds} and \ref{thm: psh properties of minimal solutions on qck manifolds with possiblt singular phi, special metrics}). The key here is a Hörmander type of theorem with singular weights (Theorem \ref{thm: singular Hormander}), which we obtain using regularization of quasi-plurisubharmonic functions (or currents) on compact manifolds, due to Demailly (see Theorem 16.3 in \cite{D}).

The Hilbert bundles that we are going to consider will be ''function spaces'', or ''functional spaces'', equipped with hermitian metrics given by integrating against some measure. Our formulation will be somewhat abstract, but the reader may have settings like those mentioned above (in the variation of Bergman kernels and its vector bundles analogue), and in Example 3.2 from \cite{Tai1}, in mind for more concrete examples. See also Theorem \ref{thm: Nakano semipositivity of quasi-complete Kahler, smooth case}. Let $X$ be a measure space with measure $d\mu_{X}$, and let $\mathcal{F}$ be a family of maps on $X$ such that $H$ is some subfamily of $\mathcal{F}$. Let $Y$ be the $m$-dimensional complex manifold over which we have our Hilbert bundles $H$ and $H_0$, and suppose that $f=\{f_{y}\}_{y\in Y}$ is a collection of maps on $H\times H$ such that for each pair of elements $u$ and $v$ in $H$, $f_{y}(u,v)$ is a complex-valued function on $X$ that can be integrated over $X$ against $d\mu_{X}$. Assume moreover that this defines a smooth hermitian metric $h=h^{f}$ on $H$. That is, explicitly, $h:=\{h_{y}\}_{y\in Y}$ is a smooth hermitian metric on $H$ given by

\begin{align}
h(u,v)(y)&=\int_{X}f_{y}(u,v)\;d\mu_{X}.
\end{align}Let next $T$ be a linear operator that can act on complex-valued functions, and assume that $T$ extends to an operator, also denoted $T$, on (local) sections of $H$, and even on (local) $H$-valued differential forms on $Y$. We assume moreover that if $\chi$ is a smooth function in the domain of definition of $T$, then $T(\chi)$ is, or at least can be identified with, a (collection of) (local) vector-valued map (we will also say a vector field) on $X$ in $\set{C}^{n}$. Here is an example illustrating what we mean by this: Consider the case that $T=\dbar^{\Omega}$, where $\dbar^{\Omega}$ is the $\dbar$-operator on $\Omega$, and $\Omega$ is some domain in $\set{C}^{n}$. Then if $\chi$ is a function on $\Omega$, and we choose $z:=(z_1,\ldots,z_n)$ as generic (global) coordinates on $\Omega$, we may identify $\dbar^{\Omega}\chi$ with the vector-valued map (we will also say the vector field) 
\begin{align}
a\mapsto \begin{pmatrix}
\del{\chi}{z_{1}}(a)\\
\vdots\\
\del{\chi}{z_{n}}(a)
\end{pmatrix}
\end{align} on $\Omega$. We can also do something similar in the case that $\Omega$ is a more general $n$-dimensional complex manifold, but then the identifications are in general only local. In this case, we therefore identify $T(\chi)$ instead with a \ita{collection} of (local) vector fields on $Y$ in $\set{C}^{n}$. Finally, we suppose that $\text{ker}(T)\sub H$ is closed in $H$ with respect to $h$, and take $H_0:=\text{ker}(T)$. We choose generic (local) holomorphic coordinates $t:=(t_1,\ldots, t_m)$ on $Y$, write with respect to these, 
\begin{align}
\Theta&=\s{j,k=1}{m}dt_{j}\wedge d\bar{t}_{k}\wedge \Theta_{j\bar{k}},
\end{align}and assume that for each $x\in X$, $[\Theta_{j\bar{k}}(\cdot, x)]$ defines a (collection of local) matrix field(s) on $Y$. With the above set-up specified, we may now proceed to state our main results. The following definition, whose motivation and relevance will be elaborated on in section \ref{sec: motivations}, is used:

\begin{definition}\label{def: Hörmander estimate condition}
With the above notation (and set-up), suppose that $v\in \text{Im}(T)$ (the image of $T$), and let $u\in H$ be a solution to the equation $T(\cdot)=v$. Let for each $x\in X$, $\vect{F}_{x}$ and $A_x$ be (a collection of local) matrix fields of sizes $n\times 1$ and $n\times n$ respectively, on $Y$. Suppose also that for all $y\in Y$ and each $x\in X$, $A_{x}(y)$ is semi positive definite hermitian. Then we shall say that (the triple) $\{u, \{\vect{F}_{x}\}, \{A_{x}\}\}$ satisfies the \bt{Hörmander $h$-estimate for the equation $T(\cdot)=v$} if
\begin{align}
\norm{u}^2_{h}\leq \lim_{\epsilon \downarrow 0}\int_{X}\norm{\vect{F}_{x}}^2_{(A_{x}+\epsilon I_{n})^{-1}}\;d\mu_{X}(x),
\end{align}assuming $\norm{\vect{F}_x}^2_{(A_x+\epsilon I_{n})^{-1}}$, as $x$ varies, can be integrated over $X$ with respect to $d\mu_{X}$ for all sufficiently small $\epsilon>0$, where $I_p$ denotes  the  identity matrix of  size $p\times p$, and where
\begin{align}
\norm{\vect{F}_x}^2_{(A_x+\epsilon I_{n})^{-1}}&=\vect{F}_x^*(A_x+\epsilon I_n)^{-1}\vect{F}_x,
\end{align}with $\vect{F}_x^*$ denoting the adjoint (that is, the complex conjugate of the transpose) of $\vect{F}_x$.
\end{definition} 
Our first main result concerns (Nakano) (semi)positivity of $\Theta_0$, and is the  following theorem:

\begin{theorem}\label{thm: first main theorem, Nakano semipositivity of subbundle under Hormander condition}
With the above notation and set-up, suppose that $i\Theta\geq_{G}0$, let $M:=M(y,x)$ be the matrix whose $(j,k)$th entry is $\Theta_{j\bar{k}}(y,x)$, and put $M_{x}:=M(\cdot,x)$. Suppose that, for each $x\in X$, $A_{x}$ and $B_{x}$ are collections of (local) matrix fields on $Y$ of sizes $n\times n$ and $n\times m$ respectively, and that each $A_{x}$ semi positive definite (at each point $y\in Y$). Let $B$ be defined by $B(y,x):=B_{x}(y)$, and let $\vect{b}_{j}$ denote the $j$th column of $B$. Assume that for all $\epsilon >0$
\begin{align}
(A_{x}+\epsilon I_{n})-B_{x}(M_{x}+\epsilon I_{m})^{-1}(B_{x})^*\geq 0,
\end{align}for each $x\in X$, and that for all $u_1,\ldots, u_m\in H_0$ there exists a solution $u\in H$ such that $\{u, \{\vect{F}_x\}, \{A_{x}\}\}$ satisfies the Hörmander $h$-estimate for the equation $T(\cdot)=v$, where

\begin{align}
v&=T\left(\s{j=1}{m}\delta_{j} u_{j}\right),
\end{align}
\begin{align}
\vect{F}&=\s{j=1}{m}\sqrt{f(u_j, u_j)}\vect{b}_{j},
\end{align}and $\vect{F}_{x}:=\vect{F}(\cdot,x)$.
Then $i\Theta_{0}\geq_{N}0$.
\end{theorem}
The following setting is of particular  interest: Consider a \ita{quasi-complete Kähler  manifold} $(X,\omega)$ $=\{(X_j, \omega_j, \chi_{j})\}$ of (complex) dimension $n$, and let $(L,e^{-\phi})$ be a hermitian holomorphic line bundle over $X$. The precise definition of \ita{quasi-complete Kähler manifolds} will come later (see Definition \ref{def: quasi-complete Kahler}), and we here simply mention that in particular all complete Kähler manifolds, and also the Zariski open subsets in these, are quasi-complete Kähler manifolds. Denote by $h^{0}=h^{0,\phi}$ the $L^2$-type of hermitian metric (inner product) on $L$-valued $(n,q)$-forms induced by the metric on $L$ and $\omega$, and by $\fancy{D}^{n,q}(L)$ the space of smooth $L$-valued $(n,q)$-forms with compact support in $X$. That is, explicitly, $h^{0}$ is given by
\begin{align}
h^{0}(u,v)&=\int_{X}(u,v)_{\omega}e^{-\phi}\omega_{n},
\end{align}where $(\cdot,\cdot)_{\omega}$ denotes the $\omega$-metric on differential forms, and $\omega_{n}:=\fr{\omega^{n}}{n!}$ (the product here is of course the exterior product). As is common, we shall usually omit (writing) $\omega_{n}$ in the integral. Let $H$ be the completion of $\fancy{D}^{n,q}(L)$ with respect to $h^0$. We view $H$ as a trivial Hilbert bundle over $B$, and $h^0$ as a hermitian metric of zero variation on $H$. Suppose that $\psi=\{\psi\}$ is a (collection of local) function(s) on $Y\times X$ such that multiplication by $e^{-\psi}$ defines a smooth weight for $h^0$, and such that $e^{-(\phi+\psi)}$ is a hermitian metric on $L$. Let $h$ denote the smooth weighted hermitian metric induced from $h^0$ with the smooth weight given by multiplication by $e^{-\psi}$. We put 
\begin{align}
\theta:=\phi+\psi,
\end{align} denote by $\theta^{y}$, for each $y\in Y$, the restriction of $\theta$ to $\{y\}\times X$, and let $H_0:=H\cap \text{ker}(\dbar^{X})$, where $\dbar^{X}$ is the $\dbar$-operator on $X$. In this setting, the following theorem generalizes in the case of trivial fibrations the positivity of direct images from \cite{Bo1} (to $(n,q)$-forms with $q\geq 0$, and more general Kähler manifolds):

\begin{theorem}\label{thm: Nakano semipositivity of quasi-complete Kahler, smooth case}
With the same above notation and set-up, suppose that $\theta$ is plurisubharmonic. Then $i\Theta_{0}\geq_{N}0$.
\end{theorem}
Theorem \ref{thm: Nakano semipositivity of quasi-complete Kahler, smooth case} is a special case of Theorem \ref{thm: first main theorem, Nakano semipositivity of subbundle under Hormander condition}. Now, by considering (a) special $\psi$ (and $Y$) in Theorem \ref{thm: Nakano semipositivity of quasi-complete Kahler, smooth case}, we may even allow $\phi$ to be possibly singular. More precisely, we have the following variant of Theorem \ref{thm: Nakano semipositivity of quasi-complete Kahler, smooth case} where $\phi$ may be singular:

\begin{theorem}\label{thm: Nakano semipositivity for singular metrics of special type on qck Kahler manifolds}
With the same notation and set-up as in Theorem \ref{thm: Nakano semipositivity of quasi-complete Kahler, smooth case}, suppose that $Y=\set{C}$, that $\phi$ is pseudoeffective (that is, $\phi$ is plurisubharmonic; in particular, $\phi$ may be singular), and that $G$ is a non-positive function on $X$. Let $\ld>0$, and let $\chi$ be a smooth convex function on $\set{R}$ such that $0\leq \chi'\leq \ld$, and $\chi(s)=0$ for all $s\leq 0$. Consider that $\theta$ is defined by
\begin{align}
\theta^{y}:=\phi+\chi(G-\text{Re}(y)),
\end{align}for each $y\in Y$, and suppose that
\begin{align}
i\partial_{z}\dbar_{z}(\phi+\ld G)&\geq 0,
\end{align}where the $\partial$ and $\dbar$ operators here are those on $X$. Then $i\Theta_{0}\geq_{N}0$.
\end{theorem}
The main step in proving both Theorems \ref{thm: Nakano semipositivity of quasi-complete Kahler, smooth case} and \ref{thm: Nakano semipositivity for singular metrics of special type on qck Kahler manifolds}, is to verify that the conditions of Theorem \ref{thm: first main theorem, Nakano semipositivity of subbundle under Hormander condition}, with $T:=\dbar_z$, are satisfied in the setting in which these theorems take place. In particular, we need to verify that the Hörmander $h$-estimate for a certain $\dbar_{z}$-equation, as formulated in Definition \ref{def: Hörmander estimate condition}, is satisfied. The main novelty here is the possibly non-smoothness of $\phi$ in Theorem \ref{thm: Nakano semipositivity for singular metrics of special type on qck Kahler manifolds} (and also Theorem \ref{thm: psh properties of minimal solutions on qck manifolds with possiblt singular phi, special metrics} later), and for this, we shall use the following Hörmander type of theorem for special weights with possibly singular parts, on quasi-complete Kähler manifolds:

\begin{theorem}\label{thm: singular Hormander}
With the same notation and set-up as in Theorem \ref{thm: Nakano semipositivity for singular metrics of special type on qck Kahler manifolds}, fix $y\in \set{C}$, put $t:=\text{Re}(y)$, and let $u_0\in H_0$ be given such that
\begin{align}
\norm{u_0}^2_{h_{y}}:=\norm{u_0}^2_t=\int_{X}\abs{u_0}_{\omega}^2e^{-\theta^{y}}<\infty,
\end{align}where $\abs{\cdot}^2_{\omega}=(\cdot,\cdot)_{\omega}$. Then there exists a solution $u\in H$ to the equation 
\begin{align}
\dbar_{z}(\cdot)=-\dbar_{z}\theta_{t}\wedge u_{0}
\end{align}
 satisfying
\begin{align}
\norm{u}^2_{t}=\int_{X}\abs{u}^2_{\omega}e^{-\theta^{y}}&\leq \norm{\sqrt{\theta_{t\bar{t}}}u_0}^2_{t}=\int_{X}(\theta_{tt}u_0,u_0)_{\omega}e^{-\theta^{y}},
\end{align}where we use the notation $\theta_{t}$ for $\del{\theta}{t}$, and $\theta_{tt}$ for $\del{^2\theta}{t^2}$.
\end{theorem}

 Theorem \ref{thm: singular Hormander} may be seen as a main result in itself. As mentioned above, we shall use regularization of quasi-plurisubharmonic functions (currents) on compact manifolds, due to Demailly (see Theorem 16.3 in \cite{D}), to prove it. Note that due to the $\chi(G-t)$-term we cannot apply the regularization directly, so there \ita{is} some work that needs to be done and shown. We will refrain from going into details here, and simply mention that the main idea is to apply the regularization to two (collections of) (quasi-)plurisubharmonic functions simultaneously and separately.

Our next main result is the following theorem on plurisubharmonicity properties of holomorphic sections of the quotient bundle $Q:=H/H_0$:

\begin{theorem}\label{thm: second main result on psh properties of holomorphic sections of Q under certain assumptions}
With the above notation, suppose that $i\Theta\geq_{G}0$, and let $[u]$ be a holomorphic (local) section of $Q$. Denote by $m_{l}$ the minimal lifting operator with respect to $h$, and by $h^{Q}$ the hermitian metric on $Q$ induced by $h$. Let $A_x,B_x$, and $M_x$, for each $x\in X$, be as in Theorem \ref{thm: first main theorem, Nakano semipositivity of subbundle under Hormander condition}. Assume that given any $\xi:=(\xi_1,\ldots, \xi_{m})\in \set{C}^{m}$, and any $\epsilon>0$,

\begin{align}
(A_{x}+\epsilon I_{n})-B_{x}(M_{x}+\epsilon I_{m})^{-1}(B_x)^*\geq 0,
\end{align}and that there exists $U\in H$ such that $\{U, \{\vect{F}_x\}, \{A_{x}\}\}$ satisfies the Hörmander $h$-estimate for the equation $T(\cdot)=v$, where

\begin{align}
v&=T\left(\s{j=1}{m}\xi_{j}\delta_{j}m_{l}[u]\right),
\end{align}
\begin{align}
\vect{F}&=\s{j=1}{m}\sqrt{f(m_{l}[u], m_{l}[u])}\abs{\xi_j}\vect{b}_{j},
\end{align}$\vect{b}_{j}$ denotes the $j$th column if $B$, defined also as in Theorem \ref{thm: first main theorem, Nakano semipositivity of subbundle under Hormander condition} (that is, $B(\cdot,x)=B_{x}$ for each $x\in X$), and $\vect{F}_{x}:=F(\cdot,x)$.
Then for all strictly increasing and concave twice differentiable functions $\tilde{f}$ on $[0,\infty)$, $-\tilde{f}(\norm{[u]}^2_{h^{Q}})$ is plurisubharmonic.
\end{theorem}
In particular, if the conditions of the theorem are met, then $-\norm{[u]}^2_{h^Q}$ and $-\log(\norm{[u]}^2_{h^{Q}})$ are plurisubharmonic. A special case of Theorem \ref{thm: second main result on psh properties of holomorphic sections of Q under certain assumptions} in the setting of Theorem \ref{thm: Nakano semipositivity of quasi-complete Kahler, smooth case}, is the following theorem; we use the same notational scheme for operators with subscripts as in \cite{Tai1}:

\begin{theorem}\label{thm: special case of psh for quotient sections, psh of minimal solutions in the setting of quasi-complete kahler manifolds, smooth case}
In the same setting as Theorem \ref{thm: Nakano semipositivity of quasi-complete Kahler, smooth case}, and using the same notation as in there, suppose that $\theta$ is plurisubharmonic. Let $v$ be an $L$-valued $(n,q+1)$-form, and suppose that 
\begin{align}
\partial_{t}\theta=0
\end{align} on the support of $v$, choosing $t$ as generic (local) holomorphic coordinates on $Y$. Choose also $z$ as generic (local) holomorphic coordinates on $X$, and suppose that $u\in H$ satisfies $\dbar_{z}u=v$. Let $u_{m}$ denote the $h$-minimal solution to $\dbar_z(\cdot)=v$. Then for all strictly increasing and concave twice differentiable functions $\tilde{f}$ on $[0,\infty)$, $-\tilde{f}(\norm{u_m}_{h}^2)$ is plurisubharmonic.
\end{theorem}
Theorem \ref{thm: special case of psh for quotient sections, psh of minimal solutions in the setting of quasi-complete kahler manifolds, smooth case} gives plurisubharmonicity properties of minimal solutions (to certain inhomogeneous $\dbar_z$-equations), and may (also) be viewed as a main result in itself. It might be related to variations of so-called generalized Green energies (\cite{XuDiag}), but the details are far from clear to us. We hope that this might be a future story for another time. Finally, similar to Theorem \ref{thm: Nakano semipositivity for singular metrics of special type on qck Kahler manifolds}, by considering (a) special $\psi$ (and $Y$), we may allow for $\phi$ in Theorem \ref{thm: special case of psh for quotient sections, psh of minimal solutions in the setting of quasi-complete kahler manifolds, smooth case} to be possibly singular. More precisely, we have the following version of Theorem \ref{thm: special case of psh for quotient sections, psh of minimal solutions in the setting of quasi-complete kahler manifolds, smooth case},  analogous to Theorem \ref{thm: Nakano semipositivity for singular metrics of special type on qck Kahler manifolds},  where $\phi$ may be singular:

\begin{theorem}\label{thm: psh properties of minimal solutions on qck manifolds with possiblt singular phi, special metrics}
With the same notation and assumptions as in Theorem \ref{thm: special case of psh for quotient sections, psh of minimal solutions in the setting of quasi-complete kahler manifolds, smooth case}, suppose that the metric $\theta$ is of the same form as in Theorem \ref{thm: Nakano semipositivity for singular metrics of special type on qck Kahler manifolds}. That is, for each $y\in \set{C}$, $t:=\text{Re}(y)$, $\theta(y,\cdot)=\theta^{y}$ is given by

\begin{align}
\theta^{y}&=\phi+\chi(G-t),
\end{align}where $\phi$ is pseudoeffective (in particular, it may be singular), $G\leq 0$ on $X$ and satisfies that
\begin{align}
i\partial_z\dbar_z(\phi+\ld G)&\geq 0
\end{align}for some $\ld>0$, and $\chi$ is a smooth convex function on $\set{R}$ which vanishes identically on $(-\infty, 0]$ and which satisfies that $0\leq \chi'\leq \ld$. Then the conclusion of Theorem \ref{thm: special case of psh for quotient sections, psh of minimal solutions in the setting of quasi-complete kahler manifolds, smooth case} holds. That is, $-\tilde{f}(\norm{u_m}^2_{h})$ is plurisubharmonic for all strictly increasing and concave twice differentiable functions $\tilde{f}$ on $[0,\infty)$.
\end{theorem}
As in the case for Theorem \ref{thm: Nakano semipositivity for singular metrics of special type on qck Kahler manifolds}, the key to proving Theorem \ref{thm: psh properties of minimal solutions on qck manifolds with possiblt singular phi, special metrics} (from Theorem \ref{thm: special case of psh for quotient sections, psh of minimal solutions in the setting of quasi-complete kahler manifolds, smooth case}) is Theorem \ref{thm: singular Hormander}. 

This concludes our introduction. In the remainder of the paper we shall detail the proofs of the above theorems, and give the motivation behind, and relevance of, Definition \ref{def: Hörmander estimate condition}. We will also give the precise definition of \ita{quasi-complete Kähler manifolds}, which we have yet to do. We end the paper by briefly commenting on the theme of the sequel of this paper (\cite{TaiXuHilb}).

\section{A linear algebraic result and motivations}
In this section we discuss a few preliminary linear algebraic results that we shall have use for when proving our main results from the introduction. We also explain in this section the motivation behind, and the relevance of, the Hörmander estimate in Definition \ref{def: Hörmander estimate condition}. As one might expect,  the latter is naturally related to Hörmander's $L^2$-estimates for solutions of the $\dbar$-equation (see, for example, \cite{Bo1}, section 3). Actually, the motivation for our upcoming linear algebraic discussion comes from there as well.

\subsection{Some linear algebra on norms induced by matrices}\label{section: linear algebra}
We shall view elements in $\set{C}^{n}$ as column vectors. Let $A$ be a positive hermitian matrix of size $n\times n$. The matrix $A$ defines an inner product $(\cdot,\cdot)_{A}$ on $\set{C}^{n}$ given simply by 
\begin{align}
(\vect{v}, \vect{w})_{A}:=\vect{w}^*A\vect{v},
\end{align} for $\vect{v}$ and $\vect{w}$ elements in $\set{C}^{n}$, and $\vect{w}^*$ the adjoint of $\vect{w}$ (that is, the complex conjugate of its transpose). The converse is of course also true. We will denote the norm induced from the inner product $(\cdot,\cdot)_{A}$ by $\norm{\cdot}_{A}$. Note that we may do similar things also in the case that $A$ is only semi positive definite, but in this case, $\norm{\cdot}_{A}$ is in general only a semi-norm. Let $\vect{F}\in \set{C}^{n}$. We may view $\vect{F}$ also as a complex-antilinear functional on $\set{C}^{n}$ given simply by 
\begin{align}
\set{C}^{n}\to \set{C}:\vect{v}\mapsto \ip{\vect{F}}{\vect{v}}:=\vect{v}^*\vect{F}.
\end{align} Let us denote the operator norm of $\vect{F}$, when viewed as such a functional, with respect to the norm $\norm{\cdot}_{A}$, by $\norm{\vect{F}}_{*,A}$. By definition,
\begin{align}
\norm{\vect{F}}^2_{*,A}&=\sup_{\vect{v}\in \set{C}^{n}\backslash\{\vect{0}\}}\fr{\abs{\ip{\vect{F}}{\vect{v}}}^2}{\norm{\vect{v}}^2_{A}}.
\end{align}
We may also consider the norm of $\vect{F}$ viewed as a vector in $\set{C}^{n}$ with respect to $A^{-1}$, the inverse of $A$. It follows by the Riesz representation theorem (or simply by a direct computation) that these two norms are the same. That is, 
\begin{align}
\norm{\vect{F}}_{*,A}=\norm{\vect{F}}_{A^{-1}}.
\end{align} Observe that, for $\vect{v}\in \set{C}^{n}$, $\abs{\ip{\vect{F}}{\vect{v}}}^2$ is equal to the semi norm of $\vect{v}$ with respect to a certain semi positive definite matrix $M$. Namely, $M=\vect{F}\vect{F}^*$. Indeed, since $\ip{\vect{F}}{\vect{v}}$ is a number,
\begin{align}
\abs{\ip{\vect{F}}{\vect{v}}}^2&=\ip{\vect{F}}{\vect{v}}\ip{\vect{F}}{\vect{v}}^*=(\vect{v}^*\vect{F})(\vect{v}^*\vect{F})^*=\vect{v}^*(\vect{F}\vect{F}^*)\vect{v}.
\end{align}
From this we immediately get the following lemma:

\begin{lemma}\label{lem: useful linear algebra lemma}
With the above notation, suppose that there is some $\eta>0$ such that 
\begin{align}
\vect{F}\vect{F}^*\leq \eta A.
\end{align} Then \begin{align}
\norm{\vect{F}}^2_{A^{-1}}\leq \eta.
\end{align}
\end{lemma}
Consider now the special case that $\vect{F}$ is of the form 
\begin{align}
\vect{F}=\s{j=1}{m}\vect{b}_{j}
\end{align} where each $\vect{b}_{j}\in \set{C}^{n}$. In this situation we may try and compare $\vect{F}\vect{F}^*$ with the matrix $BB^*$, where $B$ is the matrix of size $n\times m$ whose $j$th column is $\vect{b}_{j}$. Let $I_{p}$ denote the identity matrix of size $p\times p$. We may of course also write $BB^*=BI_{m}^{-1}B^*$. We generalize this somewhat by replacing $I_{m}$ with a diagonal matrix $\Ld$ consisting of positive entries, and compare $\vect{F}\vect{F}^*$ with $B\Ld^{-1}B^*$. Granted, this may seem a bit unmotivated at the moment, but as with everything else, we shall in due time explain our motivation(s) for this. Hopefully things become clearer as we go along, and in the end, all come together nicely. A pertinent example where the above is of natural interest, albeit still somewhat unmotivated, is the following: Consider a block matrix 
\begin{align}
M=\begin{pmatrix}
\Ld&B^*\\
B&A
\end{pmatrix}.
\end{align} Then $M\geq 0$ if and only if 
\begin{align}
A-B\Ld^{-1}B^*\geq 0;
\end{align} the interesting ''part'' here is of course the ''term'' $B\Ld^{-1}B^*$. This is a special case of the following more general \ita{Schur complement theorem} from matrix analysis (see \cite{MatrixAnalysis}):

\begin{theorem}\label{thm: Schurs complement}
Let $M$ be a hermitian matrix of the block form 
\begin{align}
M=\begin{pmatrix}
\Sigma &B^*\\
B&C
\end{pmatrix},
\end{align} where $\Sigma$ is non-singular. Then $M\geq 0$ if and only if  $\Sigma\geq 0$
\begin{align}
C-B\Sigma^{-1}B^*\geq 0.
 \end{align}
\end{theorem}
Our main linear algebraic result is the following theorem on the comparison of $\vect{F}\vect{F}^*$ and $B\Ld^{-1}B^*$:

\begin{theorem}\label{thm: main linear algebraic result}
With the above notation and set-up, let 
\begin{align}
\Ld:=\text{diag}(\ld_1,\ldots, \ld_{m})=\begin{pmatrix}
\ld_1\\
&\ld_2\\
&&\ddots&\\
&&&\ld_m
\end{pmatrix}
\end{align} and 
\begin{align}\ld:=\s{j=1}{m}\ld_{j}.
\end{align} Then 
\begin{align}
\vect{F}\vect{F}^*\leq \ld B\Ld^{-1}B^*.
\end{align}
\end{theorem}
\begin{proof}
 Let $\Sigma:=(\Ld^{-1})^{1/2}$ and $B:=[\vect{b}_1\cdots \vect{b}_m]$. Then we have 
\begin{align}
B\Ld^{-1}B^*&=(B\Sigma)(B\Sigma)^*,
\end{align}
\begin{align}
B\Sigma&=\left[B\begin{pmatrix}
\ld_1^{-1/2}\\
0\\
\vdots\\
0
\end{pmatrix}\cdots B\begin{pmatrix}
0\\\
0\\
\vdots\\
\ld_{m}^{-1/2}
\end{pmatrix}\right],
\end{align}and
\begin{align}
B\begin{pmatrix}
0\\
\vdots\\
\ld_{j}^{-1/2}\\
\vdots\\
0
\end{pmatrix}=\fr{1}{\sqrt{\ld_j}}\vect{b}_{j}.
\end{align}Hence, 
\begin{align}
B\Sigma&=\left[\fr{1}{\sqrt{\ld_1}}\vect{b}_1\cdots \fr{1}{\sqrt{\ld_m}}\vect{b}_m\right],
\end{align}so
\begin{align}
(B\Sigma)^*&=\begin{bmatrix}
\fr{1}{\sqrt{\ld_1}}\vect{b}_1^*\\
\vdots\\
\fr{1}{\sqrt{\ld_m}}\vect{b}_m^*
\end{bmatrix}.
\end{align}Let $\vect{v}\in \set{C}^{n}\backslash \{\vect{0}\}$. We then get 
\begin{align}
(B\Sigma)^*\vect{v}&=\begin{bmatrix}
\fr{1}{\sqrt{\ld_1}}\vect{b}_1^*\vect{v}\\
\vdots\\
\fr{1}{\sqrt{\ld_m}}\vect{b}_m^*\vect{v}
\end{bmatrix}.
\end{align}Thus, 
\begin{align}
\vect{v}^*(B\Ld^{-1}B)\vect{v}&=\vect{v}^*B\Sigma (B\Sigma)^*\vect{v}=((B\Sigma)^*\vect{v})^*((B\Sigma)^*\vect{v})\notag\\
&=\norm{(B\Sigma)^*\vect{v}}^2_e=\fr{1}{\ld_1}\abs{\vect{b}_1^*\vect{v}}^2+\cdots+\fr{1}{\ld_m}\abs{\vect{b}_m^*\vect{v}}^2.\label{eq: first in linear algebraic result}
\end{align}Since $\vect{F}=\s{j=1}{m}\vect{b}_{j}$, we also directly get, where $\norm{\cdot}_e$ here denotes the usual euclidean norm,
\begin{align}
\vect{v}^*(\vect{F}\vect{F}^*)\vect{v}&=\norm{\vect{F}^*\vect{v}}^2_{e}=\norm{\s{j=1}{m}\vect{b}_{j}^*\vect{v}}^2_e=\abs{\vect{b}_1^*\vect{v}+\cdots+\vect{b}_m^*\vect{v}}^2.\label{eq: second in linear algebraic result}
\end{align}Multiplying \eqref{eq: first in linear algebraic result} by $\ld$ and subtracting \eqref{eq: second in linear algebraic result}, gives

\begin{align}
\vect{v}^*(\ld B\Ld^{-1}B^*-\vect{F}\vect{F}^*)\vect{v}&=\s{j=1}{m}\fr{\ld}{\ld_{j}}\abs{\vect{b}^*_{j}\vect{v}}^2-\abs{\s{j=1}{m}\vect{b}_{j}^*\vect{v}}^2\notag\\
&\s{j=1}{m}\fr{\ld_j}{\ld}\abs{\fr{\ld}{\ld_j}\vect{b}_{j}^*\vect{v}}^2-\abs{\s{j=1}{m}\fr{\ld_j}{\ld}\left(\fr{\ld}{\ld_j}\vect{b}_j^*\vect{v}\right)}^2.\label{eq: last in linear algebraic result}
\end{align}Hence, as $\s{j=1}{m}\fr{\ld_j}{\ld}=1$, it follows by convexity of $\abs{\cdot}^2$ that the right-hand side in \eqref{eq: last in linear algebraic result}, thus the left-hand side, is no less than 0. That is:
\begin{align}
\vect{v}^*(\ld B\Ld^{-1}B^*-\vect{F}\vect{F}^*)\geq 0,
\end{align}which is what we wanted to prove.
\end{proof}
We have the following important corollary to Theorem \ref{thm: main linear algebraic result}:

\begin{corollary}\label{cor: corollary to main linear algebraic result}
Let $A$ and $M$ be positive definite hermitian matrices of sizes $n\times n$ and $m\times m$ respectively, and let $B$ be a matrix of size $n\times m$. Let $\vect{b}_{j}$ denote its $j$th column, and define $\vect{F}:=\s{j=1}{m}\vect{b}_{j}$. Suppose that 
\begin{align}
A-BM^{-1}B^*\geq 0.
\end{align}Then 
\begin{align}
\norm{\vect{F}}^2_{A^{-1}}&\leq T_{r}(M),
\end{align}where $T_{r}(M)$ denotes the trace of $M$.
\end{corollary}

\begin{proof}
Let $\Sigma$ be the orthogonal diagonalization of $M$. Then
\begin{align}
BM^{-1}B^*&\geq B\Sigma^{-1}B^*.
\end{align}The sum of entries of $\Sigma$ is $T_{r}(M)$, so it follows by Theorem \ref{thm: main linear algebraic result} and Lemma \ref{lem: useful linear algebra lemma} that
\begin{align}
\norm{\vect{F}}^2_{A^{-1}}&\leq T_{r}(M).
\end{align}
\end{proof}

\subsection{Motivations}\label{sec: motivations}

Having finished section \ref{section: linear algebra}, we are now in a position to prove both Theorems \ref{thm: first main theorem, Nakano semipositivity of subbundle under Hormander condition} and \ref{thm: second main result on psh properties of holomorphic sections of Q under certain assumptions} from the introduction. Before we do this, however, let us first discuss the motivation(s) behind everything up to this point. 

Consider a domain $U\times \Omega$ in $\set{C}^{m}\times \set{C}^{n}$, where $U$ is open and $\Omega$ is bounded. Consider next a strictly plurisubharmonic function $\phi$ on $U\times \Omega$ which is smooth up to the boundary, and let for each $t\in U$, $\phi^{t}$ denote the restriction of $\phi$ to $\{t\}\times \Omega$. The set-up that we shall now consider is precisely that in the introduction in \cite{Tai1}. We denote the holomorphic functions on $\Omega$ by $\fancy{O}(\Omega)$, and let
\begin{align}
H&:=\cbrac{f\in L^2_{\text{loc}}(\Omega):\int_{\Omega}\abs{f}^2<\infty}\quad \quad \text{ and }\\
H_0&:=H\cap \fancy{O}(\Omega),
\end{align}where we integrate with respect to the Lebesgue measure. Then $H$ and $H_0$ are trivial Hilbert bundles over $U$ with respect to the usual $L^2$-metric, and $H_0\leq H$. We consider the $L^2$-metric as a hermitian metric $h^0$ of zero variation on $H$ and $H_0$. Multiplication by $e^{-\phi}$ defines a smooth weight for $h^0$ (on $H$(!)), and we denote the induced smooth weighted hermitian metric by $h$. That is, $h$ is defined by 
\begin{align}
h_{t}(\cdot,\cdot):=h^{0}(e^{-\phi^{t}}\cdot,\cdot),
\end{align} for each $t\in U$. We have seen in \cite{Tai1} that the Chern connection on $H$ and $H_0$ with respect to $h$ exist. Let use denote these by $D$ and $D_0$, respectively, and let us write for their respective Chern curvatures, $\Theta$ and $\Theta_0$. We choose $t:=(t_1,\ldots, t_m)$ and $z:=(z_1,\ldots,z_n)$ as generic coordinates on $U$ and $\Omega$, respectively. We may consider the $\partial$ and $\dbar$ operators on both $U$ and $\Omega$, and will then as usual use subscripts with $t$ and $z$ to distinguish whether we consider the operators on $U$ or $\Omega$. Thus, $\dbar_t$ denotes the $\dbar$-operator on $U$, $\partial_z$ denotes the $\partial$-operator on $\Omega$, and so on. In the situation at hand, we have (seen that we have)
\begin{align}
\delta=\partial_t-\partial_t \phi\wedge
\end{align} and 
\begin{align}\Theta=\partial_t\dbar_t\phi\wedge.
\end{align} Let us consider (Nakano) (semi)positivity of $\Theta_0$. For simplicity, and to fix the ideas, we will here stick to the case that $m=1$. Let $S:=D-D_0$ be the second fundamental form (associated with $D$ and $D_0$). Our starting point is the subbundle curvature formula
\begin{align}
\Theta_0&=\Theta-S^*S,\label{eq: subbundle curvature formula}
\end{align}from \cite{Tai1}. Since $\phi$ is strictly plurisubharmonic, we have $i\Theta>_{G}0$, and following Berndtsson (see \cite{Bo0}, \cite{Bo1}), the idea is to use $L^2$-methods to control the $S^*S$ term in \eqref{eq: subbundle curvature formula}. We first repeat Berndtsson's argument. This will (hopefully) make the motivation behind our discussion in section \ref{section: linear algebra}, and our interest in norms like $\norm{\vect{F}}_{A^{-1}}$ in there, apparent. Let $u\in H_0$. Then $Su\perp H_0$ (in fact this holds even if $u\in H$), which means that $Su$ is, slightly imprecisely spoken, the $h$-minimal solution to the equation 
\begin{align}
\dbar_{z}(\cdot)=\dbar_{z}(Su)=\dbar_{z}(\delta u)=-\dbar_{z}\partial_{t}\phi\wedge u.
\end{align} Under the additional assumption that $\Omega$ is pseudoconvex, we may therefore (try and) use Hörmander's $L^2$-estimates for (solutions of) the $\dbar_{(z)}$-equation to control the $h$-norm of $Su$. The precise estimates that we shall need says the following (see, for example, section 2 in \cite{Bo0}, section 3 in \cite{Bo1}): Let 
\begin{align}
f=\s{j=1}{n}f_{j}d\bar{z}_{j}
\end{align} be a $\dbar_z$-closed (0,1)-form on $\Omega$. Then there exists a solution $v$ to the equation $\dbar_z(\cdot)=f$ (in the sense of currents) such that
\begin{align}
\int_{\Omega}\abs{v}^2e^{-\phi^{t}}\leq \int_{\Omega}\s{j,k=1}{n}(\phi^{t})^{j\bar{k}}f_{j}\overline{f}_{k}e^{-\phi^{t}}\label{eq: Hormander estimate local setting}
\end{align}for each $t\in U$, where $(\phi^{t})^{j\bar{k}}$ denotes the $(j,k)$th entry of the inverse of the complex Hessian matrix $[\phi^{t}]$ of $\phi^{t}$. Consider $f:=-\dbar_{z}\partial_{t}\phi\wedge u$ (''divided by'', or without, $dt$). The point, of course, is that since $Su$ (again, ''divided by'', or without, $dt$) is the \ita{minimal solution} to this equation, it too must satisfy the estimate in \eqref{eq: Hormander estimate local setting}. Indeed, we have (still somewhat imprecisely),
\begin{align}
\norm{Su}^2_{h_{t}}\leq \displaystyle \int_{\Omega}\abs{v}^2e^{-\phi^{t}}
\end{align}
(and more precisely, $i\norm{Su}^2_{h_{t}}\leq i\displaystyle\int_{\Omega}\abs{v}^2e^{-\phi^{t}} dt\wedge d\bar{t}$) if $v$ is any solution to the equation. To show that 
\begin{align}
ih(\Theta u,u)-i\norm{Su}^2_{h}\geq 0,
\end{align} which is what we want to do according to \eqref{eq: subbundle curvature formula}, it therefore suffices to show that $h(\Theta u,u)(t)$ (''divided by'', or without, $dt$) is no smaller than the right-hand side in \eqref{eq: Hormander estimate local setting}, where $f$ is chosen as above to be 
\begin{align}
f=-\left(\dbar_{z}\phi^{t}_{t}\right)\wedge u,
\end{align} and where we here use the notation $\phi_{t}$ for $\del{\phi}{t}$ and $\phi_{t\bar{t}}$ for $\del{^2\phi}{t\partial \bar{t}}$. Similarly, we write, $\phi_{\bar{z}_{j}}$ for $\del{\phi}{\bar{z}_{j}}$, $\phi_{\bar{z}_{j}t}$ for $\del{^2\phi}{\bar{z}_{j}t}$, and so on. Unravelling everything, we are then interested in showing that
\begin{align}
\int_{\Omega}\left(\phi_{t\bar{t}}-\s{j,k=1}{n}(\phi^{t})^{j\bar{k}}\phi_{\bar{z}_{j}t}\overline{\phi_{\bar{z}_{k}t}}\right)\abs{u}^2e^{-\phi^{t}}\geq 0.\label{eq: determinant full hessian divided by determinant hessian phit}
\end{align}The key now, as observed by Berndtsson (see \cite{Bo0}), is that the expression inside the parenthesis in \eqref{eq: determinant full hessian divided by determinant hessian phit} is the determinant of the full complex Hessian of $\phi$ divided by the determinant of the complex Hessian of $\phi^{t}$. Thus \eqref{eq: determinant full hessian divided by determinant hessian phit} is indeed true (under our assumptions of strict plurisubharmonicity, we even have strict inequality). The case $m>1$ is analogous, but somewhat ''messier''; see \cite{Bo1}. The connection to our linear algebraic discussion in section \ref{section: linear algebra} is the following: Let $[\phi^{t}]:=A$ (that is $A$ is the complex Hessian of $\phi$ with respect to the $z$-variable), and identify $f$, as we may, with the vector field 
\begin{align}
\vect{f}:=\begin{pmatrix}
f_1\\
\vdots\\
f_n
\end{pmatrix}.
\end{align}Then the right-hand side in \eqref{eq: Hormander estimate local setting} is equal to
\begin{align}\displaystyle \int_{\Omega}\norm{\vect{f}}^2_{A^{-1}} e^{-\phi^{t}}.
\end{align} This is precisely where our motivation to study norms like $\norm{\vect{F}}_{A^{-1}}$, from earlier, comes from. We may next consider a similar example as above, but now in the global setting of complex manifolds; see also example in \cite{Tai1}. In this setting, $u$ above is a holomorphic $L$-valued section of $K_{X}$, the canonical bundle of $X$, an $n$-dimensional complex manifold that plays the role of $\Omega$. The role of the weight function $e^{-\phi}$ is played by a smooth hermitian metric on $L$, and we also replace $U$ with an $m$-dimensional manifold $Y$. The hermitian metric $h^0$ of zero variation is in this case defined by \begin{align}
h^0(u',v')=\displaystyle \int_{X}i^{n^2}u'\wedge \bar{v'},
\end{align} for $L$-valued sections $u'$ and $v'$ of $K_{X}$.  In this setting, under appropriate assumptions on $X$ and $\phi$, we still have analogous $L^2$-estimates for the $\dbar_{(z)}$-equation as above, and the (more or less)  same arguments as before give (Nakano) (semi)positivity of the subbundle. Since the complex dimension of $X$ is $n$, $u$ above in this setting may also be interpreted as an $L$-valued $(n,q)$ form with $q=0$, and as mentioned in the introduction, we can generalize this further by also allowing for $q\geq 1$. In the case that $q\neq 0$, we can of course not directly use the same definition as before for $h^0$, but there are natural analogues if we add to the picture some hermitian metric on $X$, that is, on its holomorphic tangent bundle. In these situations we still have Hörmander $L^2$-estimates for the $\dbar_{z}$-equation, but these estimates differ somewhat from those above and the case $q=0$. Fortunately, they do not differ too much, and with some more work, we are able to relate these to certain estimates where the right-hand side upper bound is of the form (an integral of) $\norm{\vect{F}}^2_{A^{-1}}$ like before, for suitable  local vector and matrix fields $\vect{F}$ and $A$. More specifically, we are able to show that a Hörmander $h$-estimate for a certain $\dbar_{z}$-equation, as defined in Definition \ref{def: Hörmander estimate condition}, holds. From here, we apply Theorem \ref{thm: main linear algebraic result}. As we shall show, this provides an alternative way, which differs somewhat from the argument given earlier, of showing that the right-hand side \eqref{eq: subbundle curvature formula} is positive. In particular, there will be no mentioning of any determinants, and the (slightly messy) equation \eqref{eq: determinant full hessian divided by determinant hessian phit} above will also come into play, at least not as directly as earlier. Moreover, there is virtually no difference between the case $m=1$ and the case $m>1$, and the arguments that we give are in either of these cases more or less identical. Finally, Theorem \ref{thm: main linear algebraic result} will also be very useful when we later consider $\phi$ to be possibly singular; that is, when we prove Theorem \ref{thm: singular Hormander}.

We now explain this (first part) a bit more, without requiring that $m=1$. For simplicity, we stick, however, still to the the case that $q=0$. We want to prove that $i\Theta_{0}\geq_{N}0$, and we want to use what we have discussed in section \ref{section: linear algebra}. Let $u_1,\ldots, u_{m}\in H$. According to the subbundle curvature formula, that is \eqref{eq: subbundle curvature formula}, we need to show, point-wisely, that
\begin{align}
\s{j,k=1}{m}h(\Theta_{j\bar{k}}u_{j}, u_{k})-\norm{\s{j=1}{m}s_{j}u_{j}}^2_{h},
\end{align}where we write $S:=\s{j=1}{m}s_j dt_{j}$. Let 
\begin{align}
f:=\dbar_{z}\left(\s{j=1}{m}s_{j}u_{j}\right).
\end{align} A similar computation as before shows that 
\begin{align}
f=-\s{j=1}{m}\dbar_{z}(\partial_{t_{j}}\phi)\wedge u_{j},
\end{align} where we as usual write $\phi_{t_{j}}$ for $\del{\phi}{t_{j}}$ and so on. Let 
\begin{align}
\vect{b}_{j}:=\begin{pmatrix}
\phi_{t_{j}\bar{z}_1}\\
\vdots\\
\phi_{t_{j}\bar{z}_{n}}
\end{pmatrix}
\end{align} and 
\begin{align}
\vect{F}:=\s{j=1}{m}(-u_{j})\vect{b}_{j}.
\end{align} We let as before $[\phi^{t}]:=A$. The $L^2$-estimates of Hörmander, \eqref{eq: determinant full hessian divided by determinant hessian phit}, now says that if $v$ is any solution to $\dbar_{z}(\cdot)=f$, then for each $t\in U$, 
\begin{align}
\norm{v}^2_{h}(t)&\leq \int_{\Omega}\norm{\vect{F}}^2_{A^{-1}}e^{-\phi^{t}}(t).
\end{align}The solution $\s{j=1}{m}s_{j}u_{j}$ to the equation is the $h$-minimal solution, so we of course also have 
\begin{align}
\norm{\s{j=1}{m}s_{j}u_{j}}^2_{h}(t)&\leq \int_{\Omega}\norm{\vect{F}}^2_{A^{-1}}e^{-\phi^{t}}(t).
\end{align}To prove that $i\Theta_0\geq_{N}0$, it therefore suffices to show that at each $t\in U$,
\begin{align}
\s{j,k=1}{m}h(\Theta_{j\bar{k}}u_{j},u_{k})(t)-\int_{\Omega}\norm{\vect{F}}^2_{A^{-1}}e^{-\phi^{t}}(t).
\end{align}That is,
\begin{align}
\int_{\Omega}\left(\s{j,k=1}{m}\Theta_{j\bar{k}}u_{j}\bar{u}_{k}-\norm{\vect{F}}^2_{A^{-1}}\right)e^{-\phi^{t}}(t)\geq 0.\label{eq: subbundle curvature, clean version}
\end{align}Fix $z\in \Omega$. It suffices to show that the integrand in \eqref{eq: subbundle curvature, clean version} is non-negative at $z$. Having fixed $z$, we may choose the coordinates $t$ such that the matrix $\Theta_{j\bar{k}}(t,z)$ is a diagonal matrix $\Ld$. We denote the (positive) entries of $\Ld$ by $\ld_1,\ldots, \ld_m$. Then, want we want to show is that; for simplicity we drop evaluations at points:
\begin{align}
\s{j=1}{m}\ld_{j}\abs{u_{j}}^2-\norm{\vect{F}}^2_{A^{-1}}\geq 0.
\end{align}By replacing $\vect{b}_{j}$ with $u_{j}\vect{b}_{j}$, and $\ld_{j}$ with $\abs{u_{j}}^2\ld_{j}$, this now follows by (the proof of) Theorem \ref{thm: main linear algebraic result} and Corollary \ref{cor: corollary to main linear algebraic result} if we can show that $A\geq B\Ld^{-1}B^*$, where $B$ is the matrix whose $j$th column is $\vect{b}_{j}$. Consider the (full) complex Hessian matrix of $\phi$ (at the point $(t,z)$). It is of the block form 
\begin{align}
[\phi]=\begin{pmatrix}
\Ld&B^*\\
B&A
\end{pmatrix},
\end{align} and from Theorem \ref{thm: Schurs complement}, since $[\phi]\geq 0$, we know that \begin{align}
A-B\Ld^{-1}B^*\geq 0.
\end{align} Thus we get $i\Theta_0\geq_{N}0$ as required. 

We see from our discussion where the motivation to include $\Ld^{-1}$ in $B^*\Ld^{-1}B$ comes from. That is, why we earlier have insisted on comparing $\vect{F}\vect{F}^*$ with $B\Ld^{-1}B^*$ instead of just $BB^*$. We also see that the argument just given is indeed a little bit different from Berndtsson's original argument that we previously recapitulated. Of course, the key idea
 is still owing to Berndtsson, namely to use $L^2$-methods to control the second fundamental form term in the subbundle curvature formula. Our contribution is instead Theorem \ref{thm: first main theorem, Nakano semipositivity of subbundle under Hormander condition}, which is more general than what we have discussed, and which also applies to the case that $q\geq 1$ (see Theorems \ref{thm: Nakano semipositivity of quasi-complete Kahler, smooth case} and \ref{thm: Nakano semipositivity for singular metrics of special type on qck Kahler manifolds}). It remains to say a few words about the motivation behind Definition \ref{def: Hörmander estimate condition}. This will be more apparent from the proof of Theorem \ref{thm: first main theorem, Nakano semipositivity of subbundle under Hormander condition} that we give below, but the main idea is of course that it replaces the $L^2$-estimates of Hörmander above, \eqref{eq: Hormander estimate local setting}. To access these estimates, we usually need to impose certain pseudoconvexity assumptions on the domain, and plurisubharmonicity properties on the weight function. Definition \ref{def: Hörmander estimate condition} gives us the estimates somehow for free without imposing such conditions, allowing for some general statements. In practice, we must of course verify that we do have these estimates, and in these cases, the former conditions on the domain and weight function again must be imposed. As mentioned above, and in the introduction, there are analogous situations, or settings, in which we also have Hörmander $L^2$-estimates, but the required conditions, and also the directly given estimates, may vary somewhat. For example, in the global setting of complex manifolds, there is usually some kind of compact or complete Kähler condition on the underlying manifold (\cite{Bodbar}). The idea is that by using Definition \ref{def: Hörmander estimate condition}, Theorem \ref{thm: first main theorem, Nakano semipositivity of subbundle under Hormander condition} somehow places all these analogous examples under the same umbrella. Finally, taking the approach through Theorem \ref{thm: first main theorem, Nakano semipositivity of subbundle under Hormander condition}, has also led us to discover Theorem \ref{thm: second main result on psh properties of holomorphic sections of Q under certain assumptions} (and its subsequent special cases; Theorems \ref{thm: special case of psh for quotient sections, psh of minimal solutions in the setting of quasi-complete kahler manifolds, smooth case} and \ref{thm: psh properties of minimal solutions on qck manifolds with possiblt singular phi, special metrics}).
 
This concludes our discussion on linear algebraic preliminaries, and the motivation behind, and relevance of, Definition \ref{def: Hörmander estimate condition}. In the remainder of the paper, starting with the next section, we give the proofs of the theorems in the introduction (which are not from \cite{Tai1}). We begin with the proofs of Theorem \ref{thm: first main theorem, Nakano semipositivity of subbundle under Hormander condition} and Theorem \ref{thm: second main result on psh properties of holomorphic sections of Q under certain assumptions}, as we are going to use theorems in the proofs of the remaining theorems.

\section{Proofs of Theorems 1.4 and 1.8}
In this section we prove Theorems \ref{thm: first main theorem, Nakano semipositivity of subbundle under Hormander condition} and \ref{thm: second main result on psh properties of holomorphic sections of Q under certain assumptions} from the introduction. Since much of its proof has  actually already been discussed in the previous section, we start with Theorem \ref{thm: first main theorem, Nakano semipositivity of subbundle under Hormander condition}. We will use similar notation as above, choosing generic (local) holomorphic coordinates $t:=(t_1,\ldots, t_m)$ on $Y$, and writing with respect to these:
\begin{align}
\Theta_{(0)}&=\s{j,k=1}{m}dt_{j}\wedge d\bar{t}_{j}\Theta_{j\bar{k}}^{(0)},
\end{align}and
\begin{align}
\delta_{(0)}&=\s{j=1}{m}dt_{j}\wedge \delta^{(0)}_{j}.
\end{align}

\subsection{Nakano (semi)positivity of subbundles: Proof of Theorem 1.4}
\begin{proof}
 Let $u_1,\ldots, u_{m}\in H_0$. We want to show that at each $y\in Y$,
\begin{align}
\s{j,k=1}{m}h(\Theta^{0}_{j\bar{k}}u_{j}, u_{k})&\geq 0.
\end{align} Let $y\in Y$, $\epsilon>0$, and put $s_{j}:=\delta_{j}-\delta_{j}^{0}$. By the subbundle curvature, \eqref{eq: subbundle curvature formula}, it suffices to show that
\begin{align}
\left(\s{j,k=1}{m}h(\Theta_{j\bar{k}}u_{j}, u_{k})-\norm{\s{j=1}{m}s_{j}u_{j}}^2_{h}\right)(y)\geq 0.
\end{align}Note that $T(s_{j}u_{j})=T(\delta_{j}u_{j})$, so linearity gives
\begin{align}
v&=T\left(\s{j=1}{m}s_{j}u_{j}\right).
\end{align}It therefore follows by orthogonality properties of $s_{j}$ that $\s{j=1}{m}s_{j}u_{j}$ is the $h$-minimal solution to the equation $T(\cdot)=v$. By assumption, there is $u\in H$ such that $\{u, \{\vect{F}_{x}\},\{A_x\}\}$ satisfies the Hörmander $h$-estimate for the equation $T(\cdot)=v$. In particular, by ($h$)-minimality of $\s{j=1}{m}s_{j}u_{j}$, we have
\begin{align}
\norm{\s{j=1}{m}s_{j}u_{j}}^2_{h}&\leq \norm{u}^2_{h}.
\end{align}Hence, it suffices to show that
\begin{align}
\left(\s{j,k=1}{m}h(\Theta_{j\bar{k}}u_{j}, u_{k})-\norm{u}^2_{h}\right)(y)&\geq 0.
\end{align}By assumption, 
\begin{align}
\norm{u}^2_{h}&\leq \lim_{\epsilon \downarrow 0}\int_{X}\norm{\vect{F}_x}^2_{(A_x+\epsilon I_{n})^{-1}}\;d\mu_{X}(x).
\end{align}Writing out the definition of $h$, it therefore suffices to show that 
\begin{align}
\int_{X}\left(\s{j=1}{m}f(\Theta_{j\bar{k}}u_{j}, u_{k})-\lim_{\epsilon \downarrow }\norm{\vect{F}_x}^2_{(A_x+\epsilon I_{n})^{-1}}\right)\;d\mu_{X}(x)\left(y\right)&\geq 0
\end{align}For this, it suffices to verify that the integrand is non-negative at each $x\in X$. Hence, fix $x\in X$. Having fixed $x$, we may choose $t$-coordinates near $y$ such that $\Theta_{j\bar{k}}(y,x)$ is diagonal. Let $\Ld$ denote the associated diagonal matrix, and let $\ld_{j}$ denote its $j$th entry (on the main diagonal). Thus, it suffices to show that
\begin{align}
\left(\s{j=1}{m}f(u_j, u_j)\ld_{j}-\lim_{\epsilon \downarrow 0}\norm{\vect{F}}^2_{(A+\epsilon I_{n})^{-1}}\right)(y,x)\geq 0.\label{eq: inequality to be proved, with epsilon, in first main result}
\end{align}In the above sum, we may assume that each $f(u_{j}, u_{j})(y,x)\neq 0$, else we reduce the sum. By assumption, at $y$,
\begin{align}
(A_{x}+\epsilon I_{n})-B_{x}(\Ld+\epsilon I_{m})^{-1}(B_x)^*\geq 0.
\end{align}Let $\tilde{B}$ be the (collection of local) matrix (fields) whose $j$th column is $\sqrt{f(u_j, u_j)}\vect{b}_{j}$, and let $\tilde{\Ld}_{\epsilon}$ be the diagonal matrix whose $j$th entry is $f(u_j, u_j)(\ld_j+\epsilon)$. Then,
\begin{align*}
B^*_x(\Ld+\epsilon I_{m})^{-1}(B_x)^*&=\tilde{B}_x(\tilde{\Ld}_{\epsilon})^{-1}(\tilde{B}_x)^{*},
\end{align*}so it follows by (the proof of) Corollary \ref{cor: corollary to main linear algebraic result} that 
\begin{align}
\norm{\vect{F}_{x}}^2_{(A_x+\epsilon I_{n})^{-1}}(y)&\leq \s{j=1}{m}f(u_j, u_j)(\ld_j+\epsilon).
\end{align}Substituted into the above, this shows that \eqref{eq: inequality to be proved, with epsilon, in first main result} holds, which completes the proof.
\end{proof}

\subsection{Variations of holomorphic sections of quotient bundles: Proof of Theorem 1.8}
Next, we prove Theorem \ref{thm: second main result on psh properties of holomorphic sections of Q under certain assumptions}. The proof is analogous to that of Theorem \ref{thm: first main theorem, Nakano semipositivity of subbundle under Hormander condition}, but we shall additionally need from \cite{Tai1}, the quotient bundle curvature formula (Theorem 4.7 in \cite{Tai1}) and a variational formula (equation (4.60) in \cite{Tai1}).

\begin{proof}[Proof of Theorem \ref{thm: second main result on psh properties of holomorphic sections of Q under certain assumptions}]
Let $u$ be a holomorphic (local) section of $H$. By the variational formula (4.60) in \cite{Tai1}, we have 
\begin{align}
\fancy{V}^{1,\bar{1}}_{\tilde{f}(u)}&=\tilde{f}''(\norm{u}^2_h)\abs{h(\delta u, u)}^2+\tilde{f}'(\norm{u}^2_h)\left(-h(i\Theta u,u)+\norm{\delta u}^2_h\right).
\end{align}Now, replace in the above $H$ with $Q$, $u$ with $[u]$, $h$ with $h^{Q}$, and $\Theta$ with $\Theta^{Q}$, the latter being the Chern curvature of the Chern connection on $Q$ with respect to $h^{Q}$. Multiplying the previous variational formula by $(-i)$ and using the quotient bundle curvature formula (Theorem 4.7 in \cite{Tai1}) from \cite{Tai1}, since $\tilde{f}$ is concave, we then find
\begin{align}
-i\fancy{V}^{1,\bar{1}}_{\tilde{f}(h^{Q})}([u])&\geq i\tilde{f}'(\norm{[u]}^2_{h^{Q}})\left(h(\Theta m_{l}[u], m_{l}[u])-\norm{\delta m_{l}[u]}^2_{h}\right).
\end{align}Since $\tilde{f}$ is strictly increasing, it therefore suffices to prove that
\begin{align}
i\left(h(\Theta m_{l}[u], m_{l}[u])-\norm{\delta m_{l}[u]}^2_h\right)&\geq 0.
\end{align}That is, for all $\xi:=(\xi_1,\ldots, \xi_m)$, that at each $y\in Y$,
\begin{align}
\s{j,k=1}{m}h(\Theta_{j\bar{k}}m_{l}[u], m_{l}[u])\xi_{j}\bar{\xi}_{k}-\s{j,k=1}{m}h(\delta_{j}m_{l}[u], \delta_{k}m_{l}[u])\xi_{j}\bar{\xi}_{k}\geq 0.\label{eq: ineqality we want to prove in second main result}
\end{align}For simplicity we put $u_{m}:=m_{l}[u]$. Note that since for all $v_0\in H_0$, $h(\delta_{j}u_{m}, v_0)=0$, it follows that $\s{j=1}{m}\delta_{j}\xi_{j}u_{m}$ is the $h$-minimal solution to the equation $T(\cdot)=v$. The rest of the proof now follows in exactly the same manner as the proof of Theorem \ref{thm: first main theorem, Nakano semipositivity of subbundle under Hormander condition}, taking $u_{j}:=\xi_{j}u_{m}$ for each $j$.
\end{proof}

\section{Proof of Theorems 1.5, 1.6, 1.9, and 1.10}
In this section we prove the remaining theorems from the introduction which we have yet to prove, starting with Theorem \ref{thm: Nakano semipositivity of quasi-complete Kahler, smooth case}. 

\subsection{Quasi-complete Kähler manifolds and Hörmander's $L^2$ theorem}

The first step in proving either of the remaining theorems is of course to give the definition of \ita{quasi-complete Kähler}, which appears in the statement of the theorems, but which we have yet to give. The notion of \ita{quasi-complete Kähler (manifolds)} generalizes that of complete Kähler (manifolds), and is due to Xu Wang and Bo Yong Chen (\cite{XuDiag}). The definition is as follows: 

\begin{definition}\label{def: quasi-complete Kahler}
A Kähler manifold $(X,\omega)$ is said to be \bt{(a) quasi-complete (Kähler manifold)} if there exists a family of Kähler manifolds $\{(X_{j}, \omega_{j})\}_{j}$ and a family of smooth $[0,1]$-valued functions $\{\chi_{j}\}_{j}$ on $X$ such that the following properties are satisfied:
\begin{enumerate}[(i)]
\item For each $j$, $X_{j}\sub X_{j+1}$, $X_{j}$ is open in $X$, and $\bigcup_{j}X_{j}=X$. 
\item For each $j$, $\omega_{j}\geq \omega$ on $X_{j}$, and for each compact subset $K$ of $X$, 
\begin{align}
\lim_{j\to \infty}\sup_{K}\abs{\omega_{j}-\omega}_{\omega}=0.
\end{align}
\item Each $\chi_{j}$ has compact support in $X_{j}$, satisfies
\begin{align}
\lim_{j\to \infty}\sup_{X_{j}}\abs{\dbar \chi_{j}}_{\omega_{j}}=0,
\end{align}and that for each compact subset $K$ of $X$, there is $j=j(K)$ such that $\chi_{j}|_{K}\equiv 1$ for all $j\geq j(K)$. 
\end{enumerate}
In the case that $(X,\omega)$ is quasi-complete Kähler, and $\{(X_{j},\omega_{j})\}$ and $\{\chi_{j}\}$ are as above, satisfying the needed criteria,  we will refer to $\{(X_{j},\omega_{j},\chi_{j})\}_{j}$ as an \bt{approximation family for $(X,\omega)$}. In this case, we will also frequently write $(X,\omega)=\{(X_{j},\omega_{j},\chi_{j})\}_{j}$.
\end{definition}
All complete Kähler manifolds are quasi-complete Kähler manifolds. In fact, if $X$ admits a complete Kähler metric, then $(X,\omega)$ is quasi-complete Kähler. Indeed, suppose that $\hat{\omega}$ is a complete Kähler metric on $X$. It suffices to observe that $\{(X,\omega_{j},\chi_{j})\}$ is an approximation family for $(X,\omega)$ where, $\omega_{j}:=\fr{1}{j^{k}}\hat{\omega}+\omega$, $\chi_{j}:=\chi(\fr{\rho}{j^{k+1}})$, for $k\geq 1$, $\chi$ is a smooth $[0,1]$-valued function on $\set{R}$ satisfying that $\chi|_{(-\infty, 1]}\equiv 1, \chi_{[a,\infty)}=0$ for some $a>1$, $\abs{\chi'}\leq 1$, and $\rho$ is a smooth exhaustion function for $X$ such that $\abs{d\rho}_{\hat{\omega}}\leq 1$ (which exists since $\hat{\omega}$ is complete). One may further show that if $Y$ is a complete Kähler manifold, $\omega$ is a Kähler metric on $Y$, and $S$ is an analytic subset of $Y$, then $(Y\backslash S,\omega)$ is quasi-complete Kähler (\cite{XuDiag}). In fact, there is the following result:

\begin{lemma}\label{le:qck} If $(X,\omega)$ is quasi-complete Kähler, then $(X\backslash S, \omega)$ is also quasi-complete Kähler for every analytic subset $S$ in $X$. 
\end{lemma}
We include a proof here, courtesy of Xu Wang, since we shall have use for it later (see the proof of Theorem \ref{thm: singular Hormander} below):

\begin{proof} By Lemma 2.12 in \cite{Demailly12}, we know that there exists $\psi<-1$ on $X$ such that
$$
S=\{x\in X: \psi(x) =-\infty\}, \ \ \psi\in C^\infty(X\setminus S)
$$
and $i\partial\dbar \psi+ \theta \geq 0$ for some real smooth $(1,1)$-form $\theta$ on $X$. Note that
\begin{equation*}
    i\partial\dbar(-\log-\psi) \geq -A+i\partial\log(-\psi)\wedge\dbar\log(-\psi), \ \ \ A:=\max\{0, \theta\}.
\end{equation*}
Choose a decreasing sequence of sufficiently small positive numbers $a_k$ such that
$$
\lim_{k\to \infty} 2^k a_k^2  =0
$$
and $ 5^ka_k^2 A<\omega$ on $X_k$ (note that one may take an approximation family $(X_j, \omega_j, \chi_j)$ for $(X, \omega)$ such that each $X_j$ is relatively compact). Take $\kappa\in C^{\infty}(\mathbb R,[0,1])$ such that $\kappa\equiv1$ on $(-\infty,1/2)$ and $\kappa\equiv0$ on $(1,\infty)$. Put
\begin{equation*}
    w_{j,k}:=\omega_j+2^{-k}\omega+2^ka_k^2\cdot i\partial\dbar(-\log-\psi)
\end{equation*}
and
\begin{equation*}
    \chi_{j,k}:=\chi_j \cdot \kappa(a_k\log-\psi).
\end{equation*}
Then $(X_j\setminus S, ~ w_{j,j}, ~ \chi_{j,j})$ is an approximation family for $(X\setminus S, ~ \omega)$.
\end{proof}

The next step is that we have a Hörmander $L^2$-type of theorem on quasi-complete Kähler manifolds for the $\dbar$-equation. We now discuss this; it is known that we have this on complete Kähler manifolds (see \cite{Bodbar}). To set the stage, let $(X,\omega):=\{(X_j,\omega_j, \chi_{j})\}_{j}$ be an $n$-dimensional quasi-complete Kähler manifold, $(L,e^{-\phi})$ be a positive holomorphic line bundle over $X$, and let $\ip{\cdot}{\cdot}$ denote the $L^2$-type of inner product on $L$-valued $(n,q)$-forms on $X$ induced by $\omega$ and $e^{-\phi}$, given by
\begin{align}
\ip{u}{v}&:=\int_{X}(u,v)_{\omega}e^{-\phi}
\end{align}for all $L$-valued $(n,q)$-forms $u$ and $v$ on $X$. Now, in general, if $u$ is an $L$-valued $(n,q)$-form, the hard Lefschetz theorem provides the existence of a unique $L$-valued (and $\omega$-primitive) $(n-q,0)$-form $\gamma_u$, such that
\begin{align}
u&=\gamma_{u}\wedge \omega_{q}.
\end{align}Let $\Theta^{L}$ denote the Chern curvature of the Chern connection on $L$ with respect to the metric $e^{-\phi}$. We define what we shall call the ''$B$-operator'', denoted by $B$, on ($L$-valued) $(n,q)$-forms by:
\begin{align}
Bu:=i\Theta^{L}(\gamma_{u}\wedge \omega_{q-1}).
\end{align}One may then check that if we write (locally) with respect to a local orthonormal frame $e_1,\ldots, e_{n}$ for $\Ld^{1,0}(T^*X)$, the holomorphic cotangent bundle of $X$, $i\Theta^{L}$ as
\begin{align}
i\Theta^{L}&=\s{j=1}{n}\ld_{j}e_{j}\wedge \bar{e}_{j}, 
\end{align}then
\begin{align}
B(e\wedge \bar{e}_{J})&=\left(\s{j\in J}{}\ld_{j}\right)e\wedge \bar{e}_{J},
\end{align}where $e:=e_1\wedge\cdots\wedge e_{n}$ and $\bar{e}_{J}:=\bar{e}_{j_1}\wedge\cdots\wedge \bar{e}_{j_{q}}$ for $J=(j_1,\ldots, j_q)$ a strictly increasing multiindex of length $q$. In particular it follows that $B^{-1}$, the inverse operator of $B$, exists when $(L,e^{-\phi})$ is positive (and that it is given by, using the same notation as above, $e\wedge \bar{e}_{J}\mapsto \left(\s{j\in J}{}\ld_{j}\right)^{-1}e\wedge \bar{e}_{J}$). Let us choose $z:=(z_1,\ldots, z_n)$ as generic (local) holomorphic coordinates on $X$, put $\dbar_{z}:=\dbar$, and denote the formal adjoint of $\dbar$ with respect to the metric $\ip{\cdot}{\cdot}$ by $\dbar^*$. Let us denote by $\norm{\cdot}$ the norm induced by $\ip{\cdot}{\cdot}$. Then we have the following Hörmander $L^2$-type of theorem for the $\dbar$-equation, on quasi-complete Kähler manifolds under positive holomorphic line bundles:

\begin{theorem}\label{thm: Hormander on quasi-complete Kahler manofolds, smooth case}
With the above notation and set-up, let $v$ be a $\dbar$-closed $L$-valued $(n,q+1)$-form such that $\ip{B^{-1}v}{v}<\infty$. Then there exists an $L$-valued $(n,q)$-form $u$ such that $\dbar u=v$ in the sense of currents, and
\begin{align}
\norm{u}^2&\leq \ip{B^{-1}v}{v}.
\end{align}
\end{theorem}
As mentioned earlier, the theorem is known in the case that $(X,\omega)$ is complete Kähler. The quasi-complete Kähler case for $q=0$ was taught to us by Xu Wang (\cite{XuDiag}). The following proof is therefore essentially due to him:

\begin{proof}[Proof of Theorem \ref{thm: Hormander on quasi-complete Kahler manofolds, smooth case}]
Let $\square_{j}$ denote the $\dbar$-Laplace operator on $X_j$ with respect to the metric $\ip{\cdot}{\cdot}_{j}$, where the latter is defined like $\ip{\cdot}{\cdot}$, except that we replace $X$ with $X_{j}$ and $\omega$ with $\omega_{j}$. Let similarly $B_{j}$ denote the ''$B$-operator'' with respect to $\ip{\cdot}{\cdot}_{j}$. We first solve the ($\dbar$-Laplace) equation $\square_{j}(\cdot)=v$ in the sense of currents and obtain a solution $u_{j}$ satisfying
\begin{align}
\norm{\dbar u_{j}}^2_{j}+\norm{\dbar^*u_{j}}_{j}&\leq \ip{B^{-1}_{j}v}{v},\label{eq: apriori estimate in proof of solving j-dbar Laplace equation}
\end{align}where we denote by $\norm{\cdot}_{j}$ the norm induced by $\ip{\cdot}{\cdot}_{j}$, and $B^{-1}_{j}=(B_{j})^{-1}$. Since $\dbar v=0$, we get that 
\begin{align}
\dbar \square_{j}u_{j}&=\square_{j}\dbar u_{j}=\dbar\dbar^*\dbar u_{j}=0.
\end{align}Hence, by elliptic regularity, we may assume that $\dbar u_{j}$ is smooth. Using smoothness and compact support we can then write
\begin{align}
\norm{\chi_{j}\dbar^*\dbar u_{j}}^2_{j}&=\ip{\dbar (\chi_{j}^2\dbar u_{j})}{\dbar^*\dbar u_{j}}_{j}.
\end{align}We use the Leibniz rule to compute $\dbar(\chi^2\dbar u_j)$. The Cauchy-Schwarz inequality (and \eqref{eq: apriori estimate in proof of solving j-dbar Laplace equation}) then give(s)
\begin{align}
\norm{\chi_{j}\dbar^*\dbar u_{j}}^2_{j}&\leq 2\sup_{X_{j}}\abs{\dbar\chi_{j}}_{\omega_{j}}\sqrt{\ip{B_{j}^{-1}v}{v}}.
\end{align}
Now fix $X_{0}\Subset X$. We will denote restricted norms and inner products to $X_0$ using subscripts. Letting $j\to  \infty$, we find that
\begin{align}
\lim_{j\to \infty}\norm{\dbar^*\dbar u_{j}}^2_{j,X_0}=0.
\end{align}Thus on $X_0$, $\dbar(\dbar^*u_{j})\to v$ as a current as $j\to \infty$. Note that $\norm{\dbar^*u_{j}}^2_{j,X_0}\leq \ip{B^{-1}v}{v}_{X_0}$, so we get that $\{\dbar^* u_{j}|_{X_0}\}_{j}$ is a uniformly bounded sequence with respect to the $\norm{\cdot}_{X_0}$-norm. Hence it admits a convergent subsequence. Let us denote by $u_{X_0}$ the weak limit of this convergent subsequence. Then $\dbar u_{X_0}=v$ on $X_0$, and $\norm{u_{X_0}}_{X_0}^2\leq \ip{B^{-1}v}{v}_{X_0}$. We now use a diagonal argument. Let $\{U_{k}\}_{k}$ be an increasing sequence of relatively compact subsets of $X$. To each $U_{k}$ we have by the previous argument a $u_{k}$ such that $\dbar u_{k}=v$ in the sense of currents on $U_{k}$ and $\norm{u_{k}}^2_{U_{k}}\leq \ip{B^{-1}v}{v}_{U_k}$. Consider first the sequence $\{u^{(1)}_{j}\}_{j}$ given by restricting, for large $j$, each $u_{j}$ to $U_{1}$. Then let $\{u^{(2)}_{j}\}_{j}$ be the sequence given by restricting each $u^{(1)}_{j}$, again for large $j$, to $U_{2}$, and so on. Finally, let $u$ be the weak limit of the diagonal sequence $\{u^{(j)}_{j}\}_{j}$. Then $u$ satisfies that $\dbar u=v$ in the sense of currents on $X$ and $\norm{u}^2\leq \ip{B^{-1}v}{v}$.
\end{proof}
We shall need a variant of Theorem \ref{thm: Hormander on quasi-complete Kahler manofolds, smooth case} where $(L, e^{-\phi})$ is semi positive. In this case, $B^{-1}$ may not make sense, and the condition $\ip{B^{-1}v}{v}<\infty$ needs to be replaced with the condition 
\begin{align}
\lim_{\epsilon \downarrow 0}\ip{(B+\epsilon (q+1))^{-1}v}{v}<\infty.
\end{align} Since the proof is analogous to that of Theorem \ref{thm: Hormander on quasi-complete Kahler manofolds, smooth case}, we shall content ourselves here with simply stating the result:

\begin{theorem}\label{thm: Hörmander's theorem on quasi-complete Kahler manifolds, smooth case, semipositive case}
Let $(X, \omega)$ be an $n$-dimensional quasi-complete Kähler manifold, and let $(L,e^{-\phi})$ be an ample holomorphic line bundle over $X$. Let $\Theta$ denote the Chern curvature of $L$, and let $B$ be defined by 
\begin{align}
B(\omega_{q}\wedge \gamma_u):=i\Theta \gamma_{u}\wedge \omega_{q-1}
\end{align}for all $L$-valued $(n,q)$-forms $u=\gamma_u\wedge \omega_{q}$. Let $v$ be a $\dbar$-closed $L$-valued $(n,q+1)$-form such that $\lim_{\epsilon \to 0}\ip{(B+(q+1)\epsilon )^{-1}v}{v}<\infty$, where $\ip{\cdot}{\cdot}$ is the $L^2$-type of inner product on $L$-valued differential forms induced by $e^{-\phi}$ and $\omega$. Then there exists an $L$-valued $(n,q)$-form $u$ such that $\dbar u=v$ in the sense of currents, and 
\begin{align}
\norm{u}^2<\lim_{\epsilon \downarrow 0}\ip{(B+\epsilon (q+1))^{-1}v}{v},
\end{align}where $\norm{\cdot}$ denotes the norm induced by $\ip{\cdot}{\cdot}$.
\end{theorem}

\subsection{Nakano (semi)positivity for quasi-complete Kähler manifolds and (semi)positive line bundles: Proof of Theorem 1.5}
We are now ready for the proof of Theorem \ref{thm: Nakano semipositivity of quasi-complete Kahler, smooth case}, for which we shall use Theorem \ref{thm: Hörmander's theorem on quasi-complete Kahler manifolds, smooth case, semipositive case}:

\begin{proof}[Proof of Theorem \ref{thm: Nakano semipositivity of quasi-complete Kahler, smooth case}]
Let $u_1,\ldots, u_m\in H_0$ be given, and let $\epsilon>0$. The idea is to use Theorem \ref{thm: first main theorem, Nakano semipositivity of subbundle under Hormander condition} with $T:=\dbar_z$, taking as usual $z:=(z_1,\ldots, z_n)$ as generic (local) holomorphic coordinates on $X$. We take, also as usual, $t:=(t_1,\ldots, t_m)$ as generic (local) holomorphic coordinates on $Y$. For each $x\in X$, let $A_{x}$ denote the matrix field associated to $\partial_{z}\dbar_{z}\theta(\cdot,x)$. We write $\theta_{t_{j}}$ for $\del{\theta}{t_{j}}$, $\partial_{t_{j}}$ for $\del{}{t_{j}}$, $\theta_{t_{j}\bar{t}_{k}}$ for $\del{^2\theta}{t_{j}\bar{t}_{k}}$,  $\theta_{t_{j}\bar{z}_{k}}$ for $\del{}{\bar{z}_{k}}\theta_{t_{j}}$, and so on. We have seen that in our situation, 

\begin{align}
\delta_{j}&=\partial_{t_{j}}-\theta_{t_{j}},
\end{align}so that 
\begin{align}
\dbar_{z}(\delta_{j}u_{j})&=-\dbar_{z}(\theta_{t_{j}})\wedge u_{j}.
\end{align}We put $v:=\s{j=1}{m}\dbar_{z}(\delta_{j}u_{j})$, and let, as in Theorem \ref{thm: first main theorem, Nakano semipositivity of subbundle under Hormander condition}, $M_{x}$ be the matrix field associated with $\Theta(\cdot,x)$. Thus, letting $M$ be defined by $M(\cdot,x):=M_{x}$, the $(j,k)$th entry of $M$ is $\Theta_{j\bar{k}}$. We have seen that we have (and this also follows from the expression of $\delta$ above) 
\begin{align}
\Theta_{j\bar{k}}=\theta_{t_{j}\bar{t}_{k}}.
\end{align} 
Let 
\begin{align}
b_{s}:=\dbar_{z}\theta_{t_{s}},
\end{align} and let $\vect{b}_{s}$ be the (collection of) (local) vector fields identified with $b_{s}$. It follows by the plurisubharmonicity of $\theta$ that, at each $x\in X$,
\begin{align}
 (A_{x}+\epsilon I_{n})-B_{x}(M_{x}+\epsilon I_{m})^{-1}(B_{x})^*\geq 0.\label{eq: linear alg condition in proof of Nakano positivity on qck , smooth case}
 \end{align}
Let $B$ be the ''$B$-operator'' with respect to the metric $h$, and let $e_1,\ldots, e_n$ be a (local) orthonormal frame for $\Ld^{1,0}(T^*X)$ with respect to which we write
\begin{align}
\partial_{z}\dbar_{z}\theta&=\s{j=1}{n}\eta_{j} e_{j}\wedge \bar{e}_{j}\\
u_{j}&=\s{|J_{j}|=q}{'}(-1)\alpha^{J_{j}}e\wedge \bar{e}_{J_{j}}\\
\dbar_{z}\theta_{t_{j}}&=\s{l=1}{n}\beta^{l}_{j}\bar{e}_{l}.
\end{align}Then we get
\begin{align}
v&=\s{j=1}{m}\s{|J_j|=q}{'}\s{l=1}{n}\alpha^{J_{j}}\beta^{l}_{j}e\wedge \bar{e}_{l}\wedge \bar{e}_{J_{j}},
\end{align}
and we may assume that $l\not\in J_{j}$ in the above sum. Computing, we now find
\begin{align}
((B+(q+1)\epsilon)^{-1}v,v)_{\omega}&=\left((B+(q+1)\epsilon)^{-1}\left(\s{j,J_j, l}{}\alpha^{J_j}\beta^{l}_{j}e\wedge \bar{e}_{l}\wedge \bar{e}_{J_{j}}\right),\s{k,J_{k},p}{}\alpha^{J_{k}}\beta^{p}_{k}e\wedge \bar{e}_{p}\wedge \bar{e}_{J_{k}}\right)_{\omega}\notag\\
&=\left(\s{j, J_{j},l}{}\alpha^{J_j}\beta^{l}_{j}\left(\s{r\in \{l\}\cup J_{j}}{}(\eta_{r}+\epsilon)\right)^{-1}e\wedge \bar{e}_{l}\wedge \bar{e}_{J_{j}},\s{k,J_{k},p}{}\alpha^{J_{k}}\beta^{p}_{k}e\wedge \bar{e}_{p}\wedge \bar{e}_{J_{k}}\right)_{\omega}\notag\\
&\leq \left(\s{j,J_j,l}{}\fr{1}{\eta_{l}+\epsilon}\alpha^{J_j}\beta^{l}_{j}e\wedge \bar{e}_{l}\wedge \bar{e}_{J_{j}},\s{k,J_{k},p}{}\alpha^{J_{k}}\beta^{p}_{k}e\wedge \bar{e}_{p}\wedge \bar{e}_{J_{k}}\right)_{\omega}\notag\\
&=\s{j,k}{}(e,e)_{\omega}\left(\s{J_j}{'}\alpha^{J_j}\bar{e}_{J_j},\s{J_k}{'}\alpha^{J_k}\bar{e}_{J_k}\right)_{\omega}\left(\s{l}{}\fr{1}{\eta_{l}+\epsilon}\beta^{l}_{j}\bar{e}_{l},\s{p}{} \beta^{p}_{k}\bar{e}_{p}\right)_{\omega}\notag\\
&\leq \s{j,k}{}\abs{\s{J_j}{'}\alpha^{J_j}\bar{e}_{J_j}}_{\omega}\abs{\s{J_k}{'}\alpha^{J_k}\bar{e}_{J_k}}_{\omega}\left(\s{l}{}\fr{1}{\eta_{l}+\epsilon}\beta^{l}_{j}\bar{e}_{l}, b_{k}\right)_{\omega}\notag\\
&=\s{j,k}{}\left(\s{l}{}\fr{1}{\eta_{l}+\epsilon}\beta^{l}_{j}\bar{e}^{l}\abs{u_{j}}_{\omega}, b_{k}\abs{u_{k}}_{\omega}\right)_{\omega}\notag\\
&=\left((B+\epsilon)^{-1}\s{j=1}{m}b_{j}\abs{u_{j}}_{\omega}, \s{k=1}{m}b_{k}\abs{u_k}_{\omega}\right)_{\omega}.\label{eq: computation of B inverse in proof of Nakano on qck}
\end{align}Let $\vect{F}:=\s{s=1}{m}\abs{u_s}_{\omega}\vect{b}_{s}$. Then the (last) right-hand side in \eqref{eq: computation of B inverse in proof of Nakano on qck} is equal to $\norm{\vect{F}}_{(A+\epsilon I_{n})^{-1}}^2$. That is, we have (locally)
\begin{align}
((B+(q+1)\epsilon)^{-1}v,v)_{\omega}&\leq \norm{\vect{F}}_{(A+\epsilon I_{n})^{-1}}^2.\label{eq: bound on B plus q plus 1 epsilon, by norm of F in proof of Nakanos emipostitiv for cqk, smooth case}
\end{align}
Integrating this and applying Corollary \ref{cor: corollary to main linear algebraic result}, we find that $\lim_{\epsilon \downarrow 0}\ip{(B+(q+1)\epsilon)^{-1}v}{v}<\infty$. Hence, we may apply Theorem \ref{thm: Hörmander's theorem on quasi-complete Kahler manifolds, smooth case, semipositive case} and find a $U\in H$ such that $\dbar_z(U)=v$ and 
\begin{align}
\norm{U}^2_{h}&\leq \lim_{\epsilon \downarrow 0}\ip{(B+(q+1)\epsilon )^{-1}v}{v}.
\end{align}Incidentally, \eqref{eq: bound on B plus q plus 1 epsilon, by norm of F in proof of Nakanos emipostitiv for cqk, smooth case} also shows that $\{U, \{\vect{F}_x\}, \{A_x\}\}$ satisfies the Hörmander $h$-estimate for the equation $\dbar_z(\cdot)=v$, so the assertion follows by Theorem \ref{thm: first main theorem, Nakano semipositivity of subbundle under Hormander condition}.
\end{proof}

\subsection{A Hörmander $L^2$-type of theorem with singular weights: Proof of Theorem 1.7 }
Having proved Theorem \ref{thm: Nakano semipositivity of quasi-complete Kahler, smooth case}, we proceed to prove the variant of it where $\phi$ may possibly be singular. That is, Theorem \ref{thm: Nakano semipositivity for singular metrics of special type on qck Kahler manifolds}. The first step in proving Theorem \ref{thm: Nakano semipositivity for singular metrics of special type on qck Kahler manifolds}, is to prove a singular variant of Theorem \ref{thm: Hormander on quasi-complete Kahler manofolds, smooth case} in this setting. That is, Theorem \ref{thm: singular Hormander}. We now do this. For parts of the proof it might be helpful to refer back to the proof of Theorem \ref{thm: Hormander on quasi-complete Kahler manofolds, smooth case}. 

\begin{proof}[Proof of Theorem \ref{thm: singular Hormander}]
Since $\theta^{y}$ is independent of $t$ when $t\geq 0$, we may assume that $t<0$. Fix $X_{j}\Subset X$ and $t<0$; we are of course later going to let $j\to \infty$. Let $\phi_0$ be a fixed smooth metric on $L$. The idea is to apply the regularization of quasi-plurisubharmonic currents (see Theorem 16.3 in \cite{D}) to $\theta-\phi_0$ (or rather $\theta^{t}-\phi_0)$), but we cannot do this directly. Instead, note that if we define $\psi:=\phi+\ld G$, then $\phi$ and $\psi$ are plurisubharmonic by assumption. We shall apply the regularization to $\phi-\phi_0$ and $\psi-\phi_0$ simultaneously and separately. The idea is that since 
\begin{align}
G=\fr{\phi-\psi}{\ld}=\fr{\phi\pm \phi_0-(\psi\pm \phi_0)}{\ld},
\end{align} approximations of $\phi-\phi_0$ and $\psi-\phi_0$, give approximations for $G$, and hence $\theta-\phi_0$. Let $\epsilon_{j}>0$. Applying the regularization to $\phi-\phi_0$ and $\psi-\phi_0$, we find decreasing sequences $\{\phi^{j}_{\nu}\}_\nu$ and $\{\psi^{j}_{\nu}\}_\nu$ converging to $\phi-\phi_0$ and $\psi-\phi_0$ respectively, and such that
\begin{align}
i\partial\dbar(\phi^{j}_{\nu}+\phi_{0}), i\partial\dbar (\psi_{\nu}^{j}+\phi_{0})\geq -\epsilon_{j}\omega.\label{eq: estimate of regularizations, quasi-plurisubharmonicity}
\end{align}The $\partial$ and $\dbar$ operators are here with respect to $X$, the corresponding operators for $Y(=\set{R})$ will have no play in the proof. We put 
\begin{align}
G_{\nu}^{j}:=\fr{\phi^{j}_{\nu}-\psi^{j}_{\nu}}{\ld}
\end{align}
 and 
 \begin{align}
 \theta^{j}_{\nu}:=\phi^{j}_{\nu}+\chi\left(G_{\nu}^{j}-t\right)
 \end{align} (as with $\theta$, $\theta^{j}_{\nu}$ really depends on $t$, but we will neglect to write this to avoid too messy a notation). A direct computation now gives
\begin{align}
i\partial\dbar (\theta^{j}_{\nu}+\phi_0)&=i\partial\dbar(\phi^{j}_{\nu}+\phi_0)\left(1-\fr{\chi'(G_{\nu}^{j}-t)}{\ld}\right)+\fr{\chi'(G_{\nu}^{j}-t)}{\ld}i\partial\dbar(\psi^{j}_{\nu}+\phi_0)\notag\\
&+\chi''(G_{\nu}^{j}-t)i\partial G_{\nu}^{j}\wedge \dbar G_{\nu}^{j}.
\end{align}Since $0\leq \chi'\leq \ld$, \eqref{eq: estimate of regularizations, quasi-plurisubharmonicity}, gives
\begin{align}
i\partial\dbar(\theta^{j}_{\nu}+\phi_0)&\geq -2\epsilon_{j}\omega+\chi''(G_{\nu}^{j}-t)i\partial G_{\nu}^{j}\wedge \dbar G_{\nu}^{j}.\label{eq: inequality/estimate for theta j nu}
\end{align}Our approximations $\phi^{j}_{\nu}$ (and $\psi^{j}_{\nu}$) are smooth outside analytic subsets $Z^{j}_{\nu}:=\{\phi^{j}_{\nu}=-\infty\}\sub X_{j}$. We are going to solve a perturbed $\dbar$-Laplace equation outside $Z^{j}_{\nu}$, and for this we will mimic the proof of Lemma \ref{le:qck}. By choosing $\phi_0$ sufficiently large, we may assume that
\begin{align}
\phi^{j}_{\nu}<-1+\ld t,
\end{align} which gives $\theta^{j}_{\nu}<-1$. With this, let $A:=\max\{i\partial\dbar \phi_{0}+\omega\}$, and let $\{a_{k}\}$ be a decreasing sequence of sufficiently small positive real numbers such that 
\begin{align}
2^{k}a_{k}\to 0
\end{align} as $k\to \infty$, and 
\begin{align}
5^{k}a_{k}^2 A<\omega
\end{align} for large $k$. Let $\kappa$ be a smooth $[0,1]$-valued function on $\set{R}$ such that $\kappa|_{(-\infty, 1/2)}\equiv 1$ and $\kappa_{(1,\infty)}\equiv 0$. Define
\begin{align}
\omega_{j,\nu,k}&:=\omega_{j}+2^{-k}\omega+2^{k}a_{k}^2\partial\dbar(-\ln(-\theta^{j}_{\nu}))\\
\chi_{j,\nu,k}&:=\chi_{j}\kappa\left(a_{k}\ln(-\theta^{j}_{\nu})\right).
\end{align}We now mimic the first part of the proof of Theorem \ref{thm: Hormander on quasi-complete Kahler manofolds, smooth case} (actually the proof of Theorem \ref{thm: Hörmander's theorem on quasi-complete Kahler manifolds, smooth case, semipositive case}) for the family $(X_{j}\backslash Z^{j}_{\nu}, \omega_{j,\nu,k}, \chi_{j,\nu,k})$, for fixed $j$. Let $\square_{j,\nu,k}$ be the $\dbar$-Laplace operator associated with $X^{j}\backslash Z^{j}_{\nu}$ and $\omega_{j,\nu,k}$. We find a solution $v_{j,\nu,k}$ such that
\begin{align}
(\square_{j,\nu,k}+2\epsilon_{j}(q+1))v_{j,\nu,k}&=\chi''(G_{\nu}^{j}-t)\dbar G_{\nu}^{j}\wedge u:=c_{j,\nu,k}\label{eq: perturbed Laplace j,nu k equation}
\end{align}and
\begin{align}
\norm{\dbar^*v_{j,\nu,k}}_{j,\nu,k}^2+\norm{\dbar v_{j,\nu,k}}^2_{j,\nu,k}+2\epsilon_{j}(q+1)\norm{v_{j,\nu,k}}^2_{j,\nu,k}&\leq \ip{(B^{\epsilon}_{j,\nu, k})^{-1}c_{j,\nu,k}}{c_{j,\nu,k}}_{j,\nu,k},
\end{align}where $B^{\epsilon_{j}}_{j,\nu,k}$ is the ''$B$-operator'' associated with $\ip{\cdot}{\cdot}_{j,\nu,k}$ and the positive ''curvature form'' $i\partial\dbar(\theta^{j}_{\nu}+\phi_0)+2\epsilon_{j}\omega_{j,\nu,k}$, and where $\ip{\cdot}{\cdot}_{j,\nu,k}$ is defined similar as $h$, except that we replace $\omega$ with $\omega_{j,\nu,k}$, $X$ with $X_{j}\backslash Z^{j}_{\nu}$, and $\theta$ with $\theta^{j}_{\nu}+\phi_0$, Let us denote by $B_{j,\nu,k}$ the ''$B$-operator'' associated with $i\partial\dbar(\theta^{j}_{\nu}+\phi_0)$. Then $B_{j,\nu,k}^{\epsilon_{j}}=B_{j,\nu,k}+2\epsilon_{j}(q+1)$. Let $b:=\dbar G_{\nu}^{j}$ and let $\vect{b}$ be the vector (field) identified with $b$. Let $\vect{F}:=\chi''(G_{\nu}^{j}-t)\abs{u}_{\omega_{j,\nu,k}}\vect{b}$, and let $A$ be the matrix (field) identified with $\partial\dbar(\theta^{j}_{\nu}+\phi_{0})+2\epsilon_{j}\omega_{j,\nu,k}$. By \eqref{eq: inequality/estimate for theta j nu} we have
\begin{align}
A\geq \vect{b}\chi''(G_{\nu}^{j}-t)\vect{b}^*.
\end{align}Hence, as in the proof of Theorem \ref{thm: Nakano semipositivity of quasi-complete Kahler, smooth case}, it follows by Corollary \ref{cor: corollary to main linear algebraic result} (in our case now, $m=1$, which simplifies things)
\begin{align}
((B^{\epsilon_{j}}_{j,\nu,k})^{-1}c_{j,\nu,k},c_{j,\nu,k})_{\omega_{j,\nu,k}}&\leq \norm{\vect{F}}^{2}_{A^{-1}}\leq \chi''(G_{\nu}^{j}-t)\abs{u}^2_{\omega_{j,\nu,k}}.
\end{align}Integrating we therefore get that $v_{j,\nu,k}$ satisfies
\begin{align}
\norm{\dbar^*v_{j,\nu,k}}_{j,\nu,k}^2+\norm{\dbar v_{j,\nu,k}}^2_{j,\nu,k}+2\epsilon_{j}(q+1)\norm{v_{j,\nu,k}}^2_{j,\nu,k}&\leq \ip{\chi''(G_{\nu}^{j}-t)u}{u}_{j,\nu,k}.\label{eq: partial final estimate of v j nu k}
\end{align}Now, $\dbar c_{j,\nu,k}=0$, so applying $\dbar$ to \eqref{eq: perturbed Laplace j,nu k equation} we get
\begin{align}
0&=(\square_{j,\nu,k}+2\epsilon_{j}(q+1))\dbar v_{j,\nu,k}.
\end{align}Hence we may assume that $\dbar v_{j,\nu,k}$ is smooth on $X^{j}\backslash Z^{j}_{\nu}$. Thus
\begin{align}
0&=\ip{\chi_{j,\nu,k}^2\dbar v_{j,\nu,k}}{(\dbar\dbar^*\dbar+2\epsilon_{j}(q+1)\dbar)v_{j,\nu,k}}_{j,\nu,k},
\end{align}which gives
\begin{align}
\norm{\chi_{j,\nu,k}\dbar^*\dbar v_{j,\nu,k}}^2_{j,\nu,k}&=\ip{\dbar v_{j,\nu,k}}{2\chi_{j,\nu,k}\wedge \dbar \chi_{j,\nu,k}\wedge \dbar^*\dbar v_{j,\nu,k}}_{j,\nu,k}-2\epsilon_{j}(q+1)\norm{\chi_{j,\nu,k}\dbar v_{j,\nu,k}}^2_{j,\nu,k}.
\end{align}Hence, 
\begin{align}
\norm{\chi_{j,\nu,k}\dbar^*\dbar v_{j,\nu,k}}_{j,\nu,k}&\leq 2\norm{\dbar \chi_{j,\nu,k}}_{j,\nu,k}\norm{\dbar v_{j,\nu,k}}_{j,\nu,k}.
\end{align}Put $C_{j,\nu,k}:=\ip{\chi''(G_{\nu}^{j}-t)u}{u}_{j,\nu,k}$. Then $\norm{\dbar v_{j,\nu,k}}_{j,\nu,k}\leq \sqrt{C_{j,\nu,k}}$. Computing we have
\begin{align}
\dbar \chi_{j,\nu,k}&=\dbar \chi_{j}\kappa(a_{k}\ln(-\theta^{j}_{\nu}))+\chi_{j}\kappa'(a_{k}\ln(-\theta^{j}_{\nu}))a_{k}\dbar\ln(-\theta^{j}_{\nu}).
\end{align}Hence there exist a constant $M^{j}_{\nu}$ such that 
\begin{align}
\norm{\dbar \chi_{j,\nu,k}}_{j,\nu,k}&\leq \norm{\dbar \chi_{j}}_{j,\nu}+a_{k}M^{j}_{\nu},
\end{align}where $\norm{\cdot}_{j,\nu}:=\lim_{k\to \infty}\norm{\cdot}_{j,\nu,k}$. Put $C_{j,\nu}:=\lim_{k\to \infty}C_{j,\nu,k}$. We now let $k\to \infty$, for fixed $j$ and $\nu$. By \eqref{eq: partial final estimate of v j nu k}, $\{v_{j,\nu,k}\}_{k}$ is uniformly bounded with respect to the $L^2$-metric on $X^{j}\backslash Z^{j}_{\nu}$ determined by $\omega_{j}$ and $\theta^{j}_{\nu}+\phi_0$. Let $v_{j,\nu}$ denote the weak limit of a converging subsequence. Then 
\begin{align}
\norm{\chi_{j}\dbar^*\dbar v_{j,\nu}}_{j,\nu}\leq \norm{\dbar \chi_{j}}_{j,\nu}\sqrt{C_{j,\nu}}
\end{align}We have that $v_{j,\nu}$ satisfies a perturbed $\dbar$-Laplace equation outside the analytic subset $Z^{j}_{\nu}$. By standard results (see \cite{Demailly12}), it satisfies the equation on all of $X_{j}$. We now let $\nu\to \infty$. Let us define $C_{j}:=\ip{\chi''(G-t)u}{u}_{X_{j}}$, where the subscript notation here means that we restrict integration to over $X_{j}$. The weight that we are integrating with respect to here is $\theta$ (or rather $\theta^{t}$). Letting $\nu\to \infty$, and letting similarly as above $v_{j}$ be the weak limit of a convergent subsequence of $\{v_{j,\nu}\}_{\nu}$, we have
\begin{align}
(\dbar\dbar^*+\dbar^*\dbar+2\epsilon_{j}(q+1))v_{j}&=\chi''(G-t)\dbar G\wedge u
\end{align}on $X_{j}$, that
\begin{align}
\norm{\dbar v_{j}}_{X_j}^2+\norm{\dbar^*v_{j}}^2_{X_{j}}+2\epsilon_{j}(q+1)\norm{v_{j}}^2_{X_j}&\leq C_{j},
\end{align}and moreover that
\begin{align}
\norm{\chi_{j}\dbar^*\dbar v_j}_{X_{j}}&\leq \norm{\dbar \chi_{j}}_{X_{j}}\sqrt{C_{j}}.
\end{align}We may assume that $\epsilon_{j}\to 0$ as $j\to \infty$. From here we mimic the last part of the proof of Theorem \ref{thm: Hormander on quasi-complete Kahler manofolds, smooth case}: Let $U\Subset X$. Letting $j\to \infty$, we find
\begin{align}
\lim_{j\to \infty}\norm{\dbar^*\dbar v_{j}}^2_{j, U}=0.
\end{align}Let then $v_{U}$ denote the weak limit of a convergent subsequence of $\{\dbar^*v_{j}|_{U}\}_{j}$, consider an exhaustion of $X$ by relatively compact subsets, say $\{U_{l}\}_{l}$, and let finally $v$ be the weak limit a diagonal (sub)sequence constructed from (subsequences of) $\{v_{U_{l}}\}_{l}$. Then 
\begin{align}
\dbar v&=c
\end{align}in the sense of currents on $X$, where $c:=\chi''(G-t)\dbar G\wedge u$, and
\begin{align}
\norm{v}^2_{t}&\leq h_{t}(\chi''(G-t)u,u).
\end{align}Noting that $\chi''(G-t)=\theta_{tt}(t,\cdot)$ and $c=-(\dbar \theta_{t})\wedge u$, this completes the proof.
\end{proof}

\subsection{Nakano (semi)positivity for quasi-complete Kähler manifolds and pseudoeffective line bundles: Proof of Theorem 1.6}
We may now prove Theorem \ref{thm: Nakano semipositivity for singular metrics of special type on qck Kahler manifolds}. Apart from a slight technicality on $e^{-\theta+\phi}$ being a valid smooth weight for the hermitian metric of zero variation $h^0$, we get the theorem quite directly from Theorem \ref{thm: singular Hormander}. That $e^{-\theta+\phi}$ is indeed a valid smooth weight is explained below. Assuming this to be true for now, the proof of the theorem can be given as follows:

\begin{proof}[Proof of Theorem \ref{thm: Nakano semipositivity for singular metrics of special type on qck Kahler manifolds}]
Let $y\in Y$. We use the same notation as above, and let $t$ generically be the real part of $y$. We also use $t$ as generic coordinates on $Y$; hopefully, this will not cause any confusion. Let $u\in H_0$, and let $dt\wedge s$ be second fundamental form associated with $D$ and $D_0$, where $D$ is the Chern connection of $H$, and $D_0$ is the Chern connection of $H_0$, both with respect to $h$. Then $su$ is the $h$-minimal solution to $\dbar(su)=-\dbar\theta_{t}\wedge u$, where the $\dbar$ operator here is the one on $X$. Using Theorem \ref{thm: singular Hormander}, let $v\in H$ be such that $\dbar v=\dbar \theta_{t}\wedge u$ and such that $\norm{v}^2_{t}\leq h(\theta_{tt}u,u)(y)$. Let $\Theta$ be the Chern curvature of $D$. Then $\Theta=\theta_{tt}\wedge$. Hence we have
\begin{align}
h(\theta_{tt}u,u)(y)\geq \norm{v}^2_{h}(y)\geq \norm{su}^2_{h}(y),
\end{align}which gives that
\begin{align}
h(\theta_{tt}u,u)(t)-\norm{su}^2_{h}(t)\geq 0.
\end{align}Hence $i\Theta\geq_{N}0$ by the subbundle curvature formula, \eqref{eq: subbundle curvature formula}.\\

Alternatively, we can use Theorem \ref{thm: first main theorem, Nakano semipositivity of subbundle under Hormander condition} to give a proof. Let $u\in H_0$ be given, and let for each $x\in X$, $A_x$ be the matrix (field) associated with $\partial_z\dbar_z\theta(\cdot,x)$, and $M_{x}:=\theta_{tt}(\cdot,x)$. For each $x\in X$, we also let $B_{x}$ be the vector (field) identified with $\dbar_{z}\theta_{t}(\cdot,x)$, and put $\vect{F}_{x}:=\abs{u}_{\omega}(x)B_{x}$. Let $\epsilon>0$. Since $\theta$ is plurisubharmonic, it follows that we have, for each $x\in X$,
\begin{align}
(A_{x}+\epsilon I_{n})-B_x(M_{x}+\epsilon I_{m})^{-1}B_{x}^*\geq 0.
\end{align}Let $v:=\dbar_z(su)$. By Theorem \ref{thm: singular Hormander}, we can find $u\in H$ such that $\{u, \{\vect{F}_x\}, \{A_x\}\}$ satisfies the Hörmander $h$-equation $\dbar_z(\cdot)=v$. Hence the assertion follows by Theorem \ref{thm: first main theorem, Nakano semipositivity of subbundle under Hormander condition}.
\end{proof}

\subsubsection{Some special weights}\label{sec: special weights}
We now deal with the aforementioned technicality that $h$ is a valid smooth weighted hermitian metric as defined in \cite{Tai1}, and which in the above proof was taken for granted. For this we do not really need all the assumptions in the theorem:

Let $X$ be a complex $n$-dimensional manifold with a hermitian metric (on its complex tangent bundle), and let $\omega$ denote the fundamental form of said metric. Suppose that $dV_{X}$ is a given volume form on $X$, and also that $L$ is a given holomorphic line bundle over $X$ with a smooth hermitian metric $e^{-\phi}$. Let $G\leq 0$ be a negative function on $X$, and for each $t\in (-\infty, 0]$, let $\psi^{t}:=\chi(G-t)$, where $\chi$ is a non-negative smooth real-valued function on $\set{R}$ which vanishes on $(-\infty, 0]$. Let $H$ be the Hilbert space consisting of $L$-valued $(n,q)$-forms $u$ on $X$ such that 
\begin{align}
\norm{u}_{\phi}^2&:=\int_{X}\abs{u}^2_{\omega}e^{-\phi}dV_{X}<\infty.
\end{align}
We first show that multiplication by $e^{-\psi^{t}}$, for each $t$, is an element in $\text{Aut}(H)$, the space of complex-linear automorphisms $H\to H$. That is, so that $f:=e^{-\psi}:=t\mapsto e^{-\psi^{t}}$ is a map $\set{R}\to \text{Aut}(H)$.

\begin{proposition}
With the above notation, $f$ is a map $\set{R}\to \text{Aut}(H)$.
\end{proposition}

\begin{proof} It is clear that $f(t):=f^{t}$ is linear. We need to show that $f^{t}u\in H$ for $u\in H$, and that $f^{t}$ is continuous. Let $u\in H$. We look at $\norm{f^{t}u}^2_{\phi}$. For simplicity we shall write $dV$ for $e^{-\phi}dV_{X}$. We have 
\begin{align*}
\norm{f^{t}u}^2_{\phi}&=\int_{X}\abs{f^{t}u}^2_{\omega}dV=\int_{X}\abs{e^{-\chi(G-t)}u}^2_{\omega}dV=\int_{X}\abs{e^{-\chi(G-t)}}^2\abs{u}^2_{\omega}dV.
\end{align*}
Since $\chi\geq 0$ we obviously have $e^{-\chi}\leq 1$, so
\begin{align*}
\norm{f^{t}u}^2_{\phi}\leq \int_{X}\abs{u}^2_{\omega}dV=\norm{u}^2_{\phi}.
\end{align*}Thus, the proposition is true.
\end{proof}
Next, we shall show that for all $n\in \set{N}$, $f^{(n)}(t)\in \text{Aut}(H)$. 

\begin{lemma}
For each $t\set{R}$, $f^{(n)}(t)\in \text{Aut}(H)$.
\end{lemma}

\begin{proof}
We may assume that $t<0$, and again, linearity is clear. Let $u\in H$. We shall study $\norm{f^{(n)}(t)u}^2_{\phi}$, and we write as in the previous proof $e^{-\phi}dV_{X}$ simply as $dV$. We find
\begin{align*}
\norm{f^{(n)}(t)u}^2_{\phi}&=\int_{X}\abs{f^{(n)}(t)}\abs{u}^2_{\omega}dV=\int_{X}\abs{\fr{\partial ^{n}}{\partial t^{n}}e^{-\chi(G-t)}}\abs{u}^2_\omega dV.
\end{align*}The key observation now is that since $\chi|_{(-\infty, 0]}=0$, we have that $\fr{\partial ^{n}}{
\partial t^{n}}\left(e^{-\chi(G-t)}\right)$ when $G-t\leq 0$. Hence, 
\begin{align*}
\norm{f^{(n)}(t)u}^2_{\phi}&=\int_{G>t}\abs{\fr{\partial^n}{\partial t^{n}}e^{-\chi(G-t)}}\abs{u}^2_{\omega}dV.
\end{align*}Note that since $G\leq 0$, on $G>t$, it follows that $0\geq G>t$ so that $-t>G-t>0$. Thus, for any $\epsilon>0$, $G-t\in [-\epsilon,-t+\epsilon]$, a compact subset of $\set{R}$. Since $\fr{\partial^{n}}{\partial t^{n}}e^{-\chi(G-t)}$ is smooth, it is bounded on this compact. Thus, there exists a constant $M=M(t)$, such that 
\begin{align*}
\norm{f^{(n)}(t)u}^2_{\phi}&\leq M(t)\int_{G>t}\abs{u}^2_{\omega}dV\leq M\int_{X}\abs{u}^2_{\omega}dV=M\norm{u}^2_{\phi}.
\end{align*}This proves the lemma.
\end{proof}
Finally, we shall prove that $f:\set{R}\to \text{Aut}(H)$ is smooth, its $n$th derivative being given by $f^{(n)}$.

\begin{theorem}
$f:\set{R}\to \text{Aut}(H)$ is smooth, with $f^{(n)}$ as its $n$th order derivative. That is, more precisely, for $t, \Delta t_1,\ldots,\Delta t_n\in \set{R}$ and for all $u\in H$,
\begin{align*}
f^{(n)}(t)(\Delta t_1,\ldots, \Delta t_n)u:=\prod_{j=1}^{n}\Delta t_{j}f^{(n)}(t)u,
\end{align*} where we of course view $f^{(n)}(t)$ as a multilinear continuous map.
\end{theorem}

\begin{proof}
We may as in the previous proof assume $t<0$. The proof will be by induction on $n$, and we begin with the case that $n=1$. By definition of the derivative, what we need to show is that
\begin{align}
\lim_{h\to 0}\fr{\norm{f(t+h)-f(t)-f'(t)(h)}_{H^*}}{\abs{h}}=0.\label{eq: identity 1 in proof}
\end{align}By definition of the norm on $H^*$, it suffices to show that for all $u\in H\backslash \{0\}$, 
\begin{align}
\lim_{h\to 0}\fr{\norm{(f(t+h)-f(t)-f'(t)(h))(u)}^2_{\phi}}{\abs{h}^2}=0.\label{eq: id 2 in proof}
\end{align}Indeed, given any $F\in \text{Aut}(H)$, by definition, 
\begin{align*}
\norm{F}^2_{H^*}&=\sup_{u\in H\backslash\{0\}}\fr{\norm{F(u)}^2_{\phi}}{\norm{u}^2_{\phi}},
\end{align*}so taking the supremum in \eqref{eq: id 2 in proof} we get \eqref{eq: identity 1 in proof}.
Hence, we first take a closer look at the numerator in \eqref{eq: id 2 in proof}, namely,  $\norm{(f(t+h)-f(t)-f'(t)(h))u}^2_{\phi}$. This is equal to 
\begin{align*}
\int_{X}\abs{f(t+h)-f(t)-f'(t)h}^2\abs{u}^2_{\omega}dV,
\end{align*}where, as before, we write $e^{-\phi}dV_{X}$ simply as $dV$. We split the integral over the regions $G\leq t+h$ and $G>t+h$, assuming here that $h$ is negative. In the case that $h$ is positive, we split instead the integral over $G\leq t$ and $G>t$. The proof is similar in either case from this point on, and we stick here to the case that $h$ is negative (in the case that $h$ is positive, we replace everywhere in what follows $G>t+h$ with $G>t$). Hence, on $G\leq t+h$, we have $G-t\leq h<0$. That means $f(t+h)-f(t)-f'(t)h=0$ here. Hence, we have 
\begin{align*}
\norm{(f(t+h)-f(t)-f'(t)h)u}^2_{\phi}&=\int_{G>t+h}\abs{f(t+h)-f(t)-f'(t)h}^2\abs{u}^2_{\omega}dV.
\end{align*}We now divide by $\abs{h}$ and get 

\begin{align*}
\fr{\norm{(f(t+h)-f(t)-f'(t)h)u}^2_{\phi}}{\abs{h}}&=\int_{G>t+h}\abs{\fr{f(t+h)-f(t)-f'(t)h}{h}}^2\abs{u}^2_{\omega}dV\notag\\
&=\int_{G>t+h}\abs{\fr{f(t+h)-f(t)}{h}-f'(t)}^2\abs{u}^2_{\omega}dV.
\end{align*}
The idea now is to let $h\to 0$, and we will then need to show that $\fr{f(t+h)-f(t)}{h}\to f'(t)$ uniformly for all $x\in \{G>t+h\}$. That is, writing $G=G(x)$,  $x\in X$,
\begin{align}
\fr{e^{-\chi(G(x)-t-h)}-e^{-\chi(G(x)-t)}}{h}\to \chi'(G(x)-t)e^{-\chi(G(x)-t)}=(-1)(e^{-\chi(G(x)-t)})'
\end{align}uniformly in $x$. Let $g:=g(r):=e^{-\chi(r)}$. It suffices to show that given any $\epsilon>0$ we can find a $\delta=\delta(\epsilon)$ such that whenever $|s-r|<\delta$, 
\begin{align*}
\abs{\fr{g(r)-g(s)}{r-s}-g'(s)}<\epsilon.
\end{align*}Indeed, then we can take $r:=G(x)-t-h$ and $s:=G(x)-t$, so that $r-s=-h$. Now, by the mean-value theorem, there is $p\in (s,r)$ such that
\begin{align}
\fr{g(r)-g(s)}{r-s}&=g'(p).
\end{align}
We can write $p=\eta s+(1-\eta)r$ for some $\eta\in (0,1)$. We can write $r=(r-s)+s$. Substituting this we find
\begin{align*}
p&=\eta s+(1-\eta)((r-s)+s)=\eta s+(1-\eta)s+(1-\eta)(r-s)=s+(1-\eta)(r-s).
\end{align*}We now get 
\begin{align*}
\abs{\fr{g(r)-g(s)}{r-s}-g'(s)}&=\abs{g'(p)-g'(s)}.
\end{align*}Now comes the crux of the argument. Since $g'$ is smooth, in particular continuous, it is \ita{uniformly continuous} on compact subsets. Recall that we are considering $G=G(x)$ for $x$ in the set $\{G>t+h\}$. That is $-t>G-t>h$, so there is a compact, say $K\Subset \set{R}$, such that we are considering $g'$ on (a subset of) $K$. Hence give $\epsilon>0$ we can find $\delta=\delta(\epsilon)$ such that 
\begin{align*}
\abs{p-s}<\delta\implies \abs{g'(p)-g'(s)}<\epsilon.
\end{align*}Now, $\abs{p-s}=(1-\eta)\abs{r-s}$, and so if we let $\delta':=\fr{\delta}{1-\eta}$, it follows that for $\abs{r-s}<\delta'$, we have that $\abs{\fr{g(r)-g(s)}{r-s}-g'(s)}<\epsilon$. It follows that given $\epsilon>0$, we can choose sufficiently small such that 
\begin{align*}
\fr{\norm{(f(t+h)-f(t)-f'(t)h)u}^2_{\phi}}{\abs{h}}&\leq \epsilon \int_{G>t+h}\abs{u}^2_{\omega}dV\leq \epsilon \norm{u}^2_{\phi},
\end{align*}This suffices to show that $f'(t)$ is the derivative of $f=f(t)$ as a map $U\ni t\mapsto f(t)\in \text{Aut}(H)$, and proves the case that $n=1$. For general case of $n$, the exact same argument applies, except that we replace $f$ above with $f^{(n-1)}$; as before we may focus on integrating over $\{G>t+h\}$ for $h$ small. The same ''trick'' with the mean-value theorem and uniform continuity applies for the conclusion.
\end{proof}

\subsection{Plurisubharmonicity properties of minimal solutions associated with quasi-complete Kähler manifolds and (semi)positive line bundles: Proof of Theorem 1.9}We have now come quite far, and it remains only to prove Theorems \ref{thm: special case of psh for quotient sections, psh of minimal solutions in the setting of quasi-complete kahler manifolds, smooth case} and Theorem \ref{thm: psh properties of minimal solutions on qck manifolds with possiblt singular phi, special metrics}. We shall begin, as is natural, with the proof of Theorem \ref{thm: special case of psh for quotient sections, psh of minimal solutions in the setting of quasi-complete kahler manifolds, smooth case}. This is analogous to the proof of Theorem \ref{thm: Nakano semipositivity of quasi-complete Kahler, smooth case}, except that we use Theorem \ref{thm: second main result on psh properties of holomorphic sections of Q under certain assumptions} instead of Theorem \ref{thm: first main theorem, Nakano semipositivity of subbundle under Hormander condition}.

\begin{proof}[Proof of Theorem \ref{thm: special case of psh for quotient sections, psh of minimal solutions in the setting of quasi-complete kahler manifolds, smooth case}]
We use notation similar to before. Let $[u]$ be a holomorphic (local) section of $Q:=H/H_0$ such that $u$ is a (smooth) representative of $[u]$. Then $u_{m}=m_{l}[u]$, where $m_{l}$ denotes the minimal lifting operator with respect to $h$. Thus, we want to prove that $-\tilde{f}(\norm{[u]}^2_{h^{Q}})$ is plurisubharmonic, $h^{Q}$ being the hermitian metric on $Q$ induced by $h$, and may therefore use Theorem \ref{thm: second main result on psh properties of holomorphic sections of Q under certain assumptions}.  Let, as in the proof of Theorem \ref{thm: Nakano semipositivity of quasi-complete Kahler, smooth case}, $A_{x}:=\partial_z\dbar_z\theta(\cdot,x)$, and $ M_{x}:=[\Theta_{j\bar{k}}](\cdot,x)$, for each $x\in X$. Let $\xi:=(\xi_1,\ldots, \xi_m)\in \set{C}^{m}$, and $\epsilon>0$, be given. We put $V:=\dbar_z\left(\s{j=1}{m}\delta_{j}\xi_{j}u_{m}\right)$, write $u=u_m+u_0$, where $u_0\in H$, and find, since $\delta_{j}=\partial_{t_{j}}-\partial_{t_{j}}\theta\wedge$, that
\begin{align}
V&=\s{j=1}{m}\xi_{j}\dbar_z\left(-\partial_{t_{j}}\theta\wedge u_m\right)=\s{j=1}{m}\xi_{j}\left(-\dbar_z\partial_{t_j}\theta\wedge u_{m}-\partial_{t_{j}}\wedge v\right).
\end{align}By our assumption, $\partial_{t}\theta=0$ on the support of $v$, so it follows that 
\begin{align}
V&=-\s{j=1}{m}\dbar_{z}\partial_{t_{j}}\theta\wedge (\xi_{j}u_{m}).
\end{align}Let $\vect{b}_{s}$ be the vector (field) identified with $b_{s}:=\dbar_z\partial_{t_s}\theta$, and let $\vect{F}:=\s{j=1}{m}\abs{\xi_{j}u_{m}}_{\omega}\vect{b}_{j}$. Let, for each $x\in X$, $\vect{F}_{x}:=\vect{F}(\cdot,x)$, and $B_{x}$ have $j$th column equal to $\vect{b}_{j}(\cdot,x)$. Considering $u_{j}:=\xi_{j}u_{m}$, the proof of Theorem \ref{thm: Nakano semipositivity of quasi-complete Kahler, smooth case} lets us apply Theorem \ref{thm: Hörmander's theorem on quasi-complete Kahler manifolds, smooth case, semipositive case} to find a $U\in H$ such that $\{U, \{\vect{F}_{x}\}, \{A_x\}\}$ satisfies the Hörmander $h$-estimate for the $\dbar_z(\cdot)=V$ equation. Also the plurisubharmonicity of $\theta$ gives that at each $x\in X$, 
\begin{align*}
(A_x+\epsilon I_{n})-B_x(M_x+\epsilon I_{m})^{-1}B_x^*\geq 0.
\end{align*}Hence the assertion follows by Theorem \ref{thm: second main result on psh properties of holomorphic sections of Q under certain assumptions}.
\end{proof}

\subsection{Plurisubharmonicity properties of minimal solutions associated with quasi-complete Kähler manifolds and pseudoeffective line bundles: Proof of Theorem 1.10}
Finally, we prove Theorem \ref{thm: psh properties of minimal solutions on qck manifolds with possiblt singular phi, special metrics}. Thanks to our discussion in section \ref{sec: special weights}, the proof of Theorem \ref{thm: special case of psh for quotient sections, psh of minimal solutions in the setting of quasi-complete kahler manifolds, smooth case}., and Theorem \ref{thm: singular Hormander}, the proof is going to be very short.

\begin{proof}[Proof of Theorem \ref{thm: psh properties of minimal solutions on qck manifolds with possiblt singular phi, special metrics}]
By our discussion in section \ref{sec: special weights}, $y\mapsto e^{-\chi(G-t)}$ is a valid smooth weight for $h^0$ of $h$.
The remaining proof is analogous to that of Theorem \ref{thm: special case of psh for quotient sections, psh of minimal solutions in the setting of quasi-complete kahler manifolds, smooth case}: We use Theorem \ref{thm: second main result on psh properties of holomorphic sections of Q under certain assumptions}. The condition that $(A_{x}+\epsilon I_{n})-B_x(M_{x}+\epsilon I_{m})^{-1}B_x^*\geq 0$ follows as before from $\theta$ being plurisubharmonicity, and the Hörmander $h$-estimate part follows as in the proof of Theorem \ref{thm: special case of psh for quotient sections, psh of minimal solutions in the setting of quasi-complete kahler manifolds, smooth case}, except that we use Theorem \ref{thm: singular Hormander} instead of Theorem \ref{thm: Hörmander's theorem on quasi-complete Kahler manifolds, smooth case, semipositive case}.
\end{proof}

\section{End and preview}
We end by commenting briefly on the theme of the sequel to this paper (\cite{TaiXuHilb}). In that paper we turn to applications in $L^2$-holomorphic extension of Ohsawa-Takegosghi type and  (strong) openness. The paper will be joint work with Xu Wang, and among our main results in there are a sharp Ohsawa-Takegoshi extension theorem on quasi-complete Kähler manifolds, and a new, complex Brunn-Minkowski theoretical, proof of the strong openness conjecture.

%
%
%
%
%
%
%
%
%
%
%
%
%
%
%
%
%
%
%
%
%
%

\chapter{Paper 3}

\begin{center}
\normalsize{\textbf{A HILBERT BUNDLES APPROACH TO COMPLEX BRUNN-MINKOWSKI THEORY, III}}\\[0.5cm]
\small{TAI TERJE HUU NGUYEN AND XU WANG}

\begin{abstract} This is a sequel to \cite{Tai2} and an expansion of \cite{wang-soc}. We consider applications of the material from the previous two papers, \cite{Tai1}, \cite{Tai2}, to effective sharp estimates of $L^2$-holomorphic extensions, and to strong openness. In particular, we generalize the Berndtsson-Lempert method (see \cite{BoLempOT}) to the class of quasi-complete Kahler manifolds (see \cite{Tai2}, Definition 4.1), and discuss related applications in the Ohsawa-Takegoshi theory and singularity theory of plurisubharmonic functions. The main result in this paper is a general monotonicity theorem. We use it to obtain a Guan-Zhu type of sharp Ohsawa-Takegoshi extension theorem, and give a complex Brunn-Minkowski proof of the strong openness conjecture of Demailly, as well as a generalization of Guan's sharp strong openness theorem. 
\end{abstract}

\end{center}




\section{Introduction}

In \cite{BoLempOT}, Berndtsson and Lempert gave a ''complex Brunn-Minkowski proof'' of the sharp Ohsawa-Takegoshi extension theorem \cite{Bl0, GZ} for pseudoconvex domains in $\mathbb C^n$, and introduced the so-called \emph{Berndtsson-Lempert method} for $L^2$-holomorphic extension. The reader may refer to \cite{OT} for the original Ohsawa-Takegoshi theorem, and \cite{Chen} for a simple proof. The main observation of Berndtsson and Lempert is that the sharp effective $L^2$ estimate, in the Ohsawa-Takegoshi extension theorem, can be seen to follow from a certain decreasingness property of minimal extensions. Moreover, this monotonicity property is a direct consequence (see \cite[Corollary 4.3]{Bern20}) of the convexity theorem of Berndtsson (see \cite[Lemma 4.2]{Bern20} and \cite{Bo0, Bo1}), which is a fundamental result in complex Brunn-Minkowski theory. In this paper, which is a sequel to \cite{Tai1}, we shall generalize the Berndtsson-Lempert method from the local setting of domains in $\mathbb{C}^{n}$ to the global setting of complex manifolds. We also discuss applications to (the) strong openness (conjecture of Demailly) and singularity theory of plurisubharmonic functions. More specifically, we shall prove a general monotonicity theorem, and use this to obtain the results related to $L^2$-holomorphic extension and strong openness. The general monotonicity theorem is thus our main result, and the other results are applications of it (and its proof). Related results to, and applications of, the main theorem (and its proof), include the following: a H\"ormander proof (see Theorem \ref{th:DK}) of the Demailly-Koll\'ar semi-continuity theorem \cite{DK}; a Donnelly-Fefferman proof (see Theorem \ref{th:so-DF}) of Guan-Zhou's strong openness theorem \cite{GZ0}; a  decreasing theorem approach (see Theorem \ref{th:so} and Theorem \ref{th:so1}) to Guan's sharp strong openness theorem \cite{G}; a decreasing theorem approach (see Theorem \ref{th:BLOT1}) to Zhou-Zhu's sharp Ohsawa-Takegoshi extension theorem  \cite{ZZ}.

We next state our main result. To set the stage, let $(X,\omega)$ be an $n$-dimensional quasi-complete Kahler manifold, and let $(L, e^{-\phi})$ be a pseudoeffective holomorphic line bundle over $X$. Let $H$ be the trivial Hilbert bundle over $\mathbb{C}$ consisting of $L$-valued $(n,q)$-forms $u$ such that 
\begin{align}
||u||^2_{\phi}&:=\int_{X}|u|^2_{\omega}e^{-\phi}<\infty,
\end{align}
equipped with the inner product that induces the norm $||\cdot||_{\phi}$. We denote this inner product by $h^{0}$, and view it as a hermitian metric of zero variation on $H$. We let $H_0$ be the subbundle over $H$ consisting of elements in $H$ which vanishes under the action of the $\dbar$-operator; the $\dbar$-operator here being the one on $X$ (the one on $\mathbb{C}$ will have no play in our discussion). We write a generic point in $\mathbb{C}$ as $y$, and the real part of $y$ as $t$. Suppose that $G\leq 0$ is a function on $X$ such that 
\begin{align}
\phi+\lambda G
\end{align}is plurisubharmonic for some $\lambda>0$, and consider the smooth weight $w$ for $h^0$ given by $y\mapsto e^{-\lambda \max\{G-t\}}$; this is a limiting case of the weights in Theorems 1.6 and 1.10 in \cite{Tai2}. We denote the smooth weighted hermitian metric induced by $h^0$ with weight $w$ by $h$. Thus, explicitly, if $u,v\in H$, we have, for all $y\in \mathbb{C}$,
\begin{align}
h_{y}(u,v)&=\int_{X}(u,v)_{\omega}e^{-\phi-\lambda \max\{G-t\}}.
\end{align}
Let $||\cdot||_{t}$ denote the norm induced by $h_{y}$; it is clear that it only depends on $t$.
Let $0\leq \alpha<\lambda$, and let $S_{\alpha}$ consists of the elements of $H_0$ which are locally $L^2$-integrable against $e^{-\phi-\alpha G}$. Our main result is the following theorem:

\begin{theorem}[Main theorem]\label{thm: main theorem, monotonicity property of minimal solutions, or holomorpphic quotient sections}
With the above notation and set-up, let $u\in H_0$. Then,
\begin{align}
\mathbb{R}\ni t\mapsto e^{-\alpha t}\inf\left\{||v||^2_{t}:(v-u)\in S_{\alpha}\right\}
\end{align}is decreasing.
\end{theorem}
For the proof of Theorem \ref{thm: main theorem, monotonicity property of minimal solutions, or holomorpphic quotient sections}, we shall use Theorem 4.15 in \cite{Tai1} and Theorem 1.6 in \cite{Tai2}, together with the following Hormander type of theorem:

\begin{theorem}\label{thm: Hormander theorem for decreasingness}
With the above notation and set-up, replace $\lambda \max\{G-t\}$ with $\chi(G-t)$, where $\chi$ is a smooth convex function on $\mathbb{R}$ which vanishes identically on $(-\infty, 0]$, and which satisfies that $0\leq \chi'\leq \lambda$. Fix $t<0$ and $\alpha$ (as above), and let $F\in H_0$. Then there exists $\tilde{F}\in H_0$ such that 
\begin{align}
||F-\tilde{F}||^2_{\phi+\alpha G}\leq \frac{9\lambda (1-t)^{\lambda-\alpha+1}e^{-\alpha t}}{\lambda-\alpha}||F||^2_{t}+2e^{-\alpha t}||F||^2_{\phi}\label{eq:h11},
\end{align}and
\begin{align}
||\tilde{F}||^2_{\phi}&\leq \frac{ 4\lambda (1-t)^{\lambda-\alpha+1}e^{-\alpha t}   }{\lambda-\alpha}||F||^2_{t},\label{eq:h12}
\end{align}where $||\cdot||_{\phi+\alpha G}$ is defined like $||\cdot||_{\phi}$, except that we replace $\phi$ with $\phi+\alpha G$.
\end{theorem}
The proof of Theorem \ref{thm: Hormander theorem for decreasingness} is based on the same idea for that of Theorem 1.7 in \cite{Tai2}. For the benefit of the reader, we give the proof in the next section. This concludes our short introduction. In the remainder of the paper we give the proof of the main theorem, and discuss related results and applications as mentioned above.

\section{Proof of Theorems 1.1 and 1.2}
In this section we prove our main theorem, that is, Theorem \ref{thm: main theorem, monotonicity property of minimal solutions, or holomorpphic quotient sections}, and Theorem \ref{thm: Hormander theorem for decreasingness}. Since the proof of the former uses the latter, we begin with the latter.

\subsection{A Hormander type of theorem: Proof of Theorem 1.2}
The proof follows much of that of Theorem 1.7 in \cite{Tai2}, and we shall use the same notation, referring also to \cite{Tai2}. The main difference is that instead of solving the $\dbar$-Laplace equation with the weight $\theta^{j}_{\nu}+\phi_0$, we solve it with a weight of the form

\begin{align}
\Theta^{j}_{\nu}&:=\phi^{j}_{\nu}+\phi_0+\alpha G_{\nu}^{j}-c\ln(1-f(G_{\nu}^{j})),
\end{align}for some constant $c>0$ and some sufficiently smooth real-valued function $f$; both to be further specified. As in the proof of Theorem 1.7 in \cite{Tai2}, we have 
\begin{align}
i\partial\dbar(\phi^{j}_{\nu}+\phi_0), i\partial\dbar(\psi_{j}^{\nu}+\phi_{0})&\geq -\epsilon_{j}\omega.
\end{align}A direct computation gives, where $F$ is also a sufficiently smooth real-valued function,
\begin{align}
\partial\dbar F(1-f(G^{j}_{\nu}))&=(F''(1-f(G^{j}_{\nu}))|f'(G^{j}_{\nu})|^2-F'(1-f(G^{j}_{\nu}))f''(G^{j}_{\nu}))\partial G^{j}_{\nu}\wedge \dbar G^{j}_{\nu}\notag\\
&-F'(1-f(G^{j}_{\nu}))f'(G^{j}_{\nu})\partial\dbar G^{j}_{\nu}.
\end{align}Hence, taking $F:=-c\ln(\cdot)$, we get
\begin{align}
i\partial \dbar (-\ln(1-f(G^{j}_{\nu})))&=\left(\frac{c|f'(G^{j}_{\nu})|^2}{(1-f(G^{j}_{\nu}))^2}+\frac{cf''(G
{j}_{\nu})}{1-f(G^{j}_{\nu})}\right)i\partial G^{j}_{\nu}\wedge \dbar G^{j}_{\nu}+\frac{cf'(G^{j}_{\nu})}{1-f(G^{j}_{\nu})}i\partial\dbar G^{j}_{\nu}.
\end{align}Assume that $f(G^{j}_{\nu})\leq 0$, and that $f'(G_{\nu}^{j})\geq 0$. If $i\partial\dbar G^{j}_{\nu}\geq 0$, then 
\begin{align}
i\partial\dbar \Phi^{j}_{\nu}&\geq -\epsilon_{j}\omega+\frac{cf''(G^{j}_{\nu})}{1-f(G^{j}_{\nu})}i\partial G^{j}_{\nu}\wedge \dbar G^{j}_{\nu}.
\end{align}On the other hand, if $i\partial\dbar G^{j}_{\nu}\leq 0$, then 
\begin{align}
i\partial\dbar \Phi^{j}_{\nu}&\geq -\epsilon_{j}\omega+\frac{cf''(G^{j}_{\nu})}{1-f(G^{j}_{\nu})}i\partial G^{j}_{\nu}\wedge \dbar G^{j}_{\nu}+\left(\alpha+cf'(G^{j}_{\nu}\right)i\partial \dbar G^{j}_{\nu}.
\end{align}Assume that $0\leq f'(G^{j}_{\nu})\leq 1$, and consider $c:=\lambda-\alpha$. Recall that $G^{j}_{\nu}=\frac{(\psi^{j}_{\nu}+\phi_0)-(\phi^{j}_{\nu}+\phi_0)}{\lambda}$. Hence, in either case, whether $i\partial \dbar G^{j}_{\nu}\geq 0$, or not, we get, choosing $c:=\lambda-\alpha$ as above,
\begin{align}
i\partial\dbar \Phi^{j}_{\nu}&=-\epsilon_{j}\omega +\frac{(\lambda-\alpha)f''(G^{j}_{\nu})}{1-f(G^{j}_{\nu})}i\partial G^{j}_{\nu}\wedge \dbar G^{j}_{\nu}.
\end{align}
Consider $c^{j}_{\nu}:=\dbar\left( (1-f'(G^{j}_{\nu}))F\right)=-f''(G^{j}_{\nu}) \dbar G^{j}_{\nu}\wedge F$, and let, using the same notation as in the proof of Theorem 1.7 in \cite{Tai2}, $v_{j,\nu,k}$  satisfy
\begin{align}
\left(\square_{j,\nu,k}+(q+1)\epsilon_{j}\right)v_{j,\nu,k}&=c_{j,\nu}
\end{align}and
\begin{align}
||\dbar v_{j,\nu,k}||^2_{j,\nu,k}+||\dbar^* v_{j,\nu,k}||^2_{j,\nu,k}+\epsilon_{j}(q+1)||v_{j,\nu,k}||^2_{j,\nu,k}&\leq \left|\left| \sqrt{\frac{(1-f(G^{j}_{\nu})|f''(G^{j}_{\nu})|}{\lambda-\alpha}}F \right|\right|^2_{j,\nu,k}.
\end{align}
The expression in the square root in the right-hand side above is found in a similar way as how (4.43) in the proof of Theorem 1.7 in \cite{Tai2} is found, using Corollary 2.4 in \cite{Tai2}. The idea now is to choose $f$ such that the $t$-norm can be introduced from the $(j,\nu,k)$-norm above. Note that in the latter there is a term $e^{-\alpha G^{j}_{\nu}}$, and this can be replaced bounded by $e^{-\alpha t}$ on $\{G^{j}_{\nu}>t\}$. Note also that since $f'\geq 0$, $f$ is increasing, so on $\{G^{j}_{\nu}>t\}$, we have $1-f(G^{j}_{\nu})\leq 1-t$. To introduce the $t$-norm we can therefore consider $f$ such that $f''(s)=\chi'(s-t)e^{-\chi(s-t)}=-\frac{d}{ds}\left(e^{-\chi(s-t)}\right)$. Then the  factor $f''(G^{j}_{\nu})$ in the above $(j,\nu,k)$-norm will introduce the factor $e^{-\chi(G^{j}_{\nu}-t)}$ appearing in the (approximated) $t$-norm, and also reduce the domain of integration to over where $\{G^j_{\nu}>t\}$. Integrating 
\begin{align}
f''(s)&=\frac{d}{ds}(-e^{-\chi(s-t)}),
\end{align}we get
\begin{align}
f'(s)&=C-e^{-\chi(s-t)}
\end{align}for some constant $C$. We want $0\leq f'\leq 1$, so we consider taking $C:=1$. Integrating once more, we finally get
\begin{align}
f(s)&=\int_{0}^{s}1-e^{-\chi(x-t)}\;dx.
\end{align}There is another factor in the $(j,\nu,k)$-norm which also does not appear in the $t$-norm, namely,
\begin{align}
e^{(\lambda-\alpha)\ln\left(1-f(G^{j}_{\nu})\right)}=\left(1-f(G^{j}_{\nu})\right)^{\lambda-\alpha},
\end{align}but this can, similar to the term $1-f(G^{j}_{\nu})$ above, be bounded above by $(1-t)^{\lambda-\alpha}$ on $\{G^{j}_{\nu}>t\}$. Choosing $f$ as above, we therefore get (recall also that $0\leq \chi'\leq \lambda$)
\begin{align}
||\dbar v_{j,\nu,k}||^2_{j,\nu,k}+||\dbar^* v_{j,\nu,k}||^2_{j,\nu,k}+\epsilon_{j}(q+1)||v_{j,\nu,k}||^2_{j,\nu,k}&\leq \frac{\lambda (1-t)^{\lambda-\alpha+1}}{\lambda-\alpha}e^{-\alpha t}||F||^2_{t}.
\end{align}As in the proof of Theorem 1.7 in \cite{Tai2}, we now let $k\to \infty$ first, and then $\nu,j\to \infty$, considering as usual weak limits. Analogous to in there, we then get a solution $u$ satisfying
\begin{align}
\dbar u&=\dbar ((1-f(G))F)
\end{align}on $X$ and
\begin{align}
||u||_{\phi+\alpha G}^2&\leq ||u||^2_{\phi+\alpha G-(\lambda-\alpha)\ln(1-f(G)))}\leq \frac{\lambda (1-t)^{\lambda-\alpha+1}}{\lambda-\alpha}e^{-\alpha t}||F||^2_{t},
\end{align}the middle norm being defined similarly as the norm furthest to the left,except that we replace $\phi$ with $\phi-(\lambda-\alpha)\ln(1-f(G))$. Finally, let $\tilde{F}:=(1-f'(G))F-u$. Then $\tilde{F}$ fits our needs. Indeed, the estimate on $u$ immediately gives

\begin{align}
||\tilde F - (1-\sigma'(G)) F||^2_{\phi+\alpha G }\leq  \frac{\lambda(1-t)^{\lambda -\alpha+1} e^{-\alpha t}}{\lambda-\alpha}   ||F||_t^2. \label{eq:x30}
\end{align}
Since $G\leq 0$, the above inequality gives
\begin{align*}
||\tilde F||^2_\phi & \leq 2 ||\tilde F - (1-\sigma'(G)) F||^2_{\phi}+2||(1-\sigma'(G)) F  ||^2_\phi \\
& \leq \frac{2\lambda(1-t)^{\lambda -\alpha+1} e^{-\alpha t}}{\lambda-\alpha}   ||F||_t^2 + 2 ||(1-\sigma'(G)) F  ||^2_\phi \\
& = \frac{2\lambda(1-t)^{\lambda -\alpha+1} e^{-\alpha t}}{\lambda-\alpha}  ||F||_t^2 + 2  \int_X  i^{n^2} F \wedge \bar F  e^{-\phi -2 \chi(G-t) }   \\ 
& \leq  \left(\frac{2\lambda(1-t)^{\lambda -\alpha+1} e^{-\alpha t}}{\lambda-\alpha}  +2\right)  ||F||_t^2 \\
& \leq  \frac{4\lambda(1-t)^{\lambda -\alpha+1} e^{-\alpha t}}{\lambda-\alpha}  ||F||_t^2, 
\end{align*}
which gives \eqref{eq:h12}. To prove \eqref{eq:h11}, note that \eqref{eq:x30} gives
$$
||\tilde F - F||^2_{\phi+\alpha G, G<t} \leq  \frac{\lambda(1-t)^{\lambda -\alpha+1} e^{-\alpha t}}{\lambda-\alpha}   ||F||_t^2.
$$
Hence together with
$$
||\tilde F - F||^2_{\phi+\alpha G, G \geq t} \leq e^{-\alpha t}    ||\tilde F - F||^2_\phi \leq 2 e^{-\alpha t} \left(||\tilde F||^2_\phi +||F||^2_\phi \right),
$$
we see that \eqref{eq:h12} gives \eqref{eq:h11}.

\subsection{A general decreasingness theorem: Proof of Theorem 1.1}
Next, we give the proof of our main theorem, Theorem \ref{thm: main theorem, monotonicity property of minimal solutions, or holomorpphic quotient sections}. It follows as a special case of the following general monotonicity theorem by letting $\chi(x)$ converge to $\lambda \max\{s,0\}$:

\begin{theorem}\label{thm: general monotonicity theorem, generalization of main theorem}
With the same assumptions and notation in Theorem \ref{thm: general monotonicity theorem, generalization of main theorem}, let instead $h_{y}$ be defined by 
\begin{align}
h_{y}(u,v)&:=\int_{X}(u,v)_{\omega}e^{-\phi-\chi(G-t)},
\end{align}writing $t$ for the real part of $y$ (as usual), where $\chi$ is a smooth convex function on $\mathbb{R}$ vanishing identically on $(-\infty,0]$, and satisfying that $0\leq \chi'\leq \lambda$. Then the same conclusion as in Theorem \ref{thm: general monotonicity theorem, generalization of main theorem} holds. That is, given $u\in H_0$, 
\begin{align}
\mathbb{R}\ni t\mapsto e^{-\alpha t}\inf \{||v||^2_{t}:(u-v)\in S_{\alpha}\}
\end{align}is decreasing, where now $||\cdot||_{t}$ is the norm induced by $h_y$ above.
\end{theorem}

\begin{proof}
Let $\Theta_0$ denote the Chern curvature of the Chern connection of $H_0$ with respect to $h$. By Theorem 1.6 in \cite{Tai2}, it follows that $i\Theta_0\geq_{N}0$ since $\theta$ is plurisubharmonic under our assumptions. Let $\Theta_0^{\vee}$ denote the Chern curvature of the dual bundle $(H_0)^*:=H_0^*$ with respect to the dual hermitian metric $h^{\vee}$ of $h$. Then $i\Theta^{\vee}\leq 0$, so it follows by Theorem 4.15 in \cite{Tai1}, that 
\begin{align}
t\mapsto \ln(||f||^2_{t}),
\end{align}is convex for all holomorphic (local) sections of $H_0^*$, where the $t$-norm here is the norm induced by $h^{\vee}_{y}$. Consider $Q_{\alpha}:=H_0/S_{\alpha}$, and let $h^{Q_{\alpha}}$ denote the hermitian metric on $Q_{\alpha}$ induced by $h$. By the quotient bundle curvature formula from \cite{Tai1}, it follows that for all holomorphic (local) sections $\hat{f}$ of $(Q_{\alpha})^*$, 
\begin{align}
t\mapsto \ln(||\hat{f}||^2_{t})
\end{align}is convex, where the $t$-norm here denotes the norm induced by the dual hermitian metric of $h^{Q_{\alpha}}$. The point  of all this is of course that $\inf\{||v||^2_{t}:(v-u)\in S_{\alpha}\}$ is precisely equal to the norm of $[u]$ with respect to $h_y^{Q_{\alpha}}$, where $[u]$ denotes any holomorphic (local) section of $Q_{\alpha}$ with $u$ as a (smooth) representative. What we are trying to prove amounts therefore to certain decreasingness properties of norms of holomorphic (local) sections of $Q_{\alpha}$. That is, turning to the dual, certain increasingness properties of the norms of continuous linear functionals on $Q_{\alpha}$. Now, with this in mind, fix $t<0$, let $F\in H_0$ be such that $||[u]||^2_{t}=||F||^2_{t}$, and let $\tilde{F}$ be given as in Theorem \ref{thm: Hormander theorem for decreasingness}. That is, $(\tilde{F}-F)\in S_{\alpha}$, and
\begin{align}
||\tilde{F}||^2_{0}=||\tilde{F}||^2_{\phi}&\leq \frac{4\lambda (1-\lambda)^{\lambda-\alpha+1}e^{-\alpha t}}{\lambda-\alpha}||F||^2_{t}=\frac{4\lambda (1-\lambda)^{\lambda-\alpha+1}e^{-\alpha t}}{\lambda-\alpha}||[u]||^2_{t}.
\end{align}That $(\tilde{F}-F)\in S_{\alpha}$ means that $[\tilde{F}]=[u]$. Hence, by definition of quotient norms, we have shown that
\begin{align}
||[u]||^2_{0}&\leq \frac{4\lambda (1-\lambda)^{\lambda-\alpha+1}e^{-\alpha t}}{\lambda-\alpha}||[u]||^2_{t}.
\end{align}Hence for all $\epsilon>0$, there exists $C=C(\epsilon,\lambda, \alpha)>0$, such that
\begin{align}
||[u]||&\leq C e^{-(\alpha+\epsilon)t}||[u]||^2_{t}.
\end{align}Turning to the dual we get that for all $\hat{f}\in (Q_{\alpha})^*$, where the norm notation is as above,
\begin{align}
||\hat{f}||^2_{t}&\leq \frac{1}{C}e^{-(\alpha+\epsilon )t}||\hat{f}||^2_{0}.\label{eq: functional estimate in proof of main theorem}
\end{align}The left-hand side is by our discussion above logarithmically convex in $t$. Then  

\begin{align}
\ln(||\hat{f}||^2_{t})+(\alpha+\epsilon) t
\end{align}
is convex in $t$, and
 \eqref{eq: functional estimate in proof of main theorem} shows that it is bounded above as $t\to -\infty$.
Hence it must be increasing as a function of $t$. Hence, taking the exponential, $e^{(\alpha+\epsilon)t}||\hat{f}||^2_{t}$ is increasing. Thus, turning back to $Q_{\alpha}$, $e^{-(\alpha+\epsilon)t}||[u]||^2_{t}$ is decreasing. Letting $\epsilon \to 0$, we finally get that recall that if $t\geq 9$, then $\theta_{y}$ is independent of $t$, so the same is true for all the $t$-norms)
\begin{align}
\mathbb{R}\ni t\mapsto e^{-\alpha t}||[u]||^2_{t}
\end{align}is decreasing, which completes the proof.

\end{proof}

Note that the previous proof also gives the following general convexity theorem:

\begin{theorem}\label{thm: general convexity theorem}
With the same notation and assumptions as in Theorem \ref{thm: general monotonicity theorem, generalization of main theorem}, let $f$ be any holomorphic (local) section of $H_0^*$ (or $(Q_{\alpha})^*$) . Then 
\begin{align}
t\mapsto \ln(||f||^2_{t})
\end{align}is convex.
\end{theorem}
Similar to the case of Theorem \ref{thm: main theorem, monotonicity property of minimal solutions, or holomorpphic quotient sections}, Theorem \ref{thm: general convexity theorem} also holds with $\theta^{y}=\phi+\lambda \max\{G-t,0\}$, by letting $\chi(s)$ tend to $\lambda \max\{s,0\}$.  

This finishes the proofs of the theorem from the introduction. In the remainder of the paper we turn to related results and applications.

\section{Singularity theory of quasi-plurisubharmonic functions}
In this section we discuss applications to singularity theory of quasi-plurisubharmonic functions. In particular, we shall give a (semi) new proof, using Theorem \ref{thm: Hormander theorem for decreasingness}, of a semi-continuity result of Demailly-Kollar (\cite{DK}), and a complex Brunn-Minkowski proof, using Theorem \ref{thm: main theorem, monotonicity property of minimal solutions, or holomorpphic quotient sections}, of strong openness.

\subsection{Multiplier ideal sheaf and Demailly-Koll\'ar theory}
We begin by recalling the notions of \emph{quasi-plurisubharmonic} metrics, and \emph{multiplier ideal sheaves}:
\begin{definition}\label{de:quasi-psh} Let $L$ be a holomorphic line bundle over a complex manifold $X$. A metric $e^{-\phi}$ on $L$ is said to be \textbf{quasi-plurisubharmonic} (quasi psh) if its local potential $\phi$ is quasi-psh, i.e. $\phi$ can be locally written as the sum of a psh function and a smooth function. If $(L, e^{-\phi})$ is a line bundle with quasi-psh metric, then the sheaf $\mathcal I(\phi)$ (of  germs of holomorphic functions $f$ such that $|f|^2e^{-\phi}$ is locally integrable) is called the \textbf{multiplier ideal sheaf} associated to $\phi$.
\end{definition}

The notion of multiplier ideal sheaf comes from Nadel's paper \cite{Nad}. Nadel proved that $\mathcal I(\phi)$ is a coherent ideal sheaf of $\mathcal O_X$ (see \cite[Proposition 5.7]{D} for the details), and obtained a nice criterion for the existence of K\"ahler-Einstein metrics on certain Fano manifolds using $\mathcal I(\phi)$. Later, Demailly-Koll\'ar \cite{DK} found a simpler proof of Nadel's criterion using the following semi-continuity result (see \cite[Theorem 0.2]{DK}); we use the abbreviation \emph{psh} for plurisubharmonic:

\begin{theorem}\label{th:DK} Let $\phi$ be a psh function on the unit ball $\mathbb B$ in $\mathbb C^n$, and let $G, G_j\leq 0$ such that
\begin{align}
\text{$\phi+\lambda G_j$ and $\phi+\lambda G$}
\end{align}are psh on $\mathbb B$ for some constant $\lambda>1$. Assume that (the Lebesgue measure is omitted)
\begin{equation}\label{eq:DK0}
\int_{\mathbb B} e^{-\phi-\beta G}<\infty, \ \ \lim_{j\to\infty} \int_{\mathbb B} |G_j-G|=0,
\end{equation}
for some $1<\beta<\lambda$. Then
\begin{equation}\label{eq:DK1}
\lim_{j\to \infty} \int_{|z|<r} |e^{-\phi-G_j}-e^{-\phi-G}|  =0, \  \ \ \forall \ 0<r<1.
\end{equation}
\end{theorem}

\begin{proof} For the convenience of the reader, we shall include a proof here. Our proof depends only on the H\"ormander estimate, but the idea is already implicitly included in \cite{DK}. We shall use
\begin{equation}\label{eq:DK2}
\int  |e^{-\phi-G_j}-e^{-\phi-G}|   \leq \int_{G_j \geq t} |e^{-\phi-G_j}-e^{-\phi-G}|+ \int_{G_j < t} e^{-\phi-G} + \int_{G_j < t} e^{-\phi-G_j}
\end{equation}
By our assumption, the psh functions $\phi+G_j$ converge to the psh function $\phi+G$ in $L^1(\mathbb B)$. Since for psh functions, $L^1$ convergence implies almost everywhere  convergence, we know that $G_j$ converges to $G$ almost everywhere. Thus the Lebesgue bounded convergence theorem implies 
\begin{equation}\label{eq:DK3}
\lim_{j\to \infty}  \left(\int_{G_j \geq t} |e^{-\phi-G_j}-e^{-\phi-G}|+ \int_{G_j < t} e^{-\phi-G} \right)= \int_{G < t} e^{-\phi-G}.
\end{equation}
for every fixed $t<0$. The main difficulty is to control the last term $\int_{G_j < t} e^{-\phi-G_j}$ in \eqref{eq:DK2}, and the idea is to use Theorem \ref{thm: Hormander theorem for decreasingness} with $F=1$ and $\alpha=1$. Letting $\chi(s)$ converge to $\lambda\max\{s,0\}$ we get holomorphic functions $F_{j, t}$ on $\mathbb B$ with  (see \eqref{eq:h12})
\begin{equation}\label{eq:DK4}
\int_{\mathbb B} |F_{j, t}|^2 e^{-\phi} \leq \frac{4\lambda(1-t)^{\lambda} e^{-t}}{\lambda-1}\, \int_{\mathbb B} e^{-\phi-\lambda\max\{G_j-t, 0\}},
\end{equation}
and (by \eqref{eq:x30})
\begin{equation}\label{eq:DK5}
\int_{G_j <t} |F_{j, t} -1|^2 e^{-\phi-G_j} \leq   \frac{\lambda(1-t)^{\lambda } e^{-t}}{\lambda-1} \int_{\mathbb B} e^{-\phi-\lambda\max\{G_j-t, 0\}}.
\end{equation}
By the Lebesgue bounded convergence theorem we have
$$
\lim_{j\to \infty} \int_{\mathbb B} e^{-\phi-\lambda\max\{G_j-t, 0\}}  =  \int_{\mathbb B} e^{-\phi-\lambda\max\{G-t, 0\}}:=A(t). 
$$
Hence for every fixed $t \leq -1$, we can choose $j(t)$ such that
$$
\int_{\mathbb B} e^{-\phi-\lambda\max\{G_j-t, 0\}}  \leq 2 A(t),  \ \ \forall \ j\geq j(t).
$$
Lemma \ref{lem: lemma in semicontinuity proof} below further gives
$$
e^{-\beta t}\int_{\mathbb B} e^{-\phi-\lambda\max\{G_j-t, 0\}}  \leq \frac{4\lambda}{\lambda-\beta}  \int_{\mathbb B} e^{-\phi-\beta G},  \ \ \forall \ j\geq j(t), \ t\leq -1.
$$
Hence \eqref{eq:DK4} implies that there exists $t_0<-1$ such that
$$
\sup_{|z|<r} |F_{j, t}(z)| \leq  \frac{1}{2} , \ \ \forall \ j\geq j(t), \ t\leq t_0.
$$
Thus by \eqref{eq:DK5} (when $t_0$ is smooth enough)
$$
\int_{G_j <t, |z|<r}  e^{-\phi-G_j}  \leq e^{(\beta-1)t/2} , \ \ \forall \ j\geq j(t), \ t\leq t_0.
$$
By \eqref{eq:DK3}, we can choose $j(t)$ such that 
\begin{equation}\label{eq:DK7}
 \int_{G_j \geq t, |z|<r} |e^{-\phi-G_j}-e^{-\phi-G}|+ \int_{G_j < t, |z|<r} e^{-\phi-G} \leq 2 \int_{G < t, |z|<r} e^{-\phi-G} ,  \ \ \forall \ j\geq j(t),
\end{equation}
so \eqref{eq:DK2} gives
$$
\int_{|z|<r}  |e^{-\phi-G_j}-e^{-\phi-G}|   \leq  e^{(\beta-1)t/2}+ 2 \int_{G < t, |z|<r} e^{-\phi-G} , \ \ \forall \ j\geq j(t), \ t\leq t_0.
$$
Note that the right hand side goes to zero when $t\to-\infty$. Hence, \eqref{eq:DK1} follows.
\end{proof}
In the previous proof we used the following lemma:
\begin{lemma}\label{lem: lemma in semicontinuity proof} For every $t\leq -1$, we have 
\begin{align}
A(t) e^{-\beta t} \leq \frac{2\lambda}{\lambda-\beta}  \int_{\mathbb B} e^{-\phi-\beta G}.
\end{align}
\end{lemma}

\begin{proof} Note that
\begin{equation}\label{eq:so3}
\int_{-\infty}^0 e^{-\lambda\max\{s-t, 0\}} e^{-\beta t}\, dt= \frac{\lambda}{(\lambda-\beta)\beta}  e^{-\beta s}-\frac1\beta, \ \ \forall  \ \lambda>\beta,
\end{equation}
gives
\begin{equation}\label{eq:so4}
\int_{-\infty}^0  A(t) e^{-\beta t}\, dt   \leq  \frac{\lambda}{(\lambda-\beta)\beta}  \int_{\mathbb B} e^{-\phi-\beta G} <\infty.
\end{equation}
Since $A(t)$ is increasing in $t$, we have
$$
\int_{-\infty}^0  A(s) e^{-\beta s}\, ds  \geq  \int_{t}^{0}  A(s) e^{-\beta s}\, ds \geq A(t) 
\int_{t}^0  e^{-\beta s}\, ds     \geq  \frac{A(t) e^{-\beta t}}{2\beta}
$$
for every $t\leq -1$. Hence the lemma follows from \eqref{eq:so4}.
\end{proof}
By Berndtsson's solution of the openness conjecture (see \cite{Bern15b}, and also the next section for further developments), we know that Theorem \ref{th:DK} is also true in case that $\beta=1$. 

\subsection{The strong openness theorem (the strong openness conjecture)}
In this section, we shall show how to use the Donelly-Fefferman estimate to prove the strong openness theorem of Guan-Zhou (\cite{GZ0}), which used to be known as the strong openness conjecture of Demailly \cite[Remark 15.2.2]{Dem00}. The reader may also see \cite{JM, Bern20, Bern15b, FJ, Hiep, Lempert} and references therein for related results. We start by giving the needed details on the Donelly-Fefferman estimates; for the proof it might be helpful to also look back at the proof of Theorem \ref{thm: Hormander theorem for decreasingness}:

\subsubsection{The Donnelly-Fefferman estimate}

\begin{theorem} \label{th:DF} Let $F$ be a holomorphic function on the unit ball $\mathbb B$ in $\mathbb C^n$, and let $G 
\leq  0$ be a smooth plurisubharmonic function on 
$\mathbb B$. Let  $\chi \geq 0$ be a smooth convex increasing function on $\mathbb R$ vanishing identically on $(-\infty, 0]$, and satisfying $0\leq \chi'\leq \lambda$ for some $\lambda>0$ (by approximation, we may consider in the limit $\chi(s)=\lambda\max\{s,0\}$). Then for every  $t<0$ and $\alpha \geq 0$, we can find another holomorphic function $\tilde F $ on $\mathbb B$ such that 
\begin{equation}\label{eq:h21}
\int_{G<t} |\tilde F-F|^2 e^{-\alpha G} \leq  -4\lambda t  e^{-\alpha t}   \int_{\mathbb B} |F|^2 e^{-\chi(G-t)} 
\end{equation}
and
\begin{equation}\label{eq:h22}
\int_{\mathbb B} |\tilde F|^2 \leq (2-8\lambda t  e^{-\alpha t}  ) \int_{\mathbb B} |F|^2 e^{-\chi(G-t)}.
\end{equation}
\end{theorem}

\begin{proof} Put
$
\sigma(s):=\int_{0}^s 1- e^{-\chi(x-t)} \, dx.
$
We know that
$$
\sigma (0)= 0,  \ \  \sigma'(s)= 1-e^{-\chi(s-t)} \geq 0.
$$
Hence $0\leq \sigma' \leq 1$ and $\sigma'(s)=0$ when $s\leq t$. Thus we have
\begin{equation}\label{eq:x00}
t \leq\sigma(t) \leq \sigma(G) \leq 0, \ \ \sigma''(G)= \chi'(G-t) e^{-\chi(G-t)}  \geq 0.
\end{equation}
Now we know that
$
\psi:= -\log(-\sigma(G))
$
satisfies
$$
i\partial\dbar \psi \geq i\partial\psi \wedge \dbar \psi.
$$
By the Donnelly-Fefferman estimate, we can solve $\dbar u = \dbar ((1-\sigma'(G)) F)$ with
$$
\int_{\mathbb B} |u|^2 e^{-\alpha G} \leq 4 \int_{\mathbb B} |\sigma''(G) \dbar G|^2_{i\partial\dbar \psi } |F|^2 e^{-\alpha G}.
$$
Notice that
$$
i\partial\dbar \psi \geq  \frac{i\partial\dbar  (\sigma(G))}{-\sigma(G)} \geq \frac{\sigma''(G) i\partial G \wedge  \dbar G}{ -\sigma(G)}
$$
gives
$$
  |\sigma''(G) \dbar G|^2_{i\partial\dbar \psi } \leq  \sigma''(G) |\sigma(G)| \leq -\lambda t \, e^{-\chi(G-t)} . 
$$
Hence we have
$$
\int_{\mathbb B} |u|^2 e^{-\alpha G} \leq  -4\lambda t  \int_{G\geq t} |F|^2 e^{-\alpha G} e^{-\chi(G-t)} \leq -4\lambda t  e^{-\alpha t} \int_{G\geq t} |F|^2  e^{-\chi(G-t)} .
$$
Take $\tilde F = (1-\sigma'(G)) F -u $, we get
$$
\int_{\mathbb B} |\tilde F - (1-\sigma'(G)) F|^2 e^{-\alpha G} \leq -4\lambda t  e^{-\alpha t} \int_{G\geq t} |F|^2  e^{-\chi(G-t)}  ,
$$
which gives \eqref{eq:h21}. Moreover, 
\begin{align*}
\int_\mathbb B |\tilde F|^2 & \leq 2 \int_\mathbb B  |\tilde F - (1-\sigma'(G)) F  |^2 +2\int_\mathbb B  |(1-\sigma'(G)) F  |^2  \\
& \leq -8\lambda t  e^{-\alpha t}  \int_{\mathbb B} |F|^2 e^{-\chi(G-t)}  + 2 \int_\mathbb B  |(1-\sigma'(G)) F  |^2 \\
& =  -8\lambda t  e^{-\alpha t}  \int_{\mathbb B} |F|^2 e^{-\chi(G-t)}  + 2  \int_\mathbb B  |F|^2 e^{-2\chi(G-t)}  \\ 
& \leq (2-8\lambda t  e^{-\alpha t}  )  \int_{\mathbb B} |F|^2 e^{-\chi(G-t)} , 
\end{align*}
gives \eqref{eq:h22}.
\end{proof}

The strong openness theorem of Guan-Zhou (\cite{GZ0}) is the following; as mentioned, we shall prove it using the above Donelly-Fefferman estimate:
\begin{theorem}\label{th:so-DF}
Let $F$ be a holomorphic function on the unit ball $\mathbb B$ in $\mathbb C^n$. Let $G \leq 0$ be a psh function on 
$\mathbb B$. If $\int_\mathbb B |F|^2  \,e^{-\beta G}<\infty$ for some $\beta>0$ then there exists $\alpha>\beta$ such that $|F|^2  \,e^{-\alpha G}$ is integrable near the origin of $\mathbb B$.
\end{theorem}

\begin{proof}  Put
\begin{equation}\label{eq:sopen1}
\alpha_0:=\sup \{\alpha\geq 0: F_0 \in \mathcal I(\alpha G)_0\}.
\end{equation}
It suffices to show that $\alpha_0>\beta$. Fix $\lambda>\beta$ and put
$$
||F||_t^2:= \int_{\mathbb B} |F|^2 e^{-\lambda\max\{G-t, 0\}}, \ \ t\leq 0.
$$
If $\alpha_0 <\infty$ then the Donnelly-Fefferman Lemma below gives
$$
C:=\inf_{t<-1} \{(-t) e^{-\alpha_0 t}||F||_t^2\} >0.  
$$
But since $\int_\mathbb B |F|^2  \,e^{-\beta G}<\infty$, we know \eqref{eq:so3} gives
$$
\infty > \int_{-\infty}^0 ||F||_t^2 e^{-\beta t}\, dt \geq \int_{-\infty}^{-1} \frac{C}{(-t) e^{-\alpha_0 t}}e^{-\beta t}\, dt = C \int_1^\infty \frac{e^{(\beta-\alpha_0) s}} s \, ds.
$$
Thus we must have $\alpha_0> \beta$ (otherwise we would get $\int_1^\infty \frac1 s <\infty$).
\end{proof}
In the previous proof we used the following  Donelly-Fefferman lemma:
\noindent
\textbf{Donnelly-Fefferman Lemma.} \emph{If $\alpha_0 <\infty$ then $\inf_{t<-1} \{(-t) e^{-\alpha_0 t}||F||_t^2\} >0$.}

\begin{proof}  Put
$$
I_0:= \inf \left\lbrace\int_{\mathbb B} |\tilde F|^2:  (\tilde F-F)_0 \in \bigcup_{\alpha>\alpha_0}  \mathcal I(\alpha G)_0 \right\rbrace.
$$
If $\alpha_0 <\infty$ then the strong Noetherian property of coherent ideal sheaves gives $I_0>0$ (this is a well-known result, see \cite[Lemma 2.1]{wang-soc}). On the other hand, the Donnelly-Fefferman estimate, Theorem \ref{th:DF}, gives
$$
I_0 \leq (2-8\lambda t  e^{-\alpha_0 t}  ) ||F||_t^2 , \ \ \forall \ t<0,
$$
from which the lemma follows.
\end{proof}

The optimal estimate 
\begin{equation}\label{eq:sopen2}
\frac{\alpha_0}{\alpha_0-\beta} \geq \frac{I_0}{\int_\mathbb B |F|^2  \,e^{-\beta G}}
\end{equation}
was first proved by Guan \cite{G} using \cite[Lemma 2.1]{GZ1}; an effective but non-optimal estimate is also given in \cite{GZ1}. Another proof of Guan's optimal estimate can be given using our main theorem, Theorem \ref{thm: main theorem, monotonicity property of minimal solutions, or holomorpphic quotient sections}. In fact, Theorem \ref{thm: general monotonicity theorem, generalization of main theorem} directly implies (see the next section for details)
$$
I_0 \leq  e^{-\alpha_0 t}||F||_t^2,
$$
which gives
$$
\frac{I_0}{\alpha_0-\beta}  \geq \int_{-\infty}^0 ||F||_t^2 e^{-\beta t}\, dt  =  \frac{\lambda}{(\lambda-\beta)\beta} \int_\mathbb B |F|^2  \,e^{-\beta G}- \frac1{\beta}\int_\mathbb B |F|^2.   
$$
Letting $\lambda\to \infty$, we get
$$
\frac{I_0}{\alpha_0-\beta}  \geq \frac1\beta \int_\mathbb B |F|^2  \,e^{-\beta G}-\frac1{\beta}\int_\mathbb B |F|^2    \geq \frac1\beta \int_\mathbb B |F|^2  \,e^{-\beta G}-\frac1{\beta} I_0,
$$
from which \eqref{eq:sopen2} follows.

\subsection{A generalization of Guan's sharp strong openness theorem}

In this section, we shall show how to use the main theorem, Theorem \ref{thm: main theorem, monotonicity property of minimal solutions, or holomorpphic quotient sections}, to prove the following generalization of Guan's sharp estimate \eqref{eq:sopen2} for the strong openness theorem; below $q=0$, for which we write $H_0$ as $H^0(X,K_X+L)$:

\begin{theorem}\label{th:so} Let $(L, e^{-\phi})$ be a pseudoeffective line bundle on an $n$-dimensional quasi-complete K\"ahler manifold $X$. Let $G\leq 0$ on $X$ such that $\phi+\lambda G$ is psh for some constant $\lambda>1$. Fix $F\in H^0(X, K_X+L)$. Assume that 
$$||F||^2_{\phi+\beta G}:=\int_X i^{n^2} F\wedge \bar F  \,e^{-\phi-\beta G}<\infty$$ for some constant $0< \beta<\lambda$. For every compact subset $K$ of $X$, put
\begin{equation}\label{eq:jump-K}
\alpha_K :=\sup\left\lbrace 0\leq \alpha \leq \lambda: F_x\in \mathcal I(\phi+\alpha G), \ \forall \, x\in K\right\rbrace. 
\end{equation} 
Assume that $\alpha_K <\lambda$. Put
\begin{equation}\label{eq:so1}
I_K:=\inf \left\lbrace ||\tilde F||_\phi^2: \tilde F \in H^0(X, K_X+L), \, \tilde F_x-F_x \in \bigcup_{\alpha>\alpha_K}  \mathcal I(\phi+\alpha G)_x, \ \forall \, x\in K    \right\rbrace,
\end{equation}
then
\begin{equation}\label{eq:so2}
\alpha_K- \beta \geq     \frac{I_K}{ \frac{\lambda}{(\lambda-\beta)\beta} ||F||^2_{\phi+\beta G} - \frac1{\beta} ||F||^2_\phi }>0.
\end{equation}
\end{theorem}

\begin{proof} Since $\alpha_K<\lambda$, the proof of  \cite[Lemma 2.1]{wang-soc} gives $I_K>0$. Thus it suffices to prove the estimate in \eqref{eq:so2}. For $\alpha_K<\alpha<\lambda$, put
$$
||[F]||_t^2:=
\inf \left\lbrace ||\tilde F||_t^2 : \tilde F-F \in H^0(X, \mathcal O(K_X+L)\otimes\mathcal I(\phi+\alpha G))\right\rbrace,
$$
where $||\tilde F||_t^2:= \int_{X} i^{n^2} \tilde F \wedge \overline{\tilde F}   e^{-\phi-\lambda\max\{G-t, 0\}}$,
we have
\begin{equation}\label{eq:31}
||[F]||_0^2 \geq I_K.
\end{equation}
By our main theorem, Theorem \ref{thm: main theorem, monotonicity property of minimal solutions, or holomorpphic quotient sections}, we know that $e^{-\alpha  t} ||[F]||_t^2$ is decreasing in $t$. Hence 
$$
\int_{-\infty}^0 ||[F]||_t^2 \, e^{-t\beta} \, dt  \geq \int_{-\infty}^0 ||[F]||_0^2 \, e^{t(\alpha-\beta)} dt
$$
and \eqref{eq:31} gives
$$
\int_{-\infty}^0 ||[F]||_t^2 \, e^{-t\beta} \, dt \geq \frac{I_K}{\alpha-\beta}.
$$
Integrating 
$$
||[F]||_t^2 \geq \int_{X} i^{n^2} F \wedge \bar F  e^{-\phi-\lambda\max\{G-t, 0\}}
$$
and applying \eqref{eq:so3} we have
$$
\frac{\lambda}{(\lambda-\beta)\beta} ||F||^2_{\phi+\beta G} - \frac1{\beta} ||F||^2_\phi  \geq \frac{I_K}{\alpha-\beta}, \  \ \forall \ \alpha_K<\alpha<\lambda.
$$
Letting $\alpha$ go to $\alpha_K$, the estimate in \eqref{eq:so2} follows. 
\end{proof}
In the case that $\beta=0$, \eqref{eq:so3} degenerates to    
\begin{equation}\label{eq:32}
\int_{-\infty}^0 e^{-\lambda\max\{s-t, 0\}} \, dt = \frac1{\lambda}-s, \ \ \forall \ s\leq 0
\end{equation}
and we have
\begin{equation}\label{eq:bo-int-0}
\int_{-\infty}^0 ||F||_t^2 \, dt =  \frac1{\lambda} ||F||^2_\phi+\int_{X} i^{n^2} (-G) \, F\wedge \bar F  \,e^{-\phi},
\end{equation}
thus a similar proof gives the following result:

\begin{theorem}\label{th:so1} Let $(L, e^{-\phi})$ be a pseudoeffective line bundle on an $n$-dimensional quasi-complete K\"ahler manifold $X$. Let $G\leq 0$ on $X$ such that $\phi+\lambda G$ is psh for some constant $\lambda>1$. Fix $F\in H^0(X, K_X+L)$. Assume that 
$$
||F||^2_{\phi}:=\int_X i^{n^2} F\wedge \bar F  \,e^{-\phi}<\infty, \ \  \int_{X} i^{n^2} (-G) \, F\wedge \bar F  \,e^{-\phi}  <\infty.
$$
For every compact subset $K$ of $X$, put
\begin{equation}\label{eq:jump-K1}
\alpha_K :=\sup\left\lbrace 0\leq \alpha \leq \lambda: F_x\in \mathcal I(\phi+\alpha G), \ \forall \, x\in K\right\rbrace. 
\end{equation} 
Assume that $\alpha_K <\lambda$. Put
\begin{equation}\label{eq:so11}
I_K:=\inf \left\lbrace ||\tilde F||_\phi^2: \tilde F \in H^0(X, K_X+L), \, \tilde F_x-F_x \in \bigcup_{\alpha>\alpha_K}  \mathcal I(\phi+\alpha G)_x, \ \forall \, x\in K    \right\rbrace,
\end{equation}
then
\begin{equation}\label{eq:so21}
\alpha_K  \geq  \frac{I_K}{ \frac1{\lambda} ||F||^2_\phi+\int_{X} i^{n^2} (-G) \, F\wedge \bar F  \,e^{-\phi} } >0.
\end{equation}
\end{theorem}

\section{The Berndtsson-Lempert approach for the Ohsawa-Takegoshi theory}

\subsection{Vanishing theorems and extension of holomorphic sections}
In this section we discuss further applications, now to the Ohsawa-Takegoshi theory and sharp effective estimates in the $L^2$-holomorphic extension theorem. The starting point of the story is the following Nadel vanishing theorem (see the proof of Theorem 5.11 and Theorem 6.25 in \cite{D}):

\begin{theorem}\label{th:Nad} Let $(L, e^{-\phi})$ be a holomorphic line bundle  with smooth metric (i.e. $\phi$ is smooth) on a weakly
pseudoconvex K\"ahler manifold $(X, \omega)$. Assume that there exists an upper semi-continuous function $G$ on $X$ such that $$i\partial\dbar (\phi+G) \geq \epsilon \,\omega$$ for some continuous positive function $\varepsilon$ on X. Then 
\begin{equation}\label{eq:nad}
H^q(X, \mathcal O(K_X+L) \otimes \mathcal I(\phi+G)) =0, \ \ \forall\ q\geq1.
\end{equation}
\end{theorem}
It is known that the short exact sequence $0 \to  I_0  \to I_1  \to  I_2
 \to 0$, where
$$
 I_0:=\mathcal O(K_X+L) \otimes \mathcal I(\phi+G), \  \   I_1:=\mathcal O(K_X+L) \otimes \mathcal J, \ \  I_2:=\mathcal O(K_X+L) \otimes \mathcal J/\mathcal I(\phi+G),
$$
and $\mathcal J\supset\mathcal I(\phi+G)$ is a coherent ideal sheaf of $\mathcal O(X)$,  
induces a  long exact sequence 
$$
0 \to H^0(I_0) \to H^0(I_1) \to H^0(I_2) \to H^1(I_0) \to H^1(I_1) \to \cdots.
$$
Hence \eqref{eq:nad} implies the following extension theorem:

\begin{theorem}\label{th:Nad1}  With the notation above,
$$
H^q(X, \mathcal O(K_X+L) \otimes \mathcal J) \to H^q(X, \mathcal O(K_X+L) \otimes \mathcal J/\mathcal I(\phi+G))
$$
is surjective  for every $q\geq 0$.
\end{theorem}
In general, $H^1(I_0)=0$ is strictly stronger than the surjectivity of $H^0(I_1) \to H^0(I_2)$. Hence one may expect a weaker conditions (without assuming that $\phi+G$ is strictly psh) for surjectivity. Such a result is first obtained by Demailly \cite{Dem18} in the compact case, and later generalizations include Cao-Demailly-Matsumura \cite{CDM} (when $G$ has neat analytic singularity in Theorem \ref{th:CDMZZ}) and Zhou-Zhu \cite{ZZ}. The last-mentioned theorem is the following:

\begin{theorem}\label{th:CDMZZ} 
Let $(L, e^{-\phi})$ be a holomorphic line bundle  with quasi-psh metric (i.e. $\phi$ is quasi-psh) on a K\"ahler manifold $X$ that possesses a proper holomorphic mapping to $\mathbb C^N$. Let $G$ be locally bounded above on $X$ such that
$$
\text{$\phi+G$ and $\phi+\lambda G$ are psh} 
$$
for some positive constant $\lambda >1$ . Then 
$$
H^q(X, \mathcal O(K_X+L)\otimes \mathcal I(\phi) ) \to H^q(X, \mathcal O(K_X+L) \otimes \mathcal I(\phi)/\mathcal I(\phi+G))
$$
is surjective for every $q\geq 0$.
\end{theorem}

\noindent
\textbf{Remark.} \emph{In the case that $q=0$, Theorem \ref{th:CDMZZ} implies that every holomorphic section of  $\mathcal O(K_X+L) \otimes \mathcal I(\phi) $ over the support scheme of $\mathcal I(\phi)/\mathcal I(\phi+G)$ extends to $X$. The theory on the effective $L^2$ estimate of those extensions is usually known as the Ohsawa-Takegoshi theory \cite{OT, Ohsawa}.}

\subsection{The Ohsawa-Takegoshi theory} In this section, we discuss the Ohsawa-Takegoshi theory, and prove a sharp version of the $L^2$-holomorphic extension theorem using our main theorem. The first part of the story concerns a certain $L^2$-norm on $H^0(X, \mathcal O(K_X+L) \otimes \mathcal I(\phi)/\mathcal I(\phi+G))$. The idea is to use the following isomorphism
\begin{equation}\label{eq:iso-smooth}
H^0(X, \mathcal O(K_X+L) \otimes \mathcal I(\phi)/\mathcal I(\phi+G)) \simeq  \frac{\{u\in \mathcal A_{\phi}^{n,0}: 
\dbar u \in\mathcal A_{\phi+G}^{n,1}\} }{\{u\in \mathcal A_{\phi+G}^{n,0}: 
\dbar u \in\mathcal A_{\phi+G}^{n,1}\}},
\end{equation}
where $\mathcal A_\psi^{n,k}$ denotes the space of smooth $L$-valued $(n, k)$-forms on $X$ whose coefficients are locally square integrable with respect to $e^{-\psi}$. The above isomorphism allows us to write an element in $H^0(X, \mathcal O(K_X+L) \otimes \mathcal I(\phi)/\mathcal I(\phi+G))$ as an equivalence class $\{F\}$, where 
$$
F\in \{u\in \mathcal A_{\phi}^{n,0}: 
\dbar u \in\mathcal A_{\phi+G}^{n,1}\}.
$$ 
The surjectivity of 
$$
H^0(X, \mathcal O(K_X+L) \otimes \mathcal I(\phi)) \to H^0(X, \mathcal O(K_X+L) \otimes \mathcal I(\phi)/\mathcal I(\phi+G))
$$
is equivalent to the existence of 
$$
\tilde F \in \{u\in \mathcal A_{\phi}^{n,0}: 
\dbar u =0\} = H^0(X, \mathcal O(K_X+L) \otimes \mathcal I(\phi)) 
$$ 
such that
$$
\tilde F - F \in \mathcal A_{\phi+G}^{n,0}.
$$
That is, we have the following lemma:

\begin{lemma}\label{le:extension} The surjectivity of 
$$
H^0(X, \mathcal O(K_X+L) \otimes \mathcal I(\phi)) \to H^0(X, \mathcal O(K_X+L) \otimes \mathcal I(\phi)/\mathcal I(\phi+G))
$$
is equivalent to that for every $F\in \mathcal A_{\phi}^{n,0}$ with $\dbar F \in \mathcal A_{\phi+G}^{n,1}$, we can find
$
\tilde F \in  H^0(X, \mathcal O(K_X+L) \otimes \mathcal I(\phi)) 
$
such that $\tilde F - F \in \mathcal A_{\phi+G}^{n,0}.$
\end{lemma}
The previous lemma suggests to use H\"ormander $L^2$ estimate for $\dbar$ to study the effective extension theorem. The following general theorem is proved by Zhou-Zhu in \cite{ZZ} (see also \cite{BCP, Chan, CC} for some recent related results):

\begin{theorem}\label{th:ZZ1} Let $(L, e^{-\phi})$ be a holomorphic line bundle, with $\phi$ quasi-psh, over an $n$-dimensional weakly
pseudoconvex K\"ahler manifold $X$. Let $G\leq 0$ on $X$ such that
$$
\text{$\phi+G$ and $\phi+\lambda G$ are psh} 
$$
for some constant $\lambda >1$ on $X$. Then for every $F\in \mathcal A_{\phi}^{n,0}$ with $\dbar F \in \mathcal A_{\phi+G}^{n,1}$ and
\begin{equation}\label{eq:ZZ11}
||F||_G^2:=\sup_{K \ \text{compact in} \ X}\limsup_{t\to -\infty} \int_{\{t<G<t+1\} \cap K} i^{n^2} F\wedge \overline{F} \,e^{-(\phi+G)}<\infty,
\end{equation}
there exists $\tilde F \in  H^0(X, \mathcal O(K_X+L) \otimes \mathcal I(\phi))$ such that $\tilde F - F \in \mathcal A_{\phi+G}^{n,0}$ and
\begin{equation}\label{eq:ZZ12}
\int_{X} i^{n^2} \tilde F \wedge \overline{\tilde F} \, e^{-\phi}  \leq \frac{\lambda}{\lambda-1} ||F||^2_G.
\end{equation}
\end{theorem}

\begin{remark}
Ohsawa \cite{Ohsawa} further observed that one may use \eqref{eq:iso-smooth} to define a norm for $\{F\} \in H^0(X, \mathcal O(K_X+L) \otimes \mathcal I(\phi)/\mathcal I(\phi+G)) $. In fact,  for $F_1, F_2\in \{u\in \mathcal A_{\phi}^{n,0}: 
\dbar u \in\mathcal A_{\phi+G}^{n,1}\}$ with
$$
F_1-F_2 \in \mathcal A_{\phi+G}^{n,0},
$$
we have
\begin{equation}\label{eq:ZZ13}
\limsup_{t\to -\infty} \int_{\{t<G<t+1\} \cap K} i^{n^2} F_1\wedge \overline{F_1} \, e^{-(\phi+G)} = \limsup_{t\to -\infty} \int_{\{t<G<t+1\} \cap K} i^{n^2} F_2\wedge \overline{F_2} \, e^{-(\phi+G)} 
\end{equation}
for every compact set $K$ in $X$, which implies that
$$
||F_1||^2_G = ||F_2||^2_G, \ \ \text{if} \ \{F_1\}=\{F_2\}. 
$$
\end{remark}
The previous remark shows thus that
\begin{equation}\label{eq:ZZ14}
||\{F\}||^2:=||F||^2_G 
\end{equation}
is well-defined. We may therefore rephrase Theorem \ref{th:ZZ1} as follows (see \cite[Theorem 1.2, $R(t)=e^{-t}$]{ZZ}):

\begin{theorem}\label{ZZ2} Let $(L, e^{-\phi})$ be a holomorphic line bundle, with $\phi$ quasi-psh, over an $n$-dimensional weakly
pseudoconvex K\"ahler manifold $X$. Let $G\leq 0$ on $X$ such that
$$
\text{$\phi+G$ and $\phi+\lambda G$ are psh} 
$$
for some constant $\lambda >1$ on $X$.  Then every $\{F\} \in H^0(X, \mathcal O(K_X+L) \otimes \mathcal I(\phi)/\mathcal I(\phi+G))$ with $||\{F\}||^2<\infty$ has an extension $\tilde F \in  H^0(X, \mathcal O(K_X+L) \otimes \mathcal I(\phi))$ with
\begin{equation}\label{eq:ZZ22}
\int_{X} i^{n^2} \tilde F \wedge \overline{\tilde F} \, e^{-\phi}  \leq \frac{\lambda}{\lambda-1} ||\{F\}||^2.
\end{equation}
\end{theorem}

\subsection{The Berndtsson-Lempert method for $L^2$-holomorphic extension} A direct corollary of our main theorem, Theorem \ref{thm: main theorem, monotonicity property of minimal solutions, or holomorpphic quotient sections}, is the following theorem on; for relation with what we have discussed above, $q=0$ below, but the theorem immediately generalizes to also the case that $q\geq 0$:

\begin{theorem}\label{th:BLOT} Let $(L, e^{-\phi})$ be a pseudoeffective line bundle on an $n$-dimensional quasi-complete K\"ahler manifold $(X, \omega)$. Let $G\leq 0$ be a function on $X$ such that
$\phi+\lambda G$ is psh for some constant $\lambda>1$.
Fix
$F\in H^0(X, K_X+L)$ with
$\int_X i^{n^2} F\wedge \bar F \,e^{-\phi} <\infty$. Then there exists 
$
\tilde F \in H^0(X, K_X+L)
$
such that 
$$
\tilde F-F \in H^0(X, \mathcal O(K_X+L)\otimes\mathcal I(\phi+G))
$$ 
and
$$
\int_X i^{n^2} \tilde F \wedge \overline{ \tilde F} \, e^{-\phi}  \leq \frac{\lambda}{\lambda-1} \limsup_{s\to -\infty} e^{-s}\int_{G<s} i^{n^2} F\wedge \bar F \,e^{-\phi}.
$$
\end{theorem}

\begin{proof} By Theorem \ref{thm: main theorem, monotonicity property of minimal solutions, or holomorpphic quotient sections},
$$
t\mapsto A(t):=e^{- t} \inf\left\lbrace \int_X i^{n^2} \tilde F \wedge \overline{\tilde F} \, e^{-\phi -\lambda\max\{G-t,0\}} :   \tilde F-F \in S_1 \right\rbrace
$$
is decreasing in $t \in \mathbb R$. It suffices to show that 
\begin{equation}\label{eq:BLOT1-0}
A(0) \leq \frac{\lambda}{\lambda-1} \limsup_{s\to -\infty} e^{-s}\int_{G<s} i^{n^2} F\wedge \bar F \,e^{-\phi}.
\end{equation} 
But notice that
$$
A(t) \leq e^{-t}   \int_X i^{n^2}  F \wedge \overline{ F} \, e^{-\phi -\lambda\max\{G-t,0\}}.
$$
Hence, the decreasingness property of $A(t)$ implies
$$
A(0) \leq \lim_{t\to -\infty} A(t) \leq \liminf_{t\to -\infty} e^{-t}   \int_X i^{n^2}  F \wedge \overline{ F} \, e^{-\phi -\lambda\max\{G-t,0\}},
$$
and the calculus lemma below gives \eqref{eq:BLOT1-0}.
\end{proof}
In the previous proof we used the following calculus lemma:
\begin{lemma}\label{le:cal} Put $||F||^2_t:=\int_X i^{n^2} F\wedge \bar F \,e^{-\phi-\lambda\max\{G-t,0\}}.$ Then 
$$\limsup_{t\to -\infty} e^{-t}||F||_t^2 \leq \frac{\lambda}{\lambda-1} \limsup_{s\to -\infty} e^{-s}\int_{G<s} i^{n^2} F\wedge \bar F \,e^{-\phi}.$$
\end{lemma}

\begin{proof} Consider a positive Borel measure $\mu$ on $X$ defined by
$$
d\mu:= i^{n^2} F \wedge \bar F \, e^{-\phi}.
$$
Then 
$$
\mu(X)= \int_X i^{n^2} F \wedge \bar F \, e^{-\phi}, \ \ \ \mu(G<s)=\int_{G<s}  i^{n^2} F \wedge \bar F \, e^{-\phi}, 
$$
and we have
$$
||F||_t^2= \int_X e^{-\lambda \max\{G-t,0\}} \, d\mu= e^{\lambda t} \mu(X)- \int_{-\infty}^0  \mu(G<s) \,d \,e^{-\lambda \max\{s-t,0\}},
$$
Since $\mu(X)$ is finite and $\lambda>1$, we have
$$
\lim_{t\to -\infty} e^{-t} e^{\lambda t} \mu(X) =0,
$$
hence
\begin{align*}
\limsup_{t\to -\infty} e^{-t}||F||_t^2 & = \limsup_{t\to -\infty} \left(- e^{-t} \int_{-\infty}^0  \mu(G<s) \,d \,e^{-\lambda \max\{s-t,0\}}\right)\\
& =  \limsup_{t\to -\infty} \left(- e^{-t} \int_{t}^0  \mu(G<s) \,d \,e^{-\lambda(s-t)}\right) \\ 
& =  \limsup_{t\to -\infty} \left(\lambda \int_{t}^0  \mu(G<s) e^{-s} e^{-(\lambda-1)(s-t)}\, ds\right) \\
& =  \limsup_{t\to -\infty} \left(\lambda \int_{0}^{-t}  \mu(G<x+t) e^{-(x+t)} e^{-(\lambda-1)x}\, dx\right) \\
&  \leq   \left(\lambda \int_{0}^{\infty} e^{-(\lambda-1)x}\, dx\right) \limsup_{s\to -\infty} \left(\mu(G<s) e^{-s}\right) \\
&  =  \frac{\lambda}{\lambda-1} \limsup_{s\to -\infty} e^{-s}\int_{G<s} i^{n^2} F\wedge \bar F \,e^{-\phi},
\end{align*}
gives the lemma.
\end{proof}

Finally, we may also rephrase the Theorem \ref{th:BLOT} as a sharp Ohsawa-Takegoshi extension theorem:

\begin{theorem}\label{th:BLOT1} Let $(L, e^{-\phi})$ be a pseudoeffective line bundle on an $n$-dimensional quasi-complete K\"ahler manifold $(X, \omega)$. Let $G\leq 0$ be a function on $X$ such that
$\phi+\lambda G$ is psh for some constant $\lambda>1$.
Fix $\{F\} \in H^0(X, \mathcal O(K_X+L) \otimes \mathcal I(\phi)/\mathcal I(\phi+G))$. If $\{F\}$ has an extension  $F_1\in H^0(X, \mathcal O(K_X+L) \otimes \mathcal I(\phi))$  with
$\int_X i^{n^2} F\wedge \bar F \,e^{-\phi} <\infty$ then we can choose an extension $F_2$ such that
$$
\int_X i^{n^2} F_2 \wedge \overline{F_2} \, e^{-\phi}  \leq \frac{\lambda}{\lambda-1} \limsup_{s\to -\infty} e^{-s}\int_{G<s} i^{n^2} F_1\wedge \overline{F_1} \,  e^{-\phi}.
$$
\end{theorem}
We end our paper with the following remark:
\begin{remark}
Assume that
\begin{equation}\label{eq:BLOT11}
\limsup_{s\to -\infty} e^{-s}\int_{\{G<s\}\cap K} i^{n^2} F_1\wedge \overline{F_1} \,  e^{-\phi} = \liminf_{s\to -\infty} e^{-s}\int_{\{G<s\}\cap K} i^{n^2} F_1\wedge \overline{F_1} \,  e^{-\phi}
\end{equation}
for every compact set $K$ in $X$, then we have
$$
||F_1||_G^2  = \sup_{K \ \text{compact in} \ X} \limsup_{s\to -\infty} e^{-s}\int_{\{G<s\}\cap K} i^{n^2} F_1\wedge \overline{F_1} \,  e^{-\phi}.
$$
Hence 
\begin{equation}\label{eq:BLOT12}
||F_1||_G^2   \leq \limsup_{s\to -\infty} e^{-s}\int_{G<s} i^{n^2} F_1\wedge \overline{F_1} \,  e^{-\phi}.
\end{equation}
Thus if  $X$ is weakly pseudoconvex  and \eqref{eq:BLOT11} holds then the above theorem follows from Theorem \ref{th:ZZ1}.
But even in that special case the proofs are quite different. The main idea in our proof comes from  \cite[page 5]{BoLempOT} ($\lambda$ is denoted by $p$ there). The Berndtsson-Lempert approach  in \cite{BoLempOT} is to first let $\lambda$ go to infinity, then generalize Theorem \ref{thm: general convexity theorem} to the non-product family $
\{G<t\}$. Thus the have to assume that $\phi+ \lambda G$ is psh for all positive constant $\lambda$.  This means that both $\phi$ and $G$ are psh. Hence the method in \cite{BoLempOT} does not directly apply to the compact K\"ahler case (since there is no non-trivial psh $G$ there). In order to generalize the approach in \cite{BoLempOT} to the compact K\"ahler case, we rephrased the proof of  \cite[Theorem 3.8]{BoLempOT} in \cite[Section 10]{wang}.  Eventually we realized that  (see \cite[Section 10.4]{wang}) our new formulation of \cite[Theorem 3.8]{BoLempOT} is equivalent to Theorem \ref{thm: general convexity theorem} for a \emph{fixed} $\lambda>1$. This observation finally leads to the above generalization of the Berndtsson--Lempert theorem; see also \cite{Hosono0, TaiXu, wang21} for other recent applications of the Berndtsson--Lempert approach.
\end{remark}

\chapter{Paper 4}

\begin{center}
\normalsize{\bt{ON A REMARK BY OHSAWA RELATED TO THE BERNDTSSON-LEMPERT METHOD FOR $L^2$-HOLOMORPHIC EXTENSION}}\\[0.5cm]
\small{TAI TERJE HUU NGUYEN AND XU WANG}

\begin{abstract} In \cite[Remark 4.1]{Ohsawa17}, Ohsawa asked whether it is possible to prove Theorem 4.1 and Theorem 0.1 in \cite{Ohsawa17} using the Berndtsson-Lempert method. We  shall answer Ohsawa's question affirmatively in this paper. Our approach also suggests to introduce the Legendre-Fenchel  theory and weak psh-geodesics into the Berndtsson-Lempert method.
\end{abstract}
\end{center}

\section{Introduction}  

In \cite{Ohsawa17} Ohsawa gave a new proof of two theorems of Guan-Zhou (see \cite{GZ, YLZ, GL2} for further details and related results), Theorem  \ref{Theorem A} and \ref{Theorem B} below, and asked in a remark (remark 4.1 in \cite{Ohsawa17}) whether a proof using the Berndtsson-Lempert method in \cite{BoLempOT} can be had.

\medskip

\begin{theorem}[Theorem 4.1 in \cite{Ohsawa17} and Corollary 1.8 in \cite{GZ2}]\label{Theorem A}
{\it Let $\Omega$ be a pseudoconvex domain in $\set{C}^n$, $\Omega':=\{(z',z_n)\in \Omega:z_n=0\}$, $\phi$ a plurisubharmonic function on $\Omega$, $f\in \mathcal{O}(\Omega')$ a holomorphic function on $\Omega'$, and let $\alpha>0$. Then there exists $\tilde{f}\in \mathcal{O}(\Omega)$ a holomorphic function on $\Omega$ such that $\tilde{f}|_{\Omega'}=f$ and
\begin{align*}
\int_{\Omega}|\tilde{f}|^2e^{-\phi-\alpha|z_n|^2}\leq \fr{\pi}{\alpha}\int_{\Omega'}|f|^2e^{-\phi},
\end{align*}}where integration is with respect to the Lebesgue measure.
\end{theorem}

\begin{theorem} [Theorem 0.1 in \cite{Ohsawa17}]\label{Theorem B} {\it With the same assumptions and notation as in Theorem \ref{Theorem A}, there exists $\tilde{f}\in \mathcal{O}(\Omega)$ such that $\tilde{f}|_{\Omega'}=f$ and satisfying the estimate
\begin{align*}
\int_{\Omega}\fr{|\tilde{f}|^2e^{-\phi}}{(1+|z_n|^2)^{1+\alpha}}\leq \fr{\pi}{\alpha}\int_{\Omega'}|f|^2e^{-\phi}.
\end{align*}}
\end{theorem}
In this paper we answer Ohsawa's question affirmatively by proving an extension theorem with estimate using the Berndtsson-Lempert method. This is Theorem \ref{main main theorem} below and the main theorem of this paper. A key ingredient in our proof is the construction of a weight function $\theta$ using the Legendre-Fenchel transform and weak geodesics for plurisubharmonic functions (see the appendix in Section \ref{Appendix}) to which we can apply the Berndtsson-Lempert technique.\\  

Below we use the following terminology. Let $\sigma:U\to \set{R},\;x\mapsto \sigma(x)$ be a real-valued function where $U\sub \set{R}^n$ is some subset and we write $x:=(x_1,\ldots, x_n)$. We will say that $\sigma$ is increasing if $\sigma$ is separately increasing in each  argument. That is, more precisely, if for each $j\in \{1,2,\ldots, n\}$, the one-variable function
$
t\mapsto \sigma(x_1,\ldots, x_{j-1}, t, x_{j+1},\ldots, x_{n})
$ is increasing for all $(x_1,\ldots, x_{j-1},t,x_{j+1},\ldots, x_n)\in U$. 

\begin{theorem}[Main Theorem]\label{main main theorem}
{\it Let $\Omega$ be a pseudoconvex domain in $\set{C}^n$,  $1\leq k \leq n$ an integer,
$$
\Omega':=\{(z',z'')\in \Omega: z''=0\},  \ \ \ z'':=(z_{n-k+1}, \cdots, z_n),
$$
$\phi$ a plurisubharmonic function on $\Omega$, $f\in \mathcal{O}(\Omega')$ a holomorphic function on $\Omega'$, and $\sigma$ a convex increasing function on $\mathbb R^n$. Then there exists $\tilde{f}\in \mathcal{O}(\Omega)$ a holomorphic function on $\Omega$ such that $\tilde{f}|_{\Omega'}=f$ and 
\begin{align*}
\int_{\Omega}|\tilde{f}|^2e^{-\phi-\sigma(\ln|z''|)}\leq \int_{\Omega''} e^{-\sigma(\ln|z''|)} \int_{\Omega'}|f|^2e^{-\phi}, 
\end{align*}
where we write $\ln |z''|:=(\ln|z_{n-k+1}|, \cdots, \ln|z_n|)$, $\Omega'':=\{z''\in \mathbb C^k: \ln |z''| \in 
\Omega''_{\mathbb R}\}$, $ \Omega''_{\mathbb R}$ denotes the convex hull of $\{\ln |z''| : z\in \Omega\}$}, and integration is with respect to the Lebesgue measure.
\end{theorem}

In the case that $\sigma$ depends only on $|z''|^2$, Theorem \ref{main main theorem} reduces to \cite[Theorem 1.7]{YLZ}. It is also very likely that the main theorem in \cite{GZ} implies Theorem \ref{main main theorem}. So our contribution is not the theorem itself, but rather a method of proof using the Berndtsson-Lempert method.

\section{The Berndtsson-Lempert method in a simple setting}\label{BL method}
We start by recalling quickly the Berndtsson-Lempert method for $L^2$-holomorphic extension (\cite{BoLempOT}) in the simple setting of domains in $\set{C}^n$ as it will apply to us. Specifically, we consider the following set-up. Let $\Omega\sub \set{C}^n$ be a pseudoconvex domain and let $\phi$ be a plurisubharmonic function on $\Omega$. For $z\in \Omega$, we write $z=(z',z_n)$. Let $\Omega':=\{z\in \Omega:z_n=0\}$, suppose that $f\in \mathcal{O}(\Omega')$ is a given holomorphic function on $\Omega'$, and suppose that $|z_n|<1$ on $\Omega$. We are interested in finding a holomorphic extension $\tilde{f}\in \mathcal{O}(\Omega)$ of $f$ defined on all of $\Omega$ for which we have a good weighted $L^2$ estimate of the form (\cite{OT})
\begin{align*}
\int_{\Omega}|\tilde{f}|^2e^{-\phi}\leq C\int_{\Omega'}|f|^2e^{-\phi},
\end{align*}
where integration is with respect to the Lebesgue measure, and $C$ is some universal constant. The Berndtsson-Lempert method demonstrated in \cite{BoLempOT} is the following. Fix a positive constant $j>1$ and consider the following plurisubharmonic function 
\begin{align}\label{eq:bl0}
\theta: \mathbb C \times \Omega \to \mathbb R , \quad (\tau, z)\mapsto \theta(\tau,z):=\theta^{{\rm Re}\,\tau}(z):=\phi+ j \max\{\ln|z_n|^2-{\rm Re}\, \tau, 0\}.
\end{align}
For each $t:={\rm Re}\,\tau \leq 0$, $\theta^t$ provides a weight function for a weighted $L^2$ inner product $\ip{\cdot}{\cdot}_t$ and the corresponding induced weighted $L^2$ norm $\norm{\cdot}_t$, given by
\begin{align*}
\ip{u}{v}_{t}&:=\int_{\Omega}u\bar{v}e^{-\theta^{t}},\quad ||u||^2_t:= \ip{u}{u}_{t}.
\end{align*}
Since $|z_n|<1$ on $\Omega$, we have $\theta^{0}=\phi$. Let $f_{t}$ denote the holomorphic extension of $f$ with minimal $\norm{\cdot}_t$-norm. Then we can write  \cite[section 3]{BoLempOT}
\begin{align*}
\norm{f_t}^{2}_{t}&=\sup_{g\in C^{\infty}_{0}(\Omega')}\fr{|\xi_{g}(F)|^2}{\norm{\xi_g}^2_{t}},\quad \xi_g(F):=\int_{\Omega'} F\bar{g}e^{-\phi},
\end{align*}
where $F$ denotes \ita{any} holomorphic extension of $f$, and where $\norm{\xi_{g}}_{t}^2$ denotes the (squared) dual norm of the continuous linear functional $\xi_{g}$, given by 
\begin{align*}
\norm{\xi_{g}}_{t}^2&=\sup_{\norm{F}_{t}=1}|\xi_{g}(F)|^2.
\end{align*}
The crux of the argument is that by complex Brunn-Minkowski theory (see e.g. \cite{Bo3, Bo2, Bo0, Bo1}), under the assumption that $\theta$ is plurisubharmonic, we have that $t\mapsto \norm{\xi_g}_{t}^2$ is log-convex. Hence $t\mapsto e^{t}\norm{\xi_g}_{t}^2$ is log-convex. Suppose that this is also bounded as $t\to-\infty$. Then $t\mapsto t+\ln(\norm{\xi_g}_t^2)$ is convex and bounded as $t\to-\infty$, and hence increasing. Taking the exponential, we infer that $t\mapsto e^{t}\norm{\xi_g}_{t}^2$ is increasing, which implies that $t\mapsto e^{-t}\norm{f_t}_t^2$ is decreasing. Thus, we get the inequality
\begin{align*}
\norm{f_0}^2_{0}&\leq \lim_{t\to-\infty}e^{-t}\norm{f_t}^2_t.
\end{align*}Since $f_t$ is the extension of $f$ with \ita{minimal} $\norm{\cdot}_t$-norm, the right-hand side limit is bounded by the limit of $e^{-t}\norm{F}_{t}^2$ as $t\to-\infty$ with $F$ \ita{any} extension as above. Writing $f_0:=\tilde{f}$ and recalling that we are assuming $\theta^{0}=\phi$, we therefore find
\begin{align*}
\int_{\Omega}|\tilde{f}|^2e^{-\phi}&\leq \liminf_{t\to-\infty}e^{-t}\norm{F}_t^2=\liminf_{t\to-\infty}e^{-t}\int_{\Omega}|F|^2e^{-\theta^{t}}.
\end{align*} 
By \cite[Lemma 3.2, Lemma 3.3]{BoLempOT}, as $j\to \infty$, the right-hand side limit converges to $\pi\int_{\Omega'}|f|^2e^{-\phi}$. Thus, we get an extension $\tilde{f}$ with the following estimate:
\begin{align*}
\int_{\Omega}|\tilde{f}|^2e^{-\phi}&\leq \pi\int_{\Omega'}|f|^2e^{-\phi}.
\end{align*}
This estimate is known to be optimal (we have equality in the case $\Omega=\Omega' \times \{|z_n|<1\}$ and $\phi$ does not depend on $z_n$) and was first proved in \cite{Bl0, GZ}.

\section{A solution to Ohsawa's question}

In this section we apply the Berndtsson-Lempert method in Section \ref{BL method} to a weight function $\theta$ that is different from the one in \eqref{eq:bl0} to prove Theorem \ref{main theorem} below.
Notice that taking
$\sigma=\alpha e^{2x}$ or $\sigma=(1+\alpha)\ln(1+e^{2x})$ we get $\mathrm{L}=\pi/\alpha$ in Theorem \ref{main theorem}. Hence Theorem \ref{Theorem A} and \ref{Theorem B} follow.  Our construction of $\theta$ is based on a Legendre-Fenchel transform approach to weak geodesics for plurisubharmonic functions (see the appendix in Section \ref{Appendix}).

\begin{theorem}\label{main theorem}
Let $\Omega\sub \set{C}^n$ be pseudoconvex and write for $z\in \Omega$, $z=(z',z_n)$. Let $\phi$ be a plurisubharmonic function in $\Omega$, $\Omega':=\{z\in \Omega:z_n=0\}$, and  $f\in \mathcal{O}(\Omega')$ a holomorphic function on $\Omega'$. Put $x:=\ln|z_n|$ and let $\sigma=\sigma(x)$ be a 
convex and increasing function in $x$ with the property that $\mathrm{L}:=\displaystyle \int_{\mathbb C} e^{-\sigma(\ln |w|)} <\infty$. Then there exists $\tilde{f}\in \mathcal{O}(\Omega)$ a holomorphic function on $\Omega$ such that $\tilde{f}|_{\Omega'}=f$ and
\begin{align}
\int_{\Omega}|\tilde{f}|^2e^{-\phi-\sigma(\ln|z_n|)}\leq \mathrm{L}\int_{\Omega'}|f|^2e^{-\phi}.\label{final estimate}
\end{align}
\end{theorem}

\begin{proof}  Fix $c\in \mathbb R$ and consider the weight function $\theta$ defined by
$$
\theta(\tau, z):= \theta^{{\rm Re}\,\tau}(z) := \phi(z)+ \psi^{{\rm Re}\,\tau}(\ln|z_n|), \ \ \  t={\rm Re}\,\tau >0,
$$
where
\begin{align}
\psi^t(x):=\psi(t,x)&:=t\sigma\rbrac{\fr{x-c}{t}+c}.\label{4}
\end{align}
One may directly verify that (see also the appendix in Section \ref{Appendix} for another approach) $\psi$ is convex in $(t,x)$ and increasing with respect to $x$. Since $\phi$ is assumed to be plurisubharmnic, we know that $\theta$ defined above is plurisubharmonic in $(\tau, z)$. By the Berndtsson-Lempert method, using the same notation as in Section \ref{BL method}, it follows that
$$
\rho_c(t):=\ln ||\xi_g||^2_t-2c(t-1)
$$ 
is convex as a function of $t$. Hence for all $t\in [1/2, 1)$, 
$$
\frac{\rho_c(1)-\rho_c(t)}{1-t} \geq  \frac{\rho_c(1)-\rho_c(1/2)}{1-1/2}.
$$
Observe that $\rho_c(1)= \ln ||\xi_g||^2_1$ does not depend on $c$. Thus,
$$
\ln ||\xi_g||^2_1 \geq \rho_c(t) + 2(1-t) \left(\ln ||\xi_g||^2_1-\rho_c(1/2) \right), 
$$ for all $t\in [1/2,1)$ \ita{and} for all $c\in \set{R}$.
Now let
$$
c:= \frac{-1}{(1-t)^2},
$$
and let $t \uparrow 1$. Lemma \ref{le:1} below implies that
\begin{align}
\ln ||\xi_g||^2_1 \geq  \limsup_{t\uparrow 1}  \rho_c(t).\label{convex ineq}
\end{align}
Hence, by \eqref{convex ineq}, we have a holomorphic extension $\tilde f$ with the following estimate
\begin{align*}
\int_{\Omega}|\tilde{f}|^2e^{-\phi-\sigma(\ln|z_n|)}&\leq  \liminf_{t \uparrow 1} e^{2c(t-1)}  \int_{\Omega}  |F|^2 e^{-\phi-\psi^t(\ln|z_n|)}  \\
&= \liminf_{t \uparrow 1}  e^{2/(1-t)}  \int_{\Omega}|F|^2  e^{-\phi-t \sigma\left( \frac{\ln|z_n|+\frac1{1-t}}{t}\right)} ,
\end{align*}
where $F$ is any arbitrary holomorphic extension. From the change of variables
$$
z_n= e^{\frac{1}{t-1}} w, \ \ \Omega_t:= \{(z', w) \in \mathbb C^n: (z', e^{\frac{1}{t-1}} w) \in \Omega\}, 
$$
and the definition of $\mathrm{L}$, we find that letting $t\uparrow 1$, 
$$
e^{2/(1-t)}  \int_\Omega |F|^2  e^{-\phi-t \sigma\left( \frac{\ln|z_n|+\frac1{1-t}}{t}\right)} = \int_{\Omega_t} |F(z', e^{\frac{1}{t-1}} w)|^2 e^{-\phi\left(z', e^{\frac{1}{t-1}} w\right) - t \sigma \left(\frac{\ln |w|}{ t}\right)}
$$
converges to the right hand side of \eqref{final estimate}. Hence the theorem follows.
\end{proof}

\begin{lemma}\label{le:1} Let $c:= \frac{-1}{(1-t)^2}$. Then $
\lim_{t \uparrow 1} (1-t) \rho_c(1/2) =0$.
\end{lemma}

\begin{proof} It suffices to show that $\limsup_{c\to -\infty} \rho_c(1/2) <\infty$. That is, 
$$
 \limsup_{c\to 
 \infty}  \sup_{F \in \mathcal O(\Omega)}   \frac{\left|\displaystyle\int_{\Omega'} F\bar{g}e^{-\phi}\right|^2}{e^{-c}\displaystyle\int_{\Omega} |F|^2  e^{-\phi-\frac12 \sigma\left( 2\ln|z_n|-c\right)} } <\infty.
$$
This now follows from the change of variables $z_n= e^{c/2} w$ and the definition of $\mathrm{L}$.
\end{proof}

\textbf{Remark}: The change of variables argument in the proof of Theorem \ref{main theorem} also suggests to use another weight function, with
\begin{equation}\label{eq:new1}
\psi^t(x)=\sigma(x-t), \ \  t \in \mathbb R.
\end{equation}
In fact, as we shall see next, the proof of  Theorem \ref{main theorem} is simpler if we use this new weight function. Nevertheless, we still wish to include the proof given above since it is our first example of the use of weak geodesics in the Berndtsson-Lempert method. 

\subsection[Proof of the main theorem]{Proof of the main theorem (Theorem \ref{main main theorem})}
\begin{proof}[\ita{Since $\Omega\subset \mathbb C^{n-k} \times \Omega''$, one may assume that $\sigma=\infty$ outside $\Omega''_{\mathbb R}$}] Also we may assume that the right-hand side in the estimate in the theorem is finite, else there is nothing to prove.
Replace $\psi^t$ in \eqref{4} by
\begin{equation}\label{eq:new2}
\psi^t(x)=\sigma(x-t), \ \  t \in \mathbb R, \ \  x-t:=(x_1-t, \cdots, x_k-t).
\end{equation}
It is clear that $\psi^t(x)$ in convex in $(t,x)$ and increasing in $x$. Hence, the corresponding weight function 
$$
\theta(\tau, z)= \phi(z)+ \psi^{{\rm Re}\,\tau}(\ln|z_n|)
$$
is plurisubharmonic in $(\tau, z)$, and the Berndtsson-Lempert method applies. Thus, we know that $t\mapsto \ln \norm{\xi_g}_t^2-2kt$ is convex. Moreover, the change of variables $z''= e^t w$ gives, 
\begin{align}\label{eq:rotation}
e^{2kt} \int_{\Omega}  |F|^2 e^{-\phi-\sigma(\ln|z''|-t)} & =   \int_{\{(z', e^t w) \in \Omega\}} |F(z', e^t w)|^2 e^{-\phi(z', e^t w)-\sigma(\ln|w|)}, 
\end{align} 
where $F$, as before, is any arbitrary fixed extension of $f$. From \eqref{eq:rotation} we infer that $\ln \norm{\xi_g}_t^2-2kt$ is also bounded near $t=-\infty$. Hence, $t\mapsto \ln \norm{\xi_g}_t^2-2kt$ is increasing and there exists a holomorphic extension $\tilde f$ with 
\begin{align*}
\int_{\Omega}|\tilde{f}|^2e^{-\phi-\sigma(\ln|z''|)} & \leq 
 \liminf_{t\to -\infty}   \int_{\{(z', e^t w) \in \Omega\}} |F(z', e^t w)|^2 e^{-\phi(z', e^t w)-\sigma(\ln|w|)}   \\
  &  \leq \int_{\mathbb C^k} e^{-\sigma(\ln|w|)} \int_{\Omega'}|f|^2e^{-\phi}  \\
   \ & = \int_{\Omega''} e^{-\sigma(\ln|w|)} \int_{\Omega'}|f|^2e^{-\phi},
\end{align*} 
where in the last equality we have used that \text{$\sigma=\infty$ outside $\Omega''_{\mathbb R}$}.
This completes the proof. \end{proof}

\medskip

\textbf{Remark}: From \eqref{eq:rotation}, in the case that $\Omega$ is horizontally balanced (that is, $(z', z'')\in \Omega \Rightarrow (z', \tau z'') \in \Omega$ for all $\tau \in \mathbb C$ with $|\tau|<1$), the method above is equivalent to the one associated to the weight function
$$
\theta(\tau, z):=\phi(z', \tau z'')+\sigma(\ln |z''|).
$$
This construction has already been used in \cite[Section 2]{Bern20}. The advantage of \eqref{eq:new2} is that it also applies to general $\Omega$. 

\section{An application to Bergman kernels asymptotics}
As an interesting application of Theorem \ref{main main theorem}, we shall now give a new simple and short proof of a well-known asymptotic lower bound for the Bergman kernel.

Let $\Omega$ be a pseudoconvex domain in $\set{C}^{n}$ which contains the origin and let $\phi$ be a twice continuously differentiable plurisubharmonic function on $\Omega$ such that $i\partial\dbar \phi(0)>0$. Let $k\geq 1$ be an integer and let $K_{k}(0)$ denote the Bergman kernel on $\Omega$ with weight $e^{-k\phi}$ on the diagonal at the origin. The following is a well-known asymptotic lower bound for $K_k(0)$:

\begin{theorem}\label{intro: lb asymptotic for Bergman kernel}
With the above notation, 
\begin{align}
\liminf_{k\to \infty}\fr{K_{k}(0)e^{-k\phi(0)}}{k^{n}}&\geq \fr{(i\partial\dbar \phi)_{n}(0)}{\pi^{n}}.
\end{align}
\end{theorem}
One possible proof uses Hörmander $L^2$-estimates for the $\dbar$-equation. We shall here give a new, and in our opinion simpler and shorter, proof, which instead uses Theorem \ref{main main theorem}.  Of course, both proofs utilize non-trivial results, the first the $L^2$-estimates of Hörmander, and the second, an Ohsawa-Takegoshi type of extension theorem. The point is that the proof that we shall discuss is simpler and shorter given both these results. To prove the extension theorem that we shall employ, we used the Berndtsson-Lempert method. Our proof here is therefore complex Brunn-Minkowski theoretical.

\subsection{Proof of Theorem \ref{intro: lb asymptotic for Bergman kernel}}
Let $\ld_1,\ldots,\ld_n$ denote the eigenvalues of the complex Hessian of $\phi$ at the origin. By a unitary change of variables, we may assume that $i\partial \dbar \phi=i\s{j=1}{n}\ld_{j}dz_{j}\wedge d\bar{z}_{j}$. Here, we let $z:=(z_1,\ldots, z_n)$ be generic coordinates on $\Omega$. The key observation is the following simple identity. Let $x:=(x_1,\ldots,x_n):=\ln\abs{z}:=(\ln\abs{z_1},\ldots, \ln\abs{z_n})$, and let $\Delta(0,r)$, where $r:=(r_1,\ldots, r_n)$ be the polydisc centred at the origin with polyradius $r$. Then
\begin{align}
\int_{\Delta(0,r)}e^{-\s{j=1}{n}\ld_{j}\abs{z_{j}}^2}&=\fr{\pi^{n}}{\ld_1\cdots\ld_n}\left(1-e^{-\ld_1 r_1^2}\right)\cdots\left(1-e^{-\ld_n r_n^2}\right).\label{eq: key observation}
\end{align}Note that this tends to the inverse of the right-hand side in Theorem \ref{intro: lb asymptotic for Bergman kernel}. Let $\epsilon>0$. By continuity of $\partial \dbar \phi$ at the origin, there are $\delta_{j}>0$ such that if $\abs{z_j}<\delta_j$, for each $j$, we have $i\partial\dbar \phi(z)-i\s{j=1}{n}\ld_{j}dz_{j}\wedge d\bar{z}_{j}>-i\epsilon \partial\dbar \abs{z}^2$. We actually want the left-hand side to be positive, so we replace $\ld_j$ with some smaller $\eta_j$; the idea is that we can let $\eta_j\to \ld_j$ in the limit later when we let $\epsilon\to 0$. More precisely, let $\eta_j$, for each $j$, be such that $\ld_j-\eta_j>2\fr{\epsilon}{n}$. Then it follows that on $\Delta(0,\delta)$, where we write $\delta:=(\delta_1,\ldots, \delta_n)$, 
\begin{align}
i\partial \dbar \phi(z)-i\s{j=1}{n}\eta_{j}dz_{j}\wedge d\bar{z}_{j}>0.\label{eq: positivity of psik on small disc}
\end{align}
It follows that if we define $\sigma_{k}:=k\s{j=1}{n}\eta_{j}\abs{z_{j}}^2$ on $\Delta(0,\delta)$, then, there, 
\begin{align}
\psi_{k}&:=k\phi-\sigma_k
\end{align}is plurisubharmonic. The idea is to define $\sigma_k$ outside $\Delta(0,\delta)$ such that $\psi_k$ is plurisubharmonic on all of $\Omega$. This is really the crux of the proof and we do this as follows. We define $\sigma_k=\sigma_k(\ln\abs{z}):=f_k(x)$, to be a function of $\abs{z}$, where $x=\ln\abs{z}$ as above. What we need to define is then $f_k$ on $\{x_{j}\geq \ln\abs{\delta_j}\}$. We will define $f_k$ to be the tangent space of $x\mapsto k\s{j=1}{n}\eta_{j}e^{2x_{j}}\left(=k\s{j=1}{n}\eta_j\abs{z_j}^2\right)$. That is, we define $f_k=f_k(x)$ in general by
\begin{align}
f_{k}(x)&:=\left\{\begin{array}{cccc}
k\s{j=1}{n}\eta_j e^{2x_{j}}, &\{x_{j}<\ln\abs{\delta_j}\}\\
k\s{j=1}{n}\eta_{j}\delta_{j}^2+2k\s{j=1}{n}\eta_{j}\delta_{j}^2(x_{j}-\ln\abs{\delta_j}),&\{x_j\geq \ln\abs{\delta_j}\}.
\end{array}\right.
\end{align}
Let $\delta_{k}:=(\delta_{k,1},\ldots, \delta_{k,n})$ be such that $\delta_{k,j}\to 0$ as $k\to \infty$.  Then for sufficiently large $k$, $\Delta(0,\delta)\sub \Delta(0,\delta)$. We choose $k$ sufficiently large such that this holds and replace above $\delta$ everywhere with $\delta_k$. Hence, we define
\begin{align}
\sigma_k(z)&:=\left\{\begin{array}{cccc}
k \s{j=1}{n}\eta_j\abs{z_j}^2,&z\in \Delta(0,\delta_k)\\
k\s{j=1}{n}\eta_j\left(\delta_{k,j}^2+2\delta_{k,j}^2\left(\ln\abs{z_j}-\ln\abs{\delta_{k,j}}\right)\right),&z\in \set{C}^{n}\backslash \Delta(0,\delta_k).
\end{array}\right.\label{eq: definition of sigma k}
\end{align}Finally, we let $\psi_k:=k\phi-\sigma_k$ (defined on $\Omega$ as $\phi$ is). It is clear that $\sigma_k=\sigma_k(x)$, as a function of $x$, is increasing and convex. Moreover, by construction, $\psi_k$ is plurisubharmonic and we have $\psi_k(0)=k\phi(0)$. Hence, by Theorem \ref{main main theorem}, there exists a holomorphic function $f\in \fancy{O}(\Omega)$ on $\Omega$ such that
\begin{align}
\int_{\Omega}\abs{f}^2e^{-k\phi}&\leq \left(\int_{\set{C}^{n}}e^{-\sigma_{k}(\ln\abs{z})}\right)\abs{f(0)}^2e^{-k\phi(0)}.
\end{align}Let the parenthesis above be denoted by $L_k$ ; it is easy to see that $L_k<\infty$. By definition of $K_k(0)$ (or, if one wants, by the extremal property of the Bergman kernel), it follows then that

\begin{align}
\fr{K_k(0)e^{-k\phi(0)}}{k^{n}}&\geq \fr{1}{k^{n}L_{k}}.\label{eq: lb asymptotic part 1 inequality}
\end{align}
A simple computation shows that $L_{k}=L_{k}^{(1)}+L_{k}^{(2)}$, where $L_{k}^{(1)}$ is given from \eqref{eq: key observation} by
\begin{align}
L_{k}^{(1)}&=\int_{\Delta(0,\delta_k)}e^{-\sigma_k(\abs{z})}=\pi^{n}\fr{\left(1-e^{\eta_1 k\delta_{k,1}^2}\right)\cdots\left(1-e^{-\eta_n k\delta_{k,n}^2})\right)}{k^{n}\ld_1\cdots\ld_n},
\end{align}and
\begin{align*}
L_k^{(2)}&=\int_{\set{C}^{n}\backslash \Delta(0,\delta_k)}e^{-k\s{j=1}{n}\eta_{j}\delta_{k,j}^2-2k\s{j=1}{n}\eta_{j}\delta_{k,j}^2\left(\ln\abs{z_j}-\ln\abs{\delta_{k,j}}\right)}=\prod_{j=1}^{n}\fr{e^{-k\eta_{j}\delta_{k,j}^2}\delta_{k,j}^2}{k\left(\eta_j\delta_{k,j}^2-\fr{1}{k}\right)},
\end{align*}where we choose that $k\delta_{k,j}^2\to \infty$ as $k\to \infty$ (which we can certainly do); this is used above to choose $k$ so large that $1-2k\eta_j \delta_{k,j}^2<0$. Then for such large $k$, 
\begin{align}
k^{n}L_{k}&=\pi^{n}\prod_{j=1}^{n}\fr{1-e^{-k\eta_{j}\delta_{k,j}^2}}{\ld_{j}}+\prod_{j=1}^{n}\fr{e^{-k\eta_{j}\delta_{k,j}^2}\delta_{k,j}^2}{\eta_{j}\delta_{k,j}^2-\fr{1}{k}}=\pi^{n}\prod_{j=1}^{n}\fr{1-e^{-k\eta_{j}\delta_{k,j}^2}}{\ld_{j}}+\prod_{j=1}^{n}\fr{e^{-k\eta_{j}\delta_{k,j}^2}}{\eta_{j}-\fr{1}{k\delta_{k,j}^2}}.
\end{align}Letting $k\to \infty$, recalling that $\delta_{k,j}\to 0, k\delta_{k,j}^2\to \infty$ (consider for example $\delta_{k,j}:=1/k^{1/3}$), we can then let $\epsilon \to 0$ and \ita{then} $\eta_j\to \ld_j$, which together with \eqref{eq: lb asymptotic part 1 inequality} gives
\begin{align}
\liminf_{k\to \infty}\fr{K_k(0)e^{-k\phi(0)}}{k^{n}}&\geq \fr{\ld_1\cdots\ld_n}{\pi^{n}}=\fr{(i\partial\dbar \phi)_n(0)}{\pi^{n}}.
\end{align}This completes the proof for Theorem \ref{intro: lb asymptotic for Bergman kernel}. Note that the \ita{actual proof} (that is, without as many details as provided above) is very short (and simple): Define $\psi_k:=k\phi-\sigma_k$ with $\sigma_k$ defined as in \eqref{eq: definition of sigma k} where $\delta_{k,j}\to 0$ and $k\delta_{k,j}^2\to \infty$ as $k\to \infty$, and where $k$ is so large that \eqref{eq: positivity of psik on small disc} holds on $\Delta (0,\delta_k)$, which is possible by continuity. Since $\sigma_k$ is defined using the tangent plane, it is convex as a function of $x$ and Theorem \ref{intro: lb asymptotic for Bergman kernel} follows by \eqref{eq: lb asymptotic part 1 inequality} using Theorem \ref{main main theorem}. 

We end with a small remark that equality actually holds in Theorem \ref{intro: lb asymptotic for Bergman kernel}. The opposite inequality (which is much simpler to establish) can be shown by using the mean-value property of holomorphic functions.

\section{Appendix: weak geodesics in the space of toric plurisubharmonic functions}\label{Appendix}

By a toric plurisubharmonic function we mean a function of the following form
$$
\bm \psi(z):=\psi( \log|z_1|, 
\cdots, \log|z_n|)
$$
on $\mathbb C^n$, where $\psi$ is a convex increasing function on $\mathbb R^n$.  

\begin{definition} We call a family, say $\{\bm \psi^t\}_{0<t<1}$, of toric plurisubharmonic functions a weak geodesic \footnote{For usual weak psh geodesics on compact K\"ahler manfiolds, see \cite[Section 2]{BB}. The weak geodesic is also called generalized geodesic in \cite[Section 2.2]{Bern15}, and maximal psh segment or psh geodesic segment in \cite[page 5]{R}. See \cite{M} for the background.} if 
$$
\bm \Psi : (\tau ,z) \mapsto \bm\psi^{{\rm Re}\,\tau}(z), \ \ \tau\in \mathbb C_{(0,1)}:=\{\tau\in\mathbb C: 0<{\rm Re}\,\tau<1\} 
$$
is plurisubharmonic on $\mathbb C_{(0,1)} \times \mathbb C^n$ and $(i\partial\dbar \bm \Psi)^{n+1}=0$ on $\mathbb C_{(0,1)} \times (\mathbb C\setminus\{0\})^n$.
\end{definition}

\textbf{Remark}: Noticing that $\bm\Psi$ is locally bounded on $\mathbb C_{(0,1)} \times (\mathbb C\setminus\{0\})^n$, we know that $(i\partial\dbar \bm \Psi)^{n+1}$ is well defined. Moreover, by a local change of variables $z_j=e^{w_j}$, we know that $\{\bm \psi^t\}_{0<t<1}$ is a weak geodesic if and only if
$$
\psi(t,x):=\psi^t(x)
$$
is convex on $(0,1) \times \mathbb R^n$ and
$$
MA(\psi):=\det (D^2_{(t,x)} \psi) =0 
$$
on  $(0,1) \times \mathbb R^n$, where $D^2_{(t,x)} \psi$ denotes the real Hessian matrix of $\psi$ with respect to $(t,x)$. Hence, by Corollary 2.5 in \cite{BoConvex}, $\{\bm \psi^t\}_{0<t<1}$ is a weak geodesic if and only if the Legendre-Fenchel transform of $\psi(t,x)$ with respect to $x$ for fixed $t$, say 
$$
(\psi^{t})^{*}(\xi):= \sup_{x\in \mathbb R^n} x\cdot \xi -\psi^t(x),
$$
is an affine function in $t$. This suggest the following definition.

\begin{definition} We call a family, say $\{\bm \psi^t\}_{0 < t < 1}$, of toric plurisubharmonic functions, a weak geodesic segment if 
$$
(\psi^t)^*= t (\sigma^1)^*+(1-t) (\sigma^0)^* 
$$
with $\sigma^1, \sigma^0$ convex increasing functions on $\mathbb R^n$. 
\end{definition}

Now let us consider the case $n=1$. Let  $\sigma^1:=\sigma$ in Theorem \ref{main theorem}. For $\sigma^0$ we take
\begin{align*}
\sigma^0:=\id{(-\infty, c]}(x)=\left\{\begin{array}{ccc}
0,&x \leq c\\
\infty,&x > c.
\end{array}\right.
\end{align*}
Then we have
$$
 (\sigma^0)^*(\xi)=\sup_{x\in \mathbb R} x \xi -\sigma^0(x) = \sup_{x \leq c}  x\xi =c \xi+ \sup_{x\leq 0} x\xi= \begin{cases}
 c\xi & \xi \geq 0 \\
\infty & \xi <0.
 \end{cases}
$$
Since $\sigma^1 =\sigma$ is also increasing, we know that $(\sigma^1)^*(\xi)=\sigma^*(\xi) =\infty$ when $\xi <0$. Hence
$\psi^t$ can be written as
$$
\psi^t(x) =\sup_{\xi \in \mathbb R} x\xi-t\sigma^*(\xi)-(1-t) c\xi = t  \left(\sup_{\xi \in \mathbb R} \frac{x-(1-t)c}{t}\,\xi-\sigma^*(\xi) \right),
$$
and $\sigma^{**} =\sigma$ gives $\psi(t,x)=t\sigma\rbrac{\fr{x-c}{t}+c}$. This
is precisely \eqref{4}. Note that \eqref{eq:new1} also gives a weak geodesic since
$$
\sup_{x\in \mathbb R}  x\xi -\sigma(x-t) = t \xi +\sigma^*(\xi)
$$
is affine in $t$.

\end{document}